\documentclass[a4paper,11pt]{amsart}
\usepackage{amssymb,amsmath,amsthm,stmaryrd,bm,enumitem,hyperref,xcolor}

\def\res{\hbox{ {\vrule height .3cm}{\leaders\hrule\hskip.3cm}}\hskip5.0\mu}

\newcommand\beqn{\begin{equation}}
\newcommand\eeqn{\end{equation}}
\newcommand\beqny{\begin{eqnarray}}
\newcommand\eeqny{\end{eqnarray}}
\newcommand\beqnyn{\begin{eqnarray*}}
\newcommand\eeqnyn{\end{eqnarray*}}

\usepackage{amsmath,amsfonts,amssymb,amsthm}
\newtheorem{theorem}{Theorem}[section]
\newtheorem{lemma}[theorem]{Lemma}

\newtheorem{corollary}[theorem]{Corollary}
\newtheorem{definition}[theorem]{Definition}
\newtheorem{remark}[theorem]{Remark}

\def\b{\beta}          
\def\g{\gamma}           
\def\d{\delta}         
\def\e{\epsilon}         
\def\z{\zeta}            
\def\th{\theta}          
\def\k{\kappa}

\def\t{\tau}

\def\s{\sigma}

\def\G{\Gamma}

\def\Th{\Theta}

\numberwithin{equation}{section}

\newcommand{\op}[1]{\operatorname{\text{\rm #1}}}

\title[Analysis of singularities of area-minimizing currents]{Analysis of singularities of area-minimizing currents, Part~II: a uniform height bound, estimates away from 
branch points of rapid decay, and uniqueness of tangent cones}
\author{Brian Krummel \& Neshan Wickramasekera} 

\begin{document} 
\begin{abstract}
This paper is the second in a five-part series developing a new geometric framework to study the structure of $n$-dimensional locally area-minimizing rectifiable currents $T$ of codimension $\geq 2$. We work in ${\mathbb R}^{n+m}$ ($m \geq 2)$
for Parts~I--IV, extending the results to general Riemannian ambient spaces in Part~V. Our overarching 
program is a generalization of the classical two-dimensional theory established by White, Chang, and Micallef--White, and provides a unified approach to several basic questions: uniqueness of tangent cones at typical singular points; higher-order asymptotics at typical branch points; the size and  structure of the singular set; and criteria for the current to be a topological $n$-disk near branch points.
In particular, while relying on several ingredients of Almgren's foundational theory, this work yields a geometrically more direct proof of the main outcome of that theory: the optimal singular-set Hausdorff dimension bound $n-2.$  

Within this broader context, the present paper focuses on the fine structure of the set ${\mathcal S}$ in  
the decomposition ${\rm sing} \, T  \cap \{Z \, : \, \Theta(T, Z)=q\}= {\mathcal B} \cup {\mathcal S}$  obtained in Part~I for  density-$q$ singularities (for arbitrary integer $q \geq 2$). Part~I established: (i) at each point in ${\mathcal B}$, $T$  
decays locally uniformly to a unique tangent plane \emph{rapidly}---i.e.\ with order-of-contact with the tangent plane (or planar frequency, introduced in Part~I)
finite and bounded below by $1+\alpha$ for some universal 
$\alpha \in (0, 1);$ (ii) at each point in ${\mathcal S}$, $T$ satisfies a locally uniform weak non-planar  approximation property  (defined in Part~I). Building on this, the present paper establishes two main results:
\begin{itemize}
\item[1.] ${\mathcal H}^{n-2}$ a.e.\ point in ${\mathcal S}$ admits a unique tangent cone supported on a union of two or more planes intersecting along an $(n-2)$-dimensional subspace; consequently, possible branch points in ${\mathcal S}$ form an ${\mathcal H}^{n-2}$-null set, and $T$ has a unique tangent cone at ${\mathcal H}^{n-2}$-a.e.\ point. 
\item[2.] 
 The set ${\mathcal S}$ is locally  $(n-2)$-rectifiable (with locally finite ${\mathcal H}^{n-2}$ measure).
\end{itemize}

Consequently, rapid decay to a unique tangent plane is the generic behavior at branch points, occuring at ${\mathcal H}^{n-2}$-a.e.\ branch point.

A key ingredient in the proof of these results is a new uniform height bound:
if a finite collection of planes that are (quantitatively) disjoint in a cylinder, and the current $T$, have small tilt-excess relative to the 
base-plane of the cylinder, then in a smaller cylinder the pointwise distance of $T$
to the union of the planes is bounded 
\emph{linearly} from above by the $L^{2}$ height-excess of $T$ relative to this union.
This, together with techniques inspired by the work of L.~Simon (\cite{Sim93}) and our previous work (\cite{Wic14}, \cite{KrumWic2}),
leads to a decay estimate for $T:$ if, among other requirements, $T$ at some scale lies significantly (weakly) closer to a union of planes intersecting along an 
$(n-2)$-dimensional subspace than to any single plane, then at a fixed smaller scale the distance of $T$ to a  similar union of planes is  no more than half the distance of $T$ at the original scale to the original planes. By Part~I this decay estimate is iteratively and indefinitely applicable at ${\mathcal H}^{n-2}$-a.e.\ point in the set ${\mathcal S},$ giving the above two conclusions. 

Parts~III~and~IV analyze the current near ${\mathcal B}$ using its built-in uniform decay estimate and further new height-bounds for $T$. Part~III treats points of planar frequency 
$\neq 2$ using the monotonicity formula for planar frequency (established in Part~I) and without the use of 
Almgren center manifolds (the use of which for arbitrary branch points is a technically demanding feature of the classical  theory). Part~IV analyzes  the frequency-$2$ points---the only instance a center manifold is used in our program---crucially capitalizing on the fact that the center manifold arises canonically in this case: fixed density frequency-2 points lie locally in the critical nodal set of a \emph{single} approximating multi-valued map over a \emph{single} center manifold. Together, Parts~III~and~IV provide structure results for ${\mathcal B}$ and a frequency condition for $T$ to be a topological $n$-disk near a point in ${\mathcal B}$. 
\end{abstract}

\maketitle
\tableofcontents

\section{Introduction}
This is the second part of our study of interior singularities of $n$-dimensional locally area minimizing rectifiable currents in an open subset of ${\mathbb R}^{n+m}$  for $m \geq 2$. Extensions to general Riemannian ambient spaces are discussed in the final part of the series. In this introduction we shall focus on the  
results established in the present article and their dependence on the first part \cite{KrumWica} of the work, while describing the broader context only briefly. We refer the reader to the introduction to \cite{KrumWica} for a more detailed account of our overall approach including the main new ideas and some historical context of the problem. 

It is well known that  when the codimension $m \geq 2,$ unlike when $m=1$, a locally area minimizing rectifiable current $T$ can have (interior) branch point singularities, i.e.\ non immersed points of ${\rm spt} \, T \setminus {\rm spt} \, \partial T$ at which one tangent cone to $T$ is supported on a plane.  Non-branch-point singularities of an area minimizer can be characterized as those points where every tangent cone is translation invariant along a linear subspace of dimension at most $n-2$. By  using this tangent cone criterion and an elementary argument to bound the size of the set of non-branch-point singularities, and by developing a number of pioneering fundamental ideas to bound the size of the branch points separately, Almgren established that the Hausdorff dimension of the singular set of $T$ is $\leq n-2$. (\cite{Almgren}; see also the series of papers \cite{DeLSpa1}, \cite{DeLSpa2}, \cite{DeLSpa3} by De~Lellis and Spadaro which present Almgren's lengthy argument in a more accessible form with some technical aspects clarified and streamlined.) This is the sharp general Hausdorff dimension upper bound on the singular set. 

Apart from the question of the size of the singular set, central in the study of singularities are the questions of uniqueness of tangent cones at singular points, 
asymptotic/topological nature of the current on approach to singular points and the local structure of the singular set. While these questions for area minimizers of dimension $n=2$ have long been settled (\cite{W1}, \cite{Chang}, \cite{MicWhi}), little has been known in these directions when $n \geq 3$.  Our work is aimed at addressing these questions for general $n$ (while, as it turns out, not relying on the size bound on the singular set provided by \cite{Almgren}---in fact while providing, as a by product, a considerably more efficient proof of the dimension bound than that of \cite{Almgren}). Unlike in two-dimensional area minimizers,  singularities in higher-dimensional area minimizers need not be isolated points, and this is a key difference that makes the analysis of singularities in higher dimensions more challenging. 

The first step in our approach is to establish first order (i.e.\ tangent cone level) decay estimates for $T$ at a typical (i.e.\ ${\mathcal H}^{n-2}$ a.e.) singular 
point, which we complete in the present article building on the work in \cite{KrumWica}. 
These estimates imply, in particular, uniqueness of tangent cones to $T$ at ${\mathcal H}^{n-2}$ a.e.\ point.  They also allow us to obtain local structure results for the set of 
non-branch-point singularities, which we include in the present work (in Theorem~\ref{main-results-combined} below). Subsequently, in \cite{KrumWicc}, \cite{KrumWicd}, these first order decay estimates provide the basis for obtaining higher order estimates at typical branch points---that is, the existence, at ${\mathcal H}^{n-2}$ a.e.\ branch point $Z$, of a unique, non-zero Dirichlet energy minimizing multi-valued tangent function (blow-up) $\varphi_{Z}.$ 
The significance of these higher order estimates is two-fold: first, they readily imply local structural properties of the branch set (that is, the fact that the branch set is locally the union of finitely many pairwise disjoint locally $(n-2)$-rectifiable sets); and secondly, they serve as a step towards understanding the topology of the current near typical branch points, and indeed lead to a description of the current as a multi-valued $C^{1, \alpha}$ graph near certain branch points.

The key advantage in proceeding by first addressing the uniqueness of tangent cones question, as we do (in contrast to \cite{Almgren}), is that once the first order decay estimates needed for that purpose are in place, branch points can be analysed based on the rate of decay towards the (unique) tangent plane. For branch points where decay of the current towards a unique tangent plane is quadratic (or faster) in scale, it turns out that using  Almgren's ``center manifold'' is much simpler than for other branch points. Our approach capitalises on this fact by removing the need for center manifolds altogether in the analysis except for those branch points where decay to the tangent plane is precisely quadratic. 
Thus, as mentioned above, along the way our approach yields a considerably simpler proof of Almgren's main result that the Hausdorff dimension of the singular set is $\leq n-2$ (for which not all of our decay estimates are necessary; in particular, no result from the present article 
is necessary nor are any of the higher-order asymptotics conclusions in  \cite{KrumWicc}, \cite{KrumWicd}). See 
\cite[Corollary~1.3]{KrumWica} where these simplifications of the dimension bound is outlined, and 
\cite[Theorem~3.7]{KrumWicc}, \cite[Theorem~3.13]{KrumWicd} where they are explained in more detail.

We obtain the first order decay estimates by decomposing the singular set in a way different from the decomposition based on the tangent cone type (namely, as branch points or non branch points) considered in \cite{Almgren}. Specifically, fixing an integer $q \geq 2$, we consider the set ${\mathcal B}_{q}$ 
of points where the current upon rescaling converges \emph{rapidly} to a unique multiplicity $q$ tangent plane, and the complementary set ${\mathcal S}_{q} = {\rm sing}_{q} \, T \setminus {\mathcal B}_{q}$, where ${\rm sing}_{q} \, T$ is the set of points $Y$ with density $\Theta(T, Y) = q$; here rapid decay at $Z \in {\mathcal B}_{q}$ means decay towards the tangent plane at $Z$ at the rate $o(|X - Z|^{1+\alpha})$, $X \in {\rm spt} \, T,$ for some $\alpha \in (0, 1)$ to be fixed eventually depending only on $n$, $m$ and $q$. 

With the help of a new, intrinsic frequency function (the planar frequency function) introduced in \cite{KrumWica}, as well as certain basic results from the initial parts of Almgren's work (namely, the associated ``linear theory'' i.e.\ the theory of Dirichlet energy minimizing multi-valued functions arising as blow-ups of area minimizing currents converging to a plane, and the strong Lipschitz approximation theorem for area minimizing currents that are weakly close to a plane), it is shown in \cite[Theorem~1.1]{KrumWica} that 
${\mathcal S}_{q}$ satisfies a certain locally uniform approximation-by-non-planar-cones property. This property says that for each point $Z_{0} \in {\mathcal S}_{q}$, there is a number $\rho_{Z_{0}} >0$ such that for each $Z \in {\mathcal S}_{q} \cap {\mathbf B}_{\rho_{Z_{0}}}(Z_{0})$ and each 
scale $\sigma \in (0, \rho_{Z_{0}}),$  there is a non-planar cone ${\mathbf C}_{Z,\sigma}$ (depending on $Z$ and $\sigma$) which approximates $T$ at scale $\sigma$ \emph{significantly better} than any plane does. More precisely, \cite[Theorem~1.1]{KrumWica} asserts that 
for any given $\b \in (0, 1/2)$, there is a number $\alpha = \alpha(n, m, q, \b) \in (0, 1)$ such that if we let ${\mathcal B}_{q} = {\mathcal B}_{q}(\b)$ be the set of points $Z \in {\rm spt} \, T$ for which the rescaled current $\eta_{Z, \rho \, \#} \, T$ converges, in the $L^{2}$ sense, to a (unique) multiplicity $q$ tangent plane ${\mathbf P}_{Z}$ at a rate $o(\rho^{\alpha})$, then for every $Z_{0} \in {\mathcal S}_{q} = {\mathcal S}_{q}(\b) \equiv {\rm sing}_{q} \, T \setminus {\mathcal B}_{q}(\b)$ and some $\rho_{Z_{0}} >0,$ the following holds true: for each $Z \in {\mathcal S}_{q} \cap {\mathbf B}_{\rho_{Z_{0}}}(Z_{0})$ and each $\sigma \in (0, \rho_{Z_{0}}),$ either we have 
\begin{eqnarray*}
&&\hspace{-.3in}\sigma^{-n-2} \int_{{\mathbf B}_{\sigma}(Z)} {\rm dist}^{2} \, (X, Z + {\rm spt} \, {\mathbf C}_{Z, \sigma}) \, d\|T\|(X)\\ 
&&\hspace{.1in}+ \, \sigma^{-n-2}\int_{{\mathbf B}_{\sigma/2}(Z) \setminus \{Y \, : \, {\rm dist} \, (Y, Z + S({\mathbf C}_{Z, \sigma})) < \sigma/16\}} {\rm dist}^{2}\,  (X, Z + {\rm spt} \, T) \, d\|{\mathbf C}_{Z, \sigma}\|(X)\\ 
&&\hspace{2.5in}\leq \beta \inf_{P} \, \sigma^{-n-2}\int_{{\mathbf B}_{\sigma}(Z)} {\rm dist}^{2} \, (X, Z + P) \, d\|T\| \hspace{.5in}(\star)
\end{eqnarray*} 
for some (not necessarily area minimizing and not necessarily unique) cone ${\mathbf C}_{Z, \sigma},$ with ${\rm dim} \, S({\mathbf C}_{Z, \sigma}) = n-2,$ made up of (at least two distinct) planes, or we have 
$$\{Y \, : \, \Theta(T, Y) \geq q\} \cap {\mathbf B}_{\sigma}(Z) \subset \{Y \, : \, {\rm dist} \, (Y, Z + L) \leq \b\sigma\} \hspace{1in} (\star\star)$$ 
for some linear subspace $L \subset {\mathbb R}^{n+m}$ of dimension $\leq n-3$. Here $S({\mathbf C}_{Z, \sigma})$ denotes the spine of ${\mathbf C}_{Z, \sigma}$, i.e.\ the maximal dimensional linear subspace along which ${\mathbf C}_{Z, \sigma}$ is translation invariant. 

Of course any non-branch-point singularity with density $q$ is contained in ${\mathcal S}_{q}(\b)$ (for any $\beta \in (0, 1)$). We emphasize that we allow ${\mathcal S}_{q}(\b)$ to contain branch points also, and therein lies one of the most basic differences between our approach and that of \cite{Almgren}. 

While it follows more or less immediately from the two conditions ($\star$) and ($\star\star$) that there is $\gamma = \gamma(n, m, q, \b)$ with $\gamma \to 0^{+}$ as $\b \to 0^{+}$ such that ${\mathcal H}^{n-2+\gamma} \, ({\mathcal S}_{q}(\b))  = 0$, it is not at all clear from these two conditions alone that the set of branch points in  
${\mathcal S}_{q}(\b)$ cannot have positive $(n-2)$-dimensional Hausdorff measure. A main purpose of the present paper is to show that for an appropriate choice of the parameter $\b \in (0, 1)$ that depends only on $n$, $m$ and $q,$ this is indeed the case, i.e.\ to show that there exists $\b = \b(n, m, q) \in (0, 1)$ such that branch points that may exist in 
${\mathcal S}_{q}(\b)$ must form an ${\mathcal H}^{n-2}$ null set. In fact our analysis gives much more. It establishes that there exists $\b = \b(n, m, q) \in (0, 1)$ such that for ${\mathcal H}^{n-2}$-a.e.\ point $Z \in {\mathcal S}_{q}(\b)$, the current $T$ has a unique tangent cone ${\mathbf C}_{Z}$ with ${\rm dim} \, S({\mathbf C}_{Z}) = n-2$, and moreover that ${\mathcal S}_{q}(\b)$ is locally $(n-2)$-rectifiable (with locally finite $(n-2)$-dimensional Hausdorff measure).  These results together with the definition of ${\mathcal B}_{q}(\b)$ then lead to the following theorem, which combines the two main results (Theorem~\ref{structure thm} and Corollary~\ref{unique-tangent-and-rectifiability}) of the present work:

\begin{theorem}[Theorem~\ref{structure thm} and Corollary~\ref{unique-tangent-and-rectifiability}]\label{main-results-combined}
Let $T$ be an $n$-dimensional locally area minimizing rectifiable current in an open set $U \subset {\mathbb R}^{n+m}$, and let ${\rm sing} \, T$ denote the singular set of $T$.  
\begin{itemize}
\item [(a)] for ${\mathcal H}^{n-2}$ a.e.\ point $Z \in {\rm spt} \, T$, the current $T$ has a unique tangent cone 
${\mathbf C}_{Z}$ of the form ${\mathbf C}_{Z} = \sum_{j=1}^{p} q_{j}\llbracket P_{j} \rrbracket$ where $p$, $q_{1}, \ldots, q_{p}$ are integers $\geq 1$, and 
$P_{1}, \ldots, P_{p}$ are distinct $n$-dimensional planes (all depending on $Z$);  either $p =1 $ (i.e.\ ${\mathbf C}_{Z}$ is supported on a single plane) or there is an $(n-2)$-dimensional subspace $L$ with 
$P_{i} \cap P_{j} = L$ for every $i \neq j$; 
\item[(b)] we have that ${\rm sing} \, T = {\mathcal B} \cup {\mathcal S}$ 
where: 
\begin{enumerate}
\item [{\rm (i)}] ${\mathcal B} \cap {\mathcal S} = \emptyset$; 
\item[{\rm (ii)}] for every compact set $K \subset U$, ${\mathcal S} \cap K$ is the union of a finite number of pairwise disjoint sets each of which is locally compact and locally $(n-2)$-rectifiable (and in particular has locally finite $(n-2)$-dimensional Hausdorff measure); and
\item[{\rm (iii)}] for every point $Z \in {\mathcal B}$, $T$ has a unique tangent cone at $Z$ supported on an $n$-dimensional plane $P_{Z}$ to which the scaled current about $Z$ decays rapidly in the following sense: for every compact set $K \subset U$  there are 
numbers $\alpha_{K} = \alpha(K, T) \in (0, 1)$ and $C_{K} = C(K, T) \in (0, \infty)$ such that  for every $Z \in {\mathcal B} \,\cap \, K,$ the estimate 
\begin{eqnarray*}
&&\rho^{-n-2}\int_{{\mathbf B}_{\rho}(Z)} {\rm dist}^{2}(X, Z + P_{Z}) \, d\|T\|(X)\\ 
&&\hspace{2in}\leq C_{K}\left(\frac{\rho}{\sigma}\right)^{2\alpha_{K}}\sigma^{-n-2}\int_{{\mathbf B}_{\sigma}(Z)} {\rm dist}^{2}(X, Z + P_{Z}) \, d\|T\|(X)
\end{eqnarray*}
holds for some $\sigma_{Z}$ (depending on $Z$) and all $\rho, \sigma$ with $0 < \rho \leq \sigma \leq \sigma_{Z}.$
\end{enumerate} 
\end{itemize}
\end{theorem}

Our proof of the above theorem builds on, among other things, techniques from our previous work \cite{Wic14}, \cite{KrumWic2} which in turn were inspired by the seminal work \cite{Sim93} of L.~Simon. A fundamental new difficulty that needs to be overcome in the present setting however is the lack of a regularity theorem, analogous to Allard's regularity theorem used in \cite{Sim93} or the (inductively used) sheeting theorem established and used in \cite{Wic14}, applicable to $T$. Explicit well-known examples  show that complete regularity (as in decomposition into smooth sheets  with curvature estimates, or even decomposition into locally Lipschitz graphs) is in fact false for area minimizers of 
codimension $\ge 2$ lying close to a plane. What is needed however is not such complete regularity, but to be able to argue that the current separates as a sum of disjoint (possibly still singular) pieces whenever its $L^{2}$ distance to a union of disjoint planes (i.e.\ fine-excess) is significantly smaller than its $L^{2}$ distance to any single plane (i.e.\ optimal coarse excess). This separation property indeed holds for area minimizers, and is a direct consequence of a new height bound  (Theorem~\ref{separation thm3-intro} below) we here establish. This estimate says that whenever ${\mathbf P}$ is a sum of planes with support 
consisting of planes that are disjoint in a cylinder, subject to appropriate small-excess and mass assumptions on $T$ and a measure of disjointness of the planes making up ${\mathbf P}$ (condition~\ref{separation3 hyp2-intro} below),  the pointwise distance of $T$ to ${\rm spt} \, {\mathbf P}$ in a smaller cylinder is bounded from above \emph{linearly} by a constant times the $L^{2}$-height excess of $T$ relative to ${\rm spt} \, {\mathbf P}$. (It is a well known, easy consequence of the monotonicity formula that such a bound holds in terms of a certain sublinear expression in height excess of $T$ relative to ${\mathbf P}$, but this weaker version is inadequate for our purposes). This estimate is analogous to the well-known interior upper bound on the supremum of a solution to a uniformly elliptic equation with bounded coefficients in terms of the $L^{2}$ norm of the solution. 

\begin{theorem}[Theorem~\ref{separation thm3}] \label{separation thm3-intro}
For all integers $q, s, p\geq 1$ with $p \leq s$ and all $\gamma \in (0,1)$ and $\kappa \in (0,\infty),$ there exists $\varepsilon_0 = \varepsilon_0(n,m,q,s,\gamma,\kappa) \in (0,1)$ such that if $\mathbf{P}  = \sum_{i=1}^{p} s_{i}\llbracket P_{i} \rrbracket$  for $n$-dimensional planes $P_{i}$ with orientation $\vec P_{i},$ and integers $s_{i} \geq 1$ with 
$\sum_{i=1}^{p}s_{i} = s,$ and if $T$ is an $n$-dimensional locally area-minimizing rectifiable current in the cylinder $\mathbf{C}_1(0) = B_{1}^{n}(0) \times {\mathbb R}^{m}$ such that, writing $P_{0} = {\mathbb R}^{n} \times \{0\}$, 
\begin{gather*} 
	(\partial T) \llcorner \mathbf{C}_1(0) = 0 , \quad \pi_{\#} T = q \llbracket B_1^{n}(0) \rrbracket , \\ 
	\frac{1}{2 \omega_n} \int_{\mathbf{C}_1(0)} |\vec T - \vec P_0|^2 \,d\|T\| < \varepsilon_0^2 , \quad 
	\max_{1 \leq i \leq p} |\vec P_i - \vec P_0| < \varepsilon_0, \nonumber 
\end{gather*}
and either $p = 1$ or $p > 1$ and 
\begin{equation} \label{separation3 hyp2-intro} 
	|\vec P_i - \vec P_j| \leq \kappa \inf_{X \in P_i \cap \mathbf{C}_1(0)} \op{dist}(X, P_j) 
\end{equation}
for each $1 \leq i < j \leq p$, then 
\begin{equation*} \label{separation3 concl} 
	\sup_{X \in \op{spt} T \cap \mathbf{C}_{\gamma}(0)} \op{dist}^2(X,\op{spt} {\mathbf P}) \leq C \int_{\mathbf{C}_1(0)} \op{dist}^2(X, \op{spt} \mathbf{P}) \,d\|T\|(X) 
\end{equation*}
for some constant $C = C(n,m,q,s,\kappa,\gamma) \in (0,\infty)$. 
\end{theorem}

 If $p=1$ (i.e.\ when ${\mathbf P}$ consists of a single plane), this is a well-known result due to Allard (\cite{Allard}) and in that case it holds for any stationary 
 intergral varifold in place of $T$. Our proof of the above theorem for $p \geq 2$ builds on the case $p=1$ and relies on a uniform interior $C^{0, \alpha}$ estimate for Dirichlet energy minimizing 
 $q$-valued functions due to \cite{Almgren}.  

The paper is organized as follows: Section~\ref{height-bound} is devoted to proving Theorem~\ref{separation thm3-intro}. In Section~\ref{estimates}, using Theorem~\ref{separation thm3-intro} among other things, various a priori estimates analogous to those in \cite{Sim93}, \cite{Wic14} are established for an area minimizing current $T.$ These estimates are valid whenever $T$ satisfies, for an appropriate choice of $\b = \b(n, m, q)$, condition ($\star$) with $\sigma = 1,$ $Z = 0$ and with ${\mathbf C}_{Z, \sigma}$ equal to a some cone ${\mathbf C}$ made up of a union of planes and having ${\rm dim} \, S({\mathbf C}) = n-2$, and additionally, whenever $T$ satisfies  a certain no-large-gaps hypothesis on the singular set as in \cite{Sim93}.  (The corresponding estimates in \cite{Sim93} are valid for stationary varifolds not necessarily satisfying an area minimizing condition  but required to belong to a multiplicity 1 class; additionally, in place of condition ($\star$), in \cite{Sim93} the varifold is assumed to have small excess relative to a cylindrical cone ${\mathbf C}$ that lies close to a fixed (singular) cylindrical cone ${\mathbf C}_{0}$ in the class, and the constants in the estimates are allowed to depend on ${\mathbf C}_{0}$.  In our setting, it is important that the constants depend only on $n$, $m$ and $q$ and in particular are independent of ${\mathbf C}$.) In Section~\ref{blow-up-analysis} the estimates in Section~\ref{estimates} in conjunction with adaptations of ideas from \cite{KrumWic2} are used to establish a decay estimate for the fine blow-ups of sequences $(T_{j})$ of area minimizing currents satisfying ($\star$) with $\b  = \b_{j} \to 0^{+}$ and with $\sigma = 1,$ $Z = 0$ and ${\mathbf C}_{Z, \sigma}$ equal to some cone ${\mathbf C}_{j}$ made up of a union of planes and having ${\rm dim} \, S({\mathbf C}_{j}) = n-2$, as well as satisfying the no-large-gaps condition with the gap size tending to zero. This blow-up analysis then leads to an excess decay estimate (Theorem~\ref{excess-improvement-final-1}) for an area minimizing current $T$ satisfying, among other things, condition ($\star$) for an appropriate fixed $\b = \b(n, m, q) \in (0, 1)$ and with $\sigma = 1,$ $Z = 0$ and ${\mathbf C}_{Z, \sigma}$ equal to a some cone ${\mathbf C}$ made up of a union of planes and having ${\rm dim} \, S({\mathbf C}) = n-2$. Finally, by combining \cite[Theorem~1.1]{KrumWica} with this decay estimate, in Section~\ref{rectifiability} we obtain ${\mathcal H}^{n-2}$ a.e.\ uniqueness of tangent cones to $T$ and local $(n-2)$-rectifiability of ${\mathcal S}_{q}(\b)$, arguing exactly as in \cite[Section~5]{Sim93}.  

\noindent
\section{The work of De~Lellis, Minter and Skorobogatova: a brief comparison}

In contemporaneous, independent work \cite{DelSko1}, \cite{DelSko2}, \cite{DelMinSko}
(posted on the arXiv soon after the first posting of the present paper and \cite{KrumWica}), 
De~Lellis, Minter and Skorobogatova  have also obtained two of the main results established in our series of papers, namely, uniqueness of tangent cones  at ${\mathcal H}^{n-2}$ a.e.\ point and countable $(n-2)$-rectifiability of the singular set for $n$-dimensional area minimizing rectifiable currents $T$. 
Their approach uses the full extent of Almgren's methods and results (\cite{Almgren}) as the starting point, and therefore conceptually differs from ours in a fundamental way. The differences in the two approaches are reflected both in our  considerably simpler argument for the ${\mathcal H}^{n-2}$ a.e.\ uniqueness-of-tangent-cones conclusion and in the additional results established by our methods concerning the structure of the singular set of $T$ and the behaviour of $T$ on approach to the branch set (specifically: certain local-finiteness-of-measure conclusions for the singular set,  local uniformity of the decay estimates, uniqueness of tangent functions (blow-ups) at ${\mathcal H}^{n-2}$ a.e.\ branch points and, a topological description of the current near certain branch points). 

We refer the reader to \cite[Section~2]{KrumWica} for a detailed discussion comparing the two methods. Here we just reproduce the part of that discussion that focuses on the present paper, which addresses the similarities and differences between the present article and what may be regarded as its counterpart \cite{DelMinSko}. 
 
Between the present article and \cite{DelMinSko} there is a close similarity, in terms of the techniques vis \`a vis a certain key excess-decay lemma established in either paper, as well as a key difference in the way this decay lemma is used to prove structure results for (part of) the singular set of $T$. Because of the differences in the way a key hypothesis of the decay lemma is verified, in the present work the lemma both applies to a larger set of singularities, and yields stronger conclusions for that set, than in \cite{DelMinSko}.

This excess-decay result in the present work is given as Lemma~\ref{excess-improvement-final} and Lemma~\ref{excess-improvement-final-1}, and in \cite{DelMinSko} it is \cite[Theorem~2.5]{DelMinSko}. It gives improvement of excess for $T$ whenever condition $(\star)$ in the introduction above holds (say at scale $\sigma = 1$) for some fixed $\beta = \beta(n, m, q) \in (0, 1)$, in addition to other hypotheses including a certain 
``no-large-gaps'' assumption on the set of density $\geq q$ singular points. The close similarity between the present work and \cite{DelMinSko} lies in the fact that in either case, the proof of this lemma uses a new pointwise bound for the distance of an area minimizing current to a union of non-intersecting planes 
(Theorem~\ref{separation thm3-intro} above, and \cite[Theorem~3.2]{DelMinSko}), and proceeds by adapting techniques developed in \cite{Sim93} and \cite[Sections~14]{Wic14}.

Using this decay lemma, in the present paper an asymptotic analysis of the behaviour of $T$ about points in the set ${\mathcal S}_{q} = {\mathcal S}_{q}(\beta)$
is carried out, for a fixed choice of $\beta \in (0, 1)$ depending only on $n$, $m$, $q$ (in fact $\beta$ taken to be as in the decay lemma, Lemma~\ref{excess-improvement-final}).  Recall (from the introduction above) that 
$${\mathcal S}_{q} = {\rm sing}_{q} \, T  \setminus {\mathcal B}_{q}$$
where ${\rm sing}_{q} \, T = \{Z \in {\rm sing} \, T \, : \, \Theta \, (\|T\|, Z) = q\}$ and ${\mathcal B}_{q} = {\mathcal B}_{q}(\beta) \;\; (= {\mathcal B}_{q}^{\prime} \cup {\mathcal B}_{q}^{\prime\prime})$ is the set of branch points $Z$ where decay of $T$ towards a multiplicity $q$ plane 
holds locally uniformly at a rate $|X - Z|^{1 + \alpha}$ or faster (i.e.\ points $Z$ with planar frequency ${\mathcal N}_{T, {\rm Pl}}(Z) \geq 1 + \alpha$), where $\alpha = \alpha(n, m, q, \beta) \in (0, 1)$ (and hence $\alpha = \alpha(n, m, q)$ after we choose $\beta = \beta(n, m, q)$ as in Lemma~\ref{excess-improvement-final}).  
The starting point of this analysis is \cite[Theorem~1.1(I)]{KrumWica}, which gives that about every point in $Z \in {\mathcal S}_{q}$ and at all small scales, either $T$ satisfies condition ($\star$) or condition ($\star\star$) (in the introduction above). Of key significance here is that the validity of these conditions is locally uniform 
in the sense that a scale $\rho_{Z}>0$ can be chosen so that ($\star$) or ($\star\star$)  holds about each point of ${\mathcal S}_{q} \cap {\mathbf B}_{\rho_{Z}}(Z)$  and at each scale $\rho \leq \rho_{Z}$. 

In \cite{DelMinSko} on the other hand, the decay lemma is used to analyse the (smaller) set 
$$\Sigma_{q} = \{Z \, : \, Z\;  \mbox{is a branch point with} \, \Theta(T, Z) = q \; \mbox{and} \;  I(Z) = 1\} \cup {\mathcal S}_{q, n-2},$$ 
where ${\mathcal S}_{q, n-2}$ is the set of density $q$ non branch point singularities of $T$ (where each tangent cone has spine dimension at most $(n-2)$), and $I(Z)$ is 
the ``singularity degree'' associated to the the branch point $Z$. The singularity degree is introduced in \cite{DelSko1}, and is defined based on Almgren's methods in \cite{Almgren} including the use of a sequence of center manifolds corresponding to $Z,$ and it is in fact the rate of decay of the current towards the center manifolds on approach to $Z$. The existence of this rate (as a finite real number) is the key analytic estimate in \cite{Almgren}, and it is a consequence of the approximate monotonicity of the Almgren frequency function associated to the normal approximation maps for $T$ relative to the sequence of center manifolds. 

We have 
\emph{a posteriori} that ${\mathcal H}^{n-2}\, ({\mathcal S}_{q} \setminus \Sigma_{q})= 0$, but \emph{a priori} not much can be said about the relationship between the two sets ${\mathcal S}_{q}$, $\Sigma_{q}$ beyond the fact that $\Sigma_{q} \subset {\mathcal S}_{q}$ which follows readily from the definitions of the two sets in view of the fact that $I(Z) \geq 1 + \alpha$ for any 
point $Z \in {\mathcal B}_{q}$. A key difference  in the way the excess decay lemma is used (which influences the conclusions, summarised in the next paragraph) is that in \cite{DelMinSko}, the validity of  
($\star$) or ($\star\star$) is checked pointwise, by arguing that for each point $Z \in \Sigma_{q},$ there is a scale $\rho_{Z}$ depending on $Z$ such that for each scale $\rho < \rho_{Z}$ either  ($\star$) or ($\star\star$) holds. This is a fairly direct consequence of either the definition of ${\mathcal S}_{q, n-2}$ (if $Z$ is not a branch point) or the fact that $I(Z) = 1$ (if $Z$ is a branch point), but  because of the pointwise choice of $\rho_{Z}$ (as opposed to a locally uniform choice), the set $\Sigma_{q}$ needs to be decomposed as a union of countably many subsets (in the obvious way by setting the $j$th subset equal to points $Z \in \Sigma_{q}$ for which $\rho_{Z} = j^{-1}$) before the decay lemma is applied (to each of these subsets). The other 
significant contrast to note is that in \cite{DelMinSko}, validity of ($\star$) or ($\star\star$) along $\Sigma_{q}$ ultimately relies on center manifold constructions (through the involvement of $I(Z)$ in the definition of $\Sigma_{q}$) whereas in the present work  
the validity of these conditions for ${\mathcal S}_{q}$ is guaranteed by the more elementary result \cite[Theorem~1.1(I)]{KrumWica} which is a fairly direct consequence of the approximate monotonicity of the planar frequency function.

The main conclusions in the present paper are: (i) ${\mathcal S}_{q}$ is countably $(n-2)$-rectifiable; (ii) for ${\mathcal H}^{n-2}$ a.e.\ point in ${\mathcal S}_{q}$ the current $T$ has a unique tangent cone supported on two or more distinct planes meeting along a common $(n-2)$-dimensional subspace, and (iii) there exists an ambient open set 
$V_{q}$ with ${\rm sing}_{q} \, T \subset V_{q}$ such that $V_{q} \cap \{Z \, : \, \Theta(T, Z) \geq q\} \setminus {\mathcal B}_{q}$ is locally compact and has locally finite $(n-2)$-dimensional Hausdorff measure (Theorem~\ref{structure thm});  in particular ${\mathcal S}_{q}$ has locally finite ${\mathcal H}^{n-2}$ measure (since  ${\mathcal S}_{q} \subset V_{q} \cap \{Z \, : \, \Theta(T, Z) \geq q\} \setminus {\mathcal B}_{q}$). The main conclusions in \cite{DelMinSko} are (i) and (ii) with $\Sigma_{q}$ in place of ${\mathcal S}_{q}$.  The additional conclusions (iii) in the present paper  are consequences  
of the fact that ${\mathcal B}_{q}$ is relatively closed in $V_{q}$ and the fact (emphasized above) that ($\star$) or ($\star\star$) holds in a locally uniform way about every point in ${\mathcal S}_{q}$.

 \section{Bounding the distance of an area minimizing current to a union of non-intersecting affine planes linearly in terms of its height excess relative to the same planes}\label{height-bound}
\setcounter{equation}{0}

\subsection{Notation and preliminaries}

Throughout the paper, we shall adopt the same notation as~\cite{KrumWica}.  See Section~3 of~\cite{KrumWica} for a discussion of general notation, as well as an overview of locally area-minimizing rectifiable currents and Dirichlet energy minimizing multi-valued functions.  

Let $n,m$ be integers $\geq 2$.  We shall express each point $X \in \mathbb{R}^{n+m}$ as $X = (x,y)$ where $x \in \mathbb{R}^n$ and $y \in \mathbb{R}^m$.  For each $X_0 \in \mathbb{R}^{n+m}$ and $\rho > 0$ we let 
\begin{equation*}
	\mathbf{B}_{\rho}(X_0) = \{ X \in \mathbb{R}^{n+m} : |X - X_0| < \rho \}. 
\end{equation*}
For each $x_0 \in \mathbb{R}^{n}$ and $\rho > 0$ we let 
\begin{align*}
	B_{\rho}(x_0) &= \{ x \in \mathbb{R}^n : |x - x_0| < \rho \} , \\
	\mathbf{C}_{\rho}(x_0) &= B_{\rho}(x_0) \times \mathbb{R}^m . 
\end{align*}

We shall often write $P_0 = \mathbb{R}^n \times \{0\}$.  Throughout we will let $\pi : \mathbb{R}^{n+m} \rightarrow P_0$ denote the orthogonal projection map onto $P_0$ and $\pi^{\perp} : \mathbb{R}^{n+m} \rightarrow P_0^{\perp}$ denote the orthogonal projection map onto the orthogonal complement $P_0^{\perp}$ of $P_0$.  We shall assume that $P_0$ is oriented by the $n$-vector $\vec P_0 = e_1 \wedge e_2 \wedge \cdots \wedge e_n$, where $e_1,e_2,\ldots,e_{n+m}$ is the standard basis for $\mathbb{R}^{n+m}$.

\begin{definition} \label{sum of tilted planes defn} {\rm
For positive integers $1 \leq p \leq s$, $\Pi_{s,p}$ denote the set of all $n$-dimensional rectifiable currents $\mathbf{P}$ which can be expressed as a sum of planes 
\begin{equation} \label{sum of tilted planes form}
	\mathbf{P} = \sum_{i=1}^p s_i \llbracket P_i \rrbracket , 
\end{equation}
where $p \geq 1$ and $s_i \geq 1$ are integers such that $\sum_{i=1}^p s_i = s$ and $P_i$ are distinct $n$-dimensional planes of $\mathbb{R}^{n+m}$ oriented by $n$-vectors $\vec P_i$.  
We let $\Pi_s = \bigcup_{p=1}^s \Pi_{s,p}$. 
} \end{definition}

\begin{definition} {\rm
Let $q \geq 1$ be an integer, $x_0 \in \mathbb{R}^n$, and $\rho > 0$.  Let $T$ be an $n$-dimensional locally area-minimizing rectifiable current in $\mathbf{C}_{\rho}(x_0)$ such that $(\partial T) \llcorner \mathbf{C}_{\rho}(x_0) = 0$ and let $\mathbf{P} \in \Pi_q$.  The height excess of $T$ relative to ${\mathbf P}$ in ${\mathbf C}_{\rho}(x_{0})$ is given by 
\begin{equation*} 
	E(T, {\mathbf P}, {\mathbf C}_{\rho}(x_{0})) = \left( \frac{1}{\omega_n \rho^{n+2}} \int_{\mathbf{C}_{\rho}(x_0)} \op{dist}^2(X, \op{spt} \mathbf{P}) \,d\|T\|(X) \right)^{\frac{1}{2}} . 
\end{equation*} 
} \end{definition}

\begin{definition} {\rm  
Let $x_0 \in \mathbb{R}^n$ and $\rho > 0$.  Let $T$ be an $n$-dimensional locally area-minimizing rectifiable current in $\mathbf{C}_{\rho}(x_0)$ such that $(\partial T) \llcorner \mathbf{C}_{\rho}(x_0) = 0$ and $\sup_{X \in \op{spt} T} \op{dist}(X,P_0) < \infty$.  The (oriented) tilt excess ${\mathcal E}(T,\mathbf{C}_{\rho}(x_0))$ (of $T$ relative to $\llbracket P_{0} \rrbracket$) is given by 
\begin{equation} \label{tilt excess defn eqn}
	{\mathcal E}(T,\mathbf{C}_{\rho}(x_0))^2 = \frac{\|T\|(\mathbf{C}_{\rho}(x_0))}{\omega_n \rho^n} 
			- \frac{\|\pi_{\#} T\|(\mathbf{C}_{\rho}(x_0))}{\omega_n \rho^n} .
\end{equation}
} \end{definition}

Notice that by the constancy theorem~\cite[Theorem~26.27]{SimonGMT} 
\begin{equation} \label{projection concl} 
	\pi_{\#} T = q \llbracket B_{\rho}(x_0) \rrbracket 
\end{equation}
for some integer constant $q$.  Here $q$ can be zero, positive in the case that $\pi_{\#} T$ is oriented by $\vec P_0$, or negative in the case that $\pi_{\#} T$ is oriented by $-\vec P_0$.  Often we will assume that $q \geq 0$ since otherwise we can replace $T$ with $-T$.  Assuming $q \geq 0$, by~\cite[5.3.1]{Fed69} 
\begin{equation*} 
	{\mathcal E}(T,\mathbf{C}_{\rho}(x_0))^2 = \frac{1}{2\omega_n \rho^n} \int_{\mathbf{C}_{\rho}(x_0)} |\vec T - \vec P_0|^2 \,d\|T\| , 
\end{equation*}
where $\vec T$ is the orienting $n$-vector of $T$.  Moreover, by \eqref{tilt excess defn eqn} and \eqref{projection concl}, 
\begin{equation*} 
	{\mathcal E}(T,\mathbf{C}_{\rho}(x_0))^2 = \frac{\|T\|(\mathbf{C}_{\rho}(x_0))}{\omega_n \rho^n} - q  
\end{equation*}
or equivalently,  
\begin{equation} \label{projection eqn2}
	\|T\|(\mathbf{C}_{\rho}(x_0)) = (q + {\mathcal E}(T,\mathbf{C}_{\rho}(x_0))^2) \,\omega_n \rho^n .
\end{equation}
As a straightforward consequence of the monotonicity formula, we can often show that $q \neq 0$. 

\begin{lemma} \label{projection lemma}
Let $\gamma \in (0,1)$, $x_0 \in \mathbb{R}^n$, and $\rho > 0$.  Let $T$ be an $n$-dimensional locally area-minimizing rectifiable current of $\mathbf{C}_{\rho}(x_0)$ such that 
\begin{equation*} 
	(\partial T) \llcorner \mathbf{C}_{\rho}(x_0) = 0 , \quad\quad \sup_{X \in \op{spt} T} \op{dist}(X,P_0) < \infty, \quad\quad
	{\mathcal E}(T,\mathbf{C}_{\rho}(x_0)) < (1-\gamma)^n . 
\end{equation*}
Then either $T \llcorner \mathbf{C}_{\gamma \rho}(x_0) = 0$ or there exists a nonzero integer $q$ such that \eqref{projection concl} holds true. 
\end{lemma}

\begin{proof}
We know that \eqref{projection concl} holds true for some integer $q$.  Suppose that \eqref{projection concl} holds true with $q = 0$.  Then by \eqref{projection eqn2} and ${\mathcal E}(T,\mathbf{C}_{\rho}(x_0)) < (1-\gamma)^n$, $\|T\|(\mathbf{C}_{\rho}(x_0)) < \omega_n (1-\gamma)^n \rho^n$.  On the other hand, if there were a point $X \in \op{spt} T \cap \mathbf{C}_{\gamma\rho}(x_0)$, then by the monotonicity formula for area, $\|T\|(\mathbf{C}_{\rho}(x_0)) \geq \|T\|(\mathbf{B}_{(1-\gamma)\rho}(X)) \geq \omega_n (1-\gamma)^n \rho^n$, giving us a contradiction.  Therefore, $T \llcorner \mathbf{C}_{\gamma\rho}(x_0) = 0$.
\end{proof}

\begin{lemma} \label{projection lemma2}
Let $x_0 \in \mathbb{R}^n$ and $\rho > 0$.  Let $T$ be an $n$-dimensional locally area-minimizing rectifiable current of $\mathbf{C}_{\rho}(x_0)$ such that for some integer $q \geq 0$
\begin{equation*}
	(\partial T) \llcorner \mathbf{C}_{\rho}(x_0) = 0, \quad \sup_{X \in \op{spt} T} \op{dist}(X,P_0) < \infty, \quad 
	\pi_{\#} T = q \llbracket B_{\rho}(x_0) \rrbracket, \quad {\mathcal E}(T,\mathbf{C}_{\rho}(x_0)) < 1.
\end{equation*}
Suppose that for $i \in \{1,2,\ldots,N\}$ there are $n$-dimensional locally area-minimizing rectifiable currents $T_i$ of $\mathbf{C}_{\rho}(x_0)$ such that 
\begin{equation*}
	(\partial T_i) \llcorner \mathbf{C}_{\rho}(x_0) = 0 \text{ for all $i$,} \quad\quad T = \sum_{i=1}^N T_i, \quad\quad \|T\| = \sum_{i=1}^N \|T_i\| .
\end{equation*}
Then for each $i \in \{1,2,\ldots,N\}$ there exist an integer $q_i \geq 0$ such that $\pi_{\#} T_i = q_i \llbracket B_{\rho}(x_0) \rrbracket$ and $\sum_{i=1}^N q_i = q$.
\end{lemma}

\begin{proof}
Let $i \in \{1,2,\ldots,N\}$.  By the constancy theorem, there exist an integer $q_i$ such that $\pi_{\#} T_i = q_i \llbracket B_{\rho}(x_0) \rrbracket$.  We have that $q_i \geq 0$ since  
\begin{equation*}
	-q_i \leq \frac{\|T_i\|(\mathbf{C}_{\rho}(x_0))}{\omega_n \rho^n} - q_i 
	= \frac{1}{2\omega_n \rho^n} \int_{\mathbf{C}_{\rho}(x_0)} |\vec T - \vec P_0|^2 \,d\|T_i\| 
	\leq \frac{1}{2\omega_n \rho^n} \int_{\mathbf{C}_1(0)} |\vec T - \vec P_0|^2 \,d\|T\| < 1,
\end{equation*}
where $\vec T$ is the orientation $n$-vector of $T$ and $\vec T |_{\op{spt} T_i}$ orients $T_i$.  Clearly $\sum_{i=1}^N q_i = q$.
\end{proof}

The following elementary ``coarse'' upper bound for distance will be used in the proof of our main distance estimates. 

\begin{lemma} \label{sepmono lemma}
Let $q \geq 1$ be an integer and $\gamma \in (0,1)$.  If $\mathbf{P} \in \Pi_q$ and $T$ is an $n$-dimensional locally area-minimizing rectifiable current in $\mathbf{C}_1(0)$ with $(\partial T) \llcorner \mathbf{C}_1(0) = 0$ such that 
\begin{equation*} 
	\frac{1}{\omega_n} \int_{\mathbf{C}_1(0)} \op{dist}^2(X, \op{spt} \mathbf{P}) \,d\|T\|(X) < \left(\frac{1-\gamma}{2}\right)^{n+2} , 
\end{equation*}
then 
\begin{equation*} 
	\sup_{X \in \op{spt} T \cap \mathbf{C}_{\gamma}(0)} \op{dist}(X, \op{spt} \mathbf{P}) 
		\leq 2 \left( \frac{1}{\omega_n} \int_{\mathbf{C}_1(0)} \op{dist}^2(X, \op{spt} \mathbf{P}) \,d\|T\|(X) \right)^{\frac{1}{n+2}} . 
\end{equation*}
\end{lemma}

\begin{proof} See \cite[Lemma~3.6]{KrumWica} (with $K = \op{spt}\mathbf{P}$). 
\end{proof}

The following Lipschitz and harmonic approximation results due to Almgren (\cite{Almgren}) will play an important role in our distance estimates. 

\begin{theorem}[Almgren's strong approximation theorem] \label{lip approx thm}
For each $\gamma \in (0,1)$ there exists $\varepsilon = \varepsilon(n,m,q,\gamma) > 0$ such that the following holds true.  Let $x_0 \in \mathbb{R}^n$ and $\rho > 0$.  Let $T$ be an $n$-dimensional locally area-minimizing rectifiable current of $\mathbf{C}_{\rho}(x_0)$ such that such that 
\begin{equation*}
	(\partial T) \llcorner \mathbf{C}_{\rho}(x_0) = 0, \quad \sup_{X \in \op{spt}} \op{dist}(X,P_0) < \infty , \quad
	\pi_{\#} T = q \llbracket B_{\rho}(x_0) \rrbracket , \quad {\mathcal E} = {\mathcal E}(T,\mathbf{C}_{\rho}(x_0)) < \varepsilon , 
\end{equation*}
where $\pi$ is the orthogonal projection map onto $P_0$.  Then there exists a Lipschitz $q$-valued function $u : B_{\gamma\rho}(x_0) \rightarrow \mathcal{A}_q(\mathbb{R}^m)$ and a closed set $K \subset B_{\rho}(x_0)$ such that 
\begin{gather} 
	\label{lip approx concl} \op{Lip} u \leq C {\mathcal E}^{\alpha}, \quad T \llcorner (K \times \mathbb{R}^m) = (\op{graph} u) \llcorner (K \times \mathbb{R}^m), \\
	\mathcal{L}^n(B_{\gamma\rho}(x_0) \setminus K) + \|T\|((B_{\gamma\rho}(x_0) \setminus K) \times \mathbb{R}^m) \leq C {\mathcal E}^{2+\alpha} \rho^n, 
		\nonumber 
\end{gather} 
and 
\begin{equation} \label{lip approx tilt} 
	\left| \,\omega_n (\sigma \rho)^n \,{\mathcal E(}T,\mathbf{C}_{\sigma \rho}(x_0)) - \frac{1}{2} \int_{B_{\sigma \rho}(x_0)} |Du|^2 \right| 
	\leq C {\mathcal E}^{2+\alpha} \rho^n 
\end{equation}
for all $0 < \sigma \leq \gamma$, where $C = C(n,m,q,\gamma) \in (0,\infty)$ and $\alpha = \alpha(n,m,q) \in (0,1)$ are constants.  
\end{theorem}

\begin{proof}
See \cite[Corollary~3.29]{Almgren} or \cite[Theorem~2.4]{DeLSpa1}.
\end{proof}

\begin{theorem}[Harmonic approximation theorem] \label{harmonic thm}
For every $\eta > 0$ and $\gamma \in (0,1)$ there exists $\varepsilon = \varepsilon(n,m,q,\gamma,\eta) > 0$ such that the following holds true.  Let $x_0 \in \mathbb{R}^n$ and $\rho > 0$.  Let $T$ be an $n$-dimensional locally area-minimizing rectifiable current of $\mathbf{C}_{\rho}(x_0)$ such that 
\begin{equation*}
	(\partial T) \llcorner \mathbf{C}_{\rho}(x_0) = 0, \quad \sup_{X \in \op{spt}} \op{dist}(X,P_0) < \infty , \quad
	\pi_{\#} T = q \llbracket B_{\rho}(x_0) \rrbracket , \quad {\mathcal E} = {\mathcal E}(T,\mathbf{C}_{\rho}(x_0)) < \varepsilon , 
\end{equation*}
where $\pi$ is the orthogonal projection map onto $P_0$.  Let $u : B_{\gamma\rho}(x_0) \rightarrow \mathcal{A}_q(\mathbb{R}^m)$ be the Lipschitz approximation of $T$ as in Theorem~\ref{lip approx thm}.  Then there exists a Dirichlet energy minimizing $q$-valued function $w : B_{\gamma\rho}(x_0) \rightarrow \mathcal{A}_q(\mathbb{R}^m)$ such that 
\begin{equation} \label{harmonic concl}
	\frac{1}{\rho^{n+2}} \int_{B_{\gamma\rho}(x_0)} \mathcal{G}(u,w)^2 + \frac{1}{\rho^n} \int_{B_{\gamma\rho}(x_0)} (|Du| - |Dw|)^2 \leq \eta \,{\mathcal E}^2 . 
\end{equation} 
\end{theorem}

\begin{proof}
This is \cite[Theorem~3.33]{Almgren} or ~\cite[Theorem~2.6]{DeLSpa1} with obvious modifications.
\end{proof}

We will also need the following estimate that bounds (oriented) tilt excess from above by the height excess. 

\begin{lemma} \label{energy est lemma}
Let $\gamma \in (0,1)$, $x_0 \in \mathbb{R}^n$, and $\rho > 0$.  If $T$ is an $n$-dimensional locally area-minimizing rectifiable current of $\mathbf{C}_{\rho}(x_0)$ such that 
\begin{equation*}
	(\partial T) \llcorner \mathbf{C}_{\rho}(x_0) = 0, \quad\quad 
	\sup_{X \in \op{spt} T \cap \mathbf{C}_{\rho}(x_0)} \op{dist}(X,P_0) \leq \rho , 
\end{equation*}
then 
\begin{equation*} 
	{\mathcal E}(T,\mathbf{C}_{\gamma\rho}(x_0))
		\leq C \left( \frac{1}{\rho^{n+2}} \int_{\mathbf{C}_{\rho}(x_0)} \op{dist}^2(X,P_0) \,d\|T\|(X) \right)^{1/2}
\end{equation*} 
for some constant $C = C(n,m,q,\gamma) \in (0,\infty)$. 
\end{lemma}

\begin{proof} This is \cite[Lemma~3.7]{KrumWica} which is essentially the same as \cite[Lemma~3.2]{HardtSimon}.  
\end{proof}

\subsection{Statements of the main estimates}  
The conclusion of Lemma~\ref{sepmono lemma} is too coarse to use in practice.  More useful would be a bound, for the $L^{\infty}$-distance of $T$ to a sum of planes $\mathbf{P},$ that is linear in the height excess of $T$ with respect to $\mathbf{P},$ or linear in the tilt excess of $T$ relative to $P_0$.  We establish both of these in our main results Theorem~\ref{separation thm1} and Theorem~\ref{separation thm3}. 

\begin{theorem} \label{separation thm1} 
For each integer $q \geq 1$ and $\gamma \in (0,1)$ there exists $\varepsilon_0 = \varepsilon_0(n,m,q,\gamma) \in (0,1)$ such that if $T$ is an $n$-dimensional locally area-minimizing rectifiable current in $\mathbf{C}_1(0)$ such that 
\begin{equation} \label{separation1 hyp} 
	(\partial T) \llcorner \mathbf{C}_1(0) = 0 , \quad \sup_{X \in \op{spt} T} \op{dist}(X,P_0) < \infty, \quad
	\pi_{\#} T = q \llbracket B_1(0) \rrbracket , \quad {\mathcal E}(T,\mathbf{C}_1(0)) < \varepsilon_0 ,  
\end{equation}
then there exists $\mathbf{P} \in \Pi_q$ such that $\mathbf{P} = \sum_{i=1}^q \llbracket \mathbb{R}^n \times \{a_i\} \rrbracket$, where $a_i \in \mathbb{R}^m$ (possibly repeating) and $\mathbb{R}^n \times \{a_i\}$ is oriented by $\vec e_1 \wedge \vec e_2 \wedge\cdots\wedge \vec e_n$, and 
\begin{equation} \label{separation1 concl}
	\sup_{X \in \op{spt} T \cap \mathbf{C}_{\gamma}(0)} \op{dist}(X, \op{spt} \mathbf{P}) \leq C {\mathcal E}(T,\mathbf{C}_1(0)) , 
\end{equation}
where $C = C(n,m,q,\gamma) \in (0,\infty)$ is a constant. 
\end{theorem}

\begin{theorem} \label{separation thm3}
For all integers $q \geq 1$ and $1 \leq p \leq s$ and all $\gamma \in (0,1)$ and $\kappa \in (0,\infty),$ there exists $\varepsilon_0 = \varepsilon_0(n,m,q,s,\gamma,\kappa) \in (0,1)$ such that if $\mathbf{P} \in \Pi_{s,p}$ (as in \eqref{sum of tilted planes form}) and $T$ is an $n$-dimensional locally area-minimizing rectifiable current in $\mathbf{C}_1(0)$ such that 
\begin{gather} \label{separation3 hyp1} 
	(\partial T) \llcorner \mathbf{C}_1(0) = 0 , \quad \sup_{X \in \op{spt} T} \op{dist}(X,P_0) < \infty, \quad \pi_{\#} T = q \llbracket B_1(0) \rrbracket , \\ 
	{\mathcal E}(T,\mathbf{C}_1(0))^2 = \frac{1}{2 \omega_n} \int_{\mathbf{C}_1(0)} |\vec T - \vec P_0|^2 \,d\|T\| < \varepsilon_0^2 , \quad 
	\max_{1 \leq i \leq p} |\vec P_i - \vec P_0| < \varepsilon_0, \nonumber 
\end{gather}
and either $p = 1$ or $p > 1$ and 
\begin{equation} \label{separation3 hyp2} 
	|\vec P_i - \vec P_j| \leq \kappa \inf_{X \in P_i \cap \mathbf{C}_1(0)} \op{dist}(X, P_j) 
\end{equation}
for each $i,j \in \{1,2,\ldots,p\}$ with $i \neq j$, then 
\begin{equation} \label{separation3 concl} 
	\sup_{X \in \op{spt} T \cap \mathbf{C}_{\gamma}(0)} \op{dist}^2(X,\op{spt} \, {\mathbf P}) \leq C \int_{\mathbf{C}_1(0)} \op{dist}^2(X, \op{spt} \mathbf{P}) \,d\|T\|(X) 
\end{equation}
for some constant $C = C(n,m,q,s,\kappa,\gamma) \in (0,\infty)$. 
\end{theorem}

\subsection{$L^{\infty}$-distance estimate in terms of tilt excess relative to a plane} 
In this section we prove Theorem~\ref{separation thm1}.  
The proof is based on a controlled growth estimate for the tilt excess (Lemma~\ref{excess decay lemma}), which follows from a growth estimate for the Dirichlet energy of a Dirichlet energy minimizing multi-valued function, established in~\cite[Theorem~2.13]{Almgren}. 
We combine this with a Poincar\'e type inequality (Lemma~\ref{poincare lemma}), and iteratively apply both results 
as long as the tilt excess ${\mathcal E}(T,\mathbf{C}_{\rho}(\xi))$ remains small.  
At scales where ${\mathcal E}(T,\mathbf{C}_{\rho}(\xi))$ is no longer small, we use the coarse bound given by Lemma~\ref{sepmono lemma}.

Throughout this section and Section~\ref{sec:separation2_sec}, we consider the following class $\mathcal{P}_q$ of sums of oriented planes which are parallel to $P_0 = \mathbb{R}^n \times \{0\}$.  (C.f.~Definition~\ref{sum of tilted planes defn} and note that $\mathcal{P}_q \subset \Pi_q$.)

\begin{definition} \label{sum of planes defn} {\rm
For integers $p$, $q$ with $1 \leq p \leq q$, let $\mathcal{P}_{q,p}$ denote the set of all $n$-dimensional rectifiable currents $\mathbf{P}$ of $\mathbb{R}^{n+m}$ which can be expressed as a sum of parallel planes in the form 
\begin{equation} \label{sum of planes form}
	\mathbf{P} = \sum_{i=1}^p q_i \llbracket P_i \rrbracket , 
\end{equation}
where $q_i \geq 1$ are integers such that $\sum_{i=1}^p q_i = q$, $P_i = \mathbb{R}^n \times \{a_i\}$ for distinct $a_i \in \mathbb{R}^m$ and $P_i$ is oriented by the $n$-vector $\vec P_0 = e_1 \wedge e_2 \wedge \cdots \wedge e_n$.  Let $\mathcal{P}_q = \bigcup_{p=1}^q \mathcal{P}_{q,p}$. 
} \end{definition}

\begin{definition} \label{associated plane defn} {\rm 
We associate each $a = \sum_{i=1}^q \llbracket a_i \rrbracket \in \mathcal{A}_q(\mathbb{R}^m)$ (where $a_i \in \mathbb{R}^m,$ possibly repeating) with $\mathbf{P}_a = \sum_{i=1}^q \llbracket \mathbb{R}^n \times \{a_i\} \rrbracket \in \mathcal{P}_q$. 
} \end{definition}

\begin{definition} \label{size and width defn} {\rm
For $1 < p \leq q$ and $\mathbf{P} \in \mathcal{P}_{q,p}$ as in \eqref{sum of planes form}, we define 
\begin{equation*}
	\op{sep}(\mathbf{P}) = \min_{i \neq j} |a_i - a_j| \quad\quad \op{width}(\mathbf{P}) = \max_{i \neq j} |a_i - a_j| . 
\end{equation*}
If $\mathbf{P} \in \mathcal{P}_{q,1}$, we define $\op{sep}(\mathbf{P}) = \infty$ and $\op{width}(\mathbf{P}) = 0$. 
} \end{definition}

\begin{lemma}[Controlled growth of tilt excess] \label{excess decay lemma}
For every $\theta \in (0,1/8)$ there exists $\varepsilon_1 = \varepsilon_1(n,m,q,\theta) \in (0,1)$ such that the following holds true.  Let $x_0 \in \mathbb{R}^n$ and $\rho > 0$.  Let $T$ be an $n$-dimensional locally area-minimizing rectifiable current of $\mathbf{C}_{\rho}(x_0)$ such that 
\begin{align} \label{excess decay hyp}
	&(\partial T) \llcorner \mathbf{C}_{\rho}(x_0) = 0, \quad\quad \sup_{X \in \op{spt} T} \op{dist}(X,P_0) < \infty, \\ 
	&\pi_{\#} T = q \llbracket B_{\rho}(x_0) \rrbracket , \quad\quad {\mathcal E} = {\mathcal E}(T,\mathbf{C}_{\rho}(x_0)) < \varepsilon_1 . \nonumber 
\end{align}
Then 
\begin{equation} \label{excess decay concl}
	{\mathcal E}(T,\mathbf{C}_{\theta \rho}(x_0)) \leq C \theta^{\mu-1} \,{\mathcal E}(T,\mathbf{C}_{\rho}(x_0))
\end{equation}
for some constants $C = C(n,m,q) \in (0,\infty)$ and $\mu = \mu(n,m,q) \in (0,1)$ (independent of $\theta$). 
\end{lemma}

\begin{proof}
Without loss of generality assume that $x_0 = 0$ and $\rho = 1$.  Let $\eta = \eta(n,m,q,\theta) \in (0,1)$ and $\varepsilon_1 = \varepsilon_1(n,m,q,\theta,\eta) > 0$ to be later determined.  Let $T$ be an $n$-dimensional locally area-minimizing rectifiable current of $\mathbf{C}_1(0)$ such that \eqref{excess decay hyp} holds true.  Set ${\mathcal E} = {\mathcal E}(T,\mathbf{C}_1(0))$.  Assuming $\varepsilon_1$ is sufficiently small, 
$T$ has a Lipschitz approximation $u : B_{1/4}(0) \rightarrow \mathcal{A}_q(\mathbb{R}^m)$ as in Theorem~\ref{lip approx thm} with $\gamma = 1/4$ and a harmonic approximation $w : B_{1/4}(0) \rightarrow \mathcal{A}_q(\mathbb{R}^m)$ as in Theorems~\ref{harmonic thm} with $\gamma = 1/4$.  By \cite[Theorem~2.13]{Almgren},
there exists $\mu = \mu(n,m,q) \in (0,1)$ such that 
\begin{equation*} 
	\theta^{2-n} \int_{B_{\theta}(0)} |Dw|^2 \leq 4^{n-2+2\mu} \theta^{2\mu} \int_{B_{1/4}(0)} |Dw|^2 . 
\end{equation*}
Thus by \eqref{harmonic concl}, 
\begin{align*}
	\theta^{-n} \int_{B_{\theta}(0)} |Du|^2 
	&\leq 2 \theta^{-n} \int_{B_{\theta}(0)} |Dw|^2 + 2 \eta \theta^{-n} {\mathcal E}^2 
	\\&\leq 2 \cdot 4^{n-2+2\mu} \theta^{2\mu-2} \int_{B_{1/4}(0)} |Dw|^2 + 2 \eta \theta^{-n} {\mathcal E}^2 
	\\&\leq 4^{n-1+2\mu} \theta^{2\mu-2} \int_{B_{1/4}(0)} |Du|^2 + 4^{n-1+2\mu} \eta \theta^{2\mu-2} {\mathcal E}^2 + 2 \eta \theta^{-n} {\mathcal E}^2 . 
\end{align*}
Choose $\eta = \eta(n,m,q,\theta) \in (0,1)$ so that $(4^{n-1+2\mu} + \theta^{2-n-2\mu}) \eta < \omega_n$.  Then provided $\varepsilon_1 = \varepsilon_1(n,m,q,\theta)$ is sufficiently small
\begin{equation*} 
	\theta^{-n} \int_{B_{\theta}(0)} |Du|^2 \leq 4^{n-1+2\mu} \theta^{2\mu-2} \int_{B_{1/8}(0)} |Du|^2 + \theta^{2\mu-2} \omega_n {\mathcal E}^2. 
\end{equation*}
By \eqref{lip approx tilt}, 
\begin{align*}
	\theta^{-n} \int_{\mathbf{C}_{\theta}(0)} |\vec T - \vec P_0|^2 \,d\|T\| 
	&\leq \theta^{-n} \int_{B_{\theta}(0)} |Du|^2 + C \theta^{-n} {\mathcal E}^{2+\alpha}
	\\&\leq C\theta^{2\mu-2} \int_{B_{1/8}(0)} |Du|^2 + \theta^{2\mu-2} {\mathcal E}^2 + C \theta^{-n} {\mathcal E}^{2+\alpha} 
	\\&\leq C\theta^{2\mu-2} {\mathcal E}^2 + C \theta^{2\mu-2} {\mathcal E}^{2+\alpha} + \theta^{2\mu-2} {\mathcal E}^2 + C \theta^{-n} {\mathcal E}^{2+\alpha} 
	\\&\leq (3C+1) \,\theta^{2\mu-2} {\mathcal E}^2 ,
\end{align*}
where $C = C(n,m,q) \in (0,\infty)$ and $\alpha = \alpha(n,m,q) \in (0,1)$ are constants and we assumed $\varepsilon_1 = \varepsilon_1(n,m,q,\theta)$ is small enough that ${\mathcal E}^{\alpha} < \varepsilon_1^{\alpha} < 1$ and $\theta^{2-n-2\mu} {\mathcal E}^{\alpha} < \theta^{2-n-2\mu} \varepsilon_1^{\alpha} < 1$. 
\end{proof}

\begin{lemma}[Poincar\'e-type inequality] \label{poincare lemma}
For each $\gamma \in (0,1)$ there exists $\varepsilon_2 = \varepsilon_2(n,m,q,\gamma) > 0$ such that the following holds true.  Let $x_0 \in \mathbb{R}^n$ and $\rho > 0$.  Suppose that $T$ is an $n$-dimensional locally area-minimizing rectifiable current in $\mathbf{C}_{\rho}(x_0)$ such that 
\begin{align} \label{poincare hyp}
	&(\partial T) \llcorner \mathbf{C}_{\rho}(x_0) = 0, \quad\quad \sup_{X \in \op{spt} T} \op{dist}(X,P_0) < \infty, \\ 
	&\pi_{\#} T = q \llbracket B_{\rho}(x_0) \rrbracket , \quad\quad {\mathcal E} = {\mathcal E}(T,\mathbf{C}_{\rho}(x_0)) < \varepsilon_2 . \nonumber 
\end{align}
Let $u : B_{(1+\gamma)\rho/2}(x_0) \rightarrow \mathcal{A}_q(\mathbb{R}^m)$ and $K \subset B_{(1+\gamma)\rho/2}(x_0)$ be as in Theorem~\ref{lip approx thm} with $(1+\gamma)/2$ in place of $\gamma$.  Then there exists $a \in \mathcal{A}_q(\mathbb{R}^m)$  such that 
\begin{align} \label{poincare concl}
	\rho^{-n-2} \int_{\mathbf{C}_{\gamma\rho}(x_0)} \op{dist}^2(X,\op{spt} \mathbf{P}_a) \,d\|T\|(X) 
		+ \rho^{-n-2} \int_{B_{\gamma\rho}(x_0)} \mathcal{G}(u,a)^2 \leq C {\mathcal E}(T,\mathbf{C}_{\rho}(x_0))^2 
\end{align}
for some constant $C = C(n,m,q,\gamma) \in (0,\infty)$, where ${\mathbf P}_{a} \in {\mathcal P}_{q}$ is the sum of planes associated with $a$ as in Definition~\ref{associated plane defn}.  
\end{lemma}

\begin{proof}
Without loss of generality assume $x_0 = 0$ and $\rho = 1$.  Let $\varepsilon_2 = \varepsilon_2(n,m,q,\gamma) \in (0,1)$ be a small constant to be later determined.  Let $T$ be an $n$-dimensional locally area-minimizing rectifiable current of $\mathbf{C}_1(0)$ such that \eqref{poincare hyp} holds true.  Set ${\mathcal E} = {\mathcal E}(T,\mathbf{C}_1(0))$.  Let $u : B_{(1+\gamma)/2}(0) \rightarrow \mathcal{A}_q(\mathbb{R}^m)$ and $K \subset B_{(1+\gamma)/2}(0)$ be as in Theorem~\ref{lip approx thm} with $(1+\gamma)/2$ in place of $\gamma$.  By the Poincar\'e inequality~\cite[Proposition~2.12]{DeLSpaDirMin}, there exists $a \in \mathcal{A}_q(\mathbb{R}^m)$ such that 
\begin{equation} \label{poincare eqn1}
	\int_{B_{(1+\gamma)/2}(0)} \mathcal{G}(u,a)^2 \leq C \int_{B_{(1+\gamma)/2}(0)} |Du|^2 
\end{equation}
where $C = C(n,m,q,\gamma) \in (0,\infty)$ is a constant.  Let $\mathbf{P}_a \in \mathcal{P}_q$ be the sums of planes associated with $a$ as in Definition~\ref{associated plane defn}.  By \eqref{lip approx tilt}, 
\begin{equation} \label{poincare eqn2}
	\int_{B_{(1+\gamma)/2}(0)} |Du|^2 \leq C(n) \,{\mathcal E}^2 
\end{equation}
provided $\varepsilon_2$ is sufficiently small.  By \eqref{lip approx concl} and the area formula (see (3.29) and (3.28) of~\cite{KrumWica})
\begin{equation} \label{poincare eqn3} 
	\int_{K \times \mathbb{R}^m} \op{dist}^2(X,\op{spt} \mathbf{P}_{a}) \,d\|T\|(X) \leq 2 \int_{B_{(1+\gamma)/2}(0)} \mathcal{G}(u,a)^2 
\end{equation}
provided $\varepsilon_2$ is sufficiently small.  By \eqref{poincare eqn3}, \eqref{poincare eqn1}, and \eqref{poincare eqn2}, 
\begin{gather} 
	\label{poincare eqn4} \int_{K \times \mathbb{R}^m} \op{dist}^2(X,\op{spt} \mathbf{P}_{a}) \,d\|T\|(X) \leq C {\mathcal E}^2, \\
	\label{poincare eqn5} \int_{B_{(1+\gamma)/2}(0)} \mathcal{G}(u,a)^2 \leq C {\mathcal E}^2 
\end{gather}
for some constant $C = C(n,m,q,\gamma) \in (0,\infty)$.

Suppose that 
\begin{equation*}
	\op{dist}(Z,\op{spt} \mathbf{P}_{a}) \geq 1-\gamma
\end{equation*}
for some $Z \in \op{spt} T \cap \mathbf{C}_{\gamma}(0)$.  Then by the monotonicity formula and the fact that $\|T\|((B_{(1+\gamma)/2}(0) \setminus K) \times \mathbb{R}^m) \leq C {\mathcal E}^{2+\alpha}$ (as in \eqref{lip approx concl}), 
\begin{align*}
	\int_{K \times \mathbb{R}^m} \op{dist}^2(X,\op{spt} \mathbf{P}_{a}) \,d\|T\|(X) 
		&\geq \int_{\mathbf{B}_{(1-\gamma)/2}(Z) \cap (K \times \mathbb{R}^m)} \op{dist}^2(X,\op{spt} \mathbf{P}_{a}) \,d\|T\|(X) 
		\\&\geq \bigg(\frac{1-\gamma}{2}\bigg)^2 \|T\|(\mathbf{B}_{(1-\gamma)/2}(Z) \cap (K \times \mathbb{R}^m))
		\\&\geq \bigg(\frac{1-\gamma}{2}\bigg)^2 \big( \|T\|(\mathbf{B}_{(1-\gamma)/2}(Z)) - C\mathcal{E}^{2+\alpha} \big) 
		\\&\geq \omega_n \bigg(\frac{1-\gamma}{2}\bigg)^{n+2} - C {\mathcal E}^{2+\alpha} 
\end{align*}
which, provided $\varepsilon_2$ is sufficiently small, contradicts \eqref{poincare eqn4}.  Therefore 
\begin{equation} \label{poincare eqn8}
	\sup_{X \in \op{spt} T \cap \mathbf{C}_{\gamma}(0)} \op{dist}(X,\op{spt} \mathbf{P}_{a}) \leq 1-\gamma . 
\end{equation}
By \eqref{poincare eqn8} and the fact that $\|T\|((B_{(1+\gamma)/2}(0) \setminus K) \times \mathbb{R}^m) \leq C {\mathcal E}^{2+\alpha}$, 
\begin{equation} \label{poincare eqn9}
	\int_{(B_{\gamma}(0) \setminus K) \times \mathbb{R}^m} \op{dist}^2(X,\op{spt} \mathbf{P}_{a}) \,d\|T\|(X) \leq C {\mathcal E}^{2+\alpha} 
\end{equation}
for some constant $C = C(n,m,q,\gamma) \in (0,\infty)$.  By \eqref{poincare eqn4}, \eqref{poincare eqn9}, and \eqref{poincare eqn5}, we obtain \eqref{poincare concl}. 
\end{proof}

\begin{proof}[Proof of Theorem~\ref{separation thm1}]
Let $\varepsilon_0 = \varepsilon_0(n,m,q,\gamma) \in (0,1)$ to be later determined.  Choose $\theta = \theta(n,m,q) \in (0,1/8)$ so that $C \theta^{\mu/2} < 1$, where $\mu$ and $C$ are as in Lemma~\ref{excess decay lemma}.  Suppose that $T$ is an $n$-dimensional locally area-minimizing rectifiable current of $\mathbf{C}_1(0)$ such that \eqref{separation1 hyp} holds true.  Let $u : B_{(3+\gamma)/4}(0) \rightarrow \mathcal{A}_q(\mathbb{R}^m)$ and $K \subset B_{(3+\gamma)/4}(0)$ be the Lipschitz approximation of $T$ as in Theorem~\ref{lip approx thm} with $x_0 = 0$, $\rho = 1$, and $(3+\gamma)/4$ in place of $\gamma$.  Provided $\varepsilon_0$ is sufficiently small, by Lemma~\ref{poincare lemma} (with $(1+\gamma)/2$ in place of $\gamma$) there exists $a \in \mathcal{A}_q(\mathbb{R}^m)$ and associated sum of planes $\mathbf{P} \in \mathcal{P}_q$ as in Definition~\ref{associated plane defn} such that 
\begin{equation} \label{separation1 eqn2}
	\int_{\mathbf{C}_{(1+\gamma)/2}(0)} \op{dist}^2(X,\op{spt} \mathbf{P}) \,d\|T\|(X) 
		+ \int_{B_{(1+\gamma)/2}(0)} \mathcal{G}(u,a)^2 \leq C_1 {\mathcal E}(T,\mathbf{C}_1(0))^2 
\end{equation}
for some constant $C_1 = C_1(n,m,q,\gamma) \in (0,\infty)$.  Notice that if ${\mathcal E}(T,\mathbf{C}_1(0)) = 0$, then by \eqref{separation1 eqn2} we must have that 
\begin{equation*}
	\op{spt} T \cap \mathbf{C}_{(1+\gamma)/2}(0) \subset \op{spt} \mathbf{P} 
\end{equation*}
thereby proving the theorem.  Hence we may assume that ${\mathcal E}(T,\mathbf{C}_1(0)) > 0$. 

Let $i_0 \geq 1$ be the integer such that 
\begin{equation}\label{separation1 eqn3} 
	\theta^{i \,(\mu/2-1)} {\mathcal E}(T,\mathbf{C}_1(0)) < \varepsilon_0 
\end{equation}
for all $i \in \{0,1,2,\ldots,i_0-1\}$ and 
\begin{equation}\label{separation1 eqn4} 
	\theta^{i_0 (\mu/2-1)} {\mathcal E}(T,\mathbf{C}_1(0)) \geq \varepsilon_0 .
\end{equation}
Note that by the assumption \eqref{separation1 hyp}, \eqref{separation1 eqn3} holds true for $i = 0$.  Moreover, since ${\mathcal E}(T,\mathbf{C}_1(0)) > 0$, \eqref{separation1 eqn4} holds true for some $i_0$.  Fix any $\xi \in B_{\gamma}(0)$.  We claim that for each $i \in \{0,1,2,\ldots,i_0\}$ 
\begin{equation} \label{separation1 eqn5}
	{\mathcal E}(T,\mathbf{C}_{\theta^i (1-\gamma)/2}(\xi)) \leq \left(\frac{2}{1-\gamma}\right)^{n/2} \theta^{i (\mu/2-1)} {\mathcal E}(T,\mathbf{C}_1(0)) . 
\end{equation}
Notice that when $i = 0$, 
\begin{align*}
	{\mathcal E}(T,\mathbf{C}_{(1-\gamma)/2}(\xi)) 
		&= \frac{1}{2\omega_n ((1-\gamma)/2)^n} \int_{\mathbf{C}_{(1-\gamma)/2}(\xi)} |\vec T - \vec P_0|^2 \,d\|T\|  
		\\&\leq \frac{1}{2\omega_n ((1-\gamma)/2)^n} \int_{\mathbf{C}_1(0)} |\vec T - \vec P_0|^2 \,d\|T\| 
		= \left(\frac{2}{1-\gamma}\right)^n {\mathcal E}(T,\mathbf{C}_1(0))^2 , \nonumber 
\end{align*}
proving \eqref{separation1 eqn5} with $i = 0$.  Suppose that for some $i \in \{0,1,2,\ldots,i_0-1\}$ we know that \eqref{separation1 eqn5} holds true.  Then by \eqref{separation1 eqn3} and \eqref{separation1 eqn5}, 
\begin{equation} \label{separation1 eqn7}
	{\mathcal E}(T,\mathbf{C}_{\theta^i (1-\gamma)/2}(\xi)) \leq \left(\frac{2}{1-\gamma}\right)^{n/2} \theta^{i (\mu/2-1)} {\mathcal E}(T,\mathbf{C}_1(0)) 
		< \left(\frac{2}{1-\gamma}\right)^{n/2} \varepsilon_0 . 
\end{equation}
Hence provided $\varepsilon_0$ is sufficiently small, by Lemma~\ref{excess decay lemma} and \eqref{separation1 eqn5}, 
\begin{equation*}
	{\mathcal E}(T,\mathbf{C}_{\theta^{i+1} (1-\gamma)/2}(\xi)) \leq \theta^{\mu/2-1} {\mathcal E}(T,\mathbf{C}_{\theta^i (1-\gamma)/2}(\xi)) 
		\leq \left(\frac{2}{1-\gamma}\right)^{n/2} \theta^{(i+1) (\mu/2-1)} {\mathcal E}(T,\mathbf{C}_1(0)) ,
\end{equation*}
thereby proving that \eqref{separation1 eqn5} holds true with $i+1$ in place of $i$. 

Recall that \eqref{separation1 eqn7} holds true for all $i \in \{0,1,2,\ldots,i_0-1\}$ and observe that 
\begin{equation*} 
	{\mathcal E}(T,\mathbf{C}_{\theta^{i_0} (1-\gamma)/2}(\xi)) \leq \theta^{-n/2} {\mathcal E}(T,\mathbf{C}_{\theta^{i_0-1} (1-\gamma)/2}(\xi)) 
		\leq \theta^{-n/2} \left(\frac{2}{1-\gamma}\right)^n \varepsilon_0 .
\end{equation*}
Hence provided $\varepsilon_0$ is sufficiently small, for each $i \in \{1,2,\ldots,i_0\}$ we can let $u_i : B_{3\theta^i (1-\gamma)/8}(\xi) \rightarrow \mathcal{A}_q(\mathbb{R}^m)$ and $K_i \subset B_{3\theta^i (1-\gamma)/8}(\xi)$ be as in Theorem~\ref{lip approx thm} with $\xi$, $\theta^i (1-\gamma)/2$, $3/4$, $u_i$, and $K_i$ in place of $x_0$, $\rho$, $\gamma$, $u$, and $K$.  By Lemma~\ref{poincare lemma} (with $\gamma = 1/2$), for each $i \in \{1,2,\ldots,i_0\},$ there exists $a_i \in \mathcal{A}_q(\mathbb{R}^m)$ and associated sum of planes $\mathbf{P}_i \in \mathcal{P}_q$ as in Definition~\ref{associated plane defn} such that 
\begin{align} \label{separation1 eqn9}
	&\frac{1}{\omega_n \theta^{i (n+2)}} \int_{\mathbf{C}_{\theta^i (1-\gamma)/4}(\xi)} \op{dist}^2(X,\op{spt} \mathbf{P}_i) \,d\|T\|(X) 
	+ \frac{1}{\omega_n \theta^{i (n+2)}} \int_{B_{\theta^i (1-\gamma)/4}(\xi)} \mathcal{G}(u_i,a_i)^2 
	\\&\hspace{15mm} \leq C {\mathcal E}(T,\mathbf{C}_{\theta^i (1-\gamma)/2}(\xi))^2 \leq C \theta^{2 i (\mu/2-1)} {\mathcal E}(T,\mathbf{C}_1(0))^2 , \nonumber
\end{align}
where $C = C(n,m,q,\gamma) \in (0,\infty)$ are constants.  When $i = 0$, let $u_0 = u$ be the Lipschitz approximation of $T$ in $\mathbf{C}_{(3+\gamma)/4}(0)$ and $K = K_0$ as chosen above.  Let $a = a_0$ and $\mathbf{P}_0 = \mathbf{P}$ be as chosen above so that \eqref{separation1 eqn2} holds true.  Provided $\varepsilon_0$ is sufficiently small, by \eqref{lip approx concl} for each $i \in \{0,1,2,\ldots,i_0-1\}$ 
\begin{align} \label{separation1 eqn10}
	\mathcal{L}^n(\{ x \in B_{\theta^{i+1} (1-\gamma)/4}(\xi) : u_i(x) \neq u_{i+1}(x) \}) 
	&\leq \mathcal{L}^n(B_{\theta^{i+1} (1-\gamma)/4}(\xi) \setminus (K_i \cup K_{i+1})) 
	\\ \leq C {\mathcal E}(T,\mathbf{C}_{\theta^i  (1-\gamma)/2}(\xi))^{2+\alpha} \theta^{in} 
	&\leq C \varepsilon_0^{2+\alpha} \theta^{in} < \frac{1}{2} \,\omega_n \bigg(\frac{(1-\gamma) \,\theta^{i+1}}{4}\bigg)^n , \nonumber
\end{align}
where $C = C(n,m,q,\gamma) \in (0,\infty)$ and $\alpha = \alpha(n,m,q) \in (0,1)$ are constants.  Thus by the triangle inequality, \eqref{separation1 eqn2}, \eqref{separation1 eqn9}, and \eqref{separation1 eqn10}, 
\begin{align} \label{separation1 eqn11}
	\mathcal{G}(a_i, a_{i+1})^2 
	\leq\,& \frac{C}{\theta^{(i+1)n}} \int_{B_{\theta^{i+1} (1-\gamma)/4}(\xi) \cap \{ u_i = u_{i+1} \}} \mathcal{G}(u_{i+1},a_i)^2 
		\\&+ \frac{C}{\theta^{(i+1)n}} \int_{B_{\theta^{i+1} (1-\gamma)/4}(\xi) \cap \{ u_i = u_{i+1} \}} \mathcal{G}(u_{i+1},a_{i+1})^2 
		\nonumber
	\\ \leq\,& \frac{C}{\theta^{(i+1)n}} \int_{B_{\theta^i (1-\gamma)/4}(\xi)} \mathcal{G}(u_i,a_i)^2 
		+ \frac{C}{\theta^{(i+1)n}} \int_{B_{\theta^{i+1} (1-\gamma)/4}(\xi)} \mathcal{G}(u_{i+1},a_{i+1})^2 \nonumber
	\\ \leq\,& C \theta^{\mu i} {\mathcal E}(T,\mathbf{C}_1(0))^2 , \nonumber
\end{align}
where $C = C(n,m,q,\gamma) \in (0,\infty)$ are constants.  By applying \eqref{separation1 eqn11} using the triangle inequality, for each integer $i \in \{0,1,2,\ldots,i_0\},$ 
\begin{equation} \label{separation1 eqn12}
	\mathcal{G}(a_i, a)  
	\leq \sum_{k=0}^{i-1} \mathcal{G}(a_k, a_{k+1}) 
	\leq \sum_{k=0}^{i-1} C \theta^{k \mu/2} {\mathcal E}(T,\mathbf{C}_1(0))
	\leq C {\mathcal E}(T,\mathbf{C}_1(0)) , 
\end{equation}
where $C = C(n,m,q,\gamma) \in (0,\infty)$ are constants.  By \eqref{separation1 eqn2}, \eqref{separation1 eqn9}, and \eqref{separation1 eqn3}, for each $i \in \{0,1,2,\ldots,i_0\},$ 
\begin{equation*}
	\frac{1}{\omega_n \theta^{i (n+2)}} \int_{\mathbf{C}_{\theta^i (1-\gamma)/4}(\xi)} \op{dist}^2(X,\op{spt} \mathbf{P}_i) \,d\|T\|(X) 
	\leq C \theta^{2 i (\mu/2-1)} {\mathcal E}(T,\mathbf{C}_1(0))^2 \leq C\varepsilon_0^2 < \left(\frac{1-\gamma}{16}\right)^{n+2} 
\end{equation*}
provided $\varepsilon_0$ is sufficiently small, where $C = C(n,m,q,\gamma) \in (0,\infty)$ are constants.  Thus by Lemma~\ref{sepmono lemma} (with $1/2$, $\eta_{(\xi,0),\theta^i (1-\gamma)/4\#} T$, and $\mathbf{P}_i$ in place of $\gamma$, $T$, and $\mathbf{P}$) and \eqref{separation1 eqn9}, for each $i \in \{0,1,2,\ldots,i_0\},$ 
\begin{align} \label{separation1 eqn13}
	\sup_{X \in \op{spt} T \cap \mathbf{C}_{\theta^i (1-\gamma)/8}(\xi)} \op{dist}(X,\op{spt} \mathbf{P}_i) 
	&\leq 2 \left( \frac{1}{\omega_n} \int_{\mathbf{C}_{\theta^i (1-\gamma)/4}(\xi)} \op{dist}^2(X,\op{spt} \mathbf{P}_i) \,d\|T\|(X) \right)^{\frac{1}{n+2}} 
	\\&\leq C \theta^{i \frac{n+\mu}{n+2}} {\mathcal E}(T,\mathbf{C}_1(0))^{\frac{2}{n+2}} , \nonumber 
\end{align}
where $C = C(n,m,q,\gamma) \in (0,\infty)$ is a constant. 

Now fix $\varepsilon_0 = \varepsilon_0(n,m,q,\gamma)$ small enough for the above discussion to hold true.  By \eqref{separation1 eqn4}, 
\begin{equation*}
	\theta^{i_0} \leq \left(\frac{{\mathcal E}(T,\mathbf{C}_1(0))}{\varepsilon_0}\right)^{\frac{1}{1-\mu/2}} . 
\end{equation*}
Hence by \eqref{separation1 eqn13} 
\begin{align} \label{separation1 eqn14}
	\sup_{X \in \op{spt} T \cap \mathbf{C}_{\theta^{i_0} (1-\gamma)/8}(\xi)} \op{dist}(X,\op{spt} \mathbf{P}_{i_0}) 
	&\leq C \theta^{i_0 \frac{n+\mu}{n+2}} {\mathcal E}(T,\mathbf{C}_1(0))^{\frac{2}{n+2}} 
	\\&\leq C {\mathcal E}(T,\mathbf{C}_1(0))^{\frac{1}{1-\mu/2} \cdot \frac{n+\mu}{n+2} + \frac{2}{n+2}} = C {\mathcal E}(T,\mathbf{C}_1(0))^{\frac{1}{1-\mu/2}} , \nonumber 
\end{align}
where $C = C(n,m,q,\gamma) \in (0,\infty)$ are constants.  By the triangle inequality, \eqref{separation1 eqn12}, and \eqref{separation1 eqn14}, 
\begin{eqnarray*}
	&&\sup_{X \in \op{spt} T \cap \mathbf{C}_{\theta^{i_0} (1-\gamma)/8}(\xi)} \op{dist}(X,\op{spt} \mathbf{P})\nonumber\\ 
	&&\hspace{.5in}\leq\, \sup_{X \in \op{spt} T \cap \mathbf{C}_{\theta^{i_0} (1-\gamma)/8}(\xi)} \op{dist}(X,\op{spt} \mathbf{P}_{i_0}) 
		+ \sup_{X \in \op{spt} \mathbf{P}_{i_0}} \op{dist}(X,\op{spt} \mathbf{P}) \nonumber\\ 
	&&\hspace{1in}\leq\, \sup_{X \in \op{spt} T \cap \mathbf{C}_{\theta^{i_0} (1-\gamma)/8}(\xi)} \op{dist}(X,\op{spt} \mathbf{P}_{i_0}) 
		+ \mathcal{G}(a_{i_0}, a) \leq\, C {\mathcal E}(T,\mathbf{C}_1(0)) , \nonumber
\end{eqnarray*}
where $C = C(n,m,q,\gamma) \in (0,\infty)$ is a constant, thereby proving \eqref{separation1 concl}.
\end{proof}

\begin{remark}\label{separation thm1 rmk}{\rm
Let $\mathbf{P}$ be as in Theorem~\ref{separation thm1} and let $a \in \mathcal{A}_q(\mathbb{R}^m)$ be the point associated with $\mathbf{P}$ as in Definition~\ref{associated plane defn}.  Then by Theorem~\ref{separation thm1} 
\begin{equation} \label{separation thm1 rmk eqn1}
	\sup_{X \in \op{spt} T \cap \mathbf{C}_{\gamma}(0)} \op{dist}(X, \op{spt} \mathbf{P}) \leq C_0 {\mathcal E}(T,\mathbf{C}_1(0)) 
\end{equation}
for some constant $C_0 = C_0(n,m,q,\gamma) \in (0,\infty)$.  Let 
\begin{equation*}
	\{ X \in \mathbb{R}^{n+m} : \op{dist}(X, \op{spt} \mathbf{P}) < 2C_0 {\mathcal E}(T,\mathbf{C}_1(0)) \} = \bigcup_{i=1}^N \mathbb{R}^n \times U_i
\end{equation*}
where $\{U_i\}$ is a collection of mutually disjoint, connected, open subsets of $\mathbb{R}^m$ (and $C_0$ is as in \eqref{separation thm1 rmk eqn1}).  Then by \eqref{separation1 concl}, 
\begin{equation*}
	T \llcorner \mathbf{C}_{\gamma}(0) = \sum_{i=1}^N T_i \quad\text{where}\quad T_i = T \llcorner (B_{\gamma}(0) \times U_i) 
\end{equation*}
and $T_i$ are locally area-minimizing rectifiable currents with $(\partial T_i) \llcorner \mathbf{C}_{\gamma}(0) = 0$.  By the constancy theorem 
\begin{equation} \label{separation thm1 rmk eqn2}
	\pi_{\#} T_i = q_i \llbracket B_{\gamma}(0) \rrbracket 
\end{equation}
for some integers $q_i$ with $\sum_{i=1}^N q_i = q$.  Provided $\varepsilon_0$ is sufficiently small, by Lemma~\ref{projection lemma2} we have that $q_i \geq 0$ for each $i \in \{1,2,\ldots,N\}$.  Moreover, by Lemma~\ref{projection lemma}, $q_i = 0$ if and only if $\op{spt} T_i \cap \mathbf{C}_{\gamma/2}(0) = \emptyset$.  
We claim that provided $\mathbf{P}$ satisfies \eqref{separation1 eqn2} and \eqref{separation thm1 rmk eqn1} holds true for $C_0 = C_0(n,m,q,\gamma)$ sufficiently large, we can guarantee that $q_i > 0$ for all $i$.  
To see this, let $u : B_{(3+\gamma)/4}(0) \rightarrow \mathcal{A}_q(\mathbb{R}^m)$ and $K$ be as in Theorem~\ref{lip approx thm} with $x_0 = 0$, $\rho = 1$, and $(3+\gamma)/4$ in place of $\gamma$.  Provided $\varepsilon_0$ is sufficiently small, by \eqref{separation1 eqn2} and \eqref{lip approx concl} there exists a set $\Omega \subset K \cap B_{\gamma/2}(0)$ with $\mathcal{L}^n(\Omega) \geq \tfrac{1}{2} \,\omega_n (\gamma/2)^n$ and 
\begin{equation} \label{separation thm1 rmk eqn3}
	\mathcal{G}(u(x),a) \leq C_2 {\mathcal E}(T,\mathbf{C}_1(0)) \text{ for all } x \in \Omega,
\end{equation}
where $C_2 = \sqrt{\tfrac{4C_1}{\omega_n (\gamma/2)^n}}$ (for $C_1$ as in \eqref{separation1 eqn2}).  For each $x \in \Omega$ let $u(x) = \sum_{j=1}^q \llbracket u_j(x) \rrbracket$ where $u_j(x) \in \mathbb{R}^m$.  For each $i \in \{1,2,\ldots,N\}$ there is a plane $P_{k(i)} = \mathbb{R}^n \times \{a_{k(i)}\}$ of $\mathbf{P}$ in $\mathbb{R}^n \times U_i$.  By \eqref{separation thm1 rmk eqn3}, for each $i \in \{1,2,\ldots,N\}$ and $x \in \Omega$ there exists $j(i) \in \{1,2,\ldots,q\}$ such that $X = (x,u_{j(i)}(x)) \in \op{spt} T$ and 
\begin{equation*}
	\op{dist}(X, P_{k(i)}) \leq |u_{j(i)}(x) - a_{k(i)}| \leq \mathcal{G}(u(x),a) \leq C_2 {\mathcal E}(T,\mathbf{C}_1(0)) . 
\end{equation*}
Hence provided we take $C_0$ in \eqref{separation thm1 rmk eqn1} to be large enough that $C_2 < 2C_0$ (where $C_2$ is as in \eqref{separation thm1 rmk eqn3}), $X \in \op{spt} T_i \cap \mathbf{C}_{\gamma/2}(0)$.  Therefore, $q_i > 0$ in \eqref{separation thm1 rmk eqn2} for all $i \in \{1,2,\ldots,N\}$.
}\end{remark}

\subsection{$L^{\infty}$-distance estimate in terms of height excess relative to parallel planes}\label{sec:separation2_sec}
In the special case that $\op{spt} \mathbf{P}$ is a single plane, Theorem~\ref{separation thm3} is a well-known consequence of estimates for subharmonic functions on minimal submanifolds, first established in~\cite{Allard} (see Lemma~\ref{one plane lemma} below).

\begin{lemma} \label{one plane lemma}
If $T$ is an $n$-dimensional locally area-minimizing rectifiable current in $\mathbf{C}_1(0)$ such that $(\partial T) \llcorner \mathbf{C}_1(0) = 0$ and $T \llcorner \mathbf{C}_{\gamma}(0) \neq 0$ and if $P$ is an $n$-dimensional affine plane in $\mathbb{R}^{n+m}$, then 
\begin{equation*}
	\sup_{X \in \op{spt} T \cap \mathbf{C}_{\gamma}(0)} \op{dist}^2(X,P) \leq C \int_{\mathbf{C}_1(0)} \op{dist}^2(X,P) \,d\|T\|(X) 
\end{equation*}
for some constant $C = C(n,m,\gamma) \in (0,\infty)$. 
\end{lemma}

\begin{proof}
Without loss of generality assume $\op{spt} T \not\subseteq P$.  Choose an orthonormal basis $\nu_1,\nu_2,\ldots,\nu_m$ for the orthogonal complement $P^{\perp}$ of $P$.  For $0 < \delta < \int_{\mathbf{C}_1(0)} \op{dist}^2(X,P) \,d\|T\|(X)$, let $\phi_{\delta} : \mathbb{R} \rightarrow [0,\infty)$ be a smooth convex function such that $\phi_{\delta}(t) = 0$ if $|t| \leq \delta/2$ and $\phi_{\delta}(t) = |t| - \delta$ if $|t| \geq 2\delta$.  Noting that by~\cite[Lemma~7.5(3)]{Allard} $f_{i}(x) = \phi_{\delta}(\nu_i \cdot x)$ is subharmonic on the stationary varifold $V = |T|$ associated with $T$ (see~\cite[Section~3.3]{KrumWica}), apply~\cite[Theorem~7.5(6)]{Allard} to $f_{i}$ for $i = 1,2,\ldots,m$.
\end{proof}

Now we want to prove Theorem~\ref{separation thm3} in a special case where $\mathbf{P} \in \mathcal{P}_s$ (as in Definition~\ref{sum of planes defn}) so that $\mathbf{P}$ consists of two or more parallel planes.  We shall additional assume that $\op{sep} \mathbf{P}$ is proportional to $\op{width} \mathbf{P}$ and the $L^2$-distance of $T$ to $\mathbf{P}$ is much smaller than $\op{sep} \mathbf{P}$.

\begin{lemma} \label{separation lemma2}
For all integers $q \geq 1$ and $s \geq 2$ and for all $\gamma \in (0,1)$ and $\Lambda \in (0,\infty),$ there exists $\varepsilon = \varepsilon(n,m,q,s,\gamma,\Lambda) \in (0,1)$ and $\beta = \beta(n,m,q,s,\gamma,\Lambda) \in (0,1)$ such that if $\mathbf{P} \in \mathcal{P}_s$ and $T$ is an $n$-dimensional locally area-minimizing rectifiable current in $\mathbf{C}_1(0)$ such that 
\begin{gather} 
	\label{separation2 hyp1} (\partial T) \llcorner \mathbf{C}_1(0) = 0, \quad \sup_{X \in \op{spt} T} \op{dist}(X,P_0) < \infty, \quad 
		\pi_{\#} T = q \llbracket B_1(0) \rrbracket, \quad {\mathcal E}(T,\mathbf{C}_1(0)) < \varepsilon , \\
	\label{separation2 hyp2} \op{width} \mathbf{P} \leq \Lambda \op{sep} \mathbf{P} , \\
	\label{separation2 hyp3} \int_{\mathbf{C}_1(0)} \op{dist}^2(X, \op{spt} \mathbf{P}) \,d\|T\|(X)
		\leq \beta^2 (\op{sep} \mathbf{P})^2 , 
\end{gather}
then 
\begin{equation} \label{separation2 concl}
	\sup_{X \in \op{spt} T \cap \mathbf{C}_{\gamma}(0)} \op{dist}^2(X,\op{spt} \mathbf{P}) 
		\leq C \int_{\mathbf{C}_1(0)} \op{dist}^2(X,\op{spt} \mathbf{P}) \,d\|T\|(X) , 
\end{equation}
where $C = C(n,m,\gamma) \in (0,\infty)$ is a constant. 
\end{lemma}

Let $T$ and $\mathbf{P}$ be as in Lemma~\ref{separation lemma2}.  In Lemma~\ref{tilt vs sep lemma}, we use a blow-up argument that the tilt excess of $T$ is also much smaller than $\op{sep}\mathbf{P}$.  Thus by Theorem~\ref{separation thm1}, there is a sum-of-planes $\mathbf{Q}$ such that the $L^{\infty}$-distance of $T$ to $\mathbf{Q}$ is much smaller than $\op{sep}\mathbf{P}$.  We show that the distance of $\mathbf{P}$ to $\mathbf{Q}$ is much smaller than $\op{sep}\mathbf{P}$, and thus $T$ separates into locally area-minimizing rectifiable currents near each plane of $\mathbf{P}$.  Using Lemma~\ref{one plane lemma}, the conclusion of Lemma~\ref{separation lemma2} readily follows.

\begin{lemma} \label{tilt vs sep lemma}
For all integers $q \geq 1$ and $s \geq 2,$ and for all $\gamma \in (0,1)$, $\Lambda \in (0,\infty)$ and $\eta \in (0,1),$ there exist $\varepsilon = \varepsilon(n,m,q,s,\gamma,\Lambda,\eta) \in (0,1)$ and $\beta = \beta(n,m,q,s,\gamma,\Lambda,\eta) \in (0,1)$ such that if $\mathbf{P} \in \mathcal{P}_s$ and $T$ is an $n$-dimensional locally area-minimizing rectifiable current in $\mathbf{C}_1(0)$ such that 
\begin{gather} 
	\label{tilt vs sep hyp1} (\partial T) \llcorner \mathbf{C}_1(0) = 0, \quad \sup_{X \in \op{spt} T} \op{dist}(X,P_0) < \infty, \quad 
		\pi_{\#} T = q \llbracket B_1(0) \rrbracket, \quad {\mathcal E}(T,\mathbf{C}_1(0)) < \varepsilon , \\
	\label{tilt vs sep hyp2} \op{width} \mathbf{P} \leq \Lambda \op{sep} \mathbf{P} , \\
	\label{tilt vs sep hyp3} \int_{\mathbf{C}_1(0)} \op{dist}^2(X, \op{spt} \mathbf{P}) \,d\|T\|(X)
		\leq \beta^2 (\op{sep} \mathbf{P})^2 , 
\end{gather}
where $\pi$ is the orthogonal projection map onto $P_0$, then 
\begin{equation} \label{tilt vs sep concl1} 
	{\mathcal E}(T,\mathbf{C}_{\gamma}(0)) \leq \eta \op{sep} \mathbf{P} .
\end{equation}
\end{lemma}

\begin{proof}  
Fix $\eta \in (0,1)$ and the integers $q \geq 1$ and $s \geq 2$.  Suppose to the contrary that for $k = 1,2,3,\ldots$ there are $\varepsilon_k \rightarrow 0^+$, $\beta_k \rightarrow 0^+$, $n$-dimensional locally area-minimizing rectifiable currents $T_k$ in $\mathbf{C}_1(0)$, and $\mathbf{P}_k \in \mathcal{P}_s$ such that \eqref{tilt vs sep hyp1}, \eqref{tilt vs sep hyp2}, and \eqref{tilt vs sep hyp3} hold true with $\varepsilon_k, \beta_k, T_k, \mathbf{P}_k$ in place of $\varepsilon, \beta, T, \mathbf{P}$ but 
\begin{equation} \label{tilt vs sep eqn1} 
	{\mathcal E}(T_k,\mathbf{C}_{\gamma}(0)) > \eta \op{sep} \mathbf{P}_k . 
\end{equation}
Note that since ${\mathcal E}(T_k,\mathbf{C}_{\gamma}(0)) \rightarrow 0$, \eqref{tilt vs sep eqn1} implies that $\op{sep} \mathbf{P}_k \rightarrow 0$. 

We have that $\mathbf{P}_k = \sum_{i=1}^s \llbracket \mathbb{R}^n \times \{a_{k,i}\} \rrbracket$ for some $a_k = \sum_{i=1}^s \llbracket a_{k,i} \rrbracket \in \mathcal{A}_s(\mathbb{R}^m)$ associated with $\mathbf{P}_k$ as in Definition~\ref{associated plane defn}.  By translating, assume that $a_{k,1} = 0$.  Then by \eqref{tilt vs sep hyp2}
\begin{equation} \label{tilt vs sep eqn2}
	|a_{k,i}| \leq \Lambda \op{sep} \mathbf{P}_k 
\end{equation}
for all $k \in \{1,2,3,\ldots\}$ and $i \in \{1,2,\ldots,s\}$.  Note that $P_0 = \mathbb{R}^n \times \{a_{k,1}\}$.  By the triangle inequality, \eqref{projection eqn2}, \eqref{tilt vs sep hyp2}, and \eqref{tilt vs sep hyp3} 
\begin{align*}
	\int_{\mathbf{C}_1(0)} \op{dist}^2(X,P_0) \,d\|T\|(X) 
	\leq 2 \int_{\mathbf{C}_1(0)} \op{dist}^2(X, \op{spt} \mathbf{P}_k) \,d\|T\|(X) 
		&\\+ 2 (q+1) \omega_n (\op{width} \mathbf{P}_k)^2 
	&\leq C (\op{sep} \mathbf{P}_k)^2 
\end{align*} 
for some constant $C = C(n,q,\Lambda) \in (0,\infty)$.  Hence by Lemma~\ref{one plane lemma}, 
\begin{equation} \label{tilt vs sep eqn3} 
	\sup_{X \in \op{spt} T_k \cap \mathbf{C}_{(7+\gamma)/8}(0)} \op{dist}(X,P_0) \leq C \op{sep} \mathbf{P}_k
\end{equation}
for some constant $C = C(n,m,q) \in (0,\infty)$.  By Lemma~\ref{energy est lemma}, 
\begin{equation} \label{tilt vs sep eqn4}
	{\mathcal E}(T_k, \mathbf{C}_{(3+\gamma)/4}(0)) \leq C \op{sep} \mathbf{P}_k . 
\end{equation}
for some constant $C = C(n,m,q,\gamma) \in (0,\infty)$. 

For each sufficiently large $k$ let $u_k : B_{(1+\gamma)/2}(0) \rightarrow \mathcal{A}_q(\mathbb{R}^m)$ and $K_k \subset B_{(1+\gamma)/2}(0)$ be as in Theorem~\ref{lip approx thm} and $w_k : B_{(1+\gamma)/2}(0) \rightarrow \mathcal{A}_q(\mathbb{R}^m)$ be as in Theorems~\ref{harmonic thm} with $x_0 = 0$, $\rho = 1$, and $\eta = 1/k$ and with $\tfrac{2+2\gamma}{3+\gamma}$, $T_k$, $u_k$, $K_k$, and $w_k$ in place of $\gamma$, $T$, $u$, $K$, and $w$.  By \eqref{tilt vs sep eqn3} and truncating $u_k$ if necessary we may assume that 
\begin{equation} \label{tilt vs sep eqn5} 
	\sup_{B_{(1+\gamma)/2}(0)} |u_k| \leq C \op{sep} \mathbf{P}_k
\end{equation}
for some constant $C = C(n,m,q,\gamma) \in (0,\infty)$.  For each $x \in B_{(1+\gamma)/2}(0)$ we write $u_k(x) = \sum_{i=1}^q \llbracket u_{k,i}(x) \rrbracket$ where $u_{k,i}(x) \in \mathbb{R}^m$ and we write $w_k(x) = \sum_{i=1}^q \llbracket w_{k,i}(x) \rrbracket$ where $w_{k,i}(x) \in \mathbb{R}^m$.  We let $\op{spt} a_k = \{a_{k,1},\ldots,a_{k,s}\}$ denote the set of all values of $a_k$.  By the area formula, \eqref{lip approx concl}, \eqref{tilt vs sep eqn4}, \eqref{tilt vs sep eqn5}, and \eqref{tilt vs sep hyp3}, 
\begin{align} \label{tilt vs sep eqn7} 
	&\int_{B_{(1+\gamma)/2}(0)} \sum_{i=1}^q \op{dist}^2(u_{k,i}(x), \op{spt} a_k) \,dx
	= \int_{B_{(1+\gamma)/2}(0)} \sum_{i=1}^q \op{dist}^2((x,u_{k,i}(x)), \op{spt} \mathbf{P}_k) \,dx
	\\ \leq\,& \int_{(B_{(1+\gamma)/2}(0) \setminus K_k) \times \mathbb{R}^m} \op{dist}^2(X, \op{spt} \mathbf{P}_k) \,d\|T_k\|(X) 
		+ C (\op{sep} \mathbf{P}_k)^{4+\alpha} \nonumber 
	\\ \leq\,& \beta_k (\op{sep} \mathbf{P}_k)^2 + C (\op{sep} \mathbf{P}_k)^{4+\alpha} , \nonumber
\end{align}
where $C = C(n,m,q,\gamma) \in (0,\infty)$ and $\alpha = \alpha(n,m,q) \in (0,1)$ are constants.  
By \eqref{harmonic concl}, \eqref{tilt vs sep eqn4}, \eqref{tilt vs sep eqn5}, and \eqref{tilt vs sep eqn7}, 
\begin{gather}
	\label{tilt vs sep eqn8} \limsup_{k \rightarrow \infty} \frac{1}{(\op{sep} \mathbf{P}_k)^2} \int_{B_{(1+\gamma)/2}(0)} |w_k|^2 \leq C, \\
	\label{tilt vs sep eqn9} \lim_{k \rightarrow \infty} \frac{1}{(\op{sep} \mathbf{P}_k)^2} \int_{B_{(1+\gamma)/2}(0)} \sum_{i=1}^q 
		\op{dist}^2(w_{k,i}(x), \op{spt} a_k) \,dx = 0, 
\end{gather} 
where $C = C(n,m,q,\gamma) \in (0,\infty)$ is a constant.  By \eqref{tilt vs sep eqn8} and the compactness of Dirichlet energy minimizing $q$-valued functions, after passing to a subsequence $w_k/(\op{sep} \mathbf{P}_k) \rightarrow \widetilde{w}$ uniformly in $B_{\gamma}(0)$ for some Dirichlet energy minimizing function $\widetilde{w} : B_{\gamma}(0) \rightarrow \mathcal{A}_q(\mathbb{R}^m)$.

By \eqref{tilt vs sep eqn2}, after passing to a subsequence for each $i$ there exists $\widetilde{a}_i \in \mathbb{R}^m$ such that $a_{k,i}/(\op{sep} \mathbf{P}_k) \rightarrow \widetilde{a}_i$.  Set $\widetilde{a} = \sum_{i=1}^s \llbracket \widetilde{a}_i \rrbracket \in \mathcal{A}_s(\mathbb{R}^m)$.  By \eqref{tilt vs sep eqn9}, $\op{spt} \widetilde{w}(x) \subseteq \op{spt} \widetilde{a}$ for each $x \in B_{\gamma}(0)$, where $\op{spt} \widetilde{w}(x)$ and $\op{spt} \widetilde{a}$ denote the set of all values of $\widetilde{w}(x)$ and $\widetilde{a}$.  Since $\widetilde{w}$ is continuous on $B_{\gamma}(0)$ and $\op{spt} \widetilde{a}$ is a finite set, $\widetilde{w}$ must be a constant function on $B_{\gamma}(0)$.  By the continuity of Dirichlet energy under uniform limits of Dirichlet energy minimizing $q$-valued functions~\cite[Proposition~3.20]{DeLSpaDirMin}, since $w_k/(\op{sep} \mathbf{P}_k)$ converge uniformly to the constant function $\widetilde{w}$ in $B_{\gamma}(0)$, 
\begin{equation*} 
	\lim_{k \rightarrow \infty} \frac{1}{(\op{sep} \mathbf{P}_k)^2} \int_{B_{\gamma}(0)} |Dw_k|^2 = \int_{B_{\gamma}(0)} |D\widetilde{w}|^2 = 0 . 
\end{equation*}
By \eqref{lip approx tilt}, \eqref{harmonic concl}, and \eqref{tilt vs sep eqn4}, 
\begin{equation*}
	\lim_{k \rightarrow \infty} \frac{2\omega_n \gamma^n \,{\mathcal E}(T_k,\mathbf{C}_{\gamma}(0))}{(\op{sep} \mathbf{P}_k)^2} 
		= \lim_{k \rightarrow \infty} \frac{1}{(\op{sep} \mathbf{P}_k)^2} \int_{B_{\gamma}(0)} |Du_k|^2 
		= \lim_{k \rightarrow \infty} \frac{1}{(\op{sep} \mathbf{P}_k)^2} \int_{B_{\gamma}(0)} |Dw_k|^2 = 0, 
\end{equation*}
thereby showing that \eqref{tilt vs sep concl1} must hold true for all sufficiently large $k$, contrary to assumption. 
\end{proof}

\begin{proof}[Proof of Lemma~\ref{separation lemma2}]
Let $\mathbf{P} \in \mathcal{P}_s$ and $T$ is an $n$-dimensional locally area-minimizing rectifiable current in $\mathbf{C}_1(0)$ such that \eqref{separation2 hyp1}--\eqref{separation2 hyp3} hold true.  Let $\eta = \eta(n,m,q,s,\gamma) \in (0,1)$ to be later determined.  Provided $\varepsilon = \varepsilon(n,m,q,s,\gamma,\Lambda,\eta)$ and $\beta = \beta(n,m,q,s,\gamma,\Lambda,\eta)$ are sufficiently small, by Lemma~\ref{tilt vs sep lemma}  
\begin{equation}\label{separation2 eqn9}
	{\mathcal E}(T, \mathbf{C}_{(3+\gamma)/4}(0)) \leq \eta \op{sep} \mathbf{P} 
\end{equation}
By Theorem~\ref{separation thm1} there exists $\mathbf{Q} \in \mathcal{P}_{q}$ such that 
\begin{equation}\label{separation2 eqn10}
	\sup_{X \in \op{spt} T \cap \mathbf{C}_{(1+\gamma)/2}(0)} \op{dist}(X, \op{spt} \mathbf{Q}) \leq C_1 {\mathcal E}(T,\mathbf{C}_{(3+5\gamma)/8}(0))
\end{equation}
for some constant $C_1 = C_1(n,m,q,\gamma) \in (0,\infty)$ (which we require to be large enough that we can apply Remark~\ref{separation thm1 rmk}).  In particular, by \eqref{separation2 eqn9} 
\begin{equation}\label{separation2 eqn11}
	\sup_{X \in \op{spt} T \cap \mathbf{C}_{(1+\gamma)/2}(0)} \op{dist}(X, \op{spt} \mathbf{Q}) < C_1 \eta \op{sep} \mathbf{P} . 
\end{equation}
(where $C$ is as in \eqref{separation2 eqn10}).  Express 
\begin{equation*}
	\{ X \in \mathbb{R}^{n+m} : \op{dist}(X,\op{spt} \mathbf{Q}) < C_1 \eta \op{sep} \mathbf{P} \} 
		= \bigcup_{i=1}^{\widehat{N}} \mathbb{R}^n \times \widehat{U}_i
\end{equation*}
for some collection of mutually disjoint connected open subsets $\{\widehat{U}_i\}$ of $\mathbb{R}^m$, where $C_1$ is as in \eqref{separation2 eqn10}.  By \eqref{separation2 eqn11}, 
\begin{equation*}
	T \llcorner \mathbf{C}_{(1+\gamma)/2}(0) = \sum_{i=1}^{\widehat{N}} \widehat{T}_i 
		\quad\text{where}\quad \widehat{T}_i = T \llcorner B_{(1+\gamma)/2}(0) \times \widehat{U}_i . 
\end{equation*}
Clearly $(\partial \widehat{T}_i) \llcorner \mathbf{C}_{(1+\gamma)/2}(0) = 0$.  By \eqref{separation2 hyp1}, constancy theorem, and Lemma~\ref{projection lemma2} 
\begin{equation}\label{separation2 eqn12}
	\pi_{\#} \widehat{T}_i = \widehat{q}_i \llbracket B_{(1+\gamma)/2}(0) \rrbracket
\end{equation}
for some integers $\widehat{q}_i \geq 0$ such that $q = \sum_{i=1}^{\widehat{N}} \widehat{q}_i$.  By Remark~\ref{separation thm1 rmk}, $\widehat{q}_i > 0$ for all $i$.  Given a plane $Q_i \subset \op{spt} \mathbf{Q}$, $Q_i \subset \mathbb{R}^n \times \widehat{U}_{j(i)}$ for some $j(i) \in \{1,2,\ldots,\widehat{N}\}$ and the distance of each point $X \in \op{spt} \widehat{T}_{j(i)}$ to $Q_i$ is $\leq 2 q C_1 \eta \op{sep} \mathbf{P}$.  Hence by the triangle inequality 
\begin{equation*}
	\sup_{Y \in Q_i} \op{dist}(Y,\op{spt}\mathbf{P}) 
	\leq \op{dist}(X,\op{spt}\mathbf{P}) + \op{dist}(X,Q_i)
	\leq \op{dist}(X,\op{spt}\mathbf{P}) + 2q C_1 \eta \op{sep} \mathbf{P} 
\end{equation*}
for each $X \in \op{spt} \widehat{T}_{j(i)}$.  Integrating over $\widehat{T}_{j(i)}$ and using \eqref{separation2 hyp3}, 
\begin{align}\label{separation2 eqn13}
	&\sup_{Y \in \op{spt}\mathbf{Q}} \op{dist}(Y,\op{spt}\mathbf{P}) 
	= \max_i \sup_{Y \in Q_i} \op{dist}(Y,\op{spt}\mathbf{P}) 
	\\ \leq\,& \max_i \left( \frac{1}{q \omega_n \big(\tfrac{1+\gamma}{2}\big)^n} \int_{\mathbf{C}_{(1+\gamma)/2}(0)} \op{dist}^2(X,\op{spt}\mathbf{P}) 
		\,d\|\widehat{T}_{j(i)}\|(X) \right)^{1/2} + 2 q C_1 \eta \op{sep} \mathbf{P} \nonumber 
	\\ \leq\,& \left( \frac{1}{q \omega_n \big(\tfrac{1+\gamma}{2}\big)^n} \int_{\mathbf{C}_{(1+\gamma)/2}(0)} \op{dist}^2(X,\op{spt}\mathbf{P}) \,d\|T\|(X) 
		\right)^{1/2} + 2q C_1 \eta \op{sep} \mathbf{P} \nonumber
	\\ \leq\,& C \,(\beta + \eta) \op{sep} \mathbf{P} \nonumber
\end{align}
for some constant $C = C(n,m,q,\gamma) \in (0,\infty)$.  By the triangle inequality, \eqref{separation2 eqn11}, and \eqref{separation2 eqn13}
\begin{align*}
	\sup_{X \in \op{spt} T \cap \mathbf{C}_{(1+\gamma)/2}(0)} \op{dist}(X, \op{spt} \mathbf{P}) 
	&\leq \sup_{X \in \op{spt} T \cap \mathbf{C}_{(1+\gamma)/2}(0)} \op{dist}(X, \op{spt} \mathbf{Q}) 
		+ \sup_{X \in \op{spt} \mathbf{Q}} \op{dist}(X, \op{spt} \mathbf{P}) 
	\\&\leq C \,(\beta + \eta) \op{sep} \mathbf{P} < \tfrac{1}{3} \op{sep} \mathbf{P} ,
\end{align*}
where $C = C(n,m,q,\gamma) \in (0,\infty)$ is a constant and in the last step we assume that $\beta$ and $\eta$ are small enough that $C \,(\beta + \eta) < 1/3$.  Hence 
\begin{equation*}
	T = \sum_{i=1}^s T_i \quad\text{where}\quad 
	T_i = T \llcorner \left\{ X \in \mathbf{C}_{(1+\gamma)/2}(0) : \op{dist}(X,P_i) < \tfrac{1}{3} \op{sep} \mathbf{P} \right\} . 
\end{equation*}
We can apply Lemma~\ref{one plane lemma} with $\eta_{0,(1+\gamma)/2\#} T_i$ and $P_i$ in place of $T$ and $P$ to conclude that either $T_i \llcorner \mathbf{C}_{\gamma}(0) = 0$ or 
\begin{equation*}
	\sup_{X \in \op{spt} T_i \cap \mathbf{C}_{\gamma}(0)} \op{dist}^2(X,P_i) 
		\leq C \int_{\mathbf{C}_{(1+\gamma)/2}(0)} \op{dist}^2(X,P_i) \,d\|T_i\|(X) 
\end{equation*}
for some constant $C = C(n,m,\gamma) \in (0,\infty)$, thereby proving \eqref{separation2 concl}. 
\end{proof}

\subsection{$L^{\infty}$-distance estimate in terms of height excess relative to disjoint, not-necessarily-parallel planes}\label{sec:separation nonparallel}

In this section we prove Theorem~\ref{separation thm3} in the general case where $\mathbf{P}$ is a sum of non-intersecting affine planes which are not necessarily parallel.  

Given $1 < p \leq s$ and $\mathbf{P} \in \Pi_{s,p}$ as in \eqref{sum of tilted planes form}, we define 
\begin{align}\label{separate3 sep width}
	\op{minsep} \mathbf{P} &= \min_{1 \leq i \leq p} \,\inf_{X \in P_i \cap \mathbf{C}_1(0)} \op{dist}(X, \op{spt} \mathbf{P} \setminus P_i) , \\ 
	\op{width} \mathbf{P} &= \max_{i \neq j} \,\sup_{X \in P_i \cap \mathbf{C}_1(0)} \op{dist}(X,P_j) . \nonumber 
\end{align}
If instead $\mathbf{P} \in \Pi_{s,1}$, we define $\op{minsep} \mathbf{P} = \infty$ and $\op{width} \mathbf{P} = 0$. 

First in Lemma~\ref{separation3_1 lemma}, we prove Theorem~\ref{separation thm3} with the additional assumptions that $\op{minsep} \mathbf{P}$ is proportional to $\op{width} \mathbf{P}$ and the $L^2$-distance of $T$ to $\mathbf{P}$ is much smaller than $\op{minsep} \mathbf{P}$. We prove Lemma~\ref{separation3_1 lemma} by showing that locally in cylinders with small radii we can replace the planes of $\mathbf{P}$ with planes parallel to $P_0$ and apply Lemma~\ref{separation lemma2}.  

Note that in the proofs of Lemma~\ref{separation3_1 lemma} and Theorem~\ref{separation thm3}, since we assume that $|\vec P_i - \vec P_0| < \varepsilon_0$ for each $i \in \{1,2,\ldots,p\},$ we may assume that 
\begin{equation}\label{separation3 matrix}
	P_i = \{ (x, b_i + A_i x) : x \in \mathbb{R}^n \} , \quad \|A_i\| < 2\varepsilon_0 
\end{equation}
for some $m \times n$ matrix $A_i$ and some $b_i \in \mathbb{R}^m$.  When $p > 1$, \eqref{separation3 hyp2} gives us 
\begin{equation}\label{separation3 hyp2 matrix}
	\|A_i - A_j\| \leq 2\kappa \inf_{X \in P_i \cap \mathbf{C}_1(0)} \op{dist}(X,P_j)
\end{equation}
for all $i \neq j$.  

Let $0 < \sigma \leq \rho$ and let $P = \{ (x, b + Ax) : x \in \mathbb{R}^n \}$ be an $n$-dimensional affine plane, where $A$ is an $m \times n$ matrix and $b \in \mathbb{R}^m$.  For $i \in \{1,2,\ldots,p\}$, let $z \in B_{\sigma}(0)$ such that $Z = (z,b+Az)$ satisfies 
\begin{equation*}
	\op{dist}(Z,P_i) = \inf_{X \in P \cap \mathbf{C}_{\sigma}(0)} \op{dist}(X,P_i) . 
\end{equation*}
Then by \eqref{separation3 matrix} and the triangle inequality 
\begin{align}\label{separation3 hausdist1}
	\op{dist}_{\mathcal H}(P_i \cap \mathbf{C}_{\rho}(0), P \cap \mathbf{C}_{\rho}(0))
	&\leq \sup_{x \in B_{\rho}(0)} |(b_i + A_i x) - (b + Ax)| 
	\\&\leq |(b_i + A_i z) - (b + Az)| + 2\rho \|A_i - A\| \nonumber 
	\\&\leq 2 \inf_{X \in P \cap \mathbf{C}_{\sigma}(0)} \op{dist}(X,P_i) + 2\rho \|A_i - A\| . \nonumber 
\end{align}
In particular, if $0 < \sigma \leq 1$ and $A = A_j$ for some $j \in \{1,2,\ldots,p\}$ with $j \neq i$, then by \eqref{separation3 hyp2 matrix}
\begin{equation}\label{separation3 hausdist2}
	\op{dist}_{\mathcal H}(P_i \cap \mathbf{C}_{\rho}(0), P_j \cap \mathbf{C}_{\rho}(0))
	\leq (2 + 4\kappa\rho) \inf_{X \in P_j \cap \mathbf{C}_{\sigma}(0)} \op{dist}(X,P_i) .
\end{equation}

\begin{lemma} \label{separation3_1 lemma}
Let $q \geq 1$ and $1 < p \leq s$ be integers.  For each $\gamma \in (0,1)$ and $\kappa \in (0,1)$ there exists $\varepsilon_0(n,m,q,s,\gamma,\kappa) \in (0,1)$, $\eta(n,m,q,s,\gamma,\kappa) \in (0,1)$, and $\beta(n,m,q,s,\gamma,\kappa) \in (0,1)$ such that if $T$ is an $n$-dimensional locally area-minimizing rectifiable current in $\mathbf{C}_1(0)$ and $\mathbf{P} \in \Pi_{s,p}$ as in \eqref{sum of tilted planes form} such that \eqref{separation3 hyp1} and \eqref{separation3 hyp2} hold true and 
\begin{gather} 
	\label{separation3_1 hyp1} \inf_{1 \leq i \leq p} \frac{1}{\omega_n} \int_{\mathbf{C}_1(0)} \op{dist}^2(X, P_i) \,d\|T\|(X) < \eta^2 , \\
	\label{separation3_1 hyp2} \op{width} \mathbf{P} \leq \lambda \op{minsep} \mathbf{P} , \\
	\label{separation3_1 hyp3} \int_{\mathbf{C}_1(0)} \op{dist}^2(X, \op{spt} \mathbf{P}) \,d\|T\|(X) \leq \beta^2 (\op{minsep} \mathbf{P})^2 .
\end{gather}
Then 
\begin{equation} \label{separation3_1 concl} 
	\sup_{X \in \op{spt} T \cap \mathbf{C}_{\gamma}(0)} \op{dist}^2(X, \op{spt} \mathbf{P}) 
		\leq C \int_{\mathbf{C}_1(0)} \op{dist}^2(X, \op{spt} \mathbf{P}) \,d\|T\|(X) 
\end{equation}
for some constant $C = C(n,m,\gamma) \in (0,\infty)$. 
\end{lemma}

\begin{proof}
Without loss of generality assume that $\mathbf{P}$ is a sum of $s$ distinct multiplicity one planes and thus $p = s$.  By \eqref{separation3_1 hyp1} we may assume that 
\begin{equation} \label{separation3_1 eqn1} 
	\frac{1}{\omega_n} \int_{\mathbf{C}_1(0)} \op{dist}^2(X, P_1) \,d\|T\|(X) < \eta^2 
\end{equation}
By translating, assume that $0 \in P_1$ and thus $P_1 = \{ (x, A_1 x) : x \in \mathbb{R}^n \}$ where $A_1$ is an $m \times n$ matrix with $\|A_1\| < 2\varepsilon_0$.  Thus 
\begin{equation} \label{separation3_1 eqn2} 
	\op{dist}_{\mathcal{H}}(P_1 \cap \mathbf{C}_1(0), P_0 \cap \mathbf{C}_1(0)) \leq \|A_1\| < 2 \varepsilon_0
\end{equation}
where $\op{dist}_{\mathcal{H}}$ denotes Hausdorff distance.  By Lemma~\ref{sepmono lemma} and \eqref{separation3_1 eqn1}, 
\begin{equation} \label{separation3_1 eqn3} 
	\sup_{X \in \op{spt} T \cap \mathbf{C}_{(7+\gamma)/8}(0)} \op{dist}(X, P_1) \leq 2 \eta^{\frac{2}{n+2}} .  
\end{equation}
By \eqref{separation3_1 eqn2} and \eqref{separation3_1 eqn3}
\begin{align} \label{separation3_1 eqn4} 
	\sup_{X \in \op{spt} T \cap \mathbf{C}_{(7+\gamma)/8}(0)} \op{dist}(X, P_0) 
	&\leq \sup_{X \in \op{spt} T \cap \mathbf{C}_{(7+\gamma)/8}(0)} \op{dist}(X, P_1) 
		\\&\hspace{15mm} + \op{dist}_{\mathcal{H}}(P_1 \cap \mathbf{C}_1(0), P_0 \cap \mathbf{C}_1(0)) \nonumber
	\\&\leq 2 \eta^{\frac{2}{n+2}} + 2 \varepsilon_0 . \nonumber 
\end{align}
Let $\pi_{P_1} : \mathbb{R}^{n+m} \rightarrow P_1$ denote the orthogonal projection map onto $P_1$ and define the cylinder $\mathbf{C}_{(3+\gamma)/4}(0,P_1) = \pi_{P_1}^{-1}(P_1 \cap \mathbf{B}_{(3+\gamma)/4}(0))$.  By $(\partial T) \llcorner \mathbf{C}_1(0) = 0$, \eqref{separation3_1 eqn2}, and \eqref{separation3_1 eqn4}, 
\begin{equation*}
	(\partial T) \llcorner \mathbf{C}_{(3+\gamma)/4}(0,P_1) = 0 .
\end{equation*}
Thus by the constancy theorem, $(\pi_{P_1 \#} T) \llcorner \mathbf{B}_{(3+\gamma)/4}(0)$ is a constant integer multiple of $\llbracket P_1 \rrbracket \llcorner \mathbf{B}_{(3+\gamma)/4}(0)$.  It follows from \eqref{separation3 hyp1}, \eqref{separation3_1 eqn2}, and \eqref{separation3_1 eqn3} that $T$ is weakly close to $q \llbracket P_1 \rrbracket$ in $\mathbf{B}_{(3+\gamma)/4}(0)$ and thus by the continuity of push-forwards in the weak topology 
\begin{equation*}
	(\pi_{P_1 \#} T) \llcorner \mathbf{B}_{(3+\gamma)/4}(0) = q \llbracket P_1 \rrbracket \llcorner \mathbf{B}_{(3+\gamma)/4}(0) .
\end{equation*}
By the triangle inequality, \eqref{separation3 hyp1}, and \eqref{projection eqn2}, 
\begin{equation*}
	\int_{\mathbf{C}_{(3+\gamma)/4}(0,P_1)} |\vec T - \vec P_1|^2 \,d\|T\| 
	\leq 2 \int_{\mathbf{C}_1(0)} |\vec T - \vec P_0|^2 \,d\|T\| + 2 \,(q+1) \,\omega_n \,|\vec P_1 - \vec P_0|^2 
	\leq C \varepsilon_0^2 
\end{equation*}
for some constant $C = C(n,q) \in (0,\infty)$.  Recalling \eqref{separate3 sep width}, observe that 
\begin{gather*}
	\op{minsep}\mathbf{P} \leq \min_{1 \leq i \leq p} \,\inf_{X \in P_i \cap \mathbf{C}_{(3+\gamma)/4}(0,P_1)} \op{dist}(X,\op{spt}\mathbf{P} \setminus P_i) , \\
	\op{width}\mathbf{P} \geq \max_{i \neq j} \,\sup_{X \in P_i \cap \mathbf{C}_{(3+\gamma)/4}(0,P_1)} \op{dist}(X,P_j) . 
\end{gather*}
Thus by rotating $P_1$ slightly to $P_0$ and rescaling, we may assume that $P_1 = P_0$ and $\vec P_1 = \vec P_0$.  

Let $0 < \sigma = \sigma(n,m,q,\gamma,\kappa,\lambda) < (1-\gamma)/32$ be a constant depending to be later determined.  Let $\{ B_{\sigma}(x_k) : k = 1,2,\ldots,K \}$ be a collection of balls such that $x_k \in B_{(1+\gamma)/2}(0)$ for each $k$, $B_{(1+\gamma)/2}(0) \subseteq \bigcup_{k=1}^K B_{\sigma}(x_k)$, and $K \leq C(n,\gamma) \,\sigma^{-n}$.  Recall that  
\begin{equation*}
	\mathbf{P} = \sum_{i=1}^s \llbracket P_i \rrbracket 
\end{equation*}
for distinct oriented $n$-dimensional planes $P_i$.  For each $1 \leq k \leq K$ and $1 \leq i \leq s$, let $y_{k,i} \in \mathbb{R}^m$ such that $(x_k,y_{k,i})$ is the unique point of $P_i \cap (\{x_k\} \times \mathbb{R}^m)$ and set $\widehat{P}_{k,i} = \mathbb{R}^n \times \{y_{k,i}\}$.  For each $1 \leq k \leq K$ define $\widehat{\mathbf{P}}_k \in \mathcal{P}_q$ by 
\begin{equation*}
	\widehat{\mathbf{P}}_k = \sum_{i=1}^s \llbracket \widehat{P}_{k,i} \rrbracket . 
\end{equation*}
By \eqref{separation3 matrix}, \eqref{separation3 hyp2 matrix}, and \eqref{separation3_1 hyp2} and noting that since $P_1 = P_0$ we have $A_1 = 0$, 
\begin{align} \label{separation3_1 eqn5} 
	&\max_{1 \leq i \leq s} \op{dist}_{\mathcal{H}}(\widehat{P}_{k,i} \cap \mathbf{C}_{2\sigma}(x_k), P_i \cap \mathbf{C}_{2\sigma}(x_k)) 
	\\ \leq\,& 2\sigma \max_{1 \leq i \leq s} \|A_i\|
	= 2\sigma \max_{1 \leq i \leq s} \|A_i - A_1\| 
	\leq 4\kappa \sigma \op{width} \mathbf{P}
	\leq 4\kappa \lambda \sigma \op{minsep} \mathbf{P}. \nonumber 
\end{align}
Since $\widehat{P}_{k,i} = \mathbb{R}^n \times \{y_{k,i}\}$ and $(x_k,y_{k,i}) \in P_i$ for each $k$ and $i$, 
\begin{equation} \label{separation3_1 eqn6} 
	\op{minsep} \mathbf{P} \leq \min_{1 \leq i < j \leq p} |y_{k,i} - y_{k,j}| = \op{sep} \widehat{\mathbf{P}}_k . 
\end{equation}
By \eqref{separation3_1 hyp2} and \eqref{separation3_1 eqn6} 
\begin{equation} \label{separation3_1 eqn7} 
	\op{width} \widehat{\mathbf{P}}_k = \max_{1 \leq i < j \leq p} |y_{k,i} - y_{k,j}| \leq \op{width} \mathbf{P}
		\leq \lambda \op{minsep} \mathbf{P} \leq \lambda \op{sep} \widehat{\mathbf{P}}_k . 
\end{equation}
Provided $\eta^{\frac{2}{n+2}} < \sigma$, by \eqref{separation3_1 eqn3} and $P_1 = P_0$, 
\begin{equation*}
	\op{spt} T \cap \mathbf{C}_{2\sigma}(x_k) 
	= \op{spt} T \cap (B^n_{2\sigma}(x_k) \times B^m_{2\eta^{\frac{2}{n+2}}}(y_{k,1})) 
	\subset \op{spt} T \cap \mathbf{B}_{4\sigma}(x_k,y_{k,1}) . 
\end{equation*}
Hence by the monotonicity formula, \eqref{projection eqn2}, and ${\mathcal E}(T,\mathbf{C}_1(0)) < \varepsilon_0 < 1$, 
\begin{equation} \label{separation3_1 eqn8} 
	\|T\|(\mathbf{C}_{2\sigma}(x_k)) \leq \|T\|(\mathbf{B}_{4\sigma}(x_k,y_{k,1})) 
		\leq \omega_n (4\sigma)^n \|T\|(\mathbf{C}_1(0)) \leq (q+1) \omega_n (4\sigma)^n . 
\end{equation}
Recalling \eqref{separation3 matrix}, let $X \in \op{spt} T \cap \mathbf{C}_{2\sigma}(x_k) $ and find $Z \in \op{spt} \mathbf{P} \cap \mathbf{C}_{2\sigma}(x_k)$ such that $|X - Z| \leq 2 \op{dist}(X, \op{spt} \mathbf{P})$.  By the triangle inequality, 
\begin{equation} \label{separation3_1 eqn9} 
	\op{dist}(X, \op{spt} \widehat{\mathbf{P}}_k)
	\leq |X - Z| + \op{dist}(Z, \op{spt} \widehat{\mathbf{P}}_k)
	\leq 2 \op{dist}(X, \op{spt} \mathbf{P} ) + \op{dist}(Z, \op{spt} \widehat{\mathbf{P}}_k) . 
\end{equation}
By squaring and integrating \eqref{separation3_1 eqn9} over $\|T\|$-a.e.~$X \in \mathbf{C}_{2\sigma}(x_k)$ and using  \eqref{separation3_1 eqn8}, \eqref{separation3_1 hyp3}, and \eqref{separation3_1 eqn5}, 
\begin{align} \label{separation3_1 eqn10} 
	&\frac{1}{\omega_n (2\sigma)^n} \int_{\mathbf{C}_{2\sigma}(x_k)} \op{dist}^2(X, \op{spt} \widehat{\mathbf{P}}_k) \,d\|T\|(X)
	\\ \leq\,& \frac{8}{\omega_n (2\sigma)^n} \int_{\mathbf{C}_{2\sigma}(x_k)} \op{dist}^2(X, \op{spt} \mathbf{P}) \,d\|T\|(X) 
		+ C \sup_{Z \in \op{spt} \mathbf{P} \cap \mathbf{C}_{2\sigma}(x_k)} \op{dist}^2(Z, \op{spt} \widehat{\mathbf{P}}_k) \nonumber 
	\\ \leq\,& C \left( \frac{\beta^2}{\sigma^n} + \sigma^2 \right) (\op{minsep} \mathbf{P})^2 , \nonumber
\end{align}
where $C = C(n,q,\gamma,\kappa,\lambda) \in (0,\infty)$ are constants.  In particular, by \eqref{separation3_1 eqn6}, 
\begin{equation} \label{separation3_1 eqn11} 
	\frac{1}{\omega_n (2\sigma)^n} \int_{\mathbf{C}_{2\sigma}(x_k)} \op{dist}^2(X, \op{spt} \widehat{\mathbf{P}}_k) \,d\|T\|(X)
	\leq C \left( \frac{\beta^2}{\sigma^n} + \sigma^2 \right) (\op{sep} \widehat{\mathbf{P}}_k)^2 
\end{equation}
for some constant $C = C(n,q,\gamma,\kappa,\lambda) \in (0,\infty)$.  In light of \eqref{separation3_1 eqn7} and \eqref{separation3_1 eqn11}, we can apply Lemma~\ref{separation lemma2} to obtain 
\begin{equation} \label{separation3_1 eqn12} 
	\sup_{X \in \op{spt} T \cap \mathbf{C}_{\sigma}(x_k)} \op{dist}^2(X, \op{spt} \widehat{\mathbf{P}}_k) 
		\leq \frac{C}{\sigma^n} \int_{\mathbf{C}_{2\sigma}(x_k)} \op{dist}^2(X, \op{spt} \widehat{\mathbf{P}}_k) \,d\|T\|(X) 
\end{equation}
for some constant $C = C(n,m) \in (0,\infty)$.  By \eqref{separation3_1 eqn10} and \eqref{separation3_1 eqn12}, 
\begin{equation} \label{separation3_1 eqn13} 
	\sup_{X \in \op{spt} T \cap \mathbf{C}_{\sigma}(x_k)} \op{dist}(X, \op{spt} \widehat{\mathbf{P}}_k) 
		\leq C \left( \frac{\beta^2}{\sigma^n} + \sigma^2 \right)^{1/2} \op{minsep} \mathbf{P} 
\end{equation}
for some constant $C = C(n,m,q,\gamma,\kappa,\lambda) \in (0,\infty)$.  Choose $\sigma$ so that $4\kappa \lambda \sigma < 1/9$ and $C\sigma < 1/9$ (for $C$ is as in \eqref{separation3_1 eqn13}).  Then choose $\beta$ so that $C\beta \sigma^{-n/2} < 1/9$ (for $C$ is as in \eqref{separation3_1 eqn13}).  Hence by \eqref{separation3_1 eqn13} and \eqref{separation3_1 eqn6} 
\begin{equation*} 
	\sup_{X \in \op{spt} T \cap \mathbf{C}_{\sigma}(x_k)} \op{dist}(X, \op{spt} \widehat{\mathbf{P}}_k) 
		< \frac{2}{9} \op{minsep} \mathbf{P} \leq \frac{2}{9} \op{sep} \widehat{\mathbf{P}}_k .
\end{equation*}
Hence 
\begin{equation*}
	T_{k,i} = T \llcorner \left\{ X \in \mathbf{C}_{\sigma}(x_k) : \op{dist}(X, \widehat{P}_{k,i}) < \tfrac{2}{9} \op{minsep} \mathbf{P} \right\} 
\end{equation*}
are locally area minimizing rectifiable currents of $\mathbf{C}_{\sigma}(x_k)$ such that $(\partial T_{k,i}) \llcorner \mathbf{C}_{\sigma}(x_k) = 0$ and 
\begin{equation}\label{separation3_1 eqn14} 
	T \llcorner \mathbf{C}_{\sigma}(x_k) = \sum_{i=1}^N T_{k,i} 
\end{equation}
and 
\begin{equation} \label{separation3_1 eqn15} 
	\sup_{X \in \op{spt} T_{k,i}} \op{dist}(X, \widehat{P}_{k,i}) < \frac{2}{9} \op{minsep} \mathbf{P}
\end{equation}
for some constant $C = C(n,m,q) \in (0,\infty)$.  Moreover, for each $X \in \op{spt} T_{k,i}$ there exists $Z \in \widehat{P}_{k,i} \cap \mathbf{C}_{\sigma}(x_k)$ such that $|X - Z| \leq \op{dist}(X,\widehat{P}_{k,i})$.  Thus by the triangle inequality, \eqref{separation3_1 eqn15}, \eqref{separation3_1 eqn5} (recalling that $4\kappa \lambda \sigma < 1/9$), and \eqref{separation3_1 eqn6}, 
\begin{equation}\label{separation3_1 eqn16}
	\sup_{X \in \op{spt} T_{k,i}} \op{dist}(X, P_i) 
		\leq 2 \sup_{X \in \op{spt} T_{k,i}} \op{dist}(X, \widehat{P}_{k,i}) 
			+ \sup_{Z \in \widehat{P}_{k,i} \cap \mathbf{C}_{\sigma}(x_k)} \op{dist}(Z, P_i) 
		< \frac{1}{3} \op{minsep} \mathbf{P} .  
\end{equation}
It follows from \eqref{separation3_1 eqn14} and \eqref{separation3_1 eqn16} that 
\begin{equation*} 
	\sup_{X \in \op{spt} T \cap \mathbf{C}_{(1+\gamma)/2}(0)} \op{dist}(X, \op{spt} \mathbf{P}) < \frac{1}{3} \op{minsep} \mathbf{P} . 
\end{equation*}
Hence 
\begin{equation*} 
	T_i = T \llcorner \left\{ X \in \mathbf{C}_{(1+\gamma)/2}(0) : \op{dist}(X, P_i) < \tfrac{1}{3} \op{minsep} \mathbf{P} \right\} 
\end{equation*}
are locally area minimizing rectifiable currents of $\mathbf{C}_{(1+\gamma)/2}(0)$ such that $(\partial T_i) \llcorner \mathbf{C}_{(1+\gamma)/2}(0) = 0$ and 
\begin{equation*} 
	T \llcorner \mathbf{C}_{(1+\gamma)/2}(0) = \sum_{i=1}^N T_i , \quad\quad
	\sup_{X \in \op{spt} T_i} \op{dist}(X, P_i) \leq \frac{1}{3} \op{minsep} \mathbf{P} . 
\end{equation*}
(Note that we can characterize $T_i$ by 
\begin{equation*} 
	T_i \llcorner (\mathbf{C}_{\sigma}(x_k) \cap \mathbf{C}_{(1+\gamma)/2}(0)) 
		= T_{k,i} \llcorner (\mathbf{C}_{\sigma}(x_k) \cap \mathbf{C}_{(1+\gamma)/2}(0)) 
\end{equation*}
for each $k = 1,2,\ldots,K$.)  We can apply Lemma~\ref{one plane lemma} with $\eta_{0,(1+\gamma)/2\#} T_i$ and $\eta_{0,(1+\gamma)/2\#} P_i$ in place of $T$ and $P$ to obtain \eqref{separation3_1 concl}.
\end{proof}

We will now complete the proof of Theorem~\ref{separation thm3} using double induction on $q$ and $s$ as follows.  Using Theorem~\ref{separation thm1}, one can reduce to the case where the $L^2$-distance of $T$ to one of the planes of $\mathbf{P}$ is small.  Find two planes of $\mathbf{P}$ which are a distance $\op{minsep}\mathbf{P}$ apart and remove one of them to form $\widetilde{\mathbf{P}}$.  We may assume that the $L^2$-distance of $T$ to $\mathbf{P}$ is much smaller than $\op{minsep}\mathbf{P}$, otherwise the conclusion \eqref{separation3 concl} readily follows by induction.  Using the induction hypothesis, we show that $T$ separates into locally area-minimizing rectifiable currents $T_i$ (possibly zero) near the planes of $\widetilde{\mathbf{P}}$.  If $T$ separates into two or more non-zero locally area-minimizing rectifiable currents $T_i$, then by induction we again obtain \eqref{separation3 concl}.  (Notice that for instance $\op{spt} \mathbf{P}$ might consist of three planes $P_1,P_2,P_3$ with $T$ separating into a multiplicity two current $T_1$ near $P_1$ and a multiplicity one current $T_2$ near $P_2, P_3$.  To treat this possibility, we allow the number of planes of $\mathbf{P}$ to be $s > q$).  If on the other hand the distance of $T$ to a single plane $P_1$ of $\mathbf{P}$ is $\leq C \op{minsep}\mathbf{P}$, then we argue (in Lemma~\ref{separation3_1 lemma}) that locally in small cylinders $\mathbf{C}_{2\sigma}(x_k)$ we can rotate the planes of $\mathbf{P}$ slightly to form a union $\widehat{\mathbf P}_k$ of parallel planes.  It follows using Lemma~\ref{separation lemma2} that $T \llcorner \mathbf{C}_{\sigma}(x_k)$ separates into locally area-minimizing rectifiable currents near each plane of $\widehat{\mathbf P}_k$.  Thus $T$ separates into locally area-minimizing rectifiable currents near each plane of $\mathbf{P}$ and \eqref{separation3 concl} follows again.

\begin{proof}[Proof of Theorem~\ref{separation thm3}]
Without loss of generality we may assume that $\mathbf{P}$ is a sum of $s$ distinct multiplicity one planes and thus $p = s$.  We shall proceed by double induction on $q$ and $s$.  The base case $s = 1$ follows from Lemma~\ref{one plane lemma}.  Suppose that $q_0 \geq 1$ and $s_0 > 1$ are integers such that 
\begin{enumerate}
	\item[(H1)]  Theorem~\ref{separation thm3} holds true if $q \in \{1,2,\ldots,q_0-1\}$ and $s \in \{1,2,\ldots,s_0\}$ and 
\end{enumerate}
and either (i) $q_0 = 1$ or (ii) $q_0 > 1$ and 
\begin{enumerate}
	\item[(H2)]  Theorem~\ref{separation thm3} holds true if $q = q_0$ and $s \in \{1,2,\ldots,s_0-1\}$. 
\end{enumerate}
Let $p_0 \in \{1,2,\ldots,s_0\}$, $\mathbf{P} \in \Pi_{s_0,p_0}$, and $T$ be an $n$-dimensional locally area minimizing rectifiable current of $\mathbf{C}_1(0)$ such that \eqref{separation3 hyp1} holds true and either $p_0 = 1$ or $p_0 > 1$ and \eqref{separation3 hyp2} hold true (with $q = q_0$, $p = p_0$, and $s = s_0$).  Notice that we may assume that $\mathbf{P}$ consists of exactly $s_0$ distinct planes (so that $p_0 = s_0$), since otherwise there is a plane in $\mathcal{P}_{s_0-1}$ with the same support as $\mathbf{P}$ and thus by (H2), \eqref{separation2 concl} holds true.  

By Theorem~\ref{separation thm1} there exists $\mathbf{Q} \in \mathcal{P}_{q_0}$ such that 
\begin{equation}\label{separation3 eqn1}
	\sup_{X \in \op{spt} T \cap \mathbf{C}_{(1+\gamma)/2}(0)} \op{dist}(X, \op{spt} \mathbf{Q}) \leq C_0 {\mathcal E}(T,\mathbf{C}_1(0)) < C_0 \varepsilon_0 
\end{equation}
for some constant $C_0 = C_0(n,m,q_0,\gamma) \in (0,\infty)$ (which we require to be large enough that we can apply Remark~\ref{separation thm1 rmk} below).  Express 
\begin{equation*}
	\{ X \in \mathbb{R}^{n+m} : \op{dist}(X,\op{spt} \mathbf{Q}) \leq C_0 {\mathcal E}(T,\mathbf{C}_1(0)) \} = \bigcup_{i=1}^N \mathbb{R}^n \times U_i
\end{equation*}
for some collection of mutually disjoint connected open subsets $\{U_i\}$ of $\mathbb{R}^m$, where $C_0$ is as in \eqref{separation3 eqn1}.  By \eqref{separation3 eqn1}, 
\begin{equation*}
	T = \sum_{i=1}^N T_i \quad\text{where}\quad T_i = T \llcorner B_{(1+\gamma)/2}(0) \times U_i . 
\end{equation*}
Clearly $T_i$ are locally area minimizing rectifiable currents of $\mathbf{C}_{(1+\gamma)/2}(0)$ such that $(\partial T_i) \llcorner \mathbf{C}_{\gamma}(0) = 0$.  By \eqref{separation3 hyp1}, the constancy theorem, and Lemma~\ref{projection lemma2} 
\begin{equation}\label{separation3 eqn2}
	\pi_{\#} T_i = q_i \llbracket B_{\gamma}(0) \rrbracket
\end{equation}
for some integers $q_i \geq 0$ such that $q_0 = \sum_{i=1}^N q_i$.  By Remark~\ref{separation thm1 rmk}, $q_i > 0$ for all $i$.  Let $Q_i$ be a plane of $\mathbf{Q}$ and $Z \in Q_i \cap \mathbf{C}_{\gamma}(0)$.  We know that $Q_i \subset \mathbb{R}^n \times U_{j(i)}$ for some $j(i) \in \{1,2,\ldots,N\}$ and thus for each $X \in \op{spt} T_{j(i)} \cap (\{\pi(Z)\} \times \mathbb{R}^m)$ we have $|X-Z| \leq 2 q_0 C_0 {\mathcal E}(T,\mathbf{C}_1(0))$.  Hence by \eqref{separation3 hausdist1} (with $P = Q_i$) and the triangle inequality 
\begin{align}\label{separation3 eqn3}
	\sup_{Y \in Q_i \cap \mathbf{C}_{\gamma}(0)} \op{dist}(Y,\op{spt}\mathbf{P}) 
	&\leq 2 \op{dist}(Z,\op{spt}\mathbf{P}) + 2 \max_{1 \leq k \leq s_0} |A_k| 
	\\&\leq 2 \op{dist}(X,\op{spt}\mathbf{P}) + 2 |X - Z| + 2 \max_{1 \leq k \leq s_0} |A_k| \nonumber
	\\&\leq 2 \op{dist}(X,\op{spt}\mathbf{P}) + 4q_0 C_0 {\mathcal E}(T,\mathbf{C}_1(0)) + 2 \max_{1 \leq k \leq s_0} |A_k| . \nonumber 
\end{align}
By integrating \eqref{separation3 eqn3} over $T_{j(i)}$ and using \eqref{separation3 eqn2}, 
\begin{align*}
	\sup_{Y \in Q_i \cap \mathbf{C}_{\gamma}(0)} \op{dist}^2(Y,\op{spt}\mathbf{P}) 
	\leq\,& \frac{3}{\omega_n \gamma^n} \int_{\mathbf{C}_1(0)} \op{dist}^2(X,\op{spt}\mathbf{P}) \,d\|T_{j(i)}\|(X) 
		\\&+ 48q_0^2 C_0^2 {\mathcal E}(T,\mathbf{C}_1(0))^2 + 12 \max_{1 \leq k \leq s_0} |A_k|^2 
\end{align*}
for all $i$ and thus 
\begin{align}\label{separation3 eqn4}
	\sup_{Y \in \op{spt}\mathbf{Q} \cap \mathbf{C}_{\gamma}(0)} \op{dist}^2(Y,\op{spt}\mathbf{P}) 
	\leq\,& \frac{3}{\omega_n \gamma^n} \int_{\mathbf{C}_1(0)} \op{dist}^2(X,\op{spt}\mathbf{P}) \,d\|T_{j(i)}\|(X) 
		\\&+ 48q_0^2 C_0^2 {\mathcal E}(T,\mathbf{C}_1(0))^2 + 12 \max_{1 \leq k \leq s_0} |A_k|^2 . \nonumber 
\end{align}
By \eqref{separation3 eqn1} and \eqref{separation3 eqn4}, 
\begin{align}\label{separation3 eqn5}
	&\sup_{X \in \op{spt} T \cap \mathbf{C}_{\gamma}(0)} \op{dist}^2(X, \op{spt} \mathbf{P}) 
	\\ \leq\,& 2 \sup_{X \in \op{spt} T \cap \mathbf{C}_{\gamma}(0)} \op{dist}^2(X, \op{spt} \mathbf{Q}) 
		+ 2 \sup_{Y \in \op{spt} \mathbf{Q} \cap \mathbf{C}_{\gamma}(0)} \op{dist}^2(Y, \op{spt} \mathbf{P}) \nonumber
	\\ \leq\,& \frac{6}{q_0 \omega_n \gamma^n} \int_{\mathbf{C}_1(0)} \op{dist}^2(X,\op{spt}\mathbf{P}) \,d\|T\|(X) 
		+ 98q_0^2 C_0^2 {\mathcal E}(T,\mathbf{C}_1(0))^2 + 24 \max_{1 \leq k \leq s_0} |A_k|^2 \nonumber
\end{align}
Thus if 
\begin{equation*}
	\frac{1}{\omega_n} \int_{\mathbf{C}_1(0)} \op{dist}^2(X,\op{spt}\mathbf{P}) \,d\|T\|(X) 
	\geq {\mathcal E}(T,\mathbf{C}_1(0))^2 + \max_{1 \leq k \leq s_0} |A_k|^2  
\end{equation*}
then it follows from \eqref{separation3 eqn5} that \eqref{separation3 concl} holds true.  Hence we may assume that 
\begin{equation}\label{separation3 eqn6}
	\frac{1}{\omega_n} \int_{\mathbf{C}_1(0)} \op{dist}^2(X,\op{spt}\mathbf{P}) \,d\|T\|(X) 
	\leq {\mathcal E}(T,\mathbf{C}_1(0))^2 + \max_{1 \leq k \leq s_0} |A_k|^2 < 5\varepsilon_0^2 . 
\end{equation}

Provided $\varepsilon_0$ is sufficiently small, by Lemma~\ref{sepmono lemma} and \eqref{separation3 eqn6}, 
\begin{equation}\label{separation3 eqn7}
	\sup_{X \in \op{spt} T \cap \mathbf{C}_{(3+\gamma)/4}(0)} \op{dist}(X,\op{spt}\mathbf{P}) < 2 (5\varepsilon_0^2)^{\frac{1}{n+2}} 
	< 10\varepsilon_0^{\frac{2}{n+2}}
\end{equation}
Express 
\begin{equation*}
	\{ X \in \mathbf{C}_{(3+\gamma)/4}(0) : \op{dist}(X,\op{spt} \mathbf{P}) < 20\varepsilon_0^{\frac{2}{n+2}} \} = \bigcup_{i=1}^{\widehat{N}} \widehat{U}_i
\end{equation*}
for some collection of mutually disjoint connected open subsets $\{\widehat{U}_i\}$ of $\mathbf{C}_{(3+\gamma)/4}(0)$.  By \eqref{separation3 eqn7}, 
\begin{equation*}
	T = \sum_{i=1}^N \widehat{T}_i \quad\text{where}\quad \widehat{T}_i = \widehat{T} \llcorner \widehat{U}_i . 
\end{equation*}
Clearly $(\partial \widehat{T}_i) \llcorner \mathbf{C}_{(3+\gamma)/4}(0) = 0$.  By \eqref{separation3 hyp1}, the constancy theorem, and Lemma~\ref{projection lemma2} 
\begin{equation*}
	\pi_{\#} \widehat{T}_i = \widehat{q}_i \llbracket B_{(3+\gamma)/4}(0) \rrbracket
\end{equation*}
for some integers $\widehat{q}_i \geq 0$ such that $q_0 = \sum_{i=1}^N \widehat{q}_i$.  By Lemma~\ref{projection lemma}, whenever $\widehat{q}_i = 0$ we have $\widehat{T}_i \llcorner \mathbf{C}_{(1+\gamma)/2}(0) = 0$.  Hence if $\widehat{q}_i = 0$ for some $i$, then $\mathbf{P} \llcorner (\mathbf{C}_{(1+\gamma)/2}(0) \setminus \widehat{U}_i)$ is a sum-of-planes in $\mathbf{C}_{(1+\gamma)/2}(0)$ consisting of $s_0-1$ or fewer planes.  Hence by (H2) we can apply Theorem~\ref{separation thm3} with $\eta_{0,(1+\gamma)/2\#} T$ and $\eta_{0,(1+\gamma)/2\#} (\mathbf{P} \llcorner (\mathbf{C}_{(1+\gamma)/2}(0) \setminus \widehat{U}_i))$ in place of $T$ and $\mathbf{P}$ to obtain \eqref{separation3 concl}.  Moreover, if $\#\{ i : \widehat{q}_i > 0 \} \geq 2$, then by (H1) for each $i$ we can apply Theorem~\ref{separation thm3} with $\eta_{0,(3+\gamma)/4\#} \widehat{T}_i$ and $\eta_{0,(3+\gamma)/4\#} \mathbf{P}$ in place of $T$ and $\mathbf{P}$ to obtain \eqref{separation3 concl}.  Hence we may assume that $\widehat{N} = 1$ and $\widehat{q}_1 = q_0$.  It follows that up to reordering the planes $P_i$ of $\widehat{\mathbf{P}}$, for each $i \in \{2,3,\ldots,s_0-1\}$ there exists $j(i) \in \{1,2,\ldots,i-1\}$ such that 
\begin{equation*} 
	\inf_{X \in P_i \cap \mathbf{C}_{(3+\gamma)/4}(0)} \op{dist}(X,P_{j(i)}) \leq 20\varepsilon_0^{\frac{2}{n+2}}
\end{equation*}
Thus by \eqref{separation3 hausdist2} 
\begin{equation*}
	\op{dist}_{\mathcal{H}}(P_i \cap \mathbf{C}_1(0), P_{j(i)} \cap \mathbf{C}_1(0)) 
	\leq (2+4\kappa) \inf_{X \in P_i \cap \mathbf{C}_{(3+\gamma)/4}(0)} \op{dist}(X,P_{j(i)}) 
	\leq 40 (1+2\kappa) \varepsilon_0^{\frac{2}{n+2}} . 
\end{equation*}
Therefore 
\begin{equation}\label{separation3 eqn9} 
	\op{width} \mathbf{P} \leq 40 q_0 (1+2\kappa) \varepsilon_0^{\frac{2}{n+2}} .
\end{equation}

Express $\mathbf{P} = \sum_{i=1}^{s_0} \llbracket P_i \rrbracket$ and assume that 
\begin{equation*}
	\inf_{X \in P_{s_0} \cap \mathbf{C}_1(0)} \op{dist}(X,P_{s_0-1}) = \op{minsep} \mathbf{P}.
\end{equation*}
By \eqref{separation3 hausdist2} 
\begin{equation*}
	\op{dist}_{\mathcal H}(P_{s_0-1} \cap \mathbf{C}_2(0), P_{s_0} \cap \mathbf{C}_2(0)) 
	\leq (2+8\kappa) \inf_{X \in P_{s_0} \cap \mathbf{C}_1(0)} \op{dist}(X,P_{s_0-1}) 
	= (2+8\kappa) \op{minsep} \mathbf{P} . 
\end{equation*}
Set $\widetilde{\mathbf{P}} = \sum_{i=1}^{s_0-1} \llbracket P_i \rrbracket \in \Pi_{s_0-1}$ so that 
\begin{gather}
	\label{separation3 eqn10} \op{spt} \widetilde{\mathbf{P}} \subset \op{spt} \mathbf{P}, \\ 
	\label{separation3 eqn11} \op{dist}_{\mathcal H}(\op{spt}\widetilde{\mathbf{P}} \cap \mathbf{C}_2(0), \op{spt}\mathbf{P} \cap \mathbf{C}_2(0))
		\leq (2+8\kappa) \op{minsep} \mathbf{P} . 
\end{gather}
Set 
\begin{equation*}
	\widetilde{H} = \left( \int_{\mathbf{C}_1(0)} \op{dist}^2(X,\op{spt} \widetilde{\mathbf{P}}) \,d\|T\|(X) \right)^{1/2} . 
\end{equation*}
If $\widetilde{H} = 0$ then $\op{spt} T \subset \op{spt} \widetilde{\mathbf{P}} \subset \op{spt} \mathbf{P}$ and we have nothing further to prove.  Thus we may assume that $\widetilde{H} > 0$. 

Let $\beta = \beta(n,m,q_0,s_0,\kappa,\gamma) \in (0,1)$ to be later determined.  If  
\begin{equation*}
	\int_{\mathbf{C}_1(0)} \op{dist}^2(X, \op{spt} \mathbf{P}) \,d\|T\|(X) 
		> \beta^2 \int_{\mathbf{C}_1(0)} \op{dist}^2(X, \op{spt} \widetilde{\mathbf{P}}) \,d\|T\|(X) ,
\end{equation*}
then by (H2) we can apply Theorem~\ref{separation thm3} together with \eqref{separation3 eqn10} to obtain 
\begin{align*}
	&\sup_{X \in \op{spt} T \cap \mathbf{C}_{\gamma}(0)} \op{dist}^2(X, \op{spt} \mathbf{P})
	\leq \sup_{X \in \op{spt} T \cap \mathbf{C}_{\gamma}(0)} \op{dist}^2(X, \op{spt} \widetilde{\mathbf{P}})
	\\ \leq\,& C \int_{\mathbf{C}_1(0)} \op{dist}^2(X,\op{spt} \widetilde{\mathbf{P}}) \,d\|T\|(X) 
	\leq \frac{C}{\beta^2} \int_{\mathbf{C}_1(0)} \op{dist}^2(X,\op{spt} \mathbf{P}) \,d\|T\|(X) ,
\end{align*}
for some constant $C = C(n,m,q_0,s_0,\gamma) \in (0,\infty)$, proving \eqref{separation3 concl}.  Hence for the remainder of the proof we may assume that 
\begin{equation}\label{separation3 eqn12}
	\int_{\mathbf{C}_1(0)} \op{dist}^2(X, \op{spt} \mathbf{P}) \,d\|T\|(X) 
		\leq \beta^2 \int_{\mathbf{C}_1(0)} \op{dist}^2(X, \op{spt} \widetilde{\mathbf{P}}) \,d\|T\|(X) . 
\end{equation}
By the triangle inequality, \eqref{projection eqn2}, \eqref{separation3 eqn11}, and \eqref{separation3 eqn12} 
\begin{align*}
	&\int_{\mathbf{C}_1(0)} \op{dist}^2(X, \op{spt} \widetilde{\mathbf{P}}) \,d\|T\|(X) 
	\\ &\hspace{.5in}\leq 2 \int_{\mathbf{C}_1(0)} \op{dist}^2(X, \op{spt} \mathbf{P}) \,d\|T\|(X) 
		+ 2 (q_0+1) \omega_n (2+8\kappa)^2 (\op{minsep} \mathbf{P})^2
	\\ &\hspace{1in}\leq  2\beta^2 \int_{\mathbf{C}_1(0)} \op{dist}^2(X, \op{spt} \widetilde{\mathbf{P}}) \,d\|T\|(X) 
		+ 2 (q_0+1) \omega_n (2+8\kappa)^2 (\op{minsep} \mathbf{P})^2 ,
\end{align*}
so taking $\beta < 1/2$ we have that 
\begin{equation}\label{separation3 eqn14}
	\widetilde{H}^2 = \int_{\mathbf{C}_1(0)} \op{dist}^2(X,\op{spt}\widetilde{\mathbf{P}}) \,d\|T\|(X) 
		\leq 4 (q_0+1) \omega_n (2+8\kappa)^2 (\op{minsep} \mathbf{P})^2 . 
\end{equation}
In particular, by \eqref{separation3 eqn12} and \eqref{separation3 eqn14} 
\begin{equation}\label{separation3 eqn15}
	\int_{\mathbf{C}_1(0)} \op{dist}^2(X, \op{spt} \mathbf{P}) \,d\|T\|(X) \leq 4 (q_0+1) \omega_n (2+8\kappa)^2 \beta^2 (\op{minsep} \mathbf{P})^2 .
\end{equation}

By (H2) we can apply Theorem~\ref{separation thm3} to obtain  
\begin{equation}\label{separation3 eqn16}
	\sup_{X \in \op{spt} T \cap \mathbf{C}_{(3+\gamma)/4}(0)} \op{dist}(X,\op{spt} \widetilde{\mathbf{P}}) < C_1 \widetilde{H}
\end{equation}
for some constant $C_1 = C_1(n,m,q_0,s_0,\kappa,\gamma) \in (0,\infty)$.  Express 
\begin{equation*}
	\{ X \in \mathbf{C}_{(3+\gamma)/4}(0) : \op{dist}(X,\op{spt} \widetilde{\mathbf{P}}) < C_1 \widetilde{H} \} 
		= \bigcup_{i=1}^{\widetilde{N}} \widetilde{U}_i
\end{equation*}
for some collection of mutually disjoint connected open subsets $\{\widetilde{U}_i\}$ of $\mathbb{R}^{n+m}$, where $C_1$ is as in \eqref{separation3 eqn16}.  By \eqref{separation3 eqn16}, 
\begin{equation*}
	T = \sum_{i=1}^{\widetilde{N}} \widetilde{T}_i \quad\text{where}\quad \widetilde{T}_i = T \llcorner \widetilde{U}_i . 
\end{equation*}
Clearly $(\partial \widetilde{T}_i) \llcorner \mathbf{C}_{(3+\gamma)/4}(0) = 0$.  By \eqref{separation3 hyp1}, the constancy theorem, and Lemma~\ref{projection lemma2} 
\begin{equation*}
	\pi_{\#} \widetilde{T}_i = \widetilde{q}_i \llbracket B_{(3+\gamma)/4}(0) \rrbracket
\end{equation*}
for some integers $\widetilde{q}_i \geq 0$ such that $q_0 = \sum_{i=1}^{\widetilde{N}} \widetilde{q}_i$.  Moreover, by Lemma~\ref{projection lemma}, whenever $\widetilde{q}_i = 0$ we have $\widetilde{T}_i \llcorner \mathbf{C}_{(1+\gamma)/2}(0) = 0$.  Hence if $\widetilde{q}_i = 0$ for some $i$, then $\mathbf{P} \llcorner (\mathbf{C}_{(1+\gamma)/2}(0) \setminus \widetilde{U}_i)$ is a sum-of-planes in $\mathbf{C}_{(1+\gamma)/2}(0)$ consisting of $s_0-1$ or fewer planes.  Hence by (H2) we can apply Theorem~\ref{separation thm3} with $\eta_{0,(1+\gamma)/2\#} T$ and $\eta_{0,(1+\gamma)/2\#} (\mathbf{P} \llcorner (\mathbf{C}_{(1+\gamma)/2}(0) \setminus \widetilde{U}_i))$ in place of $T$ and $\mathbf{P}$ to obtain \eqref{separation3 concl}.  Moreover, if $\#\{ i : \widetilde{q}_i > 0 \} \geq 2$, then by (H1) we can apply Theorem~\ref{separation thm3} with $\eta_{0,(3+\gamma)/4\#} \widetilde{T}_i$ and $\eta_{0,(3+\gamma)/4\#} \mathbf{P}$ in place of $T$ and $\mathbf{P}$ to obtain \eqref{separation3 concl}.  Hence we may assume that $\widetilde{N} = 1$ and $\widetilde{q}_1 = q_0$.  It follows that up to reordering the planes $P_i$ of $\widetilde{\mathbf{P}}$, assuming $s_0 > 2$ for each $i \in \{2,3,\ldots,s_0-1\}$ there exists $j(i) \in \{1,2,\ldots,i-1\}$ such that 
\begin{equation*} 
	\inf_{X \in P_i \cap \mathbf{C}_{(3+\gamma)/4}(0)} \op{dist}(X,P_{j(i)}) \leq 2C_1 \widetilde{H} , 
\end{equation*}
where $C_1$ is as in \eqref{separation3 eqn16}.  Thus by \eqref{separation3 hausdist2}
\begin{equation*} 
	\op{dist}_{\mathcal{H}}(P_i \cap \mathbf{C}_1(0), P_{j(i)} \cap \mathbf{C}_1(0)) 
	\leq (2+4\kappa) \inf_{X \in P_i \cap \mathbf{C}_{(3+\gamma)/4}(0)} \op{dist}(X,P_{j(i)}) 
	\leq 2 (1+2\kappa) C_1 \widetilde{H} .
\end{equation*}
Hence 
\begin{equation}
	\label{separation3 eqn18} \op{width} \widetilde{\mathbf{P}} \leq 4q_0 (1+2\kappa) C_1 \widetilde{H} ,
\end{equation}
where $C_1$ is as in \eqref{separation3 eqn16}.  (Note that in the case $s_0 = 2$, $\widetilde{\mathbf{P}}$ has exactly one plane and thus $\op{width} \widetilde{\mathbf{P}} = 0$.)  Hence by the definition of $\widetilde{\mathbf{P}}$, \eqref{separation3 eqn11}, \eqref{separation3 eqn14}, and \eqref{separation3 eqn18}
\begin{equation}\label{separation3 eqn19}
	\op{width} \mathbf{P} 
	\leq \op{width} \widetilde{\mathbf{P}} + \op{dist}_{\mathcal H}(\op{spt}\widetilde{\mathbf{P}} \cap \mathbf{C}_1(0), \op{spt}\mathbf{P} \cap \mathbf{C}_1(0))
	\leq C \op{minsep} \mathbf{P}
\end{equation}
for some constant $C = C(n,m,q_0,s_0,\kappa,\gamma) \in (0,\infty)$.  Now provided $\beta$ is sufficiently small, by \eqref{separation3 eqn7}, \eqref{separation3 eqn9}, \eqref{separation3 eqn19}, and \eqref{separation3 eqn15} we can apply Lemma~\ref{separation3_1 lemma} with $\eta_{0,(3+\gamma)/4\#} T$ and $\eta_{0,(3+\gamma)/4\#} \mathbf{P}$ in place of $T$ and $\mathbf{P}$ to prove \eqref{separation3 concl}.
\end{proof}

\section{Estimates for area minimizing currents significantly closer to a union of planes meeting along an $(n-2)$-dimensional subspace than to any single plane}\label{estimates}
\setcounter{equation}{0}

Let $\mathbf{C}$ be an $n$-dimensional rectifiable current of $\mathbb{R}^{n+m}$ whose support is a union of $p$ ($p \geq 2$) distinct $n$-dimensional oriented planes $P_i$ intersecting along $\{0\} \times \mathbb{R}^{n-2}$.  Let $T$ be an $n$-dimensional locally area-minimizing rectifiable current in $\mathbf{B}_1(0)$ with $\partial T \llcorner {\mathbf B}_{1}(0)= 0$ which is weakly close to $\mathbf{C}$.  Assume that $T$ is significantly closer to $\mathbf{C}$ than to any sum-of-planes supported on fewer distinct planes than $p$ (in the sense of Hypothesis~$(\star\star)$ below). In analogy with the situations considered in \cite{Sim93} and \cite{Wic14}, we wish to express $T$, away from the singular axis $\{0\} \times {\mathbb R}^{n-2}$ of ${\mathbf C}$, as the graph of an appropriate function over the planes of $\mathbf{C}$. In ~\cite{Sim93}, the key hypothesis that allows one to do this is that $T$ (which is a stationary varifold not assumed to be area minimizing) belongs to a  ``multiplicity 1 class;''  this makes it possible to apply Allard's regularity theorem away from the singular axis of ${\mathbf C}$ (subject only to the assumption that $T$ is sufficiently weakly close to ${\mathbf C}$, with Hypothesis~$(\star\star)$ being vacuous in that setting). The situation considered in \cite{Wic14} allows higher multiplicity, but still there is a ``sheeting theorem'' applicable which guarantees complete $C^{1, \alpha}$ regularity of $T$  (which is a stable, stationary codimension 1 varifold in that setting) away from the singular axis of ${\mathbf C}$. 

In contrast to either of these settings, in the present circumstances there is no regularity theory applicable to $T$ that would provide complete regularity of $T$ away from the axis $\{0\} \times \mathbb{R}^{n-2}$.  The basic result we use as a substitute for such regularity is our height estimate, Theorem~\ref{separation thm3}. 
We use  Theorem~\ref{separation thm3} to show (in Section~\ref{notation-and-graphical} and Section~\ref{proofs-graphical} below) that away from $\{0\} \times \mathbb{R}^{n-2}$, $T$ separates as the sum of locally area minimizing rectifiable currents $T_i$ such that $T_{i}$ is close to $P_i$ for each $i$.  Applying Almgren's Strong Lipschitz Approximation Theorem (Theorem~\ref{lip approx thm}), we then approximate $T_{i}$ (in Theorem~\ref{graphrep thm}) by the graph of a Lipschitz multi-valued function $u_i$ over an appropriate domain in $P_i$.  

In Sections~\ref{sec:keyest sec}-\ref{sec:apriori nonplanar} we establish 
a number of key estimates, analogous to those in \cite{Sim93}, \cite{Wic14},  for area minimizing currents $T$ satisfying appropriate hypothesis including Hypothesis~($\star\star$). These results will allow us to produce (in Section~\ref{blow-up-analysis}) ``fine blow-ups'' of sequences of area-minimizing currents $(T_k)$ relative to sequences of sums-of-planes $(\mathbf{C}_k)$ (of the type ${\mathbf C}$ as above), with $T_{k}$, ${\mathbf C}_{k}$ satisfying appropriate hypotheses including Hypothesis~($\star\star$), and study the asymptotic behaviour of the fine blow-ups which ultimately leads to the main excess decay result of the present work, Theorem~\ref{excess-improvement}.

\subsection{Notation and statement of graphical representation results}\label{notation-and-graphical}

Given an $n$-dimensional plane $P \subset \mathbb{R}^{n+m}$, we let $\pi_P : \mathbb{R}^{n+m} \rightarrow P$ denote the orthogonal projection map onto $P$ and we let $\pi_{P^{\perp}} : \mathbb{R}^{n+m} \rightarrow P$ denote the orthogonal projection map onto the orthogonal complement $P^{\perp}$.  For each $X_0 \in \mathbb{R}^{n+m}$ and $\rho > 0$ we define 
\begin{align*}
	B_{\rho}(X_0,P) &= \{ X_0 + X : X \in P \text{ with } |X| < \rho \}, \\
	\mathbf{C}_{\rho}(X_0,P) &= \{ X_0 + X + Y : X \in P \text{ with } |X| < \rho \text{ and } Y \in P^{\perp} \} . 
\end{align*}
We define the classes $\mathcal{C}_{q,p}$ of sums-of-planes as follows: 

\begin{definition} \label{sums of planes defn} {\rm 
Given integers $p$, $q$ with $1 \leq p \leq q$, let $\mathcal{C}_{q,p}$ be the set of all $n$-dimensional rectifiable currents $\mathbf{C}$ of $\mathbb{R}^{m+n}$ of the form 
\begin{equation}\label{cone sums of planes form}
	\mathbf{C} = \sum_{i=1}^p q_i \llbracket P_i \rrbracket
\end{equation}
where $q_1,\ldots,q_p$ are positive integers such that $\sum_{i=1}^p q_i = q$ and $P_i$ are distinct $n$-dimensional planes such that $\{0\} \times \mathbb{R}^{n-2} \subset P_1$ if $p=1$, and $P_{i} \cap P_{j} = \{0\} \times {\mathbb R}^{n-2}$ for $i \neq j$ if $p \geq 2$. Each plane $P_i$ is oriented by the unit simple $n$-vector denoted $\vec P_i$. 
} \end{definition}

\begin{remark} {\rm 
We do not assume that every $\mathbf{C} \in \mathcal{C}_{q,p}$ is area minimizing. 
} \end{remark}

For $1 < p \leq q$ and $\mathbf{C} = \sum_{i=1}^p q_i \llbracket P_i \rrbracket \in \mathcal{C}_{q,p},$ we let 
\begin{align}\label{sums of planes sep}
	\op{minsep} \mathbf{C} = \min_{1 \leq i \leq p} \,\inf_{X \in P_i \cap (\mathbb{S}^{m+1} \times \mathbb{R}^{n-2})} 
		\op{dist}(X, \op{spt} \mathbf{C} \setminus P_i), \\
	\op{maxsep} \mathbf{C} = \min_{1 \leq i \leq p} \,\sup_{X \in P_i \cap (\mathbb{S}^{m+1} \times \mathbb{R}^{n-2})} 
		\op{dist}(X, \op{spt} \mathbf{C} \setminus P_i). \nonumber 
\end{align}
If $\mathbf{C} \in \mathcal{C}_{q,1}$, we define $\op{minsep} \mathbf{C} = \infty$ and $\op{maxsep} \mathbf{C} = \infty$.  $\op{minsep} \mathbf{C}$ is proportional to the least distance between a pair of points on different planes of $\mathbf{C}$ in $\mathbb{S}^{m+1} \times \mathbb{R}^{n-2}$ and quantifies how close any pair of planes of $\mathbf{C}$ are to intersecting away from $\{0\} \times \mathbb{R}^{n-2}$.  $\op{maxsep} \mathbf{C}$ is proportional to the least Hausdorff distance between a pair of distinct planes of $\mathbf{C}$ in $\mathbb{S}^{m+1} \times \mathbb{R}^{n-2}$ and quantifies how close any two or more planes of $\mathbf{C}$ are to coinciding.  Clearly $\op{minsep} \mathbf{C} \leq \op{maxsep} \mathbf{C}$.  However, in contrast with~\cite{Wic14}, which considered half-planes meeting along an $(n-1)$-dimensional linear subspace, in the present setting $\op{minsep} \mathbf{C}$ and $\op{maxsep} \mathbf{C}$ need not be equal.

Throughout the paper, we shall use the following notation:  for each $\rho > 0$, $n$-dimensional rectifiable current $T$ of $\mathbf{B}_{\rho}(0)$, and $\mathbf{C} \in \bigcup_{p'=1}^q \mathcal{C}_{q,p'}$, we define 
\begin{align*}
	E(T,\mathbf{C},\mathbf{B}_{\rho}(0)) = \Bigg( &\frac{1}{\omega_n \rho^{n+2}} \int_{\mathbf{B}_{\rho}(0)} 
		\op{dist}^2(X, \op{spt} \mathbf{C}) \,d\|T\|(X) \Bigg)^{1/2}, \\
	Q(T,\mathbf{C},\mathbf{B}_{\rho}(0)) = \Bigg( &\frac{1}{\omega_n \rho^{n+2}} \int_{\mathbf{B}_{\rho}(0)} \op{dist}^2(X, \op{spt} \mathbf{C}) \,d\|T\|(X) 
		\\& + \frac{1}{\omega_n \rho^{n+2}} \int_{\mathbf{B}_{\rho/2}(0) \cap \{ r > \rho/16 \}} \op{dist}^2(X, \op{spt} T) \,d\|\mathbf{C}\|(X) \Bigg)^{1/2} , 
\end{align*}
where $r = r(X) = \op{dist}(X, \{0\} \times \mathbb{R}^{n-2})$ for each $X \in \mathbb{R}^n$.  In Theorem~\ref{graphrep thm} and a number of other results in subsequent sections, we shall assume the first or both of the following hypotheses for appropriate choices of small constants $\varepsilon_0 \in (0,1)$ and $\beta_0 \in (0,1)$: 

\noindent\textbf{Hypothesis~$(\star)$.}  $2 \leq p \leq q$ are integers, $\mathbf{C} = \sum_{i=1}^p q_i \llbracket P_i \rrbracket \in \mathcal{C}_{q,p}$, and $T$ is an $n$-dimensional locally area-minimizing rectifiable current in $\mathbf{B}_1(0)$ such that 
\begin{gather}
	\label{main hyp eqn1} (\partial T) \llcorner \mathbf{B}_1(0) = 0, \quad \Theta(T,0) \geq q, \quad \|T\|(\mathbf{B}_1(0)) \leq (q + 1/2) \,\omega_n , \\ 
	\label{main hyp eqn2} E(T, \mathbf{C},\mathbf{B}_1(0)) < \varepsilon_0 . 
\end{gather}

\noindent\textbf{Hypothesis~$(\star\star)$.}  $2 \leq p \leq q$ are integers, $\mathbf{C} \in \mathcal{C}_{q,p}$, and $T$ is an $n$-dimensional locally area-minimizing rectifiable current in $\mathbf{B}_1(0)$ such that 
\begin{equation}\label{main hyp eqn3} 
	Q(T,\mathbf{C},\mathbf{B}_1(0)) \leq \beta_0 \inf_{\mathbf{C}' \in \bigcup_{p'=1}^{p-1} \mathcal{C}_{q,p'}} Q(T,\mathbf{C}',\mathbf{B}_1(0)) . 
\end{equation}

\begin{remark}\label{cone-multiplicity}{\rm 
Suppose that $\mathbf{C} \in \mathcal{C}_{q,p}$ and $T$ satisfy Hypothesis~$(\star)$ and Hypothesis~$(\star\star)$ for some $\varepsilon_0,\beta_0 \in (0,1)$.  If $\mathbf{C}^* \in \mathcal{C}_{q,p}$ is any other cone with $\op{spt} \mathbf{C}^* = \op{spt} \mathbf{C}$, then Hypothesis~$(\star)$ and Hypothesis~$(\star\star)$ continue to be satisfied with $\mathbf{C}^*$, $\sqrt{q}\beta_0$ in place of $\mathbf{C}$, $\beta_0$ respectively.
}\end{remark}

\begin{theorem} \label{graphrep thm}
Given integers $2 \leq p \leq q$ and $0 < \tau < \gamma < 1$ there exists $\varepsilon_0 = \varepsilon_0(n,m,q,p,\gamma,\tau) \in (0,1)$ and $\beta_0 = \beta_0(n,m,q,p,\gamma,\tau) \in (0,1)$ such that if $\mathbf{C}$ and $T$ satisfy Hypothesis~$(\star)$ and Hypothesis~$(\star\star)$, then: 
\begin{enumerate} \setlength{\itemsep}{5pt}
	\item[(a)] $\mathbf{C}$ satisfies 
	\begin{gather}
		\label{graphrep concl a1} \op{minsep} \mathbf{C} 
			\geq c \inf_{\mathbf{C}' \in \bigcup_{p'=1}^{p-1} \mathcal{C}_{q,p'}} Q(T,\mathbf{C}',\mathbf{B}_1(0)) , \\
		\label{graphrep concl a2} \op{dist}_{\mathcal{H}}(P_i \cap \mathbf{B}_1(0),P_j \cap \mathbf{B}_1(0)) 
			\leq C \inf_{X \in P_i \cap (\mathbb{S}^{m+1} \times \mathbb{R}^{n-2})} \op{dist}(X,P_j) \text{ for all $1 \leq i,j \leq p$,}
	\end{gather}
	where $c = c(n,m,q,p) > 0$ and $C = C(n,m,q,p) \in (0,\infty)$ are constants; 
	
	\item[(b)] after replacing $\mathbf{C}$ with a cone with the same support (also denoted by $\mathbf{C}$) there exist $n$-dimensional locally area minimizing rectifiable currents $T_i$ in $\mathbf{B}_{(3+\gamma)/4}(0) \cap \{r > \tau/4\}$ for which  
	\begin{gather}
		\label{graphrep concl b1} T \llcorner \mathbf{B}_{(3+\gamma)/4}(0) \cap \{r > \tau/4\} = \sum_{i=1}^p T_i , \\
		\label{graphrep concl b2} (\partial T_i) \llcorner \mathbf{B}_{(3+\gamma)/4}(0) \cap \{r > \tau/4\} = 0, \\ 
		\label{graphrep concl b3} (\pi_{P_i \#} T_i) \llcorner \mathbf{B}_{(1+\gamma)/2}(0) \cap \{r > \tau/2\} 
			= q_i \llbracket P_i \rrbracket \llcorner \mathbf{B}_{(1+\gamma)/2}(0) \cap \{r > \tau/2\} , \\
		\label{graphrep concl b4} \sup_{X \in \op{spt} T_i \cap \{r > \sigma\}} \op{dist}(X, P_i) \leq C_{\sigma} E \text{ for all } \sigma \in [\tau/2,1/2] , 
	\end{gather}
	where $r(X) = \op{dist}(X, \{0\} \times \mathbb{R}^{n-2})$, $E = E(T,\mathbf{C},\mathbf{B}_1(0))$, and $C_{\sigma} = C_{\sigma}(n,m,q,p,\gamma,\sigma) \in (0,\infty)$ are constants; 
	
	\item[(c)] for each $i \in \{1,2,\ldots,p\}$ there exists Lipschitz $q_i$-valued functions $u_i : B_{\gamma}(0,P_i) \cap \{r > \tau\} \rightarrow \mathcal{A}_{q_i}(P_i^{\perp})$ and closed sets $K_i \subseteq B_{\gamma}(0,P_i) \cap \{r > \tau\}$ such that 
	\begin{gather}
		\label{graphrep concl c1} T_i \llcorner \pi_{P_i}^{-1}(K_i) = (\op{graph} u_i) \llcorner \pi_{P_i}^{-1}(K_i) , \\
		\label{graphrep concl c2} \mathcal{H}^n(B_{\gamma}(0,P_i) \cap \{r > \sigma\} \setminus K_i)
			+ \|T_i\|(\pi_{P_i}^{-1}(B_{\gamma}(0,P_i) \cap \{r > \sigma\} \setminus K_i)) \leq C_{\sigma} E^{2+\alpha} , \\
		\label{graphrep concl c3} \sup_{B_{\gamma}(0,P_i) \cap \{r > \sigma\}} |u_i| \leq C_{\sigma} E , \quad 
			\sup_{B_{\gamma}(0,P_i) \cap \{r > \sigma\}} |\nabla u_i| \leq C_{\sigma} E^{\alpha}
	\end{gather}
	for all $\sigma \in [\tau,1/2]$, where again $E = E(T,\mathbf{C},\mathbf{B}_1(0))$ and $\alpha = \alpha(n,m,q) \in (0,1)$, $C_{\sigma} = C_{\sigma}(n,m,q,p,\gamma,\sigma) \in (0,\infty)$ are constants.
\end{enumerate}
\end{theorem}

Theorem~\ref{graphrep thm} will follow from Lemma~\ref{weak graphrep lemma} and Lemma~\ref{graphrep annuli lemma}, which will provide a local graphical representation of locally area-minimizing rectifiable currents $T$ in annuli.  Lemma~\ref{weak graphrep lemma} and Lemma~\ref{graphrep annuli lemma} will also play an important role in the proof of Theorem~\ref{keyest thm} in the next section.  Here and subsequently, we shall use the following notation:  for each $\gamma \in (0,1)$, $\zeta \in \mathbb{R}^{n-2}$, $\rho > 0$, and $\kappa \in (0,2]$, we let $A_{\rho,\kappa}(\zeta) \subset \mathbb{R}^n$ and $\mathbf{A}_{\rho,\kappa}(\zeta) \subset \mathbb{R}^{n+m}$ be annuli given by  
\begin{align*} 
	A_{\rho,\kappa}(\zeta) &= \{ (x,y) \in \mathbb{R}^2 \times \mathbb{R}^{n-2} : 
		(|x| - \rho)^2 + |y - \zeta|^2 < \kappa^2 (1-\gamma)^2 \rho^2 /64 \} , \\
	\mathbf{A}_{\rho,\kappa}(\zeta) &= \{ (x,y) \in \mathbb{R}^{2+m} \times \mathbb{R}^{n-2} : 
		(|x| - \rho)^2 + |y - \zeta|^2 < \kappa^2 (1-\gamma)^2 \rho^2 /64 \} .
\end{align*}
For each $n$-dimensional rectifiable current $T$ of $\mathbf{A}_{\rho,\kappa}(\zeta)$ and each $\mathbf{C} \in \bigcup_{p=1}^q \mathcal{C}_{q,p}$, we define 
\begin{align*}
	E(T,\mathbf{C},\mathbf{A}_{\rho,\kappa}(\zeta)) = \Bigg( &\frac{1}{\rho^{n+2}} \int_{\mathbf{A}_{\rho,\kappa}(\zeta)} 
		\op{dist}^2(X,  \op{spt} \mathbf{C}) \,d\|T\|(X) \Bigg)^{1/2}, \\
	Q(T,\mathbf{C},\mathbf{A}_{\rho,\kappa}(\zeta)) = \Bigg( &\frac{1}{\rho^{n+2}} \int_{\mathbf{A}_{\rho,\kappa}(\zeta)} 
		\op{dist}^2(X, \op{spt} \mathbf{C}) \,d\|T\|(X) \\& + \frac{1}{\rho^{n+2}} \int_{\mathbf{A}_{\rho,\kappa/2}(0)} 
		\op{dist}^2(X, \op{spt} T) \,d\|\mathbf{C}\|(X) \Bigg)^{1/2} . 
\end{align*}

\begin{lemma} \label{weak graphrep lemma}
Given integers $1 \leq p \leq q$ and $\gamma, \kappa \in (0,1)$ and $\mu \in (1,\infty)$ there exists $\overline{\varepsilon} = \overline{\varepsilon}(n,m,q,p,\gamma,\kappa,\mu) \in (0,1)$ and $\overline{\beta} = \overline{\beta}(n,m,q,p,\gamma,\kappa,\mu) \in (0,1)$ such that the following holds true.  Let $\mathbf{C} = \sum_{i=1}^p q_i \llbracket P_i \rrbracket \in \mathcal{C}_{q,p}$  and $T$ be an $n$-dimensional locally area-minimizing rectifiable current in $\mathbf{A}_{1,1}(0)$ such that 
\begin{gather}
	\label{weak graphrep hyp1} (\partial T) \llcorner \mathbf{A}_{1,1}(0) = 0, \\ 
	\label{weak graphrep hyp3} \|T\|(\mathbf{A}_{1,1}(0)) \leq (q+1/2) \,\mathcal{L}^n(A_{1,1}(0)), \\ 
	\label{weak graphrep hyp4} E(T,\mathbf{C},\mathbf{A}_{1,1}(0)) < \overline{\varepsilon} , 
\end{gather}
and either: 
\begin{enumerate}
	\item[(i)]  $p = 1$ or 
	\item[(ii)]  $p > 1,$  
	\begin{gather}
		\label{weak graphrep hyp5} E(T,\mathbf{C},\mathbf{A}_{1,1}(0)) 
			\leq \overline{\beta} \inf_{\mathbf{C}' \in \bigcup_{p'=1}^{p-1} \mathcal{C}_{q,p'}} Q(T,\mathbf{C}',\mathbf{A}_{1,1}(0)),\;\; \mbox{and} \\
		\label{weak graphrep hyp6} \op{dist}_{\mathcal{H}}(P_i \cap \mathbf{B}_1(0),P_j \cap \mathbf{B}_1(0)) 
			\leq \mu \inf_{X \in P_i \cap (\mathbb{S}^{m+1} \times \mathbb{R}^{n-2})} \op{dist}(X,P_j) \text{ for all $i \neq j$}.
	\end{gather}
\end{enumerate}
Then:
\begin{enumerate} \setlength{\itemsep}{5pt}
	\item[(a)]  up to reversing the orientation of $P_i$, there exist (possibly zero) integers $\widehat{q}_i \geq 0$ with $\sum_{i=1}^p \widehat{q}_i \leq q$ and $n$-dimensional locally area minimizing rectifiable currents $T_i$ in $\mathbf{A}_{1,(3+\kappa)/4}(0)$ such that  
\begin{gather}
	\label{graphrep annuli concl b1} T \llcorner \mathbf{A}_{1,(3+\kappa)/4}(0) = \sum_{i=1}^p T_i , \\ 
	\label{graphrep annuli concl b2} (\partial T_i) \llcorner \mathbf{A}_{1,(3+\kappa)/4}(0) = 0, \\ 
	\label{graphrep annuli concl b3} (\pi_{P_i \#} T_i) \llcorner \mathbf{A}_{1,(1+\kappa)/2}(0) 
		= \widehat{q}_i \llbracket P_i \rrbracket \llcorner \mathbf{A}_{1,(1+\kappa)/2}(0) , \\
	\label{graphrep annuli concl b4} \sup_{X \in \op{spt} T_i} \op{dist}(X, P_i) \leq C E , 
\end{gather}
where $E = E(T,\mathbf{C},\mathbf{A}_{1,1}(0))$ and $C = C(n,m,q,p,\gamma,\kappa, \mu) \in (0,\infty)$ is a constant;

	\item[(b)] for each $i \in \{1,2,\ldots,p\}$ with $\widehat{q}_i > 0,$ there exist a Lipschitz $\widehat{q}_i$-valued function $u_i : P_i \cap \mathbf{A}_{1,\kappa}(0) \rightarrow \mathcal{A}_{\widehat{q}_i}(P_i^{\perp})$ and a closed set $K_i \subseteq P_i \cap \mathbf{A}_{1,\kappa}(0)$ such that 
	\begin{gather}\label{graphrep annuli concl c} 
		T_i \llcorner \pi_{P_i}^{-1}(K_i) = (\op{graph} u_i) \llcorner \pi_{P_i}^{-1}(K_i) , \\
		\mathcal{H}^n(P_i \cap \mathbf{A}_{1,\kappa}(0) \setminus K_i)
			+ \|T_i\|(\pi_{P_i}^{-1}(P_i \cap \mathbf{A}_{1,\kappa}(0) \setminus K_i)) \leq C E^{2+\alpha} , \nonumber \\
		\sup_{P_i \cap \mathbf{A}_{1,\kappa}(0)} |u_i| \leq C E , \quad 
			\sup_{P_i \cap \mathbf{A}_{1,\kappa}(0)} |\nabla u_i| \leq C E^{\alpha} , \nonumber
	\end{gather}
	where again $E = E(T,\mathbf{C},\mathbf{A}_{1,1}(0))$ and $\alpha = \alpha(n,m,q) \in (0,1),$ $C = C(n,m,q,p,\gamma,\kappa,\mu) \in (0,\infty)$ are constants.
\end{enumerate}
\end{lemma}

In Lemma~\ref{graphrep annuli lemma}, we replace hypothesis \eqref{weak graphrep hyp5} of Lemma~\ref{weak graphrep lemma} with the stronger assumption \eqref{graphrep annuli hyp5}, and prove that \eqref{weak graphrep hyp6} in fact follows as a conclusion, and that the integers $\widehat{q}_{i}$ (as the conclusion of Lemma~\ref{weak graphrep lemma}) are all positive. 

\begin{lemma} \label{graphrep annuli lemma}
Given integers $1 \leq p \leq q$ and $\gamma, \kappa \in (0,1)$ there exists $\overline{\varepsilon} = \overline{\varepsilon}(n,m,q,p,\gamma,\kappa) \in (0,1)$ and $\overline{\beta} = \overline{\beta}(n,m,q,p,\gamma,\kappa) \in (0,1)$ such that the following holds true.  Let $\mathbf{C} = \sum_{i=1}^p q_i \llbracket P_i \rrbracket \in \mathcal{C}_{q,p}$  and $T$ be an $n$-dimensional locally area-minimizing rectifiable current in $\mathbf{A}_{1,1}(0)$ such that $\op{spt} T \cap \mathbf{A}_{1,1/2}(0) \neq \emptyset$ and \eqref{weak graphrep hyp1}, \eqref{weak graphrep hyp3}, and \eqref{weak graphrep hyp4} hold true. Suppose also that either 
\begin{enumerate}
	\item[(i)]  $p = 1$ or 
	\item[(ii)]  $p > 1$ and 
	\begin{equation}\label{graphrep annuli hyp5} 
		Q(T,\mathbf{C},\mathbf{A}_{1,1}(0)) 
			\leq \overline{\beta} \inf_{\mathbf{C}' \in \bigcup_{p'=1}^{p-1} \mathcal{C}_{q,p'}} Q(T,\mathbf{C}',\mathbf{A}_{1,1}(0)) . 
	\end{equation}
\end{enumerate}
Then:
\begin{enumerate} \setlength{\itemsep}{5pt}
	\item[(a)]  when $p > 1$, 
	\begin{gather}
		\label{graphrep annuli concl a1} \op{minsep} \mathbf{C} 
			\geq c \inf_{\mathbf{C}' \in \bigcup_{p'=1}^{p-1} \mathcal{C}_{q,p'}} Q(T,\mathbf{C}',\mathbf{A}_{1,1}(0)) , \\
		\label{graphrep annuli concl a2} \op{dist}_{\mathcal{H}}(P_i \cap \mathbf{B}_1(0),P_j \cap \mathbf{B}_1(0)) 
			\leq C \inf_{X \in P_i \cap (\mathbb{S}^{m+1} \times \mathbb{R}^{n-2})} \op{dist}(X,P_j) \text{ for all $i \neq j$,}
	\end{gather}
	where $c = c(n,m,q,p,\gamma) > 0$ and $C = C(n,m,q,p,\gamma) \in (0,\infty)$ are constants;
	
	\item[(b)]  up to reversing the orientation of $P_i$, there exist (non-zero) integers $\widehat{q}_i > 0$ with $\sum_{i=1}^p \widehat{q}_i \leq q$ and $n$-dimensional locally area minimizing rectifiable currents $T_i$ of $\mathbf{A}_{1,(3+\kappa)/4}(0)$ such that \eqref{graphrep annuli concl b1}--\eqref{graphrep annuli concl b4} hold true for some constant $C = C(n,m,\gamma,\kappa) \in (0,\infty)$; 

	\item[(c)]  for each $i \in \{1,2,\ldots,p\}$ there exists Lipschitz $\widehat{q}_i$-valued functions $u_i : P_i \cap \mathbf{A}_{1,\kappa}(0) \rightarrow \mathcal{A}_{\widehat{q}_i}(P_i^{\perp})$ and closed sets $K_i \subseteq P_i \cap \mathbf{A}_{1,\kappa}(0)$ such that \eqref{graphrep annuli concl c} holds true for some constants $\alpha = \alpha(n,m,q) \in (0,1)$ and $C = C(n,m,q,\gamma,\kappa) \in (0,\infty)$.
\end{enumerate}
\end{lemma}

\begin{remark} \label{graphrep rmk} 
Let $U = \mathbf{A}_{1,1}(0)$ or $U = \mathbf{B}_1(0)$.  There exists $\beta = \beta(q) \in (0,1)$ such that if $p \in \{2,3,\ldots,q\}$, $\mathbf{C} = \sum_{i=1}^p q_i \llbracket P_i \rrbracket \in \mathcal{C}_{q,p}$, and $T$ is an $n$-dimensional locally area minimizing rectifiable current of $U$ such that 
\begin{equation}\label{graphrep rmk1 eqn1} 
	Q(T, \mathbf{C}, U) \leq \beta \inf_{\mathbf{C}' \in \bigcup_{p'=1}^{p-1} \mathcal{C}_{q,p'}} Q(T, \mathbf{C}', U) 
\end{equation}
(as in \eqref{main hyp eqn3} and \eqref{graphrep annuli hyp5}), then 
\begin{equation}\label{graphrep rmk1 eqn2} 
	\inf_{\mathbf{C}' \in \bigcup_{p'=1}^{p-1} \mathcal{C}_{q,p'}} Q(T, \mathbf{C}', U) \leq C \op{maxsep} \mathbf{C} ,
\end{equation}
where $C = 4 \,\|T\|(U)^{1/2}$.

{\rm Notice that \eqref{graphrep concl a1} and \eqref{graphrep annuli concl a1} are stronger conclusions than \eqref{graphrep rmk1 eqn2}.  To see \eqref{graphrep rmk1 eqn2}, let us consider the case $U = \mathbf{A}_{1,1}(0).$ The case $U = \mathbf{B}_1(0)$ is similar.  Let $i \in \{1,2,\ldots,q\}$ and let $\mathbf{C}'_i \in \mathcal{C}_{q,p-1}$ be such that $\op{spt} \mathbf{C}'_i = \op{spt} \mathbf{C} \setminus P_i$ (e.g.~if $i = 1$ let $\mathbf{C}'_{i} = q_1 \llbracket P_2 \rrbracket + \sum_{j=2}^p q_j \llbracket P_j \rrbracket$).  Let $X \in \op{spt} T \cap \mathbf{A}_{1,1}(0)$ and suppose that the closest point to $X$ in $\op{spt} \mathbf{C}$ lies on $P_i$.  Since $\mathbf{A}_{1,1}(0) \subset \mathbf{B}_2(0)$, the closest point to $X$ in $\op{spt} \mathbf{C}$ lies on $P_i \cap \mathbf{B}_2(0)$.  Thus by the triangle inequality 
\begin{align*}
	\op{dist}(X, \op{spt} \mathbf{C}'_i) 
		&\leq \op{dist}(X, P_i) + \sup_{Y \in P_i \cap \mathbf{B}_2(0)} \op{dist}(Y, \op{spt} \mathbf{C}'_i) 
		\\&= \op{dist}(X, \op{spt} \mathbf{C}) + 2 \sup_{Y \in P_i \cap (\mathbb{S}^{m+1} \times \mathbb{R}^{n-2})} 
			\op{dist}(Y, \op{spt} \mathbf{C} \setminus P_i) .
\end{align*}
If on the other hand $X \in \op{spt} T \cap \mathbf{A}_{1,1}(0)$ and the closest point to $X$ in $\op{spt} \mathbf{C}$ does not lie on $P_i$, then $\op{dist}(X, \op{spt} \mathbf{C}'_i) = \op{dist}(X, \op{spt} \mathbf{C})$.  Hence 
\begin{align}\label{graphrep rmk1 eqn3} 
	\int_{\mathbf{A}_{1,1}(0)} \op{dist}^2(X, \op{spt} \mathbf{C}'_i) \,d\|T\|(X)
	\leq\,& 2 \int_{\mathbf{A}_{1,1}(0)} \op{dist}^2(X, \op{spt} \mathbf{C}) \,d\|T\|(X) \\&+ 8 \,\|T\|(\mathbf{A}_{1,1}(0)) 
		\sup_{X \in P_i \cap (\mathbb{S}^{m+1} \times \mathbb{R}^{n-2})} \op{dist}^2(X, \op{spt} \mathbf{C} \setminus P_i) \nonumber 
\end{align}
Since $\op{spt} \mathbf{C}'_{i} \subset \op{spt} \mathbf{C}$, 
\begin{align}\label{graphrep rmk1 eqn4} 
	&\int_{\mathbf{A}_{1,1/2}(0)} \op{dist}^2(X, \op{spt} T) \,d\|\mathbf{C}'_i\|(X)
	= \sum_{j=1}^p |q'_j| \int_{P_j \cap \mathbf{A}_{1,1/2}(0)} \op{dist}^2(X, \op{spt} T) \,d\mathcal{H}^n(X)
	\\ \leq\,& \sum_{j=1}^p qq_j \int_{P_j \cap \mathbf{A}_{1,1/2}(0)} \op{dist}^2(X, \op{spt} T) \,d\mathcal{H}^n(X)
	= q \int_{\mathbf{A}_{1,1/2}(0)} \op{dist}^2(X, \op{spt} T) \,d\|\mathbf{C}\|(X), \nonumber 
\end{align}
where $\mathbf{C}'_i = \sum_{j=1}^p q'_j \llbracket P_j \rrbracket$ for some integers $q'_j$ with $\sum_{j=1}^p |q'_j| = q$ and we use $|q'_j| \leq q \leq qq_j$ for each $j$.  Adding \eqref{graphrep rmk1 eqn3} and \eqref{graphrep rmk1 eqn4} and then using \eqref{graphrep rmk1 eqn1}, 
\begin{align*}
	&\hspace{-.5in}Q(T,\mathbf{C}'_i,\mathbf{A}_{1,1}(0))^2 
	\\ \leq\,& 2q \,Q(T,\mathbf{C},\mathbf{A}_{1,1}(0))^2 + 8 \,\|T\|(\mathbf{A}_{1,1}(0)) 
		\sup_{X \in P_i \cap (\mathbb{S}^{m+1} \times \mathbb{R}^{n-2})} \op{dist}^2(X, \op{spt} \mathbf{C} \setminus P_i)
	\\ \leq\,& 2q \beta^2 \,Q(T,\mathbf{C}'_i,\mathbf{A}_{1,1}(0))^2 + 8 \,\|T\|(\mathbf{A}_{1,1}(0)) 
		\sup_{X \in P_i \cap (\mathbb{S}^{m+1} \times \mathbb{R}^{n-2})} \op{dist}^2(X, \op{spt} \mathbf{C} \setminus P_i) . \nonumber 
\end{align*}
Therefore, taking $\beta \leq \tfrac{1}{2\sqrt{q}}$ and taking the infimum over all $i \in \{1,2,\ldots,p\}$ gives us \eqref{graphrep rmk1 eqn2}. 
}\end{remark}

\subsection{Proofs of the graphical representation results}\label{proofs-graphical}

\begin{proof}[Proof of Lemma~\ref{weak graphrep lemma}]
Without loss of generality assume that $\kappa \in [\kappa_0,1)$ where $\kappa_0 = \kappa_0(n,m,q,\gamma) \in (0,1)$ such that 
\begin{equation}\label{graphrep annuli kappa hyp}
	\mathcal{L}^n(A_{1,(1+\kappa_0)/2}(0)) \geq \frac{q+1/2}{q+3/4} \,\mathcal{L}^n(A_{1,1}(0)) . 
\end{equation}
Moreover, if $E(T,\mathbf{C},\mathbf{A}_{1,1}(0)) = 0$, then $\op{spt} T \subseteq \op{spt} \mathbf{C}$ and thus the conclusion of Lemma~\ref{weak graphrep lemma} clearly holds true with $T_i = T \llcorner (P_i \cap \mathbf{A}_{1,(3+\kappa)/4}(0))$, $u_i = \widehat{q}_{i}\llbracket 0 \rrbracket$ and $K_i = P_i \cap \mathbf{A}_{1,\kappa}(0)$ where $\widehat{q}_i$ are integers.  Hence we may assume that $E(T,\mathbf{C},\mathbf{A}_{1,1}(0)) > 0$. 

For $k = 1,2,3,\ldots$ let $\varepsilon_k \rightarrow 0^+$, $\beta_k \rightarrow 0^+$, $\mathbf{C}_k \in \mathcal{C}_{q,p}$, and $T_k$ be an $n$-dimensional locally area minimizing rectifiable current of $\mathbf{A}_{1,1}(0)$ such that \eqref{weak graphrep hyp1}, \eqref{weak graphrep hyp3}, and \eqref{weak graphrep hyp4} hold true with $\varepsilon_k, \mathbf{C}_k, T_k$ in place of $\overline{\varepsilon}, \mathbf{C}, T$ and either $p = 1$ or $p > 1$ and \eqref{weak graphrep hyp5} and \eqref{weak graphrep hyp6} hold true with $\beta_k, \mathbf{C}_k, T_k$ in place of $\overline{\beta}, \mathbf{C}, T$.  In view of the arbitrary choice of sequences $(\mathbf{C}_k)$ and $(T_k)$, it suffices to show that conclusion~(a) and (b) both hold true for infinitely many $k$.

By \eqref{weak graphrep hyp1}, \eqref{weak graphrep hyp3}, the Federer-Fleming compactness theorem, and~\cite[Theorem~34.5]{SimonGMT}, after passing to a subsequence there is an $n$-dimensional locally area minimizing rectifiable current $T_{\infty}$ of $\mathbf{A}_{1,1}(0)$ such that 
\begin{equation}\label{weak graphrep eqn1} 
	T_k \rightarrow T_{\infty} \text{ weakly in $\mathbf{A}_{1,1}(0)$.}
\end{equation}
By the monotonicity formula and the fact that $\|T_{k}\| \to \|T_{\infty}\|$ as Radon measures, 
\begin{equation}\label{weak graphrep eqn2} 
	\sup_{X \in \op{spt} T_k \cap \mathbf{A}_{1,\kappa'}(0)} \op{dist}(X, \op{spt} T_{\infty}) 
		+ \sup_{X \in \op{spt} T_{\infty} \cap \mathbf{A}_{1,\kappa'}(0)} \op{dist}(X, \op{spt} T_k) \rightarrow 0 
\end{equation}
for all $\kappa' \in (0,1)$.  Let $\mathbf{C}_k = \sum_{i=1}^p q^{(k)}_i \llbracket P^{(k)}_i \rrbracket$ for some integers $q^{(k)}_i \geq 1$ with $\sum_{i=1}^p q^{(k)}_i = q$ where $P^{(k)}_i$ are $n$-dimensional oriented planes with $\{0\} \times \mathbb{R}^{n-2} \subset P^{(k)}_i$ and with orienting $n$-vector $\vec P^{(k)}_i$.  After passing to a subsequence, there are integers $q^{(\infty)}_i \geq 1$, $n$-dimensional linear planes $P^{(\infty)}_i$, and orientation $n$-vectors $\vec P^{(\infty)}_i$ of $P^{(\infty)}_i$ such that 
\begin{equation}\label{weak graphrep eqn3}
	q^{(k)}_i \rightarrow q^{(\infty)}_i, \quad 
	\op{dist}_{\mathcal{H}}(P^{(k)}_i \cap \mathbf{B}_1(0), P^{(\infty)}_i \cap \mathbf{B}_1(0)) \rightarrow 0, \quad
	\vec P^{(k)}_i \rightarrow \vec P^{(\infty)}_i
\end{equation}
for each $i$.  After possibly reversing the orientations of the planes $P^{(k)}_i$, we may assume that for each $i, j \in \{1,2,\ldots,p\}$ if $P^{(\infty)}_i = P^{(\infty)}_j$ then $\vec P^{(\infty)}_i = \vec P^{(\infty)}_j$.  Thus 
\begin{equation*}
	\mathbf{C}_k \rightarrow \mathbf{C}_{\infty} = \sum_{i=1}^p q^{(\infty)}_i \llbracket P^{(\infty)}_i \rrbracket \text{ weakly in $\mathbf{A}_{1,1}(0)$.}
\end{equation*} 
By \eqref{weak graphrep hyp6} and \eqref{weak graphrep eqn3},  
\begin{equation*}
	\op{dist}_{\mathcal{H}}(P^{(\infty)}_i \cap \mathbf{B}_1(0),P^{(\infty)}_j \cap \mathbf{B}_1(0)) 
		\leq \mu \inf_{X \in P^{(\infty)}_i \cap (\mathbb{S}^{m+1} \times \mathbb{R}^{n-2})} \op{dist}(X,P^{(\infty)}_j)
\end{equation*} 
for each $i \neq j$.  Thus either $P^{(\infty)}_i = P^{(\infty)}_j$ or $P^{(\infty)}_i \cap P^{(\infty)}_j = \{0\} \times \mathbb{R}^{n-2}$ for each $i \neq j$.  
It follows from \eqref{weak graphrep hyp4} and monotonicity formula (as in Lemma~\ref{sepmono lemma}) that for each $\kappa' \in (0,1)$ and sufficiently large $k$ 
\begin{equation}\label{weak graphrep eqn4}
	\sup_{X \in \op{spt} T_k \cap \mathbf{A}_{1,\kappa'}(0)} \op{dist}(X, \op{spt} \mathbf{C}_k) \leq 2 \omega_n^{\frac{-1}{n+2}} \varepsilon_k^{\frac{2}{n+2}}.  
\end{equation}
Letting $k \rightarrow \infty$ in \eqref{weak graphrep eqn4} using \eqref{weak graphrep eqn2} and \eqref{weak graphrep eqn3} gives that $\op{spt} T_{\infty} \subseteq \op{spt} \mathbf{C}_{\infty}$.  In particular, by the constancy theorem $\op{spt} T_{\infty}$ is a union of $n$-dimensional planes contained in $\op{spt} \mathbf{C}_{\infty}$ in $\mathbf{A}_{1,1}(0)$ and $T_{\infty}$ has constant multiplicity on each plane in its support.

For each $i \in \{1,2,\ldots,p\}$, set  
\begin{equation*}
	V_i = \big\{ X \in \mathbf{A}_{1,(15+\kappa)/16}(0) : \op{dist}(X, \op{spt} P^{(\infty)}_i) < \tfrac{1}{3} \op{minsep} \mathbf{C}_{\infty} \big\} . 
\end{equation*}
Note that by \eqref{weak graphrep eqn3}, $\mathbf{C}_k \llcorner V_i$ is a sum of finitely many $n$-dimensional oriented planes with integer multiplicity.  By \eqref{weak graphrep hyp1}, \eqref{weak graphrep eqn2} and $\op{spt} T_{\infty} \subseteq \op{spt} \mathbf{C}_{\infty}$, $T_k \llcorner V_i$ is a locally area minimizing rectifiable current of $V_i$ with $(\partial (T_k \llcorner V_i)) \llcorner \mathbf{A}_{1,(15+\kappa)/16}(0) = 0$.  
By the constancy theorem, there is an integer $\widetilde{q}^{(k)}_i$ such that 
\begin{equation}\label{weak graphrep eqn5}
	\pi_{P^{(\infty)}_i \, \#} \, (T_{k} \llcorner V_{i}) \llcorner \mathbf{A}_{1,(7+\kappa)/8}(0) 
		= \widetilde{q}^{(k)}_i \llbracket P^{(\infty)}_i \rrbracket \llcorner \mathbf{A}_{1,(7+\kappa)/8}(0) .
\end{equation}
After possibly reversing the orientation of each plane of $\mathbf{C}_k$ converging to $P^{(\infty)}_i$, we may assume that $\widetilde{q}^{(k)}_i \geq 0$.  By \eqref{weak graphrep hyp3} and \eqref{graphrep annuli kappa hyp} 
\begin{align*}
	\widetilde{q}^{(k)}_i \mathcal{L}^n(A_{1,(7+\kappa)/8}(0)) \leq\,& \|T_{k}\|(V_{i} \cap \mathbf{A}_{1,(7+\kappa)/8}(0)) 
		\leq \|T_{k}\|(\mathbf{A}_{1,1}(0)) \\ \leq\,& (q+1/2) \,\mathcal{L}^n(A_{1,1}(0)) \leq (q+3/4) \,\mathcal{L}^n(A_{1,(7+\kappa)/8}(0))
\end{align*}
and thus $\widetilde{q}^{(k)}_i \leq q$.  By \eqref{weak graphrep eqn2} and Lemma~\ref{energy est lemma}, $\int_{V_{i} \cap \mathbf{A}_{1,(7+\kappa)/8}(0)} |\vec T_{k} - \vec P^{(\infty)}_i|^2 \,d\|T_k\| \rightarrow 0$, where $\vec T_{k}$ is the orientation $n$-vector of $T_k$.  
In view of these facts as well as \eqref{weak graphrep eqn3} and \eqref{weak graphrep hyp6}, we can apply Theorem~\ref{separation thm3} with $P^{(\infty)}_i$, $\mathbf{C}_k \llcorner V_i$ and $T_k \llcorner V_i$ in place of ${\mathbb R}^{n} \times \{0\}$, ${\mathbf P}$ and $T$ to obtain 
\begin{equation*}
	\sup_{\op{spt} T_k \cap V_i \cap \mathbf{A}_{1,(3+\kappa)/4}(0)} \op{dist}(X, \op{spt} \mathbf{C}_k \cap V_i) 
		\leq C E(T_k, \mathbf{C}_k, \mathbf{A}_{1,1}(0))
\end{equation*}
for some constant $C = C(n,m,q,\gamma,\kappa,\mu) \in (0,\infty)$.  In other words, 
\begin{equation}\label{weak graphrep eqn6}
	\sup_{\op{spt} T_k \cap \mathbf{A}_{1,(3+\kappa)/4}(0)} \op{dist}(X, \op{spt} \mathbf{C}_k) \leq C E(T_k, \mathbf{C}_k, \mathbf{A}_{1,1}(0))
\end{equation}
for some constant $C = C(n,m,q,\gamma,\kappa,\mu) \in (0,\infty)$.  Recall that $\mathbf{C}_k = \sum_{i=1}^p q^{(k)}_i \llbracket P^{(k)}_i \rrbracket$.  Let 
\begin{equation*}
	\big\{ X \in \mathbf{A}_{1,(3+\kappa)/4}(0) : \op{dist}(X, \op{spt} P^{(k)}_i) < 2 C E(T_k, \mathbf{C}_k, \mathbf{A}_{1,1}(0)) \big\} 
	= \bigcup_{i=1}^{N_k} U^{(k)}_i
\end{equation*}
as a union of connected components $U^{(k)}_i$, where $C$ is as in \eqref{weak graphrep eqn6}.  Thus $U^{(k)}_i$ are mutually disjoint, connected, open subsets of $\mathbf{A}_{1,(3+\kappa)/4}(0)$.  For each $i \in \{1,2,\ldots,N_k\}$, select a plane $P^{(k)}_{j(i)}$ of $\mathbf{C}_{k}$ such that $P^{(k)}_{j(i)} \cap \mathbf{A}_{1,(3+\kappa)/4}(0) \subset U^{(k)}_i$.  After passing to a subsequence we can take $N = N_k$ to be independent of $k$ and assume that $j(i)$ is independent of $k$ for each $i \in \{1,2,\ldots,N\}$.  Set 
\begin{equation*}
	T^{(k)}_{j(i)} = T_k \llcorner U^{(k)}_i 
\end{equation*}
for each $k$ and each $i \in \{1,2,\ldots,N\}$.  Set $T^{(k)}_j = 0$ if $j \not\in \{j(1),j(2),\ldots,j(N)\}$.  Clearly \eqref{graphrep annuli concl b1} and \eqref{graphrep annuli concl b2} hold true with $T_k$ and $T^{(k)}_i$ in place of $T$ and $T_i$.  By the above construction,  
\begin{equation*}
	\sup_{\op{spt} T \cap \mathbf{A}_{1,(3+\kappa)/4}(0)} \op{dist}( X, \op{spt}\mathbf{C}_k) \leq 4q C E(T_k, \mathbf{C}_k, \mathbf{A}_{1,1}(0)) 
\end{equation*}
where $C$ is as in \eqref{weak graphrep eqn6}, proving \eqref{graphrep annuli concl b4} with $\mathbf{C}_k, P^{(k)}_i, T_k, T^{(k)}_i$ in place of $\mathbf{C}, P_i, T, T_i$.  By the constancy theorem there exists integers $\widehat{q}^{(k)}_i$ such that \eqref{graphrep annuli concl b3} holds true with $\widehat{q}^{(k)}_i, P^{(k)}_i, T^{(k)}_i$ in place of $\widehat{q}_i, P_i, T_i$.  Since \eqref{weak graphrep eqn5} holds true with $\widetilde{q}^{(k)}_i \geq 0$, by Lemma~\ref{projection lemma2} and the continuity of push-forwards in the weak topology we have that $\widehat{q}^{(k)}_i \geq 0$ for all $k$ and $i$.  By \eqref{graphrep annuli concl b3}, \eqref{graphrep annuli concl b1}, \eqref{weak graphrep hyp3}, and \eqref{graphrep annuli kappa hyp}, 
\begin{align*}
	\sum_{i=1}^p \widehat{q}^{(k)}_i \,\mathcal{L}^n(A_{1,(1+\kappa)/2}(0)) 
	&\leq \sum_{i=1}^p \|T^{(k)}_i\|(\mathbf{A}_{1,(1+\kappa)/2}(0)) 
	= \|T_k\|(\mathbf{A}_{1,(1+\kappa)/2}(0)) 
	\\&\leq (q+1/2) \,\mathcal{L}^n(A_{1,1}(0))
	\leq (q+3/4) \,\mathcal{L}^n(A_{1,(1+\kappa)/2}(0))
\end{align*}
and thus $\sum_{i=1}^p \widehat{q}^{(k)}_i \leq q$.  Therefore, Lemma~\ref{weak graphrep lemma}(a) holds true.  By applying Almgren's Strong Lipschitz Approximation Theorem (Theorem~\ref{lip approx thm}) using a partition of unity argument (see the proof of~\cite[Theorem~3.8]{KrumWica}) and Lemma~\ref{energy est lemma}, Lemma~\ref{weak graphrep lemma}(b) holds true.
\end{proof}

The proof of Lemma~\ref{graphrep annuli lemma} will proceed by induction on $p$, assuming for $p_0 \in \{2,3,\ldots,q\}$ that 
\begin{enumerate}
	\item[(H3)]  Lemma~\ref{graphrep annuli lemma} holds true for all $p \in \{1, 2,\ldots,p_0-1\}$ . 
\end{enumerate}
Before proceeding with the proof of Lemma~\ref{graphrep annuli lemma} we observe the following. 

\begin{remark}\label{tildeC rmk} {\rm
(1)  Let $\widetilde{\beta} \in (0,1)$ be an arbitrary constant.  Suppose that $2 \leq p_0 \leq q$, $\mathbf{C} \in \mathcal{C}_{q,p_0}$, and $T$ is an $n$-dimensional rectifiable current of $\mathbf{A}_{1,1}(0)$.  Choose an integer $\widetilde{p} \in \{1,2,\ldots,p_0-1\}$ such that 
\begin{equation*}
	\inf_{\mathbf{C}' \in \bigcup_{p'=1}^k \mathcal{C}_{q,p'}} Q(T, \mathbf{C}', \mathbf{A}_{1,1}(0))
		\geq \frac{\widetilde{\beta}}{2} \inf_{\mathbf{C}' \in \bigcup_{p'=1}^{k-1} \mathcal{C}_{q,p'}} Q(T, \mathbf{C}', \mathbf{A}_{1,1}(0)) 
\end{equation*}
whenever $k \in \{\widetilde{p}+1,\widetilde{p}+2,\ldots,p_0-1\}$, and such that either $\widetilde{p} = 1$ or $\widetilde{p} > 1$ and 
\begin{equation*}
	\inf_{\mathbf{C}' \in \bigcup_{p'=1}^{\widetilde{p}} \mathcal{C}_{q,p'}} Q(T, \mathbf{C}', \mathbf{A}_{1,1}(0))
		< \frac{\widetilde{\beta}}{2} \inf_{\mathbf{C}' \in \bigcup_{p'=1}^{\widetilde{p}-1} \mathcal{C}_{q,p'}} Q(T, \mathbf{C}', \mathbf{A}_{1,1}(0)) .
\end{equation*}
Choose $\widetilde{\mathbf{C}} \in \mathcal{C}_{q,\widetilde{p}}$ such that 	
\begin{equation*}
	Q(T, \widetilde{\mathbf{C}}, \mathbf{A}_{1,1}(0)) 
		\leq 2 \inf_{\mathbf{C}' \in \bigcup_{p'=1}^{\widetilde{p}} \mathcal{C}_{q,p'}} Q(T, \mathbf{C}', \mathbf{A}_{1,1}(0)) .
\end{equation*}
Hence we choose $\widetilde{p} \in \{1,2,\ldots,p_0-1\}$ and $\widetilde{\mathbf{C}} \in \mathcal{C}_{q,\widetilde{p}}$ which satisfy  
\begin{equation}\label{tildeC eqn4} 
	Q(T, \widetilde{\mathbf{C}}, \mathbf{A}_{1,1}(0)) \leq 2^{p_0-1} \widetilde{\beta}^{2-p_0} 
		\inf_{\mathbf{C}' \in \bigcup_{p'=1}^{p_0-1} \mathcal{C}_{q,p'}} Q(T, \mathbf{C}', \mathbf{A}_{1,1}(0)) , 
\end{equation}
and either $\widetilde{p} = 1$ or $\widetilde{p} > 1$ and 
\begin{equation}\label{tildeC eqn5} 
	Q(T, \widetilde{\mathbf{C}}, \mathbf{A}_{1,1}(0)) 
		\leq \widetilde{\beta} \inf_{\mathbf{C}' \in \bigcup_{p'=1}^{\widetilde{p}-1} \mathcal{C}_{q,p'}} Q(T, \mathbf{C}', \mathbf{A}_{1,1}(0)) 
\end{equation}
(as in \eqref{graphrep annuli hyp5} with $\widetilde{\mathbf{C}}$ in place of $\mathbf{C}$).  

(2) Suppose that (H3) holds true and $T$ satisfies hypotheses of Lemma~\ref{graphrep annuli lemma}.   Let $\delta = \delta(n,m,q,p_0,\gamma) \in (0,1)$ and $\widetilde{\beta} = \widetilde{\beta}(n,m,q,p_0,\gamma) \in (0,1)$ be suitably small constants 
and suppose that 
\begin{equation}\label{tildeC eqn6} 
	\inf_{\mathbf{C}' \in \bigcup_{p'=1}^{p_0-1} \mathcal{C}_{q,p'}} Q(T, \mathbf{C}', \mathbf{A}_{1,1}(0)) < \delta .
\end{equation}
Let $\widetilde{p} \in \{1,2,\ldots,p_0-1\}$ and $\widetilde{\mathbf{C}} \in \mathcal{C}_{q,\widetilde{p}}$ such that \eqref{tildeC eqn4} holds true and either $\widetilde{p} = 1$ or $\widetilde{p} > 1$ and \eqref{tildeC eqn5} holds true.  Let $\op{spt} \widetilde{\mathbf{C}} = \bigcup_{i=1}^{\widetilde{p}} \widetilde{P}_i$ for $n$-dimensional oriented planes $\widetilde{P}_i$ such that $\{0\} \times \mathbb{R}^{n-2} \subset \widetilde{P}_i$.  By (H3), \eqref{weak graphrep hyp1}, \eqref{weak graphrep hyp3}, \eqref{tildeC eqn4}, \eqref{tildeC eqn5}, and \eqref{tildeC eqn6} we can apply Lemma~\ref{graphrep annuli lemma} with $\widetilde{\mathbf{C}}$ in place of $\mathbf{C}$ to deduce the following.  By Lemma~\ref{graphrep annuli lemma}(a) either $\widetilde{p} = 1$ or $\widetilde{p} > 1$ and 
\begin{gather}
	\label{graphrep tildeC eqn7-1} \op{minsep} \widetilde{\mathbf{C}} 
		\geq c \inf_{\mathbf{C}' \in \bigcup_{p'=1}^{\widetilde{p}-1} \mathcal{C}_{q,p'}} Q(T, \mathbf{C}', \mathbf{A}_{1,1}(0)) , \\
	\label{graphrep tildeC eqn8-1} \op{dist}_{\mathcal{H}}(\widetilde{P}_i \cap \mathbf{B}_1(0), \widetilde{P}_{i'} \cap \mathbf{B}_1(0))
		\leq C \inf_{X \in \widetilde{P}_i \cap (\mathbb{S}^{m+1} \times \mathbb{R}^{n-2})} \op{dist}(X, \widetilde{P}_{i'}) 
		\text{ for all $1 \leq i, i' \leq \widetilde{p}$}
\end{gather}
for some constants $c = c(n,m,q,p_0,\gamma) > 0$ and $C = C(n,m,q,p_0,\gamma) \in (0,\infty)$.  By Lemma~\ref{graphrep annuli lemma}(b), there exists integers $\widetilde{q}_i > 0$ with $\sum_{i=1}^{\widetilde{p}} \widetilde{q}_i \leq q$ and $n$-dimensional locally area minimizing rectifiable currents $\widetilde{T}_i$ of $\mathbf{A}_{1,(15+\kappa)/16}(0)$ such that 
\begin{gather} 
	\label{graphrep tildeC eqn9-1} T \llcorner \mathbf{A}_{1,(15+\kappa)/16}(0) = \sum_{i=1}^{\widetilde{p}} \widetilde{T}_i , \\
	\label{graphrep tildeC eqn10-1} (\partial \widetilde{T}_i) \llcorner \mathbf{A}_{1,(15+\kappa)/16}(0) = 0, \\ 
	\label{graphrep tildeC eqn11-1} (\pi_{\widetilde{P}_i \#} \widetilde{T}_i) \llcorner \mathbf{A}_{1,(7+\kappa)/8}(0) 
		= \widetilde{q}_i \llbracket \widetilde{P}_i \rrbracket \llcorner \mathbf{A}_{1,(7+\kappa)/8}(0) , \\
	\label{graphrep tildeC eqn12-1} \sup_{X \in \op{spt} \widetilde{T}_i} \op{dist}(X, \widetilde{P}_i) 
		\leq C \widetilde{Q} , 
\end{gather}
where $\widetilde{Q} = Q(T, \widetilde{\mathbf{C}}, \mathbf{A}_{1,1}(0))$ and $C = C(n,m,q,p_0,\gamma,\kappa) \in (0,\infty)$ is a constant.  
By Lemma~\ref{graphrep annuli lemma}(c), for each $i \in \{1,2,\ldots,\widetilde{p}\}$ there exist a Lipschitz $\widetilde{q}_i$-valued function $\widetilde{u}_i : \widetilde{P}_i \cap \mathbf{A}_{1,(3+\kappa)/4}(0) \rightarrow \mathcal{A}_{\widetilde{q}_i}(\widetilde{P}_i^{\perp})$ and a closed set $\widetilde{K}_i \subseteq \widetilde{P}_i \cap \mathbf{A}_{1,(3+\kappa)/4}(0)$ such that 
\begin{gather}\label{graphrep tildeC eqn010-1}
	\widetilde{T}_i \llcorner \pi_{\widetilde{P}_i}^{-1}(\widetilde{K}_i) = (\op{graph} \widetilde{u}_i) \llcorner \pi_{\widetilde{P}_i}^{-1}(\widetilde{K}_i) , \\
	\mathcal{H}^n(\widetilde{P}_i \cap \mathbf{A}_{1,(3+\kappa)/4}(0) \setminus \widetilde{K}_i)
		+ \|\widetilde{T}_i\|(\pi_{\widetilde{P}_i}^{-1}(\widetilde{P}_i \cap \mathbf{A}_{1,(3+\kappa)/4}(0) \setminus \widetilde{K}_i)) 
		\leq C \widetilde{Q}^{2+\alpha} , \nonumber \\
	\sup_{\widetilde{P}_i \cap \mathbf{A}_{1,(3+\kappa)/4}(0)} |\widetilde{u}_i| \leq C \widetilde{Q} , \quad 
		\sup_{\widetilde{P}_i \cap \mathbf{A}_{1,(3+\kappa)/4}(0)} |\nabla \widetilde{u}_i| \leq C \widetilde{Q}^{\alpha} , \nonumber
\end{gather}
where  and $\alpha = \alpha(n,m,q) \in (0,1)$ and $C = C(n,m,q,p_0,\gamma,\kappa) \in (0,\infty)$ are constants. 
	
(3)  Suppose that $\widetilde{T}_i$ are as in \eqref{graphrep tildeC eqn9-1}--\eqref{graphrep tildeC eqn12-1} with $\kappa = 1/2$ and $\widetilde{u}_i$ are as in \eqref{graphrep tildeC eqn010-1} with $\kappa = 1/2$.  Let $X \in \op{spt} \mathbf{C} \cap \mathbf{A}_{1,1/2}(0)$.  By \eqref{tildeC eqn4}, \eqref{tildeC eqn6} and \eqref{graphrep tildeC eqn12-1}, $\op{dist}(X, \op{spt} T) < (1-\gamma)/32$.  Thus the closest point to $X$ in $\op{spt} T$ lies in $\op{spt} T \cap \mathbf{A}_{1,3/4}(0)$.  In other words, by \eqref{graphrep tildeC eqn9-1}, the closest point to $X$ in $\op{spt} T$ lies on $\op{spt} \widetilde{T}_i$ for some $i$.  By the triangle inequality and \eqref{graphrep tildeC eqn12-1}, 
\begin{align}\label{graphrep tildeC eqn13-1}
	\op{dist}(X, \op{spt} \widetilde{\mathbf{C}}) 
	&\leq \op{dist}(X,\widetilde{P}_i) 
	\leq \op{dist}(X, \op{spt} \widetilde{T}_i) + \sup_{Y \in \op{spt} \widetilde{T}_i} \op{dist}(X, \widetilde{P}_i) 
	\\&\leq \op{dist}(X, \op{spt} T) + C \widetilde{Q} , \nonumber 
\end{align}
where $C = C(n,m,q,p_0,\gamma) \in (0,\infty)$ are constants.  Integrating \eqref{graphrep tildeC eqn13-1} over $X \in \op{spt} \mathbf{C} \cap \mathbf{A}_{1,1/2}(0)$ and using \eqref{weak graphrep hyp3} and \eqref{graphrep annuli hyp5}, 
\begin{equation}\label{graphrep tildeC eqn014-1}
	\int_{\mathbf{A}_{1,1/2}(0)} \op{dist}^2(X, \op{spt} \widetilde{\mathbf{C}}) \,d\|\mathbf{C}\|(X) 
	\leq 2 \int_{\mathbf{A}_{1,1/2}(0)} \op{dist}^2(X, \op{spt} T) \,d\|\mathbf{C}\|(X) + C \widetilde{Q}^2 \leq C \widetilde{Q}^2 , 
\end{equation}
where $C = C(n,m,q,p_0,\gamma) \in (0,\infty)$ are constants.  On the other hand, by the triangle inequality and \eqref{graphrep tildeC eqn010-1} 
\begin{align}\label{graphrep tildeC eqn015-1}
	\op{dist}(X, \op{spt} \mathbf{C}) 
	&\leq \op{dist}(X + \widetilde{u}_{i,j}(X), \op{spt} \mathbf{C}) + |\widetilde{u}_{i,j}(X)|
	\\&\leq \op{dist}(X + \widetilde{u}_{i,j}(X), \op{spt} \mathbf{C}) + C \widetilde{Q} \nonumber
\end{align}
for each $X \in \widetilde{P}_i \cap \mathbf{A}_{1/2,1}(0)$ and $j \in \{1,2,\ldots,\widetilde{q}_i\}$, where $\widetilde{u}_i(X) = \sum_{j=1}^{\widetilde{q}_i} \llbracket \widetilde{u}_{i,j}(X) \rrbracket$ for some $\widetilde{u}_{i,j}(X) \in \widetilde{P}_i^{\perp}$ and $C = C(n,m,q,p_0,\gamma) \in (0,\infty)$ is a constant.  Integrating \eqref{graphrep tildeC eqn015-1} over $X \in \widetilde{P}_i \cap \mathbf{A}_{1/2,1}(0)$ and using the area formula, \eqref{graphrep tildeC eqn010-1}, and \eqref{graphrep annuli hyp5} 
\begin{align}\label{graphrep tildeC eqn016-1}
	&\sum_{i=1}^{\widetilde{p}} \widetilde{q}_i \int_{\widetilde{P}_i \cap \mathbf{A}_{1/2,1}(0)} \op{dist}^2(X, \op{spt} \mathbf{C}) \,d\mathcal{H}^n(X)
	\\&\hspace{15mm} \leq 2 \int_{\mathbf{A}_{3/4,1}(0)} \op{dist}^2(X, \op{spt} \mathbf{C}) \,d\|T\|(X) + C \widetilde{Q}^2 \leq C \widetilde{Q}^2 . \nonumber 
\end{align}
Since $\mathbf{C} \in \mathcal{C}_{q,p}$ and $\widetilde{\mathbf{C}} \in \mathcal{C}_{q,\widetilde{p}}$, it follows from \eqref{graphrep tildeC eqn014-1} and \eqref{graphrep tildeC eqn016-1} that 
\begin{equation}\label{graphrep tildeC eqn14-1}
	\op{dist}_{\mathcal H}(\op{spt} \mathbf{C} \cap \mathbf{B}_1(0), \op{spt} \widetilde{\mathbf{C}} \cap \mathbf{B}_1(0)) \leq C \widetilde{Q} 
\end{equation}
for some constant $C = C(n,m,q,p_0,\gamma) \in (0,\infty)$. 

(4)  By \eqref{tildeC eqn5} and \eqref{graphrep tildeC eqn7-1}, $Q(T, \widetilde{\mathbf{C}}, \mathbf{A}_{1,1}(0)) \leq C\widetilde{\beta} \op{minsep} \widetilde{\mathbf{C}}$ for some constant $C = C(n,m,q,p_0,\gamma) \in (0,\infty)$.  (Recall that if $\widetilde{p} = 1$ then $\op{minsep} \widetilde{\mathbf{C}} = \infty$.)  Thus assuming $\widetilde{\beta}$ is sufficient small, by \eqref{graphrep tildeC eqn14-1} 
\begin{equation}\label{graphrep tildeC eqn15}
	\mathbf{C} = \sum_{i=1}^{\widetilde{p}} \sum_{j=1}^{s_i} q_{i,j} \llbracket P_{i,j} \rrbracket 
\end{equation} 
for some integers $s_i \geq 1$ and $q_{i,j} \geq 1$ such that $\sum_{i=1}^{\widetilde{p}} s_i = p_0$ and $\sum_{i=1}^{\widetilde{p}} \sum_{j=1}^{s_i} q_{i,j} = q$ and some distinct $n$-dimensional oriented planes $P_{i,j}$ such that $\{0\} \times \mathbb{R}^{n-2} \subset P_{i,j}$ and 
\begin{equation}\label{graphrep tildeC eqn16} 
	\op{dist}_{\mathcal{H}}(P_{i,j} \cap \mathbf{B}_1(0), \widetilde{P}_i \cap \mathbf{B}_1(0)) \leq C \,Q(T, \widetilde{\mathbf{C}}, \mathbf{A}_{1,1}(0)) ,
\end{equation} 
where $C = C(n,m,q,p_0,\gamma) \in (0,\infty)$ is a constant.  In light of \eqref{graphrep tildeC eqn16} we may assume  (by reversing the orientation of $\widetilde{P}_{i}$ if necessary) that the orientation of $P_{i,j}$ is close to the orientation of $\widetilde{P}_i$ as $n$-vectors.  When $\widetilde{p} > 1$, since $p_0 > \widetilde{p}$, there exists $i$ and $j \neq j'$ such that $P_{i,j}$ and $P_{i,j'}$ are distinct planes close to $\widetilde{P}_i$ and thus by \eqref{graphrep tildeC eqn16} 
\begin{equation}\label{graphrep tildeC eqn17} 
	\op{maxsep} \mathbf{C} \leq C \,Q(T, \widetilde{\mathbf{C}}, \mathbf{A}_{1,1}(0))
\end{equation}
for some constant $C = C(n,m,q,p_0,\gamma) \in (0,\infty)$.  
}\end{remark}

\begin{proof}[Proof of Lemma~\ref{graphrep annuli lemma}]
We shall proceed by induction on $p$.  Let $\gamma,\kappa \in (0,1)$.  Let us look at the base case $p = 1$.  For $p = 1$ we do not need to prove Lemma~\ref{graphrep annuli lemma}(a).  By \eqref{weak graphrep hyp1}, \eqref{weak graphrep hyp3}, and \eqref{weak graphrep hyp4} we can apply Lemma~\ref{one plane lemma} and Almgren's Strong Lipschitz Approximation Theorem (Theorem~\ref{lip approx thm}) to obtain Lemma~\ref{graphrep annuli lemma}(b)(c) provided we also show that $\widehat{q}_1 > 0$.  (In particular, Lemma~\ref{graphrep annuli lemma}(b) holds true with $T_1 = T \llcorner \mathbf{A}_{1,(3+\kappa)/4}(0)$.)  By \eqref{graphrep annuli concl b4}, \eqref{weak graphrep hyp4}, and Lemma~\ref{energy est lemma}, $\int_{\mathbf{A}_{1,(1+\kappa)/2}(0)} |\vec T - \vec P_1|^2 \,d\|T\|(X) \leq C\overline{\varepsilon}^2$, where $\vec T$ is the orientation $n$-vector of $T$ and $C = C(n,m,q,\gamma,\kappa) \in (0,\infty)$ is a constant.  Thus since we assumed $\op{spt} T \cap \mathbf{A}_{1,1/2}(0) \neq \emptyset$, by Lemma~\ref{projection lemma} we must have that $\widehat{q}_1 > 0$.

Suppose now that $p_0 \in \{2,3,\ldots,q\}$ and (H3) holds. We want to show Lemma~\ref{graphrep annuli lemma} holds true when $p = p_0$.

\noindent\textit{Proof of Lemma~\ref{graphrep annuli lemma}(a).}  Let $\mathbf{C}$ and $T$ satisfy the hypotheses of Lemma~\ref{graphrep annuli lemma}.  We will prove Lemma~\ref{graphrep annuli lemma}(a) by separately considering the cases 
\begin{align*}
	\text{(I)} &\quad \inf_{\mathbf{C}' \in \bigcup_{p'=1}^{p_0-1} \mathcal{C}_{q,p'}} Q(T, \mathbf{C}', \mathbf{A}_{1,1}(0)) < \delta; \\
	\text{(II)} &\quad \inf_{\mathbf{C}' \in \bigcup_{p'=1}^{p_0-1} \mathcal{C}_{q,p'}} Q(T, \mathbf{C}', \mathbf{A}_{1,1}(0)) \geq \delta;
\end{align*}
where $\delta = \delta(n,m,q,p_0,\gamma) \in (0,1)$ will be chosen so that Lemma~\ref{graphrep annuli lemma}(a) holds true in Case~(I).  

\noindent\textit{Case~(I).}  Fix $\widetilde{\beta} = \widetilde{\beta}(n,m,q,p_0,\gamma,\kappa) \in (0,1)$ are small enough that we can apply Remark~\ref{tildeC rmk}(2)(3)(4).  
By Remark~\ref{tildeC rmk}(1), there exists $1 \leq \widetilde{p} < p_0$ and $\widetilde{\mathbf{C}} \in \mathcal{C}_{q,\widetilde{p}}$ such that \eqref{tildeC eqn4} holds true and either $\widetilde{p} = 1$ or $\widetilde{p} > 1$ and \eqref{tildeC eqn5} holds true.  
Let $\op{spt} \widetilde{\mathbf{C}} = \bigcup_{i=1}^{\widetilde{p}} \widetilde{P}_i$ for $n$-dimensional oriented planes $\widetilde{P}_i$ such that $\{0\} \times \mathbb{R}^{n-2} \subset \widetilde{P}_i$.  

It suffices to first show \eqref{graphrep annuli concl a1}.  Then \eqref{graphrep annuli concl a2} will follow by observing that by \eqref{graphrep annuli concl a1}, \eqref{tildeC eqn4}, and \eqref{graphrep tildeC eqn16} 
\begin{align}\label{graphrep tildeC eqn18} 
	\op{dist}_{\mathcal{H}}(P_{i,j} \cap \mathbf{B}_1(0), P_{i,j'} \cap \mathbf{B}_1(0))
	\leq\,& C \inf_{{\mathbf C}^{\prime} \in \cup_{p^{\prime} = 1}^{p_{0} - 1} {\mathcal C}_{q, p^{\prime}}} \, Q(T, {\mathbf C}^{\prime}, {\mathbf B}_{1}(0)) \\
	\leq\,& C \op{minsep} {\mathbf C} \leq C \inf_{X \in P_{i, j} \cap ({\mathbf S}^{m-1} \times {\mathbb R}^{n-2})} \, {\rm dist} \, (X, P_{i, j^{\prime}}) \nonumber
\end{align}
for all $1 \leq i \leq p$ and $j \neq j^{\prime}$.  Moreover, by \eqref{graphrep tildeC eqn16}, \eqref{tildeC eqn5} and \eqref{graphrep tildeC eqn7-1} we have 
${\rm dist}_{\mathcal H} \, (P_{i, j} \cap {\mathbf B}_{1}(0), \widetilde{P}_{i} \cap {\mathbf B}_{1}(0)) \leq C\widetilde{\b}  \op{minsep} \widetilde{\mathbf C},$
which together with \eqref{graphrep tildeC eqn8-1} implies that 
\begin{align*} 
	\op{dist}_{\mathcal{H}}(P_{i,j} \cap \mathbf{B}_1(0), P_{i',j'} \cap \mathbf{B}_1(0))
	\leq\,& 2 \op{dist}_{\mathcal{H}}(\widetilde{P}_i \cap \mathbf{B}_1(0), \widetilde{P}_{i'} \cap \mathbf{B}_1(0)) 
	\\ \leq\,& 2C \inf_{X \in \widetilde{P}_i \cap (\mathbb{S}^{m+1} \times \mathbb{R}^{n-2})} \op{dist}(X, \widetilde{P}_{i'}) \nonumber
	\\ \leq\,& 4C \inf_{X \in P_{i,j} \cap (\mathbb{S}^{m+1} \times \mathbb{R}^{n-2})} \op{dist}(X, P_{i',j'}) \nonumber
\end{align*}
for all $i \neq i'$, $1 \leq j \leq s_i$, and $1 \leq j' \leq s_{i'}$, where $C = C(n,m,q,p_0,\gamma) \in (0,\infty)$ is a constant; this together with \eqref{graphrep tildeC eqn18} proves \eqref{graphrep annuli concl a2}.

We claim that provided $\delta$ is sufficiently small, \eqref{graphrep annuli concl a1} holds true.  Suppose to the contrary that for $k = 1,2,3,\ldots$ there are $\varepsilon_k \rightarrow 0^+$, $\beta_k \rightarrow 0^+$, $\delta_k \rightarrow 0$, $\mathbf{C}_k \in \mathcal{C}_{q,p_0}$, and $n$-dimensional locally area minimizing rectifiable currents $T_k$ of $\mathbf{A}_{1,1}(0)$ such that \eqref{weak graphrep hyp1}, \eqref{weak graphrep hyp3}, \eqref{weak graphrep hyp4}, and \eqref{graphrep annuli hyp5} hold true with $\varepsilon_k, \beta_k, \mathbf{C}_k, T_k$ in place of $\overline{\varepsilon}, \overline{\beta}, \mathbf{C}, T$ and 
\begin{equation}
	\inf_{\mathbf{C}' \in \bigcup_{p'=1}^{p_0-1} \mathcal{C}_{q,p'}} Q(T_k, \mathbf{C}', \mathbf{A}_{1,1}(0)) < \delta_k 
\end{equation}
but 
\begin{equation}\label{graphrep sep1 eqn1}
	\op{minsep}  \mathbf{C}_k \leq \frac{1}{k} \,Q(T, \widetilde{\mathbf{C}}_k, \mathbf{A}_{1,1}(0)) .
\end{equation}
Let $1 \leq \widetilde{p} < p$ and $\widetilde{\mathbf{C}}_k \in \mathcal{C}_{q,\widetilde{p}}$ such that \eqref{tildeC eqn4} holds true with $\delta_k,T_k,\widetilde{\mathbf{C}}_k$ in place of $\delta,T,\widetilde{\mathbf{C}}$ and either $\widetilde{p} = 1$ or $\widetilde{p} > 1$ and \eqref{tildeC eqn5} holds true with $T_k$ and $\widetilde{\mathbf{C}}_k$ in place of $T$ and $\widetilde{\mathbf{C}}$.  Note that after passing to a subsequence we assume that $\widetilde{p}$ is independent of $k$.  Let 
\begin{equation*}
	\op{spt} \widetilde{\mathbf{C}}_k = \bigcup_{i=1}^{\widetilde{p}} \widetilde{P}^{(k)}_i
\end{equation*} 
for some distinct $n$-dimensional oriented planes $\widetilde{P}^{(k)}_i$ such that $\{0\} \times \mathbb{R}^{n-2} \subset \widetilde{P}^{(k)}_i$.  Recall that there exists integers $\widetilde{q}_i$ and $n$-dimensional locally area minimizing rectifiable currents $\widetilde{T}^{(k)}_i$ such that 
\begin{equation*}
	T_k \llcorner \mathbf{A}_{1,31/32}(0) = \sum_{i=1}^{\widetilde{p}} \widetilde{T}^{(k)}_i
\end{equation*} 
(as in \eqref{graphrep tildeC eqn9-1}) and \eqref{graphrep tildeC eqn10-1}--\eqref{graphrep tildeC eqn12-1} hold true with $\kappa = 1/2$ and with $\widetilde{\mathbf C}_k, \widetilde{P}^{(k)}_i, T_k, \widetilde{T}^{(k)}_i$ in place of $\widetilde{\mathbf C}, \widetilde{P}_i, T, \widetilde{T}_i$.   Again, after passing to a subsequence we assume that $\widetilde{q}_i$ are independent of $k$.  Further recall that 
there exists Lipschitz $\widetilde{q}_i$-valued functions $\widetilde{u}^{(k)}_i : \widetilde{P}^{(k)}_i \cap \mathbf{A}_{1,1/2}(0) \rightarrow \mathcal{A}_{\widetilde{q}_i}((\widetilde{P}^{(k)}_i)^{\perp})$ and sets $\widetilde{K}^{(k)}_i \subseteq \widetilde{P}^{(k)}_i \cap \mathbf{A}_{1,1/2}(0)$ such that 
\begin{gather}
	\label{graphrep sep1 eqn2} \widetilde{T}^{(k)}_i \llcorner \pi_{\widetilde{P}^{(k)}_i}^{-1}(\widetilde{K}^{(k)}_i) 
		= (\op{graph} \widetilde{u}^{(k)}_i) \llcorner \pi_{\widetilde{P}^{(k)}_i}^{-1}(\widetilde{K}^{(k)}_i) ,  \\
	\mathcal{H}^n(\widetilde{P}^{(k)}_i \cap \mathbf{A}_{1,1/2}(0) \setminus \widetilde{K}^{(k)}_i) 
		+ \|\widetilde{T}^{(k)}_i\|(\pi_{\widetilde{P}^{(k)}_i}^{-1}(\widetilde{P}^{(k)}_i) \cap \mathbf{A}_{1,1/2}(0) \setminus \widetilde{K}^{(k)}_i)) 
		\leq C \widetilde{Q}_k^{2+\alpha}, \nonumber \\
	\sup_{\widetilde{P}^{(k)}_i \cap \mathbf{A}_{1,1/2}(0)} |\widetilde{u}^{(k)}_i| \leq C \widetilde{Q}_k , \quad 
	\op{Lip} \widetilde{u}^{(k)}_i \leq C \widetilde{Q}_k^{\alpha} \nonumber
\end{gather}
(as in \eqref{graphrep tildeC eqn010-1}), where $\widetilde{Q}_k = Q(T_k, \widetilde{\mathbf{C}}_k, \mathbf{A}_{1,1}(0))$ and where $\alpha = \alpha(n,m,q) \in (0,1)$ and $C = C(n,m,q,p_0,\gamma,\kappa) \in (0,\infty)$ are constants.  
Express $\mathbf{C}_k$ as 
\begin{equation*}
	\mathbf{C}_k = \sum_{i=1}^{\widetilde{p}} \sum_{j=1}^{s_i} q_{i,j} \llbracket P^{(k)}_{i,j} \rrbracket 
\end{equation*}
(as in \eqref{graphrep tildeC eqn15}) where $s_i \geq 1$ and $q_{i,j} \geq 1$ are integers with $\sum_{i=1}^{\widetilde{p}} s_i = p_0$ and $\sum_{i=1}^{\widetilde{p}} \sum_{j=1}^{s_i} q_{i,j} = q$ and $P^{(k)}_{i,j}$ are distinct $n$-dimensional oriented planes such that $\{0\} \times \mathbb{R}^{n-2} \subset P^{(k)}_{i,j}$ and \eqref{graphrep tildeC eqn16} holds true with $\widetilde{P}^{(k)}_i, P^{(k)}_{i,j}, \widetilde{C}_k, T_k$ in place of $\widetilde{P}_i, P_{i,j}, \widetilde{C}, T$.  After passing to a subsequence we assume that $s_i$ and $q_{i,j}$ are independent of $k$.  By \eqref{graphrep tildeC eqn16}, each $P^{(k)}_{i,j}$ is the graph of a linear single-valued function $\phi^{(k)}_{i,j} : \widetilde{P}^{(k)}_i \rightarrow (\widetilde{P}^{(k)}_i)^{\perp}$ such that  
\begin{equation}\label{graphrep sep1 eqn3}
	\|\phi^{(k)}_{i,j}\|_{L^{\infty}(B_1(0,\widetilde{P}^{(k)}_i))} \leq C \widetilde{Q}_k 
\end{equation} 
for some constant $C = C(n,m,q,p_0,\gamma) \in (0,\infty)$.  For each $k$ and $i$, set $\phi^{(k)}_i = \sum_{j=1}^{s_i} q_{i,j} \llbracket \phi^{(k)}_{i,j} \rrbracket$ as a multi-valued function.  

Let $\mathfrak{q}^{(k)}_i : \mathbb{R}^{n+m} \rightarrow \mathbb{R}^{n+m}$ be an orthogonal linear transformation such that $\mathfrak{q}^{(k)}_i(\{0\} \times \mathbb{R}^{n-2}) = \{0\} \times \mathbb{R}^{n-2}$ and $\mathfrak{q}^{(k)}_i(\widetilde{P}^{(k)}_i) = \{0\} \times \mathbb{R}^n$.  It follows that $\mathfrak{q}^{(k)}_i((\widetilde{P}^{(k)}_i)^{\perp}) = \mathbb{R}^m \times \{0\}$.  By the area formula and \eqref{graphrep sep1 eqn2},  
\begin{align*}
	\int_{\widetilde{P}_i \cap \mathbf{A}_{1,1/2}(0)} |D\widetilde{u}^{(k)}_i|^2 
	\leq \frac{1}{2} \int_{\mathbf{A}_{1,5/8}(0)} |\vec T^{(k)}_i - \vec P^{(k)}_i|^2 \,d\|\widetilde{T}^{(k)}_i\| + C \widetilde{Q}_k^{2+\alpha}
	\leq C \widetilde{Q}_k^2 , 
\end{align*}
where $\vec T^{(k)}_i$ and $\vec P^{(k)}_i$ are the orientation $n$-vectors of $\widetilde{T}^{(k)}_i$ and $\widetilde{P}^{(k)}_i$ respectively and where $C = C(n,m,q,p_0,\gamma) \in (0,\infty)$ are constants.  Thus by~\cite[Theorem~2.19]{Almgren} (or \cite[Proposition~2.11 and Proposition~3.20]{DeLSpaDirMin}) after passing to a subsequence there is a locally Dirichlet energy minimizing $\widetilde{q}_i$-valued function $w_i \in W^{1,2}_{\rm loc}(\mathbf{A}_{1,1/2}(0),\mathcal{A}_{\widetilde{q}_i}(\mathbb{R}^m))$ such that $(\mathfrak{q}^{(k)}_i \circ \widetilde{u}^{(k)}_i \circ (\mathfrak{q}^{(k)}_i)^{-1})/E_k \rightarrow w_i$ in $L^2(\mathbf{A}_{1,1/2}(0),\mathcal{A}_{\widetilde{q}_i}(\mathbb{R}^m))$.  By \eqref{graphrep sep1 eqn3}, after passing to a subsequence there are linear single-valued functions $\psi_{i,j} : \mathbb{R}^n \rightarrow \mathbb{R}^m$ such that $(\mathfrak{q}^{(k)}_i \circ \phi^{(k)}_{i,j} \circ (\mathfrak{q}^{(k)}_i)^{-1})/\widetilde{Q}_k \rightarrow \psi_{i,j}$ uniformly on $B_1(0)$.  For each $i$, set $\psi_i = \sum_{j=1}^{s_i} q_{i,j} \llbracket \psi_{i,j} \rrbracket$ as a multi-valued function such that $(\mathfrak{q}^{(k)}_i \circ \phi^{(k)}_i \circ (\mathfrak{q}^{(k)}_i)^{-1})/\widetilde{Q}_k \rightarrow \psi_i$ uniformly on $\mathbf{A}_{1,1/2}(0)$.

For each $X \in \mathbf{A}_{1,1/2}(0)$ let $\widetilde{u}^{(k)}_i(X) = \sum_{j=1}^{\widetilde{q}_i} \llbracket \widetilde{u}^{(k)}_{i,j}(X) \rrbracket$ for some $\widetilde{u}^{(k)}_{i,j}(X) \in \mathbb{R}^m$.  Let $\op{spt} \widetilde{u}^{(k)}_i(X) = \{\widetilde{u}^{(k)}_{i,1}(X),\ldots,\widetilde{u}^{(k)}_{i,\widetilde{q}_i}(X)\}$ be the set of all values of $\widetilde{u}^{(k)}_i(X)$ and similarly let $\op{spt} \phi^{(k)}_i(X)$, $\op{spt} w_i(X)$, and $\op{spt} \psi_i(X)$ denote the set of all values of $\phi^{(k)}_i(X)$, $w_i(X)$, and $\psi_i(X)$ respectively.  Since by \eqref{graphrep sep1 eqn3} $|\nabla \phi^{(k)}_{i,j}| \leq C E$ is small, 
\begin{equation*}
	\op{dist}(\widetilde{u}^{(k)}_{i,j}(X), \op{spt} \phi^{(k)}_i(X)) \leq 2 \op{dist}(X + \widetilde{u}^{(k)}_{i,j}(X), \op{spt} \mathbf{C}_k)
\end{equation*}
for each $X \in \mathbf{A}_{1,1/2}(0)$ and sufficiently large $k$.  Thus by the area formula, \eqref{graphrep sep1 eqn2}, and \eqref{graphrep annuli hyp5},  
\begin{align}\label{graphrep sep1 eqn5}
	&\sum_{i=1}^{\widetilde{p}} \int_{\widetilde{P}^{(k)}_i \cap \mathbf{A}_{1,1/2}(0)} \sum_{j=1}^{\widetilde{q}_i} 
		\op{dist}^2(\widetilde{u}^{(k)}_{i,j}(X), \op{spt} \phi^{(k)}_i(X)) \,d\mathcal{H}^n(X) 
	\\ \leq\,& 4 \int_{\mathbf{A}_{1,1}(0)} \op{dist}^2(X, \op{spt} \mathbf{C}_k) \,d\|T_k\|(X) + C \widetilde{Q}_k^{2+\alpha} \nonumber 
	\\ \leq\,& 4 \beta_k^2 \widetilde{Q}_k^2 + C \widetilde{Q}_k^{2+\alpha}, \nonumber
\end{align}
where $\alpha = \alpha(n,m,q) \in (0,1)$ and $C = C(n,m,q,p_0,\gamma) \in (0,\infty)$ are constants.  Dividing \eqref{graphrep sep1 eqn5} by $\widetilde{Q}_k^2$ and letting $k \rightarrow \infty$ gives us that $\op{spt} w_i(X) \subseteq \op{spt} \psi_i(X)$ for $\mathcal{L}^n$-a.e.~$X \in A_{1,1/2}(0)$.  On the other hand, by a similar argument \eqref{graphrep sep1 eqn2} and \eqref{graphrep annuli hyp5} also give us 
\begin{equation}\label{graphrep sep1 eqn6}
	\sum_{i=1}^{\widetilde{p}} \sum_{j=1}^{s_i} q_{i,j} \int_{\widetilde{P}^{(k)}_i \cap \mathbf{A}_{1,1/2}(0)} 
		\op{dist}^2(\phi^{(k)}_{i,j}(X), \op{spt} \widetilde{u}^{(k)}_i(X)) \,d\mathcal{H}^n(X) 
	\leq 4 \beta_k^2 \widetilde{Q}_k^2 + C \widetilde{Q}_k^{2+\alpha}, 
\end{equation}
where $\alpha = \alpha(n,m,q) \in (0,1)$ and $C = C(n,m,q,p_0,\gamma) \in (0,\infty)$ are constants.  Dividing \eqref{graphrep sep1 eqn6} by $\widetilde{Q}_k^2$ and letting $k \rightarrow \infty$ gives us $\op{spt} w_i(X) = \op{spt} \psi_i(X)$ for $\mathcal{L}^n$-a.e.~$X \in A_{1,1/2}(0)$.  In fact since $w_i$ is Dirichlet energy minimizing and therefore continuous in $A_{1,1/2}(0)$, $\op{spt} w_i(X) = \op{spt} \psi_i(X)$ for all $X \in A_{1,1/2}(0)$. 

By \eqref{weak graphrep hyp3} and \eqref{tildeC eqn5} we can apply Remark~\ref{graphrep rmk} to obtain  
\begin{equation*}
	\op{maxsep}  \mathbf{C}_k \geq c \inf_{\mathbf{C}' \in \bigcup_{p'=1}^{p_0-1} \mathcal{C}_{q,p'}} Q(T_k, \mathbf{C}', \mathbf{A}_{1,1}(0))
		\geq 2^{1-p_0} \widetilde{\beta}^{p_0-2} c \,Q(T_k, \widetilde{\mathbf{C}}_k, \mathbf{A}_{1,1}(0)) , 
\end{equation*}
where $c = c(n,m,q,\gamma) > 0$ is a constant and the last step follows from \eqref{tildeC eqn4}.  Hence, since each $P^{(k)}_{i,j}$ is graph of a linear single-valued $\phi^{(k)}_{i,j}$, 
\begin{equation}\label{graphrep sep1 eqn8}
	\sup_{X \in \widetilde{P}^{(k)}_i \cap (\mathbb{S}^{m+1} \times \mathbb{R}^{n-2})} |\phi^{(k)}_{i,j}(X) - \phi^{(k)}_{i,j'}(X)| \geq c \,\widetilde{Q}_k 
\end{equation}
for all $k$, $1 \leq i \leq p$, and $1 \leq j < j' \leq s_i$ and for some constant $c = c(n,m,q,p_0,\gamma) > 0$.  Dividing \eqref{graphrep sep1 eqn8} by $\widetilde{Q}_k$ and letting $k \rightarrow \infty$ gives us that $\psi_{i,j} \not\equiv \psi_{i,j'}$ for all $1 \leq i \leq p$ and $1 \leq j < j' \leq s_i$.  Similarly, by \eqref{graphrep sep1 eqn1} for each sufficiently large $k$ there exists $1 \leq i \leq p$ and $1 \leq j < j' \leq s_i$ such that 
\begin{equation}\label{graphrep sep1 eqn9}
	\inf_{X \in \widetilde{P}^{(k)}_i \cap (\mathbb{S}^{m+1} \times \mathbb{R}^{n-2})} |\phi^{(k)}_{i,j}(X) - \phi^{(k)}_{i,j'}(X)| \leq \frac{2}{k} \,\widetilde{Q}_k . 
\end{equation}
After passing to a subsequence, we may take $i$, $j$, and $j'$ to be independent of $k$.  Dividing \eqref{graphrep sep1 eqn9} by $\widetilde{Q}_k$ and letting $k \rightarrow \infty$ gives us that $\psi_{i,j}(X) = \psi_{i,j'}(X)$ for some $X \in \mathbb{S}^1 \times \mathbb{R}^{n-2}$.  Hence since $\psi_{i,j} = \psi_{i,j'} = 0$ on $\{0\} \times \mathbb{R}^{n-2}$, $\{\psi_{i,j} = \psi_{i,j'}\}$ is an $(n-1)$-dimensional linear subspace of $\mathbb{R}^n$.  Since $\op{spt} w_i(X) = \op{spt} \psi_i(X)$ for all $X \in A_{1,1/2}(0)$, the singular set of $w_i$ contains the $(n-1)$-dimensional linear subspace $\{\psi_{i,j} = \psi_{i,j'}\}$, contradicting $w_i$ being Dirichlet energy minimizing.  Therefore, \eqref{graphrep annuli concl a1} holds true.

\noindent\textit{Case~(II).}  Now fix $\delta$ such that Lemma~\ref{graphrep annuli lemma}(a) holds true in Case~(I).  Let's show that provided $\overline{\varepsilon}, \overline{\beta}$ are sufficiently small, 
\begin{equation}\label{graphrep sep2 eqn1}
	\op{minsep} \mathbf{C} \geq c
\end{equation}
for some constant $c = c(n,m,q,p_0,\gamma) \in (0,1)$.  Since 
\begin{equation*}
	\op{dist}_{\mathcal{H}}(P_i \cap \mathbf{B}_1(0), P_j \cap \mathbf{B}_1(0)) \leq 2
\end{equation*}
for each pair $P_i$ and $P_j$ of planes of $\mathbf{C}$ and by \eqref{weak graphrep hyp3} 
\begin{equation}\label{graphrep sep2 eqn2}
	\inf_{\mathbf{C}' \in \bigcup_{p'=1}^{p-1} \mathcal{C}_{q,p'}} Q(T, \mathbf{C}', \mathbf{A}_{1,1}(0)) \leq C
\end{equation}
for some constant $C = C(n,m,q,\gamma) \in (0,\infty)$, showing \eqref{graphrep sep2 eqn1} will prove Lemma~\ref{graphrep annuli lemma}(a).  Note that by Remark~\ref{graphrep rmk}, 
\begin{equation}\label{graphrep sep2 eqn3}
	\op{maxsep} \mathbf{C} \geq c \,\delta
\end{equation}
for some constant $c = c(n,m,q,p_0,\gamma) > 0$. 

To see \eqref{graphrep sep2 eqn1}, suppose to the contrary that for $k = 1,2,3,\ldots$ there are $\varepsilon_k \rightarrow 0^+$, $\beta_k \rightarrow 0^+$, $\mathbf{C}_k \in \mathcal{C}_{q,p}$, and $n$-dimensional locally area minimizing rectifiable currents $T_k$ of $\mathbf{A}_{1,1}(0)$ such that \eqref{weak graphrep hyp1}, \eqref{weak graphrep hyp3}, \eqref{weak graphrep hyp4}, and \eqref{graphrep annuli hyp5} hold true with $\varepsilon_k, \beta_k, \mathbf{C}_k, T_k$ in place of $\overline{\varepsilon}, \overline{\beta}, \mathbf{C}, T$ and 
\begin{equation}\label{graphrep sep2 eqn4}
	\inf_{\mathbf{C}' \in \bigcup_{p'=1}^{p_0-1} \mathcal{C}_{q,p'}} Q(T_k, \mathbf{C}', \mathbf{A}_{1,1}(0)) \geq \delta
\end{equation}
but 
\begin{equation}\label{graphrep sep2 eqn5}
	\op{minsep} \mathbf{C}_k \leq \frac{1}{k} . 
\end{equation}
By \eqref{weak graphrep hyp1}, \eqref{weak graphrep hyp3}, the Federer-Fleming compactness theorem, and~\cite[Theorem~34.5]{SimonGMT}, after passing to a subsequence there is an $n$-dimensional locally area minimizing rectifiable current $T_{\infty}$ such that 
\begin{equation*}
	T_k \rightarrow T_{\infty} \text{ weakly in $\mathbf{A}_{1,1}(0)$.}
\end{equation*} 
In particular, by~\cite[Theorem~34.5]{SimonGMT} and the monotonicity formula, 
\begin{gather}
	\label{graphrep sep2 eqn6} \|T_k\| \rightarrow \|T_{\infty}\| \text{ in the sense of Radon measures of $\mathbf{A}_{1,1}(0)$}, \\
	\label{graphrep sep2 eqn7} \sup_{X \in \op{spt} T_k \cap \mathbf{A}_{1,3/4}(0)} \op{dist}(X, \op{spt} T_{\infty}) 
		+ \sup_{X \in \op{spt} T_{\infty} \cap \mathbf{A}_{1,3/4}(0)} \op{dist}(X, \op{spt} T_k) \rightarrow 0 . 
\end{gather}
Let $\mathbf{C}_k = \sum_{i=1}^p q^{(k)}_i \llbracket P^{(k)}_i \rrbracket$ for some integers $q^{(k)}_i \geq 1$ with $\sum_{i=1}^p q^{(k)}_i = q$ and $P^{(k)}_i$ are $n$-dimensional oriented planes with $\{0\} \times \mathbb{R}^{n-2} \subset P^{(k)}_i$ and with orienting $n$-vector $\vec P^{(k)}_i$.  After passing to a subsequence there are integers $q^{(\infty)}_i \geq 1$, $n$-dimensional linear planes $P^{(\infty)}_i$, and orienting $n$-vectors $\vec P^{(\infty)}_i$ or $P^{(\infty)}_i$ such that 
\begin{equation}\label{graphrep sep2 eqn8}
	q^{(k)}_i \rightarrow q^{(\infty)}_i, \quad 
	\op{dist}_{\mathcal{H}}(P^{(k)}_i \cap \mathbf{A}_{1,1}(0), P^{(\infty)}_i \cap \mathbf{A}_{1,1}(0)) \rightarrow 0, \quad
	\vec P^{(k)}_i \rightarrow \vec P^{(\infty)}_i
\end{equation}
for each $i$, and thus 
\begin{equation*}
	\mathbf{C}_k \rightarrow \mathbf{C}_{\infty} = \sum_{i=1}^p q^{(\infty)}_i \llbracket P^{(\infty)}_i \rrbracket \text{ weakly in $\mathbf{A}_{1,1}(0)$.}
\end{equation*} 
Note that by \eqref{graphrep sep2 eqn3} and \eqref{graphrep sep2 eqn8}, $\op{maxsep} \mathbf{C}_{\infty} \geq c \,\delta$ (where $c$ is as in \eqref{graphrep sep2 eqn3}) and thus $\mathbf{C}_{\infty}$ consists of $p$ distinct planes $P^{(\infty)}_i$.  In particular, $\mathbf{C}_{\infty} \in \mathcal{C}_{q,p}$.  By \eqref{graphrep annuli hyp5}, \eqref{graphrep sep2 eqn2}, and the monotonicity formula, for all $\kappa' \in (0,1)$ and sufficiently large $k$ 
\begin{align*}
	&\sup_{X \in \op{spt} T_k \cap \mathbf{A}_{1,1/4}(0)} \op{dist}(X, \op{spt} \mathbf{C}_k) 
		+ \sup_{X \in \op{spt} \mathbf{C}_k \cap \mathbf{A}_{1,1/4}(0)} \op{dist}(X, \op{spt} T_k) 
	\\ \leq\,& C Q(T_k, \mathbf{C}_k, \mathbf{A}_{1,1}(0))^{\frac{2}{n+2}} 
	\leq C \beta_k^{\frac{2}{n+2}} ,
\end{align*}
where $C = C(n,m,q,\gamma) \in (0,\infty)$ are constants.  Letting $k \rightarrow \infty$ using \eqref{graphrep sep2 eqn7} and \eqref{graphrep sep2 eqn8} 
\begin{equation}\label{graphrep sep2 eqn10}
	\op{spt} T_{\infty} \cap \mathbf{A}_{1,1/4}(0) = \op{spt} \mathbf{C}_{\infty} \cap \mathbf{A}_{1,1/4}(0). 
\end{equation} 
But by \eqref{graphrep sep2 eqn5} and \eqref{graphrep sep2 eqn8}, at least two distinct planes of $P^{(\infty)}_i$ and $P^{(\infty)}_j$ of $\mathbf{C}_{\infty}$ intersect in $\mathbb{S}^{m+1} \times \mathbb{R}^{n-2}$.  Recalling $\{0\} \times \mathbb{R}^{n-2}$ is a subspace of both $P^{(\infty)}_i$ and $P^{(\infty)}_j$, $P^{(\infty)}_i$ and $P^{(\infty)}_j$ intersect along an $(n-1)$-dimensional linear subspace.  Hence by \eqref{graphrep sep2 eqn10}, $\op{spt} T_{\infty} \cap \mathbf{A}_{1,1/2}(0)$ is a union of $n$-dimensional planes, with two distinct planes intersecting along an $(n-1)$-dimensional linear subspace, contradicting $T_{\infty}$ being area-minimizing in $\mathbf{A}_{1,1}(0)$.  Therefore, \eqref{graphrep sep2 eqn1} must hold true.

\noindent\textit{Proof of Lemma~\ref{graphrep annuli lemma}(b) and (c).}  In light of \eqref{graphrep annuli concl a2}, we can apply Lemma~\ref{weak graphrep lemma}(a) to obtain 
\begin{equation*}
	\sup_{X \in \op{spt} T \cap \mathbf{A}_{1,(1+\kappa)/2}(0)} \op{dist}(X, \op{spt}\mathbf{C}) \leq C E(T,\mathbf{C},\mathbf{A}_{1,1}(0)) 
\end{equation*}
for some constant $C = C(n,m,q,p_0,\gamma,\kappa) \in (0,\infty)$.  In particular, by \eqref{graphrep annuli hyp5} and \eqref{graphrep annuli concl a1}
\begin{equation*}
	\sup_{X \in \op{spt} T \cap \mathbf{A}_{1,(1+\kappa)/2}(0)} \op{dist}(X, \op{spt}\mathbf{C}) \leq C \overline{\beta} \op{sep} \mathbf{C}
\end{equation*}
for some constant $C = C(n,m,q,p_0,\gamma,\kappa) \in (0,\infty)$.  Therefore, assuming $\overline{\beta}$ is sufficiently small and possibly reversing the orientation of $P_i$, there exists integers $\widehat{q}_i \geq 0$ and $n$-dimensional locally area minimizing rectifiable currents $T_i$ of $\mathbf{A}_{1,(3+\kappa)/4}(0)$ such that $\sum_{i=1}^{p_0} \widehat{q}_i \leq q$ and \eqref{graphrep annuli concl b1}, \eqref{graphrep annuli concl b2} and \eqref{graphrep annuli concl b3} hold true.  By Lemma~\ref{one plane lemma}, \eqref{graphrep annuli concl b4} holds true with $C = C(n,m,\gamma,\kappa)$.  To complete the proof of Lemma~\ref{graphrep annuli lemma}(b), it remains to show that $\widehat{q}_i > 0$.  Then by applying Almgren's Strong Lipschitz Approximation Theorem (Theorem~\ref{lip approx thm}) like in Lemma~\ref{weak graphrep lemma}(b), Lemma~\ref{graphrep annuli lemma}(c) holds true. 

Suppose that 
$\widehat{q}_i = 0$ for some $i \in \{1,2,\ldots,p_0\}$.  By \eqref{graphrep annuli concl b4}, \eqref{weak graphrep hyp4}, and Lemma~\ref{energy est lemma}, we have that $\int_{\mathbf{A}_{1,(1+\kappa)/2}(0)} |\vec T - \vec P_i|^2 \,d\|T\|(X) \leq C\overline{\varepsilon}^2$, where $\vec T$ is the orientation $n$-vector of $T$ and $C = C(n,m,q,\gamma,\kappa) \in (0,\infty)$ is a constant.  Hence by Lemma~\ref{projection lemma}, $\op{spt} T_i \cap \mathbf{A}_{1,(3+5\kappa)/8}(0) = \emptyset$.  
We claim that for each $X \in P_i \cap \mathbf{A}_{1,1/4}(0)$ 
\begin{equation}\label{graphrep annuli eqn2}
	\op{dist}(X, \op{spt} T) \geq c \inf_{\mathbf{C}' \in \bigcup_{p'=1}^{p_0-1} \mathcal{C}_{q,p'}} Q(T,\mathbf{C}',\mathbf{A}_{1,1}(0))
\end{equation}
for some constant $c = c(n,m,q,p_0,\gamma,\kappa) > 0$.  In the case that $\op{dist}(X, \op{spt} T) \geq (1-\gamma)/64$, \eqref{graphrep annuli eqn2} obviously holds true.  If $\op{dist}(X, \op{spt} T) < (1-\gamma)/64$, then by the triangle inequality, $\op{spt} T_i \cap \mathbf{A}_{1,3/8}(0) = \emptyset$, \eqref{graphrep annuli concl a1}, \eqref{graphrep annuli concl b4}, and \eqref{graphrep annuli hyp5} 
\begin{align*}
	\op{dist}(X, \op{spt} T) 
	&\geq \op{dist}(X, \op{spt} \mathbf{C} \setminus P_i) 
		- \sup_{Y \in \op{spt} T \cap \mathbf{A}_{1,3/8}(0)} \op{dist}(Y, \op{spt} \mathbf{C} \setminus P_i)
	\\&= \op{dist}(X, \op{spt} \mathbf{C} \setminus P_i) 
		- \sup_{Y \in \op{spt} T \cap \mathbf{A}_{1,3/8}(0)} \op{dist}(Y, \op{spt} \mathbf{C}) \nonumber 
	\\&\geq c \inf_{\mathbf{C}' \in \bigcup_{p'=1}^{p_0-1} \mathcal{C}_{q,p'}} Q(T,\mathbf{C}',\mathbf{A}_{1,1}(0)) - C \,E(T,\mathbf{C},\mathbf{A}_{1,1}(0)) 
		\nonumber 
	\\&\geq \frac{1}{2} \,c \inf_{\mathbf{C}' \in \bigcup_{p'=1}^{p_0-1} \mathcal{C}_{q,p'}} Q(T,\mathbf{C}',\mathbf{A}_{1,1}(0)) , \nonumber
\end{align*}
where $c = c(n,m,q,p_0,\gamma,\kappa) > 0$ and $C = C(n,m,q,p_0,\gamma,\kappa) \in (0,\infty)$ are constants.  Integrating \eqref{graphrep annuli eqn2} over $X \in P_i$, 
\begin{equation*}
	Q(T,\mathbf{C},\mathbf{A}_{1,1}(0))^2 \geq \int_{P_i \cap \mathbf{A}_{1,1/4}(0)} \op{dist}^2(X, \op{spt} T) \,d\mathcal{H}^n(X) 
		\geq c \inf_{\mathbf{C}' \in \bigcup_{p'=1}^{p_0-1} \mathcal{C}_{q,p'}} Q(T,\mathbf{C}',\mathbf{A}_{1,1}(0))^2
\end{equation*}
for some constant $c = c(n,m,q,p_0,\gamma,\kappa) > 0$, contradicting \eqref{graphrep annuli hyp5}.
\end{proof}

\begin{proof}[Proof of Theorem~\ref{graphrep thm}]
Suppose that $\varepsilon_0, \beta_0, \mathbf{C}, T$ satisfy the hypotheses of Theorem~\ref{graphrep thm}.  Take any $\rho \geq \tau/4$ and $\zeta \in \mathbb{R}^{n-2}$ such that $\rho^2 + |\zeta|^2 \leq (3+\gamma)^2/16$.  We claim that 
\begin{gather}
	\label{graphrep eqn1} (\partial T) \llcorner \mathbf{A}_{\rho,1}(\zeta) = 0 , \\
	\label{graphrep eqn2} E(T, \mathbf{C}, \mathbf{A}_{\rho,1}(\zeta)) 
		\leq \rho^{-(n+2)/2} \,E(T, \mathbf{C}, \mathbf{B}_1(0)) < (\tau/4)^{-(n+2)/2} \varepsilon_0 , \\
	\label{graphrep eqn3} \|T\|(\mathbf{A}_{\rho,1/2}(\zeta)) > (q-1/2) \,\mathcal{L}^n(A_{\rho,1/2}(0)) , \\
	\label{graphrep eqn4} \|T\|(\mathbf{A}_{\rho,1}(\zeta)) < (q+1/2) \,\mathcal{L}^n(A_{\rho,1}(0)) , \\
	\label{graphrep eqn5} E(T,\mathbf{C},\mathbf{A}_{\rho,1}(\zeta)) \leq C \tau^{-(n+2)/2} \beta_0 
		\inf_{\mathbf{C}' \in \bigcup_{p'=1}^{p-1} \mathcal{C}_{q,p'}} Q(T, \mathbf{C}', \mathbf{A}_{\rho,1}(\zeta)) ,
\end{gather}
where $C = C(n,m,p,q,\gamma) \in (0,\infty)$ is a constant.  Moreover, we claim that 
\begin{equation}\label{graphrep eqn6}
	Q(T,\mathbf{C},\mathbf{A}_{1/4,1}(0)) \leq C \tau^{-(n+2)/2} \beta_0 
		\inf_{\mathbf{C}' \in \bigcup_{p'=1}^{p-1} \mathcal{C}_{q,p'}} Q(T, \mathbf{C}', \mathbf{A}_{1/4,1}(0))
\end{equation}
where again $C = C(n,m,p,q,\gamma) \in (0,\infty)$ is a constant.  Since $\mathbf{A}_{\rho,1}(\zeta) \subset \mathbf{B}_1(0)$, \eqref{graphrep eqn1} follows directly from $(\partial T) \llcorner \mathbf{B}_1(0) = 0$ and \eqref{graphrep eqn2} follows from \eqref{main hyp eqn2}. 

Let us show \eqref{graphrep eqn3} holds true.  Note that the proof of \eqref{graphrep eqn4} is similar.  Suppose to the contrary that for $k = 1,2,3,\ldots$ there exists $\varepsilon_k \rightarrow 0^+$, $\mathbf{C}_k \in \mathcal{C}_{q,p}$, and $T_k$ be an $n$-dimensional locally area minimizing rectifiable current of $\mathbf{B}_1(0)$ such that \eqref{main hyp eqn1} and \eqref{main hyp eqn2} hold true with $\varepsilon_k, \mathbf{C}_k, T_k$ in place of $\varepsilon_0, \mathbf{C}, T$ but for some $\rho_k \geq \tau/4$ and $\zeta_k \in \mathbb{R}^{n-2}$ with $\rho_k^2 + |\zeta_k|^2 \leq (3+\gamma)^2/16$ 
\begin{equation}\label{graphrep eqn7}
	\|T_k\|(\mathbf{A}_{\rho_k,1/2}(\zeta_k)) \leq (q-1/2) \,\mathcal{L}^n(A_{\rho_k,1/2}(0)) . 
\end{equation}
By \eqref{main hyp eqn1}, the Federer-Fleming compactness theorem, and~\cite[Theorem~34.5]{SimonGMT}, after passing to a subsequence there exists an $n$-dimensional locally area minimizing rectifiable current $T_{\infty}$ of $\mathbf{B}_1(0)$ such that 
\begin{gather}
	T_k \rightarrow T_{\infty} \text{ weakly in $\mathbf{B}_1(0)$,} \nonumber \\
	\label{graphrep eqn8} \|T_k\| \rightarrow \|T_{\infty}\| \text{ in the sense of Radon measures of $\mathbf{B}_1(0)$.}
\end{gather}
By the monotonicity formula, 
\begin{equation}\label{graphrep eqn9} 
	\sup_{X \in \op{spt} T_{\infty} \cap \mathbf{B}_{\sigma}(0)} \op{dist}(X, \op{spt} T_k) \rightarrow 0
\end{equation}
for all $\sigma \in (0,1)$.  Let $\mathbf{C}_k = \sum_{i=1}^p q^{(k)}_i \llbracket P^{(k)}_i \rrbracket$ for some integers $q^{(k)}_i \geq 1$ such that $\sum_{i=1}^p q^{(k)}_i = q$ and for some $n$-dimensional oriented planes $P^{(k)}_i$ with $\{0\} \times \mathbb{R}^{n-2} \subset P^{(k)}_i$ and with orientation $n$-vectors $\vec P^{(k)}_i$.  After passing to a subsequence, there are integers $q^{(\infty)}_i$ and $n$-dimensional oriented planes $P^{(\infty)}_i$ with $\{0\} \times \mathbb{R}^{n-2} \subset P^{(k)}_i$ and with orientation $n$-vectors $\vec P^{(\infty)}_i$ such that 
\begin{equation}\label{graphrep eqn10}
	q^{(k)}_i \rightarrow q^{(\infty)}_i, \quad 
	\op{dist}_{\mathcal{H}}(P^{(k)}_i \cap \mathbf{B}_1(0), P^{(\infty)}_i \cap \mathbf{B}_1(0)) \rightarrow 0, \quad 
	\vec P^{(k)}_i \rightarrow \vec P^{(\infty)}_i
\end{equation}
for each $i$.  After possibly reversing the orientations of the planes $P^{(k)}_i$, we may assume that for each $i, j \in \{1,2,\ldots,p\}$ if $P^{(\infty)}_i = P^{(\infty)}_j$ then $\vec P^{(\infty)}_i = \vec P^{(\infty)}_j$.  Thus 
\begin{equation*}
	\mathbf{C}_k \rightarrow \mathbf{C}_{\infty} = \sum_{i=1}^p q^{(\infty)}_i \llbracket P^{(\infty)}_i \rrbracket \text{ weakly in $\mathbf{B}_1(0)$.}
\end{equation*} 
It follows from \eqref{main hyp eqn2} and the monotonicity formula (as in Lemma~\ref{sepmono lemma}), for all $\sigma \in (0,1)$ and sufficiently large $k$ 
\begin{equation*}
	\sup_{X \in \op{spt} T_k \cap \mathbf{B}_{\sigma}(0)} \op{dist}(X, \op{spt} \mathbf{C}_k) \leq 2 \varepsilon_k^{\frac{2}{n+2}}.
\end{equation*}
Letting $k \rightarrow \infty$ using \eqref{graphrep eqn9} and \eqref{graphrep eqn10} gives us 
\begin{equation}\label{graphrep eqn11}
	\op{spt} T_{\infty} \subseteq \op{spt} \mathbf{C}_{\infty} = \bigcup_{i=1}^p P^{(\infty)}_i . 
\end{equation}
Notice that any two distinct planes of $\mathbf{C}_{\infty}$ must intersect either precisely along $\{0\} \times \mathbb{R}^{n-2}$ or along an $(n-1)$-dimensional linear subspace.  Since $T_{\infty}$ is area-minimizing, there is no $X \in \op{sing} T_{\infty} \setminus (\{0\} \times {\mathbb R}^{n-2})$ and $\delta > 0$ such that $\op{spt} T_{\infty} \cap \mathbf{B}_{\delta}(X)$ is a union of three or more distinct $n$-dimensional half-planes meeting along an $(n-1)$-dimensional affine subspace passing through $X$.  It follows that for each $X \in \op{spt} T_{\infty} \setminus (\{0\} \times {\mathbb R}^{n-2})$ there exists $\delta > 0$ and exactly one plane $P^{(\infty)}_i$ of $\mathbf{C}_{\infty}$ such that $\op{spt} T_{\infty} \cap \mathbf{B}_{\delta}(X) = P^{(\infty)}_i \cap \mathbf{B}_{\delta}(X)$.  Hence $\op{spt} T_{\infty}$ is the union of a subcollection of distinct planes of $\mathbf{C}_{\infty}$ such that any two planes intersect precisely along $\{0\} \times \mathbb{R}^{n-2}$ and $T_{\infty}$ has constant integer multiplicity on each plane.  Consequently, $\Theta(T_{\infty},0)$ is an integer.  By the upper semi-continuity of density and $\Theta(T_k,0) \geq q$, it follows that $\Theta(T_{\infty},0) \geq q$.  But by \eqref{graphrep eqn8} and $\|T_k\|(\mathbf{B}_1(0)) \leq (q+1/2) \,\omega_n$, we must have that $\Theta(T_{\infty},0) = \omega_n^{-1} \|T_{\infty}\|(\mathbf{B}_1(0)) \leq q+1/2$.  Therefore, $\Theta(T_{\infty},0) = q$.  By \eqref{graphrep eqn8} and the fact that $\op{spt} T_{\infty}$ is the union of a subcollection of distinct planes of $\mathbf{C}_{\infty},$ 
\begin{equation*}
	\lim_{k \rightarrow \infty} \|T_k\|(\mathbf{A}_{\rho_k,1/2}(\zeta_k)) \geq \|T_{\infty}\|(\mathbf{A}_{\rho_{\infty},1/2}(\zeta_{\infty})) 
		= q \,\mathcal{L}^n(A_{\rho_{\infty},1/2}(0)) ,
\end{equation*}
contradicting \eqref{graphrep eqn7}. 

Finally, let's verify \eqref{graphrep eqn5} and \eqref{graphrep eqn6}.  Let $\rho \geq \tau/4$ and $\zeta \in \mathbb{R}^{n-2}$ with $\rho^2 + |\zeta|^2 \leq (3+\gamma)^2/16$.  By \eqref{main hyp eqn1} and \eqref{main hyp eqn3} we can apply Remark~\ref{graphrep rmk} to obtain 
\begin{equation}\label{graphrep eqn12}
	Q(T,\mathbf{C},\mathbf{B}_1(0)) \leq C \beta_0 \op{maxsep} \mathbf{C} 
\end{equation}
for some constant $C = C(n,q) \in (0,\infty)$.  Since $\mathbf{A}_{\rho,1}(\zeta) \subset \mathbf{B}_1(0)$, 
\begin{equation}\label{graphrep eqn13}
	E(T,\mathbf{C},\mathbf{A}_{\rho,1}(\zeta)) \leq C \tau^{-(n+2)/2} \beta_0 \op{maxsep} \mathbf{C} 
\end{equation}
for some constant $C = C(n,m,p,q,\gamma) \in (0,\infty)$.  Let $\delta = \delta(n,m,q,p, \gamma) \in (0,1)$ and $\widetilde{\beta} = \widetilde{\beta}(n,m,q,p, \gamma) \in (0,1)$ be small enough that we can apply Remark~\ref{tildeC rmk}(2)(3)(4).  
If 
\begin{equation*}
	\inf_{\mathbf{C}' \in \bigcup_{p'=1}^{p-1} \mathcal{C}_{q,p'}} Q(T, \mathbf{C}', \mathbf{A}_{\rho,1}(\zeta)) \geq \delta
\end{equation*}
then since $\op{minsep} \mathbf{C} \leq 2$ it follows from \eqref{graphrep eqn13} that \eqref{graphrep eqn5} holds true.  Otherwise, by Remark~\ref{tildeC rmk} there exists an integer $1 \leq \widetilde{p} < p$ and $\widetilde{\mathbf{C}} \in \mathcal{C}_{q,\widetilde{p}}$ such that 
\begin{gather}
	Q(T, \widetilde{\mathbf{C}}, \mathbf{A}_{\rho,1}(\zeta)) \leq 2^{p-1} \widetilde{\beta}^{2-p} \delta , \nonumber \\
	\label{graphrep eqn14} Q(T, \widetilde{\mathbf{C}}, \mathbf{A}_{\rho,1}(\zeta)) \leq 2^{p-1} \widetilde{\beta}^{2-p} 
		\inf_{\mathbf{C}' \in \bigcup_{p'=1}^{p-1} \mathcal{C}_{q,p'}} Q(T, \mathbf{C}', \mathbf{A}_{\rho,1}(\zeta)) 
\end{gather}
and either $\widetilde{p} = 1$ or $\widetilde{p} > 1$ and 
\begin{equation*}
	Q(T, \widetilde{\mathbf{C}}, \mathbf{A}_{\rho,1}(\zeta)) \leq \widetilde{\beta} 
		\inf_{\mathbf{C}' \in \bigcup_{p'=1}^{\widetilde{p}-1} \mathcal{C}_{q,p'}} Q(T, \mathbf{C}', \mathbf{A}_{\rho,1}(\zeta)) . 
\end{equation*}
It follows that \eqref{graphrep tildeC eqn17} holds true with $\eta_{(0,\zeta),\rho\#} T$ in place of $T$, that is 
\begin{equation}\label{graphrep eqn15} 
	\op{maxsep} \mathbf{C} \leq C \,Q(T, \widetilde{\mathbf{C}}, \mathbf{A}_{\rho,1}(\zeta))
\end{equation}
for some constant $C = C(n,m,q,p_0,\gamma) \in (0,\infty)$.  Combining \eqref{graphrep eqn13}, \eqref{graphrep eqn15}, and \eqref{graphrep eqn14} gives us \eqref{graphrep eqn5}.  When $\rho = 1/4$ and $\zeta = 0$, we also have that $\mathbf{A}_{1/4,1/2}(0) \subset \mathbf{B}_{1/2}(0) \cap \{ r > 1/16 \}$ and thus \eqref{graphrep eqn12} gives us 
\begin{equation}\label{graphrep eqn16}
	Q(T,\mathbf{C},\mathbf{A}_{1/4,1}(\zeta)) \leq C \beta_0 \op{maxsep} \mathbf{C} 
\end{equation}
for some constant $C = C(n,m,p,q,\gamma) \in (0,\infty)$.  Combining \eqref{graphrep eqn16}, \eqref{graphrep eqn15}, and \eqref{graphrep eqn14} gives us \eqref{graphrep eqn6}.  

In light of \eqref{graphrep eqn1}, \eqref{graphrep eqn2}, \eqref{graphrep eqn4}, and \eqref{graphrep eqn6}, we can apply Lemma~\ref{graphrep annuli lemma} with $\eta_{0,1/4 \#} T$ in place of $T$.  In particular, by Lemma~\ref{graphrep annuli lemma}(a) with $\eta_{0,1/4 \#} T$ in place of $T$, we deduce that Theorem~\ref{graphrep thm}(a) holds true.  Thus by \eqref{graphrep eqn1}, \eqref{graphrep eqn2}, \eqref{graphrep eqn4}, \eqref{graphrep eqn5}, and \eqref{graphrep concl a2}, we can apply Lemma~\ref{weak graphrep lemma} with $\eta_{(0,\zeta),\rho \#} T$ in place of $T$ for all $\rho \geq \tau/4$ and $\zeta \in \mathbb{R}^{n-2}$ with $\rho^2 + |\zeta|^2 \leq (3+\gamma)^2/16$.

To see conclusions~(b) and (c), observe that $\mathbf{A}_{\rho,\kappa}(\zeta) = \{ (r\omega,y) : \omega \in \mathbb{S}^{m+1}, \,(r,y) \in B^{n-1}_{(1-\gamma) \kappa \rho/8}(\rho,\zeta) \}$ for each $\kappa \in (0,1)$, $\rho > 0$, and $\zeta \in \mathbb{R}^{n-2}$.  Note that if $\mathbf{A}_{\rho,\kappa}(\zeta) \cap \mathbf{A}_{\rho',\kappa}(\zeta') \neq \emptyset$ then $|(\rho,\zeta) - (\rho',\zeta')| < \tfrac{1}{8} \,(1-\gamma) \kappa \,(\rho + \rho')$ and thus $\tfrac{8-(1-\gamma)\kappa}{8+(1-\gamma)\kappa} \,\rho < \rho' <  \tfrac{8+(1-\gamma)\kappa}{8-(1-\gamma)\kappa} \,\rho$.  By the Vitali covering lemma there is a finite collection of $\rho_j \geq \tau/4$ and $\zeta_j \in \mathbb{R}^{n-2}$ with $\rho_j^2 + |\zeta_j|^2 \leq (3+\gamma)^2/16$ (for $1 \leq j \leq N$) such that $\rho_{j+1} \leq \rho_j$ for all $1 \leq j < N$, $\{\mathbf{A}_{\rho_j,1/8}(\rho_j,\zeta_j)\}$ covers $\mathbf{B}_{(3+\gamma)/4}(0) \cap \{r > \tau/4\}$, and $\{\mathbf{A}_{\rho_j,1/40}(\rho_j,\zeta_j)\}$ is pairwise disjoint.  For each $\sigma > 0$, if $\rho_j \leq \sigma/2$ then $\mathbf{A}_{\rho_j,1}(\zeta_j) \subset \{r \leq \sigma \}$.  Set $J_{\sigma} = \max\{ j : \rho_j \geq \sigma/2 \}$.  Since $\{B^{n-1}_{(1-\gamma)\rho_j/320}(\rho_j,\zeta_j)\}$ is a collection of pairwise disjoint balls in $B^{n-1}_1(0)$, for each $\sigma > 0$ 
\begin{equation*}
	J_{\sigma} \,\omega_{n-1} \bigg(\frac{(1-\gamma)\sigma}{320}\bigg)^{n-1} 
	\leq \sum_{j=1}^{J_{\sigma}} \mathcal{L}^{n-1}(B^{n-1}_{(1-\gamma)\rho_j/320}(\rho_j,\zeta_j)) 
	\leq \mathcal{L}^{n-1}(B^{n-1}_1(0)) = \omega_{n-1} 
\end{equation*}
and thus 
\begin{equation}\label{graphrep eqn17}
	J_{\sigma} \leq C \sigma^{1-n} 
\end{equation}
where $C = ((1-\gamma)/320)^{1-n}$.  Let $\{\psi_j\}$ be a smooth partition of unity of $B_{\gamma}(0,P_i) \cap \{r > \tau\}$ subordinate to $\{\mathbf{A}_{\rho_j,1/4}(\zeta_j)\}$ such that 
\begin{equation}\label{graphrep eqn18}
	0 \leq \psi_i \leq 1, \quad \op{spt} \psi_i \subseteq \mathbf{A}_{\rho_j,1/4}(\zeta_j), \quad 
	|\nabla \psi_j| \leq \frac{C(n,\gamma)}{\rho_j}, \quad \sum_{j=1}^N \psi_i = 1 . 
\end{equation}

Recalling \eqref{graphrep eqn2}, for each $j \in \{1,2,\ldots,N\}$ and $i \in \{1,2,\ldots,p\}$ we can apply Lemma~\ref{weak graphrep lemma}(a) with $\eta_{(0,\zeta_j),\rho_j \#} T$ in place of $T$ to find integers $q_{j,i}$ with $\sum_{i=1}^p |q_{j,i}| \leq q$ and $n$-dimensional locally area minimizing rectifiable currents $T_{j,i}$ of $\mathbf{A}_{\rho_j,7/8}(\zeta_j)$ such that  
\begin{gather}
	\label{graphrep eqn19} T \llcorner \mathbf{A}_{\rho_j,7/8}(\zeta_j) = \sum_{i=1}^p T_{j,i} , \\ 
	\label{graphrep eqn20} (\partial T_{j,i}) \llcorner \mathbf{A}_{\rho_j,7/8}(\zeta_j) = 0, \\ 
	\label{graphrep eqn21} (\pi_{P_i \#} T_{j,i}) \llcorner \mathbf{A}_{\rho_j,3/4}(\zeta_j) 
		= q_{j,i} \llbracket P_i \cap \mathbf{A}_{\rho_j,3/4}(\zeta_j) \rrbracket , \\
	\label{graphrep eqn22} \sup_{X \in \op{spt} T_{j,i}} \op{dist}(X, P_i) \leq C \rho_j^{-n/2} E , 
\end{gather}
where $E = E(T,\mathbf{C},\mathbf{B}_1(0))$ and $C = C(n,m,q,p,\gamma) \in (0,\infty)$ is a constant.  By \eqref{graphrep eqn3}, \eqref{graphrep eqn19}, \eqref{graphrep eqn21}, \eqref{graphrep eqn22}, and Lemma~\ref{energy est lemma} 
\begin{align*}
	(q-1/2) \,\mathcal{L}^n(A_{\rho_j,1/2}(0)) &< \|T\|(\mathbf{A}_{\rho_j,1/2}(\zeta_j)) 
		\\&= \sum_{i=1}^p \|T_{j,i}\|(\mathbf{A}_{\rho_j,1/2}(\zeta_j)) \leq \sum_{i=1}^p |q_{j,i}| \,\mathcal{L}^n(A_{\rho_j,1/2}(0)) + C \rho_j^n E^2 ,
\end{align*}
where $C = C(n,m,q,p,\gamma,\tau) \in (0,\infty)$ is a constant.  Thus recalling $E < \tau^{-(n+2)/2} \varepsilon_0$ and assuming $\varepsilon_0$ is sufficiently small, we obtain $q-1/2 \leq \sum_{i=1}^q |q_{j,i}| + 1/4$.  In other words, since $q_{i,j}$ are integers and $\sum_{i=1}^p |q_{j,i}| \leq q$, we must have $\sum_{i=1}^p |q_{j,i}| = q$.  Notice that if $\mathbf{A}_{\rho_j,7/8}(\zeta_j) \cap \mathbf{A}_{\rho_k,7/8}(\zeta_k) \neq \emptyset$ then by \eqref{graphrep eqn19} and \eqref{graphrep eqn22} 
\begin{equation*}
	T_{j,i} = T \llcorner \big\{ X : \op{dist}(X, P_i) \leq C \,(\tau/4)^{-n/2} E \} = T_{k,i} 
		\text{ in } \mathbf{A}_{\rho_j,7/8}(\zeta_j) \cap \mathbf{A}_{\rho_k,7/8}(\zeta_k) 
\end{equation*}
for all $i \in \{1,2,\ldots,p\}$, where $C$ is as in \eqref{graphrep eqn22} and using \eqref{main hyp eqn3} and \eqref{graphrep concl a1} we assume that $\beta_0$ is small enough that $6C \,(\tau/4)^{-(n+2)/2} E < \op{minsep} \mathbf{C}$.  Similarly, if $\mathbf{A}_{\rho_j,3/4}(\zeta_j) \cap \mathbf{A}_{\rho_k,3/4}(\zeta_k) \neq \emptyset$ then by \eqref{graphrep eqn21} we have $q_{j,i} = q_{k,i}$ for all $i$.  It follows that for each $i$ there is a well-defined rectifiable current $T_i$ such that $T_i = T_{j,i}$ in $\mathbf{A}_{\rho_j,7/8}(\zeta_j)$ for all $j$.  Since $\mathbf{B}_{\gamma}(0) \cap \{r > \tau/4\}$ is connected, 
we can take $q_{j,i} = q_i$ for all $j$ and $i$.  (In particular, $T_i$ and $q_i$ are defined independent of $j$.)  For each $i$, after possibly reversing the orientation of $P_i$, we may assume that $q_i \geq 0$.  By Lemma~\ref{graphrep annuli lemma}(b) with $\eta_{0,1/4 \#} T$ in place of $T$, we obtain $q_i > 0$ (in $(\pi_{P_i \#} T_i) \llcorner \mathbf{A}_{1/4,3/4}(0) = q_i \llbracket P_i \cap \mathbf{A}_{1/4,3/4}(0) \rrbracket$).  One readily verifies from this construction and \eqref{graphrep eqn19}--\eqref{graphrep eqn22} that $T_i$ satisfies Theorem~\ref{graphrep thm}(b). 

For each $j \in \{1,2,\ldots,N\}$ and $i \in \{1,2,\ldots,p\}$ we can apply Lemma~\ref{weak graphrep lemma}(b) with $\eta_{(0,\zeta_j),\rho_j \#} T$ to find Lipschitz $q_i$-valued functions $u_{j,i} : P_i \cap \mathbf{A}_{\rho_j,1/2}(\zeta_j) \rightarrow \mathcal{A}_{q_i}(P_i^{\perp})$ and closed sets $K_{j,i} \subseteq P_i \cap \mathbf{A}_{\rho_j,1/2}(\zeta_j)$ such that 
\begin{gather} 
	\label{graphrep eqn23} T_i \llcorner \pi_{P_i}^{-1}(K_{j,i}) = (\op{graph} u_{j,i}) \llcorner \pi_{P_i}^{-1}(K_{j,i}) , \\
	\label{graphrep eqn24} \mathcal{H}^n(P_i \cap \mathbf{A}_{\rho_j,1/2}(\zeta_j) \setminus K_{j,i}) 
		+ \|T_i\|(\pi_{P_i}^{-1}(P_i \cap \mathbf{A}_{\rho_j,1/2}(\zeta_j) \setminus K_{j,i})) \leq C \rho_j^{-2-(n+2)\alpha/2} E^{2+\alpha} , \\
	\label{graphrep eqn25} \sup_{P_i \cap \mathbf{A}_{\rho_j,1/2}(\zeta_j)} |u_{j,i}| \leq C\rho_j^{-n/2} E , \quad 
		\sup_{P_i \cap \mathbf{A}_{\rho_j,1/2}(\zeta_j)} |\nabla u_{j,i}| \leq C \rho_j^{-(n+2)\alpha/2} E^{\alpha}, 
\end{gather}
where $\alpha = \alpha(n,m,q) \in (0,1)$ and $C = C(n,m,q,p,\gamma) \in (0,\infty)$ are constants.  By~\cite[Definition~1.1(6) and Theorem~1.3]{Almgren}, there exists an integer $L(q,m) \geq 1$ and bi-Lipschitz embedding $\boldsymbol{\xi} : \mathcal{A}_q(\mathbb{R}^m) \rightarrow \mathbb{R}^L$ and $\boldsymbol{\rho} : \mathbb{R}^L \rightarrow \mathcal{Q}$ such that $\op{Lip}\boldsymbol{\xi} \leq 1$, $\op{Lip}\boldsymbol{\xi}^{-1}|_{\mathcal{Q}} \leq C(m,q)$, and $\op{Lip}\boldsymbol{\rho} \leq C(m,q)$, where $\mathcal{Q} = \boldsymbol{\xi}(\mathcal{A}_q(\mathbb{R}^m))$.  For each $i \in \{1,2,\ldots,p\}$ define 
\begin{gather}
	\label{graphrep eqn26} K_i = B_{\gamma}(0,P_i) \cap \{r > \tau\} \setminus \left( \bigcup_{j=1}^N (P_i \cap \mathbf{A}_{\rho_j,1/4}(\zeta_j) \setminus K_{j,i}) 
		\right) , \\
	\label{graphrep eqn27} u_i(X) = (\boldsymbol{\xi}^{-1} \circ \boldsymbol{\rho})\left(\sum_{j=1}^N \psi_j(X) \,\boldsymbol{\xi}(u_{j,i}(X)) \right) 
		\text{ for all } X \in B_{\gamma}(0,P_i) \cap \{r > \tau\} , 
\end{gather}
where $\{\psi_j\}$ is the smooth partition of unity of $B_{\gamma}(0,P_i) \cap \{r > \tau\}$ subordinate to $\{\mathbf{A}_{\rho_j,1/4}(\zeta_j)\}$ satisfying \eqref{graphrep eqn18}.  Arguing as in~\cite[Theorem~3.8]{KrumWica}, it follows from \eqref{graphrep eqn24} that if $\mathbf{A}_{\rho_j,1/2}(\zeta_j) \cap \mathbf{A}_{\rho_k,1/2}(\zeta_k) \neq \emptyset$ then $\mathcal{H}^n(K_{j,i} \cap K_{k,i}) > 0$.  Let $Z \in K_{j,i} \cap K_{k,i}$ and note that by \eqref{graphrep eqn23} we have $u_{j,i}(Z) = u_{k,i}(Z)$.  Thus by \eqref{graphrep eqn25} and $\rho_k \leq \tfrac{65-\gamma}{63+\gamma} \,\rho_j$ 
\begin{align}\label{graphrep eqn28}
	&\sup_{P_i \cap \mathbf{A}_{\rho_j,1/4}(\zeta_j) \cap \mathbf{A}_{\rho_k,1/4}(\zeta_k)} \mathcal{G}(u_{j,i},u_{k,i}) 
	\\ \leq\,& \sup_{P_i \cap \mathbf{A}_{\rho_j,1/2}(\zeta_j)} \mathcal{G}(u_{j,i},u_{j,i}(Z)) 
		+ \sup_{P_i \cap \mathbf{A}_{\rho_k,1/2}(\zeta_j)} \mathcal{G}(u_{k,i},u_{k,i}(Z)) 
	\leq C \rho_j^{1-(n+2)\alpha/2} E^{\alpha} \nonumber 
\end{align}
where $C = C(n,m,q,p,\gamma) \in (0,\infty)$ is a constant.  By \eqref{graphrep eqn23}, \eqref{graphrep eqn26}, and \eqref{graphrep eqn27}, \eqref{graphrep concl c1} holds true.  Recalling that $\{\mathbf{A}_{\rho_j,1/4}(\zeta_j)\}$ covers $\mathbf{B}_{(3+\gamma)/4}(0) \cap \{r > \sigma\}$, $\rho_j \geq \sigma/2$ for all $j \in \{1,2,\ldots,J_{\sigma}\}$, and $J_{\sigma} \leq C(n,\gamma)\,\sigma^{1-n}$ (as in \eqref{graphrep eqn17}), it follows from \eqref{graphrep eqn24} and \eqref{graphrep eqn26} that \eqref{graphrep concl c2} holds true.  It follows from \eqref{graphrep eqn27}, \eqref{graphrep eqn18}, \eqref{graphrep eqn25}, and \eqref{graphrep eqn28} that \eqref{graphrep concl c3} holds true.
\end{proof}

\subsection{Initial a priori estimates} \label{sec:keyest sec}

Here and in Section~\ref{sec:apriori close2plane} and Section~\ref{sec:apriori nonplanar}, we establish several key integral estimates for locally area-minimizing rectifiable currents $T$ close to a sum-of-planes $\mathbf{C} \in \mathcal{C}_{q,p}$.  These estimates are inspired by the results 
of~\cite{Sim93} for stationary varifolds in a ``multiplicity 1 class'''. We note that the multiplicity 1 class hypothesis in \cite{Sim93} in particular rules out branch points a priori. In contrast to this, in the present setting higher multiplicity and branch points are permitted, and so the proofs of the estimates in the present setting require additional arguments and strategies, some of which are adaptations of arguments in \cite{Wic14}, \cite{KrumWic2}.

\begin{theorem} \label{keyest thm}
For each $\gamma \in (0,1)$ there exists $\varepsilon_0 = \varepsilon_0(n,m,q,\gamma) \in (0,1)$ and $\beta_0 = \beta_0(n,m,q,\gamma) \in (0,1)$ such that if $\mathbf{C}$ and $T$ satisfy Hypothesis~$(\star)$ and Hypothesis~$(\star\star)$, then: 
\begin{align*} 
	&(a) \quad \int_{{\mathbf B}_{\gamma}(0)} \frac{|X^{\perp}|^2}{|X|^{n+2}} \,d\|T\|(X) 
		\leq C \int_{{\mathbf B}_1(0)} \op{dist}^2(X, \op{spt} \mathbf{C}) \,d\|T\|(X) , \\
	&(b) \quad \int_{{\mathbf B}_{\gamma}(0)} \sum_{j=1}^{n-2} |e_{m+2+j}^{\perp}|^2 \,d\|T\|(X) 
		\leq C \int_{{\mathbf B}_1(0)} \op{dist}^2(X, \op{spt} \mathbf{C}) \,d\|T\|(X) , 
\end{align*}
where $(\,\cdot\,)^{\perp}$ denotes orthogonal projection onto the orthogonal complement of the approximate tangent plane to $T$ at $X$ and $C = C(n,m,q,\gamma) \in (0,\infty)$ is a constant. 
\end{theorem}

\begin{proof} 
Without loss of generality assume that $\gamma \geq 1/2$.  Express each point $X \in \mathbb{R}^{n+m}$ as $X = (x,y)$ for $x \in \mathbb{R}^{m+2}$ and $y \in \mathbb{R}^{n-2}$.  Let $r = r(X) = |x|$ and $R = R(X) = |X|$ for each $X = (x,y) \in \mathbb{R}^{n+m}$.  Let $\psi : \mathbb{R} \rightarrow [0,1]$ be a decreasing smooth function such that $\psi(t) = 1$ for all $t \leq \gamma$, $\psi(t) = 0$ for all $t \geq (1+\gamma)/2$, $|\psi'(t)| \leq 6/(1-\gamma)$, and $|\psi''(t)| \leq 36/(1-\gamma)^2$.  We have by the inequalities (2) and (3) of Lemma~3.4 of~\cite{Sim93} that 
\begin{align}
	\label{keyest eqn1} &\int_{\mathbf{B}_{\gamma}(0)} R^{-n-2} |X^{\perp}|^2 \,d\|T\|(X) 
		\leq C \left( \int_{\mathbf{B}_1(0)} \psi^2(R) \,d\|T\|(X) - \int_{\mathbf{B}_1(0)} \psi^2(R) \,d\|\mathbf{C}\|(X) \right) , \\
	\label{keyest eqn2} &\int_{\mathbf{B}_1(0)} \left( 2 + \frac{1}{2} \sum_{j=1}^{n-2} |e_{m+2+j}|^2 \right) \psi^2(R) \,d\|T\|(X) 
		\\ \leq\,& C \int_{\mathbf{B}_1(0)} |(x,0)^{\perp}|^2 \,(\psi^2(R) + (\psi'(R))^2) \,d\|T\|(X) 
		- 2 \int_{\mathbf{B}_1(0)} \frac{r^2}{R} \,\psi(R) \,\psi'(R) \,d\|T\|(X) \nonumber 
\end{align}
for some constant $C = C(n,\gamma) \in (0,\infty)$.  By the identity (6) of Lemma~3.4 of~\cite{Sim93}, 
\begin{equation} \label{keyest eqn3} 
	2 \int_{\mathbf{B}_1(0)} \psi^2(R) \,d\|\mathbf{C}\|(X) = - 2 \int_{\mathbf{B}_1(0)} \frac{r}{R^2} \,\psi(R) \,\psi'(R) \,d\|\mathbf{C}\|(X) . 
\end{equation}
Combining \eqref{keyest eqn1}, \eqref{keyest eqn2}, and \eqref{keyest eqn3} 
\begin{align} \label{keyest eqn4} 
	&\int_{\mathbf{B}_{\gamma}(0)} R^{-n-2} |X^{\perp}|^2 \,d\|T\|(X) + \int_{\mathbf{B}_{\gamma}(0)} \sum_{j=1}^{n-2} |e_{m+2+j}|^2 \,d\|T\|(X)
	\\ \leq\,& C \left( \int_{\mathbf{B}_1(0)} |(x,0)^{\perp}|^2 \,(\psi^2(R) + (\psi'(R))^2) \,d\|T\|(X) \right. \nonumber
	\\&\left. - 2 \int_{\mathbf{B}_1(0)} \frac{r^2}{R} \,\psi(R) \,\psi'(R) \,d\|T\|(X) 
		+ 2 \int_{\mathbf{B}_1(0)} \frac{r^2}{R} \,\psi(R) \,\psi'(R) \,d\|\mathbf{C}\|(X) \right) \nonumber
\end{align}
for some constant $C = C(n,\gamma) \in (0,\infty)$.  Hence to prove the theorem it suffices to bound the right-hand side of \eqref{keyest eqn4} above by $C \int_{\mathbf{B}_1(0)} \op{dist}^2(X,\op{spt} \mathbf{C}) \,d\|T\|(X)$ for some constant $C = C(n,m,q,\gamma) \in (0,\infty)$. 

Let $\delta = \delta(n,m,q,\gamma) > 0$ to be later determined.  Let $U$ be the union of all annuli $\mathbf{A}_{\rho,1/20}(\zeta)$ such that $\rho^2 + |\zeta|^2 < (3+\gamma)^2/16$ and 
\begin{equation} \label{keyest eqn5} 
	\|T\|(\mathbf{A}_{\rho,1}(\zeta)) \leq (q+1/2) \,\mathcal{L}^n(A_{\rho,1}(0)) , \quad 
	E(T, \mathbf{C}, \mathbf{A}_{\rho,2}(\zeta)) < \delta . 
\end{equation}
As we will see below, $U$ is region where we can apply Lemma~\ref{weak graphrep lemma}, Lemma~\ref{graphrep annuli lemma}, and Remark~\ref{tildeC rmk} to obtain a Lipschitz approximation of $T$ relative to $\mathbf{C}$ or some other cone in $\mathcal{C}_{q,p}$.  Clearly $U$ is open.  By \eqref{main hyp eqn2} and \eqref{graphrep eqn4} from the proof of Theorem~\ref{graphrep thm} (with $\tau = (1-\gamma)/80$), we may assume $\varepsilon_0$ is small enough that \eqref{keyest eqn5} holds true whenever $\rho > (1-\gamma)/80$ and $\rho^2 + |\zeta|^2 < (3+\gamma)^2/16$ and thus 
\begin{equation} \label{keyest eqn6} 
	\mathbf{B}_{(3+\gamma)/4}(0) \cap \{r > (1-\gamma)/80\} \subset U .
\end{equation}
Define the locally Lipschitz cutoff function $\chi : \mathbf{B}_{(1+\gamma)/2}(0) \setminus (\{0\} \times \mathbb{R}^{n-2}) \rightarrow \mathbb{R}$ by 
\begin{equation*}
	\chi(x,y) = \begin{cases} 
		1 &\text{ if } \op{dist}((x,y),\partial U) \geq |x|/2 \\
		\frac{4}{|x|} \op{dist}((x,y),\partial U) - 1 &\text{ if } |x|/4 < \op{dist}((x,y),\partial U) < |x|/2 \\
		0 &\text{ if } \op{dist}((x,y),\partial U) \leq |x|/4 .
	\end{cases}
\end{equation*}
Observe that $0 \leq \chi \leq 1$ and $|\nabla \chi(x,y)| \leq 6/|x|$.  Recall that $\mathbf{A}_{\rho,\kappa}(\zeta) = \{ (r \omega, y) : (r,y) \in B^{n-1}_{\kappa (1-\gamma)\rho/8}(\rho,\zeta), \,\omega \in \mathbb{S}^{m+1} \}$ for each $\kappa \in (0,1]$, $\rho > 0$, and $\zeta \in \mathbb{R}^{n-2}$.  Note that $\mathbf{A}_{\rho,\kappa}(\zeta) \cap \mathbf{A}_{\rho',\kappa}(\zeta') \neq \emptyset$ $\Leftrightarrow$ $|(\rho,\zeta) - (\rho',\zeta')| < \tfrac{1}{8} \,\kappa(1-\gamma)(\rho+\rho')$, in which case $\tfrac{8-\kappa(1-\gamma)}{8+\kappa(1-\gamma)} \,\rho < \rho' < \tfrac{8+\kappa(1-\gamma)}{8-\kappa(1-\gamma)}\,\rho$.  By applying the Vitali covering lemma, there exists a countable collection $\mathcal{I}$ of $(\rho,\zeta)$ with $\rho > 0$, $\zeta \in \mathbb{R}^{n-2}$, and $\rho^2 + |\zeta|^2 < (3+\gamma)^2/16$ such that \eqref{keyest eqn5} holds true for each $(\rho,\zeta) \in \mathcal{I}$, $\{\mathbf{A}_{\rho,1/20}(\zeta)\}_{(\rho,\zeta) \in \mathcal{I}}$ a collection of pairwise disjoint annuli, and 
\begin{equation*}
	U \subset \bigcup_{(\rho,\zeta) \in \mathcal{I}} \mathbf{A}_{\rho,1/4}(\zeta).
\end{equation*}
Observe that if $\mathbf{A}_{\rho,1}(\zeta) \cap \mathbf{A}_{\rho',1}(\zeta') \neq \emptyset$ then $|(\rho,\zeta) - (\rho',\zeta')| < \tfrac{1}{8} \,(1-\gamma)(\rho+\rho')$ and thus $\tfrac{7+\gamma}{9-\gamma} \,\rho < \rho' < \tfrac{9-\gamma}{7+\gamma}\,\rho$.  Hence 
\begin{equation*}
	B^{n-1}_{\tfrac{(1-\gamma)(7+\gamma)\rho}{160(9-\gamma)}}(\rho',\zeta') 
	\subset B^{n-1}_{\tfrac{(1-\gamma)\rho'}{160}}(\rho',\zeta') 
	\subset B^{n-1}_{\tfrac{(1-\gamma)\rho}{160} \Big( 20 + \tfrac{21(9-\gamma)}{7+\gamma} \Big)}(\rho,\zeta)
\end{equation*} 
which since $\{B^{n-1}_{(1-\gamma)\rho'/160}(\rho',\zeta')\}_{(\rho',\zeta') \in \mathcal{I}}$ is a pairwise disjoint collection of balls implies that 
\begin{equation*}
	\#\{ (\rho',\zeta') \in \mathcal{I} : \mathbf{A}_{\rho,1}(\zeta) \cap \mathbf{A}_{\rho',1}(\zeta') \neq \emptyset \} \leq C(n,\gamma) .
\end{equation*}
Thus there is an integer $N \leq C(n,\gamma)$ and pairwise disjoint sets $\mathcal{I}_1,\mathcal{I}_2,\ldots,\mathcal{I}_N \subset \mathcal{I}$ such that $\mathcal{I} = \bigcup_{j=1}^N \mathcal{I}_j$ and $\{\mathbf{A}_{\rho,1}(\zeta)\}_{(\rho,\zeta) \in \mathcal{I}_j}$ is a collection of pairwise disjoint annuli for each $j = 1,2,\ldots,N$.  Let $\{\phi_{(\rho,\zeta)}\}_{(\rho,\zeta) \in \mathcal{I}}$ is a smooth partition of unity subordinate to $\{\mathbf{A}_{\rho,1/2}(\zeta)\}_{(\rho,\zeta) \in \mathcal{I}}$ such that 
\begin{gather}\label{keyest eqn8} 
	\phi_{(\rho,\zeta)}(x,y) = \phi_{(\rho,\zeta)}(\widetilde{x},y) \text{ whenever } |x| = |\widetilde{x}| , \\
	\op{spt} \phi_{(\rho,\zeta)} \subset \mathbf{A}_{\rho,1/2}(\zeta), \quad
	|\nabla \widetilde{\phi}_{(\rho,\zeta)}| = |\nabla \phi_{(\rho,\zeta)}| \leq \frac{C(n,\gamma)}{\rho} , \nonumber \\
	\sum_{(\rho,\zeta) \in \mathcal{I}} \phi_{(\rho,\zeta)} = 1 \text{ on } \mathbf{B}_{(1+\gamma)/2}(0). \nonumber
\end{gather}
We claim that for each $(\rho,\zeta) \in \mathcal{I}$ 
\begin{align} \label{keyest eqn9} 
	&\int_{\mathbf{A}_{\rho,1/2}(\zeta)} |(x,0)^{\perp}|^2 \,(\psi^2(R) + (\psi'(R))^2) \,\phi_{(\rho,\zeta)} \,\chi \,d\|T\|(X) 
	\\&- 2 \int_{\mathbf{B}_1(0)} \frac{r^2}{R} \,\psi(R) \,\psi'(R) \,\phi_{(\rho,\zeta)} \,\chi \,d\|T\|(X) 
		+ 2 \int_{\mathbf{B}_1(0)} \frac{r^2}{R} \,\psi(R) \,\psi'(R) \,\phi_{(\rho,\zeta)} \,\chi \,d\|\mathbf{C}\|(X) \nonumber
	\\ \leq\,& C_0 \int_{\mathbf{A}_{\rho,1}(\zeta)} \op{dist}^2(X, \op{spt}\mathbf{C}) \,d\|T\|(X) \nonumber
\end{align}
for some constant $C_0 = C_0(n,m,q,\gamma) \in (0,\infty)$.  Then by summing \eqref{keyest eqn9} over $(\rho,\zeta) \in \mathcal{I} = \bigcup_{j=1}^N \mathcal{I}_j$ and keeping in mind that $\op{spt} \psi \subset \mathbf{B}_{(1+\gamma)/2}(0)$, we deduce that 
\begin{align} \label{keyest eqn10} 
	&\int_{\mathbf{B}_1(0)} |(x,0)^{\perp}|^2 \,(\psi^2(R) + (\psi'(R))^2) \,\chi \,d\|T\|(X) 
	\\&- 2 \int _{\mathbf{B}_1(0)}\frac{r^2}{R} \,\psi(R) \,\psi'(R) \,\chi \,d\|T\|(X) 
		+ 2 \int_{\mathbf{B}_1(0)} \frac{r^2}{R} \,\psi(R) \,\psi'(R) \,\chi \,d\|\mathbf{C}\|(X) \nonumber
	\\ \leq\,& C \int_{\mathbf{B}_1(0)} \op{dist}^2(X, \op{spt}\mathbf{C}) \,d\|T\|(X) \nonumber
\end{align}
for some constant $C = C(n,m,q,\gamma) \in (0,\infty)$.  

To see \eqref{keyest eqn9}, fix $(\rho,\zeta) \in \mathcal{I}$ and assume that $\op{spt} T \cap \mathbf{A}_{\rho,1/2}(\zeta) \neq \emptyset$.  Let $\widetilde{\beta} = \widetilde{\beta}(n,m,q,\gamma) \in (0,1)$ to be later determined.  Suppose that 
\begin{equation}\label{keyest eqn11} 
	E(T,\mathbf{C}, \mathbf{A}_{\rho,1}(\zeta)) \leq \widetilde{\beta} 
		\inf_{\mathbf{C}' \in \bigcup_{p'=1}^q \mathcal{C}_{q,p'}} Q(T,\mathbf{C}', \mathbf{A}_{\rho,1}(\zeta))
\end{equation}
Let $\mathbf{C} = \sum_{i=1}^p q_i \llbracket P_i \rrbracket$ for some integers $q_i \geq 1$ with $\sum_{i=1}^p q_i = q$ and $n$-dimensional oriented planes $P_i$ with $\{0\} \times \mathbb{R}^{n-2} \subset P_i$.  Note that by Theorem~\ref{graphrep thm}(a), \eqref{graphrep concl a2} holds true.  Thus provided $\delta$ and $\widetilde{\beta}$ are sufficiently small, by $(\partial T) \llcorner \mathbf{B}_1(0) = 0$, \eqref{keyest eqn5}, \eqref{keyest eqn11}, and \eqref{graphrep concl a2} we can apply Lemma~\ref{weak graphrep lemma} to show the following.  By Lemma~\ref{weak graphrep lemma}(a) there are integers $\widehat{q}_i$ with $\sum_{i=1}^p |\widehat{q}_i| \leq q$ and $n$-dimensional locally area minimizing rectifiable currents $T_i$ of $\mathbf{A}_{\rho,7/8}(\zeta)$ such that 
\begin{gather*}
	T \llcorner \mathbf{A}_{\rho,7/8}(\zeta) = \sum_{i=1}^p T_i , \quad 
	(\partial T_i) \llcorner \mathbf{A}_{\rho,7/8}(\zeta) = 0, \\ 
	(\pi_{P_i \#} T_i) \llcorner \mathbf{A}_{\rho,3/4}(\zeta) = \widehat{q}_i \llbracket P_i \rrbracket \llcorner \mathbf{A}_{\rho,3/4}(\zeta) , \nonumber \\
	\sup_{X \in \op{spt} T_i} \op{dist}(X, P_i) \leq C \rho E , \nonumber
\end{gather*}
where $E = E(T,\mathbf{C},\mathbf{A}_{\rho,1}(\zeta))$ and $C = C(n,m,q,\gamma) \in (0,\infty)$ is a constant.  By Lemma~\ref{weak graphrep lemma}(b), for each $i \in \{1,2,\ldots,p\}$ with $\widehat{q}_i \neq 0$ there exists Lipschitz $|\widehat{q}_i|$-valued functions $u_i : P_i \cap \mathbf{A}_{\rho,1/2}(\zeta) \rightarrow \mathcal{A}_{|\widehat{q}_i|}(P_i^{\perp})$ and closed sets $K_i \subseteq P_i \cap \mathbf{A}_{\rho,1/2}(\zeta)$ such that 
\begin{gather} 
	\label{keyest eqn13} T_i \llcorner \pi_{P_i}^{-1}(K_i) = (\op{graph} u_i) \llcorner \pi_{P_i}^{-1}(K_i) ,  \\
	\mathcal{H}^n(P_i \cap \mathbf{A}_{\rho,1/2}(\zeta) \setminus K_i) + \|T_i\|(\pi_{P_i}^{-1}(P_i \cap \mathbf{A}_{\rho,1/2}(\zeta) \setminus K_i)) 
		\leq C \rho^n E^{2+\alpha} , \nonumber \\
	\sup_{P_i \cap \mathbf{A}_{\rho,1/2}(\zeta)} |u_i| \leq C \rho E , \quad \op{Lip} u_i \leq C E^{\alpha}, \nonumber
\end{gather}
where $\alpha = \alpha(n,m,q) \in (0,1)$ and $C = C(n,m,q,\gamma) \in (0,\infty)$ are constants.  For each $X \in P_i \cap \mathbf{A}_{\rho,1/2}(\zeta)$, let $u_i(X) = \sum_{j=1}^{|\widehat{q}_i|} \llbracket u_{i,j}(X) \rrbracket$ where $u_{i,j}(X) \in \mathbb{R}^m$.  

Take any $i \in \{1,2,\ldots,p\}$ with $\widehat{q}_i \neq 0$.  By Rademacher's Theorem~\cite[Theorem~1.13]{DeLSpaDirMin}, $u_i$ is differentiable at $\mathcal{H}^n$-a.e.~$(x',y) \in K_i$ in the sense of~\cite[Definition~1.9]{DeLSpaDirMin}.  Let $u_i$ is differentiable at $(x',y) \in K_i$ and let $X = (x,y) = (x',y) + u_{i,j}(x',y)$ be a point on $\op{spt} T_i \cap \pi_{P_i}^{-1}(K_i)$.  Notice that  
\begin{equation*}
	(x,0)^{\perp} = u_{i,j}(x',y) - (\pi_X - \pi_{P_i})(x,0)
\end{equation*}
where $\pi_X$ is the orthogonal projection map onto the approximate tangent plane to $T$ at $X$ and satisfies $\|\pi_X - \pi_{P_i}\| \leq C(n,m) \,|\nabla u_{i,j}(x',y)|$.  This together with \eqref{keyest eqn13} and Lemma~\ref{energy est lemma} gives us 
\begin{align}\label{keyest eqn14} 
	&\hspace{-.3in}\int_{\mathbf{A}_{\rho,1/2}(\zeta)} |(x,0)^{\perp}|^2 \,(\psi^2(R) + (\psi'(R))^2) \,\phi_{(\rho,\zeta)}(x,y) \,d\|T_i\|(x,y) 
	\\ \leq\,& C \int_{P_i \cap \mathbf{A}_{\rho,1/4}(\zeta) \cap K_i} (|u_i|^2 + \rho^2 |\nabla u_i|^2) \,d\mathcal{H}^n(x',y) + C E^{2+\alpha} \nonumber
	\\ \leq\,& C \int_{\mathbf{A}_{\rho,1/2}(\zeta)} (\op{dist}^2(X,P_i) + \rho^2 |\vec T_i - \vec P_i|^2) \,d\|T_i\|(X) + C E^{2+\alpha} \nonumber
	\\ \leq\,& C \int_{\mathbf{A}_{\rho,7/8}(\zeta)} \op{dist}^2(X,P_i) \,d\|T_i\|(X) + C E^{2+\alpha} 
	\leq C \int_{\mathbf{A}_{\rho,1}(\zeta)} \op{dist}^2(X, \op{spt} \mathbf{C}) \,d\|T\|(X), \nonumber
\end{align}
where $C = C(n,m,q,\gamma) \in (0,\infty)$ are constants.  Next define $F : B^n_{(1+\gamma)/2}(0) \rightarrow [0,\infty)$ by 
\begin{align*}
	F(r,s,y) = -\frac{r^2 + s^2}{\sqrt{r^2 + s^2 + |y|^2}} &\,\psi\hspace{-0.8mm}\left(\sqrt{r^2 + s^2 + |y|^2}\right) 
		\psi'\hspace{-0.8mm}\left(\sqrt{r^2 + s^2 + |y|^2}\right) \\&\cdot \phi_{(\rho,\zeta)}\hspace{-0.8mm}\left(\sqrt{r^2 + s^2}, y\right) 
		\chi\hspace{-0.8mm}\left(\sqrt{r^2 + s^2}, y\right)
\end{align*}
for each $r,s \in \mathbb{R}$ and $y \in \mathbb{R}^{n-2}$ with $r^2 + s^2 + |y|^2 < (1+\gamma)^2/4$.  Here by a slight abuse of notation we let $\phi_{(\rho,\zeta)}(|x|,y)$ and $\chi(|x|,y)$ denote the values of $\phi_{(\rho,\zeta)}(x,y)$ and $\chi(x,y)$ respectively.  By the definition of $\psi$, $|\nabla \chi(x,y)| \leq 6/|x|$, and \eqref{keyest eqn8}, 
\begin{equation}\label{keyest eqn15}
	F(r,s,y) \leq C \rho^2, \quad |F(r,s,y) - F(r,0,y)| \leq C s^2 
\end{equation}
for each $r,s \in \mathbb{R}$ and $y \in \mathbb{R}^{n-2}$ with $r^2 + s^2 + |y|^2 < (1+\gamma)^2/4$, where $C = C(n,\gamma) \in (0,\infty)$ is a constant.  Hence by \eqref{keyest eqn13}, \eqref{keyest eqn15}, and Lemma~\ref{energy est lemma}, 
\begin{align}\label{keyest eqn16}
	&-2 \int_{\mathbf{A}_{\rho,1/2}(\zeta)} \frac{r^2}{R} \,\psi(R) \,\psi'(R) \,\phi_{(\rho,\zeta)} \,\chi \,d\|T_i\|(X) 
	\\&+ 2 |\widehat{q}_i| \int_{P_i \cap \mathbf{A}_{\rho,1/2}(\zeta)} \frac{r^2}{R} \,\psi(R) \,\psi'(R) \,\phi_{(\rho,\zeta)} \,\chi \,d\mathcal{H}^n(X) \nonumber 
	\\ \leq\,& 2 \int_{P_i \cap \mathbf{A}_{\rho,1/2}(\zeta)} \sum_{j=1}^{|\widehat{q}_i|} (F(r,|u_{i,j}(x,y)|,y) \,(1 + |\nabla u_{j,i}(x,y)|^2) - F(r,0,y)) 
		\,d\mathcal{H}^n(x,y) + CE^{2+\alpha} \nonumber
	\\ \leq\,& C \int_{P_i \cap \mathbf{A}_{\rho,1/2}(\zeta)} (|u_i|^2 + \rho^2 |\nabla u_i|^2) \,d\mathcal{H}^n(x,y) + CE^{2+\alpha} \nonumber
	\\ \leq\,& C \int_{\mathbf{A}_{\rho,1}(\zeta)} \op{dist}^2(X,\op{spt} \mathbf{C}) \,d\|T\|(X) . \nonumber
\end{align}
Note that by symmetry, the value of 
\begin{equation*}
	\int_{P \cap \mathbf{A}_{\rho,1/2}(\zeta)} \frac{r^2}{R} \,\psi(R) \,\psi'(R) \,\phi_{(\rho,\zeta)} \,\chi \,d\mathcal{H}^n(X)
\end{equation*}
is the same for all $n$-dimensional planes $P$ with $\{0\} \times \mathbb{R}^{n-2} \subset P$ and thus using $\sum_{i=1}^p |\widehat{q}_i| \leq q$
\begin{align}\label{keyest eqn17}
	&\int_{\mathbf{A}_{\rho,1/2}(\zeta)} \frac{r^2}{R} \,\psi(R) \,\psi'(R) \,\phi_{(\rho,\zeta)} \,\chi \,d\|\mathbf{C}\|(X) 
	\\ \leq\,& \sum_{i=1}^p |\widehat{q}_i| \int_{P_i \cap \mathbf{A}_{\rho,1/2}(\zeta)} \frac{r^2}{R} \,\psi(R) \,\psi'(R) \,\phi_{(\rho,\zeta)} 
		\,\chi \,d\mathcal{H}^n(X) . \nonumber 
\end{align}
Summing \eqref{keyest eqn14} and \eqref{keyest eqn16} over $i = 1,2,\ldots,p$ and using \eqref{keyest eqn17} gives us \eqref{keyest eqn9}. 

Suppose instead that 
\begin{equation*} 
	E(T,\mathbf{C}, \mathbf{A}_{\rho,1}(\zeta)) > \widetilde{\beta} 
		\inf_{\mathbf{C}' \in \bigcup_{p'=1}^q \mathcal{C}_{q,p'}} Q(T,\mathbf{C}', \mathbf{A}_{\rho,1}(\zeta)) . 
\end{equation*}
By Remark~\ref{tildeC rmk}, there exists an integer $1 \leq \widetilde{p} < p$ and $\widetilde{\mathbf{C}} \in \mathcal{C}_{q,\widetilde{p}}$ such that 
\begin{gather}
	Q(T, \widetilde{\mathbf{C}}, \mathbf{A}_{\rho,1}(\zeta)) < 2^{q-1} \,\widetilde{\beta}^{1-q} \delta, \nonumber \\
	\label{keyest eqn18} Q(T, \widetilde{\mathbf{C}}, \mathbf{A}_{\rho,1}(\zeta)) 
		\leq 2^{q-1} \,\widetilde{\beta}^{1-q} \,E(T,\mathbf{C}, \mathbf{A}_{\rho,1}(\zeta)) , 
\end{gather}
and either $\widetilde{p} = 1$ or $\widetilde{p} > 1$ and 
\begin{equation*} 
	Q(T, \widetilde{\mathbf{C}}, \mathbf{A}_{\rho,1}(\zeta)) \leq \widetilde{\beta} 
		\inf_{\mathbf{C}' \in \bigcup_{p'=1}^{\widetilde{p}-1} \mathcal{C}_{q,p'}} Q(T, \mathbf{C}', \mathbf{A}_{\rho,1}(\zeta)) . 
\end{equation*}
Provided $\delta$ and $\widetilde{\beta}$ are sufficiently small, we can argue as in the previous paragraph with $\widetilde{\mathbf{C}}$ in place of $\mathbf{C}$ and using Lemma~\ref{graphrep annuli lemma} in place of Lemma~\ref{weak graphrep lemma} to show that 
\begin{align}\label{keyest eqn19} 
	&\int_{\mathbf{A}_{\rho,1/2}(\zeta)} |(x,0)^{\perp}|^2 \,(\psi^2(R) + (\psi'(R))^2) \,\phi_{(\rho,\zeta)}(X) \,d\|T\|(X) 
	\\&- 2 \int_{\mathbf{A}_{\rho,1/2}(\zeta)} \frac{r^2}{R} \,\psi(R) \,\psi'(R) \,\phi_{(\rho,\zeta)}(X) \,d\|T\|(X) \nonumber
	\\&+ 2 \int_{\mathbf{A}_{\rho,1/2}(\zeta)} \frac{r^2}{R} \,\psi(R) \,\psi'(R) \,\phi_{(\rho,\zeta)}(X) \,d\|\widetilde{\mathbf{C}}\|(X) \nonumber 
	\\ \leq\,& C \int_{\mathbf{A}_{\rho,1}(\zeta)} \op{dist}^2(X,\op{spt} \widetilde{\mathbf{C}}) \,d\|T\|(X) \nonumber
\end{align}
for some constant $C = C(n,m,q,\gamma) \in (0,\infty)$.  Thus bounding the right-hand side of \eqref{keyest eqn19} using \eqref{keyest eqn18} gives us \eqref{keyest eqn9} (with $C_0 = C(n,m,q,\gamma) \,\widetilde{\beta}^{1-q}$).  Now fix $\delta$ and $\widetilde{\beta}$ small enough that \eqref{keyest eqn9} holds true. 

Observe that if $(x,y) \in \mathbf{B}_{(1+\gamma)/2}(0) \cap U$ and $(\xi,\zeta) \in \mathbf{B}_{(3+\gamma)/4}(0) \cap \partial U$ such that $|(x,y) - (\xi,\zeta)| \leq |x|/2$, then $|x| \leq 2 |\xi|$ and thus $(x,y) \in B_{2|\xi|}(0,\zeta)$.  By \eqref{keyest eqn6}, if $(x,y) \in \mathbf{B}_{(1+\gamma)/2}(0) \setminus U$, then there exists $t \geq 1$ such that $(tx,y) \in \mathbf{B}_{(3+\gamma)/4}(0) \cap \partial U$ and $(x,y) \in B_{2t|x|}(0,y)$.  Hence 
\begin{equation*}
	\{ (x,y) \in \mathbf{B}_{(1+\gamma)/2}(0) : \op{dist}((x,y), \mathbf{B}_{(3+\gamma)/4}(0) \cap \partial U) \leq |x|/2 \} 
	\subset \bigcup_{(\xi,\zeta) \in \mathbf{B}_{(3+\gamma)/4}(0) \cap \partial U} \mathbf{B}_{2|\xi|}(0,\zeta) .
\end{equation*}
By applying the Vitali covering lemma, there exists a countable collection $\mathcal{J}$ of $(\xi,\zeta) \in \mathbf{B}_{(3+\gamma)/4}(0) \cap \partial U$ such that $\{\mathbf{B}_{2|\xi|}(0,\zeta) : (\xi,\zeta) \in \mathcal{J}\}$ a collection of pairwise disjoint balls and 
\begin{equation*}
	\{ (x,y) \in \mathbf{B}_{(1+\gamma)/2}(0) : \op{dist}((x,y), \mathbf{B}_{(3+\gamma)/4}(0) \cap \partial U) \leq |x|/2 \} 
	\subset \bigcup_{(\xi,\zeta) \in \mathcal{J}} \mathbf{B}_{10|\xi|}(0,\zeta) .
\end{equation*}

Take any $(\xi,\zeta) \in \mathbf{B}_{(3+\gamma)/4}(0) \cap \partial U$.  We claim that 
\begin{equation}\label{keyest eqn20}
	\|T\|(\mathbf{A}_{|\xi|,1}(\zeta)) < (q+1/2) \,\mathcal{L}^n(A_{|\xi|,1}(0)) . 
\end{equation}
Note that by \eqref{keyest eqn6}, $|\xi| < (1-\gamma)/80$ and thus $\mathbf{A}_{|\xi|,2}(\zeta) \subset \mathbf{B}_{2|\xi|}(0,\zeta) \subset \mathbf{B}_1(0)$.  Since $(\xi,\zeta)$ does not satisfy \eqref{keyest eqn5} with $\rho = |\xi|$, it follows that   
\begin{equation}\label{keyest eqn21} 
	E(T, \mathbf{C}, \mathbf{B}_{2|\xi|}(0,\zeta)) \geq \delta . 
\end{equation}
To see \eqref{keyest eqn20}, since $(\xi,\zeta) \in \mathbf{B}_{(3+\gamma)/4}(0) \cap \partial U$ there exists $\rho' > 0$ and $\zeta' \in \mathbb{R}^{n-2}$ such that $(\rho')^2 + |\xi'|^2 < (3+\gamma)^2/16$, $\op{dist}((\xi,\zeta), \mathbf{A}_{\rho',1/20}(\zeta')) < \tfrac{1-\gamma}{40} \,|\xi|$, and \eqref{keyest eqn5} holds true with $(\rho',\zeta')$ in place of $(\rho,\zeta)$.  Hence $|(\rho',\zeta') - (|\xi|,\zeta)| <  \tfrac{1-\gamma}{160} \,\rho' +  \tfrac{1-\gamma}{40} \,|\xi|$ and consequently $\tfrac{1}{2} \,\rho' < \tfrac{159+\gamma}{164-4\gamma} \,\rho' \leq |\xi| \leq \tfrac{161-\gamma}{156+4\gamma} \,\rho' < \tfrac{5}{4} \,\rho'$.  Thus 
\begin{equation*} 
	B^{n-1}_{(1-\gamma)\,|\xi|/8}(|\xi|,\zeta) 
	\subset B^{n-1}_{(1-\gamma)\rho'/160 + 6 (1-\gamma)\,|\xi|/40}(\rho',\zeta') 
	\subset B^{n-1}_{31(1-\gamma)\rho'/160}(\rho',\zeta') ,
\end{equation*}
or equivalently $\mathbf{A}_{|\xi|,1}(\zeta) \subset \mathbf{A}_{\rho',31/20}(\zeta')$.  Moreover, by the monotonicity formula and $\tfrac{1}{2} \,\rho' < |\xi| < (1-\gamma)/80$, 
\begin{align*}
	\|T\|(\mathbf{A}_{\rho',2}(\zeta')) \leq\,& \|T\|(\mathbf{B}_{2\rho'}(0,\zeta')) \leq \bigg(\frac{40\rho'}{1-\gamma}\bigg)^n \|T\|(\mathbf{B}_{(1-\gamma)/20}(0,\zeta')) 
		\\ \leq\,& \bigg(\frac{40\rho'}{1-\gamma}\bigg)^n \|T\|(\mathbf{B}_1(0)) \leq (q+1/2) \,\omega_n \bigg(\frac{40\rho'}{1-\gamma}\bigg)^n . 
\end{align*}
Now suppose by way of contradiction that for $k = 1,2,3,\ldots$ we had $\delta_k \rightarrow 0^+$, $\rho_k, \rho'_k > 0$, $\zeta_k,\zeta'_k \in \mathbb{R}^{n-2}$, an $n$-dimensional locally area minimizing rectifiable current $T_k$ of $\mathbf{A}_{\rho'_k,2}(\zeta'_k)$, and $\mathbf{C}_k \in \mathcal{C}_{q,p}$ such that $\mathbf{A}_{\rho_k,1}(\zeta_k) \subset \mathbf{A}_{\rho'_k,31/20}(0,\zeta'_k)$ and 
\begin{gather}
	\label{keyest eqn22} (\partial T_k) \llcorner \mathbf{A}_{\rho'_k,2}(\zeta'_k) = 0 , \\ 
	\label{keyest eqn23} \|T_k\|(\mathbf{A}_{\rho'_k,2}(\zeta'_k)) \leq (q+1/2) \,\omega_n \bigg(\frac{40\rho'_k}{1-\gamma}\bigg)^n , \\ 
	\label{keyest eqn24} \|T_k\|(\mathbf{A}_{\rho'_k,1}(\zeta'_k)) \leq (q+1/2) \,\mathcal{L}^n(A_{\rho'_k,1}(0)) , \\ 
	\label{keyest eqn25} E(T_k, \mathbf{C}_k, \mathbf{A}_{\rho'_k,2}(0,\zeta'_k)) < \delta_k , \\
	\label{keyest eqn26} \|T_k\|(\mathbf{A}_{\rho_k,1}(\zeta_k)) \geq (q+1/2) \,\mathcal{L}^n(A_{\rho_k,1}(0)) . 
\end{gather}
By rescaling, we may take $\zeta'_k = 0$ and $\rho'_k = 1$.  By \eqref{keyest eqn22}, \eqref{keyest eqn23}, the Federer-Fleming compactness theorem, and~\cite[Theorem~34.5]{SimonGMT}, after passing to a subsequence there exists an $n$-dimensional locally area minimizing rectifiable current $T_{\infty}$ of $\mathbf{A}_{1,2}(0)$ such that $T_k \rightarrow T_{\infty}$ weakly in $\mathbf{A}_{1,2}(0)$ and $\|T_k\| \rightarrow \|T_{\infty}\|$ in the sense of Radon measures of $\mathbf{A}_{1,2}(0)$.  Arguing as we did to prove \eqref{graphrep eqn4}, it follows from \eqref{keyest eqn25} that $\op{spt} T_{\infty}$ is a union of $n$-dimensional planes intersecting along $\{0\} \times \mathbb{R}^{n-2}$.  By \eqref{keyest eqn24}, $\|T_{\infty}\|(\mathbf{A}_{1,1}(0)) \leq (q+1/2) \,\mathcal{L}^n(A_{1,1}(0))$ so that $T_{\infty}$ is a sum of integer multiplicity planes to total multiplicity $\leq q$, contradicting \eqref{keyest eqn26}.  Therefore, \eqref{keyest eqn20} must hold true. 

Recall that by \eqref{keyest eqn6}, $|\xi| < (1-\gamma)/80$.  Hence by the monotonicity formula, 
\begin{align}\label{keyest eqn27}
	\|T\|(\mathbf{B}_{10|\xi|}(0,\zeta)) 
	&\leq \bigg(\frac{80 |\xi|}{1-\gamma}\bigg)^n \|T\|(\mathbf{B}_{(1-\gamma)/8}(0,\zeta)) 
	\\&\leq \bigg(\frac{80 |\xi|}{1-\gamma}\bigg)^n \|T\|(\mathbf{B}_1(0)) 
	\leq (q+1/2) \,\omega_n \bigg(\frac{80 |\xi|}{1-\gamma}\bigg)^n . \nonumber 
\end{align}
By \eqref{keyest eqn21} and \eqref{keyest eqn27}, for each $(\xi,\zeta) \in \mathcal{J}$ 
\begin{align}\label{keyest eqn28}
	&\int_{\mathbf{B}_{10|\xi|}(0,\zeta)} |(x,0)^{\perp}|^2 \,(\psi^2(R) + (\psi'(R))^2) \,(1-\chi) \,d\|T\|(X) 
	\\&- 2 \int_{\mathbf{B}_{10|\xi|}(0,\zeta)} \frac{r^2}{R} \,\psi(R) \,\psi'(R) \,(1-\chi) \,d\|T\|(X) \nonumber 
	\\ \leq\,& C |\xi|^{n+2} 
	\leq \frac{C}{\delta} \int_{\mathbf{B}_{2|\xi|}(0,\zeta)} \op{dist}^2(X,\op{spt} \mathbf{C}) \,d\|T\|(X), \nonumber
\end{align}
where $C = C(n,m,q,\gamma) \in (0,\infty)$ are constants.  Summing \eqref{keyest eqn28} over $(\xi,\zeta) \in \mathcal{J}$ gives us 
\begin{align}\label{keyest eqn29}
	&\int_{\mathbf{B}_{(1+\gamma)/2}(0)} |(x,0)^{\perp}|^2 \,(\psi^2(R) + (\psi'(R))^2) \,(1-\chi) \,d\|T\|(X) 
	\\&- 2 \int_{\mathbf{B}_{(1+\gamma)/2}(0)} \frac{r^2}{R} \,\psi(R) \,\psi'(R) \,(1-\chi) \,d\|T\|(X) \nonumber 
	\\ \leq\,& C \int_{\mathbf{B}_1(0)} \op{dist}^2(X,\op{spt} \mathbf{C}) \,d\|T\|(X) \nonumber
\end{align}
for some constant $C = C(n,m,q,\gamma) \in (0,\infty)$.  Adding \eqref{keyest eqn10} and \eqref{keyest eqn29} gives us the desired upper bound on the right-hand side of \eqref{keyest eqn4}.
\end{proof}

\begin{corollary} \label{keyest cor1}
For each $\gamma \in (0,1)$ and $\sigma \in (0,1)$ there exists $\varepsilon_0 = \varepsilon_0(n,m,q,\gamma,\sigma) \in (0,1)$ and $\beta_0 = \beta_0(n,m,q,\gamma,\sigma) \in (0,1)$ if $\mathbf{C}$ and $T$ satisfy Hypothesis~$(\star)$ and Hypothesis~$(\star\star)$, then 
\begin{equation} \label{keyest cor1 concl}
	\int_{{\mathbf B}_{\gamma}(0)} \frac{\op{dist}^2(X, \op{spt} \mathbf{C})}{R^{n+2-\sigma}} \,d\|T\|(X) 
		\leq C \int_{{\mathbf B}_1(0)} \op{dist}^2(X, \op{spt} \mathbf{C}) \,d\|T\|(X) , 
\end{equation}
where $R = R(X) = |X|$ and $C = C(n,m,q,\gamma,\sigma) \in (0,\infty)$ is a constant. 
\end{corollary}

\begin{proof} 
This follows from Theorem~\ref{keyest thm}(a) exactly as in the proof of Lemma~3.4 of~\cite{Sim93}. 
\end{proof}

\subsection{A priori estimates for area-minimizers close to a plane} \label{sec:apriori close2plane}

In this section we will focus on the case where $\mathbf{C}$ and $T$ are close to an $n$-dimensional plane $P_0$.  

\noindent 
{\bf Preliminary remarks.}  After an orthogonal change of coordinates, assume that $P_0 = \{0\} \times \mathbb{R}^n$.  (This a slight change from Section~\ref{height-bound}, in which we took $P_0 = \mathbb{R}^n \times \{0\}$.)  We identify $P_0 \cong \mathbb{R}^n$ and $P_0^{\perp} \cong \mathbb{R}^m$.  Set $\mathbf{P}_0 = q \llbracket P_0 \rrbracket$ where $P_0$ is oriented by $\vec P_0 = e_{m+1} \wedge e_{m+2} \wedge \cdots \wedge e_{m+n}$, where $e_1,e_2,\ldots,e_{n+m}$ is the standard basis for $\mathbb{R}^{n+m}$.  We will represent points $X \in \mathbb{R}^{n+m}$ as $X = (z,x,y)$ where $x \in \mathbb{R}^2$, $y \in \mathbb{R}^{n-2}$, and $z \in \mathbb{R}^m$.  For each $x_0 \in \mathbb{R}^2$, $y_0 \in \mathbb{R}^{n-2}$, and $\rho > 0$ we let 
\begin{align*}
	B_{\rho}(x_0,y_0) &= \{ (x,y) \in \mathbb{R}^2 \times \mathbb{R}^{n-2} : |x-x_0|^2 + |y-y_0|^2 < \rho^2 \}, \\
	\mathbf{C}_{\rho}(x_0,y_0) &= \mathbb{R}^m \times B_{\rho}(x_0,y_0). 
\end{align*}

We shall make the following hypothesis for appropriate choices of constants $\eta_0 \in (0,1)$ and $M \in [1,\infty)$, to be chosen ultimately depending only on $n$, $m$ and $q$: 

\noindent\textbf{Hypothesis~$(\dagger)$.}  $2 \leq p \leq q$ are integers, $\mathbf{C} \in \mathcal{C}_{q,p}$, and $T$ is an $n$-dimensional locally area-minimizing 
rectifiable current in $\mathbf{B}_1(0)$ with $\partial \, T \llcorner {\mathbf B}_{1}(0) = 0$ such that 
\begin{gather}
	\label{main hyp eqn4} E(T, \mathbf{P}_0, \mathbf{B}_1(0)) < \eta_0, \\
	\label{main hyp eqn5} E(T, \mathbf{P}_0, \mathbf{B}_1(0)) \leq M \inf_{\mathbf{P} \in \mathcal{C}_{q,1}} E(T, \mathbf{P}, \mathbf{B}_1(0)) . 
\end{gather}

Suppose that $\mathbf{C}$ and $T$ satisfy Hypothesis~$(\star)$, Hypothesis~$(\star\star)$ (of Section~\ref{notation-and-graphical}) and Hypothesis~$(\dagger)$.  By \eqref{main hyp eqn4} and Lemma~\ref{one plane lemma}, 
\begin{equation} \label{close2plane eqn1}
	\sup_{X \in \op{spt} T \cap \mathbf{B}_{7/8}(0)} \op{dist}(X,P_0) \leq C E(T, P_0, \mathbf{B}_1(0)) 
\end{equation}
for some constant $C = C(n,m) \in (0,\infty)$.  By the constancy theorem, $(\pi_{P_0\#} (T \llcorner \mathbf{B}_{7/8}(0))) \llcorner \mathbf{B}_{3/4}(0)$ is an integer multiple of $\llbracket P_0 \rrbracket \llcorner \mathbf{B}_{3/4}(0)$.  We claim that in fact $(\pi_{P_0\#} (T \llcorner \mathbf{B}_{7/8}(0))) \llcorner \mathbf{B}_{3/4}(0) = \pm q \llbracket P_0 \rrbracket \llcorner \mathbf{B}_{3/4}(0)$.  Thus up to reversing the orientation of $T$, we may assume that 
\begin{equation} \label{close2plane convention1}
	(\pi_{P_0\#} (T \llcorner \mathbf{B}_{7/8}(0))) \llcorner \mathbf{B}_{3/4}(0) = q \llbracket P_0 \rrbracket \llcorner \mathbf{B}_{3/4}(0) . 
\end{equation}
To see this, for $k = 1,2,3,\ldots$ let $\varepsilon_k \rightarrow 0^+$, $\eta_k \rightarrow 0^+$, and $T_k$ be an $n$-dimensional locally area minimizing rectifiable current of $\mathbf{B}_1(0)$ satisfying \eqref{main hyp eqn1} and \eqref{main hyp eqn4} with $\varepsilon_k, \eta_k, T_k$ in place of $\varepsilon, \eta, T$.  By \eqref{main hyp eqn1} and the Federer-Fleming compactness theorem, after passing to a subsequence there exists an $n$-dimensional locally area minimizing rectifiable current $T_{\infty}$ such that $T_k \rightarrow T_{\infty}$ weakly in $\mathbf{B}_1(0)$.  By \eqref{main hyp eqn1}, \eqref{close2plane eqn1}, and the constancy theorem, $T_{\infty} \llcorner \mathbf{B}_{3/4}(0) = \pm q \llbracket P_0 \rrbracket \llcorner \mathbf{B}_{3/4}(0)$.  Up to reversing the orientation of each $T_k$, we may assume that $T_{\infty} \llcorner \mathbf{B}_{3/4}(0) = q \llbracket P_0 \rrbracket \llcorner \mathbf{B}_{3/4}(0)$.  By again applying the constancy theorem, $(\pi_{P_0\#} (T_k \llcorner \mathbf{B}_{7/8}(0))) \llcorner \mathbf{B}_{3/4}(0)$ is an integer multiple of $\llbracket P_0 \rrbracket \llcorner \mathbf{B}_{3/4}(0)$.  By continuity of push-forwards under weak limits, 
\begin{align*}
	\lim_{k \rightarrow \infty} (\pi_{P_0\#} (T_k \llcorner \mathbf{B}_{7/8}(0))) \llcorner \mathbf{B}_{3/4}(0)
	=\,& \lim_{k \rightarrow \infty} \pi_{P_0\#}(T_k \llcorner \mathbf{B}_{7/8}(0) \cap \mathbf{C}_{3/4}(0)) 
	\\=\,& \pi_{P_0\#}(T_{\infty} \llcorner \mathbf{C}_{3/4}(0)) 
	= q \llbracket P_0 \rrbracket \llcorner \mathbf{B}_{3/4}(0) 
\end{align*}
where the limits are taken in the weak topology in $\mathbf{C}_{3/4}(0)$.  Therefore, $(\pi_{P_0\#} (T_k \llcorner \mathbf{B}_{7/8}(0))) \llcorner \mathbf{B}_{3/4}(0) = q \llbracket P_0 \rrbracket \llcorner \mathbf{B}_{3/4}(0)$ for infinitely many $k$, which in view of the arbitrary choice of sequence $(T_k)$ proves \eqref{close2plane convention1}.

By \eqref{close2plane convention1}, for each $X \in P_0 \cap \mathbf{B}_{1/2}(0)$ there exists $Y \in \op{spt} T \cap \mathbf{B}_{7/8}(0)$ such that $\pi_{P_0}(Y) = X$ and thus by \eqref{close2plane eqn1}
\begin{equation*}
	\op{dist}(X, \op{spt} T) \leq |X - Y| = \op{dist}(Y, P_0) \leq C E(T, P_0, \mathbf{B}_1(0)) .
\end{equation*}
Hence 
\begin{align} \label{close2plane eqn2}
	Q(T,\mathbf{P}_0,\mathbf{B}_1)^2 &= E(T,\mathbf{P}_0,\mathbf{B}_1)^2 
			+ q \int_{P_0 \cap \mathbf{B}_{1/2}(0) \cap \{|x| > 1/16\}} \op{dist}^2(X, \op{spt} T) \,d\mathcal{H}^n(X)
		\\&\leq E(T,\mathbf{P}_0,\mathbf{B}_1)^2 + 2q \omega_n \sup_{X \in P_0 \cap \mathbf{B}_{1/2}(0)} \op{dist}^2(X, \op{spt} T) 
		< C E(T,\mathbf{P}_0,\mathbf{B}_1)^2 \nonumber 
\end{align}
for some constant $C = C(n,m,q) \in (0,\infty)$.  

By the triangle inequality, $\|T\|(\mathbf{B}_1(0)) \leq (q+1/2)\,\omega_n$, \eqref{main hyp eqn4}, \eqref{close2plane eqn1}, \eqref{main hyp eqn3}, and \eqref{close2plane eqn2}, 
\begin{align*}
	\int_{\mathbf{B}_{1/2}(0)} \op{dist}^2(X, P_0) \,d\|\mathbf{C}\|(X) 
	\leq\,& 2 \int_{\mathbf{B}_{1/2}(0)} \op{dist}^2(X, \op{spt} T) \,d\|\mathbf{C}\|(X) 
		\\&+ 2 q \omega_n \sup_{\op{spt} T \cap \mathbf{B}_{3/4}(0)} \op{dist}^2(X, P_0) 
	< C E(T,\mathbf{P}_0,\mathbf{B}_1)^2
\end{align*}
for some constant $C = C(n,m,q) \in (0,\infty)$.  Since $P_0$ is a plane and $\mathbf{C} \in \mathcal{C}_{q,p}$, 
\begin{equation} \label{close2plane eqn3}
	\op{dist}_{\mathcal{H}}(\op{spt} \mathbf{C} \cap \mathbf{B}_1(0), P_0 \cap \mathbf{B}_1(0)) < C E(T,\mathbf{P}_0,\mathbf{B}_1)
\end{equation}
for some constant $C = C(n,m,q) \in (0,\infty)$.  Hence letting $\mathbf{C} = \sum_{i=1}^p q_i \llbracket P_i \rrbracket$ as in \eqref{cone sums of planes form}, 
for each $i \in \{1,2,\ldots,p\}$ there is an $m \times 2$ matrix $A_i$ such that 
\begin{gather}
	\label{close2plane eqn4} P_i = \{ (z,x,y) \in \mathbb{R}^{n+m} : z = A_i x \} , \\
	\label{close2plane eqn5} \|A_i\| \leq C E(T,\mathbf{P}_0,\mathbf{B}_1) , 
\end{gather}
where $\|\,\cdot\,\|$ denotes the Frobenius norm and $C = C(n,m,q) \in (0,\infty)$ is a constant.  Up to reversing the orientation of $P_i$, we may assume that 
\begin{equation} \label{close2plane convention2}
	|\vec P_i - \vec P_0| < 1/2 
\end{equation}
so that $P_i$ is equipped with the induced orientation on the graph of $(x,y) \mapsto A_i x$ over the oriented plane $P_0$.  By Theorem~\ref{graphrep thm}(a), 
\begin{gather}
	\label{close2plane eqn6} \min_{1 \leq i < j \leq p} \,\inf_{x \in \mathbb{S}^1} |A_i x - A_j x| 
		\geq c \inf_{\mathbf{C}' \in \bigcup_{p'=1}^{p-1} \mathcal{C}_{q,p}} Q(T, \mathbf{C}', \mathbf{A}_{1,1}(0)) , \\
	\label{close2plane eqn7} \|A_i - A_j\| \leq C \inf_{x \in \mathbb{S}^1} |A_i x - A_j x| \text{ for all $1 \leq i < j \leq p$,} 
\end{gather}
where $c = c(n,m,q) > 0$ and $C = C(n,m,q) \in (0,\infty)$. 

Let $P$ be an $n$-dimensional oriented plane in $\mathbb{R}^{n+m}$ and let $\mathbf{P} = q \llbracket P \rrbracket$.  By the triangle inequality, $\|T\|(\mathbf{B}_1(0)) \leq (q+1/2) \,\omega_n$, \eqref{main hyp eqn5}, \eqref{main hyp eqn3}, and \eqref{close2plane eqn2}
\begin{align}\label{close2plane eqn8}
	&\frac{1}{M^2} \int_{\mathbf{B}_1(0)} \op{dist}^2(X,P_0) \,d\|T\|(X) \leq \int_{\mathbf{B}_1(0)} \op{dist}^2(X,P) \,d\|T\|(X) 
	\\ \leq\,& 2 \int_{\mathbf{B}_1(0)} \op{dist}^2(X,\op{spt}\mathbf{C}) \,d\|T\|(X) 
		+ 2 \,(q+1/2) \,\omega_n \sup_{X \in \op{spt}\mathbf{C} \cap \mathbf{B}_1(0)} \op{dist}^2(X,P) \nonumber 
	\\ \leq\,& 2C \beta_0^2 \int_{\mathbf{B}_1(0)} \op{dist}^2(X,P_0) \,d\|T\|(X) 
		+ 2 \,(q+1/2) \,\omega_n \sup_{X \in \op{spt}\mathbf{C} \cap \mathbf{B}_1(0)} \op{dist}^2(X,P) \nonumber
\end{align}
for some constant $C = C(n,m,q) \in (0,\infty)$.  Choosing $\beta_0 = \beta_0(n,m,q,M)$ small enough that $C M^2 \beta_0^2 < 1/4$ (where $C$ is as in \eqref{close2plane eqn8}), 
\begin{equation}\label{close2plane eqn9}
	E(T, \mathbf{P}_0, \mathbf{B}_1(0)) \leq C M \op{dist}_{\mathcal{H}}(\op{spt}\mathbf{C} \cap \mathbf{B}_1(0), P \cap \mathbf{B}_1(0))
\end{equation}
for some constant $C = C(n,m,q) \in (0,\infty)$.  In particular, if $P = \{ (z,x,y) \in \mathbb{R}^{n+m} : z = Ax \}$ for some $m \times 2$ matrix $A$, then 
\begin{equation}\label{close2plane eqn10}
	E(T, \mathbf{P}_0, \mathbf{B}_1(0)) \leq C M \max_{1 \leq i \leq p} \|A_i - A\| . 
\end{equation}
Setting $A = 0$ in \eqref{close2plane eqn10} gives us 
\begin{equation}\label{close2plane eqn11} 
	\max_{1 \leq i \leq p} \|A_i\| \geq \frac{c}{M} \,E(T, \mathbf{P}_0, \mathbf{B}_1(0)) 
\end{equation}
for some constant $c = c(n,m,q) > 0$.  Setting $A = A_j$ for $j \in \{1,2,\ldots,p\}$ in \eqref{close2plane eqn10} gives us 
\begin{equation}\label{close2plane eqn12} 
	\max_{1 \leq i < j \leq p} \|A_i - A_j\| \geq \frac{c}{M} \,E(T, \mathbf{P}_0, \mathbf{B}_1(0))
\end{equation}
for some constant $c = c(n,m,q) > 0$.  By \eqref{close2plane eqn7} and \eqref{close2plane eqn12}, 
\begin{equation}\label{close2plane eqn13} 
	\max_{1 \leq i < j \leq p} \inf_{x \in \mathbb{S}^1} |A_i x - A_j x| \geq \frac{c}{M} \,E(T, \mathbf{P}_0, \mathbf{B}_1(0))
\end{equation}
for some constant $c = c(n,m,q) > 0$. 

In blow-up arguments, we typically let $\mathbf{C}_k = \sum_{i=1}^p q_i \llbracket P^{(k)}_i \rrbracket$ where $P^{(k)}_i = \{ (z,x,y) \in \mathbb{R}^{n+m} : z = A^{(k)}_i x \}$ for some $m \times 2$ matrix $A^{(k)}_i$ with $\|A^{(k)}_i\| \leq C(n,m,q) \,\widehat{E}_k$, where $\widehat{E}_k = E(T_k,\mathbf{P}_0,\mathbf{B}_1(0))$ and $C = C(n,m,q) \in (0,\infty)$ is a constant.  After passing to a subsequence, there is an $m \times 2$ matrix $\Lambda_i$ such that $A^{(k)}_i/\widehat{E}_k \rightarrow \Lambda_i$ for each $i \in \{1,2,\ldots,p\}$.  By dividing \eqref{close2plane eqn5}, \eqref{close2plane eqn11}, \eqref{close2plane eqn12}, and \eqref{close2plane eqn13}, all taken with $T = T_k$ and $A_i = A^{(k)}_i,$ by $\widehat{E}_k$ and letting $k \rightarrow \infty$, 
\begin{gather}
	\label{close2plane eqn14} \|\Lambda_i\| \leq C, \quad \max_{1 \leq i \leq p} \|\Lambda_i\| \geq \frac{c}{M} , \quad 
		\max_{1 \leq i < j \leq p} \|\Lambda_i - \Lambda_j\| \geq \frac{c}{M} , \\
	\label{close2plane eqn15} \max_{1 \leq i < j \leq p} \inf_{x \in \mathbb{S}^1} |\Lambda_i x - \Lambda_j x| \geq \frac{c}{M} 
\end{gather}
for some constants $c = c(n,m,q) > 0$ and $C = C(n,m,q) \in (0,\infty)$.  We will frequently use the fact that by \eqref{close2plane eqn15} if $x \in \mathbb{R}^2$ such that $\Lambda_i x = \Lambda_j x$ for all $i,j \in \{1,2,\ldots,p\}$ then $x = 0$.

When $T$ is close to the plane $P_0$ as in \eqref{main hyp eqn4}, we can restate Theorem~\ref{graphrep thm} as follows:

\begin{theorem} \label{graphrep close2plane thm}
Given integers $2 \leq p \leq q$ and $0 < \tau < \gamma < 1$ there exists $\varepsilon_0, \beta_0, \eta_0 \in (0,1)$ depending only on $n,m,q,p,\gamma,\tau$ such that if $\mathbf{C} \in \mathcal{C}_{q,p}$ and $T$ satisfy Hypothesis~$(\star)$, Hypothesis~$(\star\star)$, \eqref{main hyp eqn4}, \eqref{close2plane convention1}, and \eqref{close2plane convention2}, then: 
\begin{enumerate} \setlength{\itemsep}{5pt}
	\item[(a)]  \eqref{graphrep concl a1} and \eqref{graphrep concl a2} hold true; 
	
	\item[(b)]  up to changing the values of the multiplicities $q_i$ of the planes of $\mathbf{C}$, there exists $n$-dimensional locally area minimizing rectifiable currents $T_i$ of $\mathbf{C}_{(3+\gamma)/4}(0) \cap \{r > \tau/4\}$ such that 
	\begin{gather}\label{graphrep close2plane concl b} 
		T \llcorner \mathbf{B}_{(7+\gamma)/8}(0) \cap \mathbf{C}_{(3+\gamma)/4}(0) \cap \{r > \tau/4\} = \sum_{i=1}^p T_i , \\
		(\partial T_i) \llcorner \mathbf{C}_{(3+\gamma)/4}(0) \cap \{r > \tau/4\} = 0, \nonumber \\ 
		(\pi_{P_0 \#} T_i) \llcorner \mathbf{C}_{(1+\gamma)/2}(0) \cap \{r > \tau/2\} 
			= q_i \llbracket P_0 \rrbracket \llcorner \mathbf{C}_{(1+\gamma)/2}(0) \cap \{r > \tau/2\} , \nonumber \\ 
		\sup_{X \in \op{spt} T_i \cap \{r > \sigma\}} \op{dist}(X, P_i) \leq C_{\sigma} E \text{ for all } \sigma \in [\tau/2,1/2] , \nonumber 
	\end{gather}
	where $r(X) = \op{dist}(X, \{0\} \times \mathbb{R}^{n-2})$, $E = E(T,\mathbf{C},\mathbf{B}_1(0))$, and $C_{\sigma} = C_{\sigma}(n,m,q,\gamma,\sigma) \in (0,\infty)$ are constants; 
	
	\item[(c)] for each $i \in \{1,2,\ldots,p\}$ there exists Lipschitz $q_i$-valued functions $u_i : B_{\gamma}(0) \cap \{r > \tau\} \rightarrow \mathcal{A}_{q_i}(\mathbb{R}^m)$ and a closed set $K \subseteq B_{\gamma}(0) \cap \{r > \tau\}$ such that 
	\begin{gather}\label{graphrep close2plane concl c} 
		T_i \llcorner (\mathbb{R}^m \times K) = \op{graph} (A_i x + u_i) \llcorner (\mathbb{R}^m \times K) , \\
		\mathcal{H}^n(B_{\gamma}(0) \cap \{r > \sigma\} \setminus K)
			+ \|T_i\|(\mathbb{R}^m \times (B_{\gamma}(0) \cap \{r > \sigma\} \setminus K)) \leq C_{\sigma} E^{2+\alpha} , \nonumber \\
		\sup_{B_{\gamma}(0) \cap \{r > \sigma\}} |u_i| \leq C_{\sigma} E , \quad 
			\sup_{B_{\gamma}(0) \cap \{r > \sigma\}} |\nabla u_i| \leq C_{\sigma} E^{\alpha} \nonumber
	\end{gather}
	for all $\sigma \in [\tau,1/2]$, where again $E = E(T,\mathbf{C},\mathbf{B}_1(0))$ and $\alpha = \alpha(n,m,q) \in (0,1)$, $C_{\sigma} = C_{\sigma}(n,m,q,\gamma,\sigma) \in (0,\infty)$ are constants.
\end{enumerate}
\end{theorem}

\begin{proof}
Conclusion (a) is the same as Theorem~\ref{graphrep thm}(a) with $(1+\gamma)/2$ and $\tau/2$ in place of $\gamma$ and $\tau$.  Conclusion (b) follows by using \eqref{main hyp eqn4}, \eqref{close2plane eqn1}, \eqref{close2plane eqn3}, and Theorem~\ref{graphrep thm}(b) with $(1+\gamma)/2$ and $\tau/2$ in place of $\gamma$ and $\tau$.  Conclusion (c) follows from Almgren's Strong Lipschitz Approximation Theorem (Theorem~\ref{lip approx thm}) and \eqref{close2plane eqn4}.
\end{proof}

We also have the following consequence of Theorem~\ref{keyest thm} and Theorem~\ref{graphrep close2plane thm}. 

\begin{corollary} \label{keyest cor2}
For each $0 < \tau < \gamma < 1$  there exists $\varepsilon_0, \beta_0, \eta_0 \in (0,1)$ depending only on $n,m,q,\gamma,\tau$ such that if $\mathbf{C}$ and $T$ satisfy Hypothesis~$(\star)$, Hypothesis~$(\star\star)$, \eqref{main hyp eqn4}, \eqref{close2plane convention1}, and \eqref{close2plane convention2}, then: 
\begin{equation} \label{keyest cor2 concl}
	\sum_{i=1}^p \int_{B_{\gamma}(0) \cap \{r > \tau\}} R^{2-n} \left| \frac{\partial (u_i/R)}{\partial R} \right|^2 
		\leq C \int_{\mathbf{B}_1(0)} \op{dist}^2(X, \op{spt} \mathbf{C}) \,d\|T\|(X) , 
\end{equation}
where $R(x,y) = |(x,y)|$, $u_i : B_{\gamma}(0) \cap \{r > \tau\} \rightarrow \mathcal{A}_{q_i}(\mathbb{R}^m)$ are as in Theorem~\ref{graphrep close2plane thm}(c), and $C = C(n,m,q,\gamma) \in (0,\infty)$ is a constant. 
\end{corollary}

\begin{proof}
Let $i \in \{1,2,\ldots,p\}$.  Throughout this proof, let $(x',y)$ denote points in $\mathbb{R}^n$, where $x' \in \mathbb{R}^2$ and $y \in \mathbb{R}^{n-2}$.  For each $(x',y) \in B_{\gamma}(0) \cap \{r > \tau\}$ let $u_i(x',y) = \sum_{j=1}^{q_i} \llbracket u_{i,j}(x',y) \rrbracket$ where $u_{i,j}(x',y) \in \mathbb{R}^m$.  By Rademacher's Theorem~\cite[Theorem~1.13]{DeLSpaDirMin}, $u_i$ is differentiable at $\mathcal{H}^n$-a.e.~$(x',y) \in K$ in the sense of~\cite[Definition~1.9]{DeLSpaDirMin}.  Let $u_i$ be differentiable at $(x',y) \in K$ and let $x = (A_i x' + u_{i,j}(x',y), x')$ for $j \in \{1,2,\ldots,q_i\}$ so that $X = (x,y) = (A_i x' + u_{i,j}(x',y), x',y)$ is a point on $\op{spt} T_i \cap (\mathbb{R}^m \times K_i)$.  Then $\tfrac{\partial}{\partial R} (u_{i,j}(x',y), x',y)$ is tangential to $\op{spt} T_i$ at $X$ and thus 
\begin{equation*}
	(x,y)^{\perp} = R^2 \Bigg(\frac{\partial}{\partial R}\Bigg(\frac{(A_i x' + u_{i,j}(x',y), x',y)}{R}\Bigg)\Bigg)^{\perp} 
		= R^2 \Bigg(\frac{\partial}{\partial R}\Bigg(\frac{(u_{i,j}(x',y), 0,0)}{R}\Bigg)\Bigg)^{\perp} , 
\end{equation*}
where $R = |(x',y)|$ and $(\,\cdot\,)^{\perp}$ denotes the orthogonal projection onto the approximate tangent plane to $T$ at $X$.  Since by Theorem~\ref{graphrep close2plane thm}(c) $\op{Lip} u_i \leq C E^{\alpha}$ is small, 
\begin{equation*}
	|(x,y)^{\perp}| \geq \frac{1}{2} \,R^2 \,\Bigg|\frac{\partial}{\partial R}\Bigg(\frac{u_{i,j}(x',y)}{R}\Bigg)\Bigg|
\end{equation*}
Now \eqref{keyest cor2 concl} follows from Theorem~\ref{keyest thm}(a) and Theorem~\ref{graphrep close2plane thm}(c). 
\end{proof}

Next we derive further a priori estimates based on~\cite[Theorem~3.1]{Sim93} and~\cite[Corollary~10.2]{Wic14}. 

\begin{theorem} \label{branchpt close2plane thm}
For all integers $2 \leq p \leq q$, $0 < \tau < \gamma < 1$, $\sigma \in (0,1)$, and $M \in [1,\infty)$ there exists $\varepsilon_0,\beta_0,\eta_0 \in (0,1)$ depending only on $n,m,q,p,\gamma,\tau,\sigma,M$ such that the following holds true.  Suppose that $\mathbf{C}$ and $T$ satisfy Hypothesis~$(\star)$, Hypothesis~$(\star\star)$, Hypothesis~$(\dagger)$, \eqref{close2plane convention1}, and \eqref{close2plane convention2}.  Let $Z = (\chi,\xi,\zeta) \in \op{spt} T \cap \mathbf{B}_{1/2}(0)$ such that $\Theta(T,Z) \geq q$.  Then 
\begin{align*} 
	&(a) \quad |\chi|^2 + \widehat{E}^2 |\xi|^2 \leq C \int_{\mathbf{B}_1(0)} \op{dist}^2(X, \op{spt} \mathbf{C}) \,d\|T\|(X) , \\
	&(b) \quad \int_{\mathbf{B}_{1/4}(Z)} \op{dist}(X-Z,\op{spt} \mathbf{C})^2 \,d\|T\|(X) 
		\leq C \int_{\mathbf{B}_1(0)} \op{dist}^2(X, \op{spt} \mathbf{C}) \,d\|T\|(X) ,
\end{align*}
where $\widehat{E} = E(T,\mathbf{P}_0,\mathbf{B}_1(0))$ and $C = C(n,m,q,p,M) \in (0,\infty)$ is a constant (independent of $\tau$).  Moreover, 
\begin{align*} 
	&(c) \quad \int_{\mathbf{B}_{\gamma}(0)} \frac{\op{dist}^2(X,\op{spt} \mathbf{C}) \,d\|T\|(X)}{|X-Z|^{n+2-\sigma}}
		\leq C \int_{\mathbf{B}_1(0)} \op{dist}^2(X, \op{spt} \mathbf{C}) \,d\|T\|(X) , \\
	&(d) \quad \int_{B_{\gamma}(0) \cap \{r > \tau\}} \frac{|u_i(x,y) - \chi + A_i \xi|^2}{|(x,y) - (\xi,\zeta)|^{n+2-\sigma}} 
		\leq C \int_{\mathbf{B}_1(0)} \op{dist}^2(X, \op{spt} \mathbf{C}) \,d\|T\|(X) \text{ for all $i$}
\end{align*}
where $A_i$ are as in \eqref{close2plane eqn4}, $u_i : B_{\gamma}(0) \cap \{r > \tau\} \rightarrow \mathcal{A}_{q_i}(\mathbb{R}^m)$ are as in Theorem~\ref{graphrep close2plane thm}(c), and $C = C(n,m,q,p,\gamma,M,\sigma) \in (0,\infty)$ is a constant (independent of $\tau$). 
\end{theorem}

The proof of Theorem~\ref{branchpt close2plane thm} is similar to that of~\cite[Corollary~10.2]{Wic14} with some minor changes.  We provide some details for completion.  We shall assume the inductive hypothesis that $p_0 \in \{2,3,\ldots,q\}$ such that either $p_0 = 2$ or $p_0 > 2$ and 
\begin{enumerate}
	\item[(H4)]  Theorem~\ref{branchpt close2plane thm} holds true for all $p \in \{2,3,\ldots,p_0-1\}$ . 
\end{enumerate}
Theorem~\ref{branchpt close2plane thm} in the case $p = p_0$ will follow from three preliminary results Lemma~\ref{branchpt close2plane lemma1}, Lemma~\ref{branchpt close2plane lemma2}, and Corollary~\ref{branchpt close2plane cor3}.

\begin{lemma} \label{branchpt close2plane lemma1}
For all integers $q$, $p$ with $2 \leq p \leq q$, any $M \in [1,\infty)$ and any $\delta > 0,$ there exist $\varepsilon,\beta,\eta \in (0,1)$ depending only on $n,m,q,M,\delta$ such that if $\mathbf{C} \in \mathcal{C}_{q,p}$ and $T$ satisfy Hypothesis~$(\star)$, Hypothesis~$(\star\star)$, and Hypothesis~$(\dagger)$ with $\varepsilon,\beta,\eta$ in place of $\varepsilon_0,\beta_0,\eta_0$ and if $Z = (\chi,\xi,\zeta) \in \op{spt} T \cap \mathbf{B}_{1/2}(0)$ is such that $\Theta(T,Z) \geq q$, then 
\begin{equation} \label{branchpt close2plane1 concl}
	|\chi|^2 + \widehat{E}^2 |\xi|^2 \leq \delta^2 \widehat{E}^2 , 
\end{equation}
where $\widehat{E} = E(T,\mathbf{P}_0,\mathbf{B}_1(0))$. 
\end{lemma}

\begin{proof}
Fix $\delta > 0$.  Suppose to the contrary that for $k = 1,2,3,\ldots$ there exists $\varepsilon_k \rightarrow 0^+$, $\beta_k \rightarrow 0^+$, $\gamma_k \rightarrow 0^+$, $\mathbf{C}_k \in \mathcal{C}_{q,p}$, and $T_k$ such that Hypothesis~$(\star)$, Hypothesis~$(\star\star)$, Hypothesis~$(\dagger)$, \eqref{close2plane convention1}, and \eqref{close2plane convention2} hold true with $\varepsilon_k, \beta_k, \eta_k, \mathbf{C}_k, T_k$ in place of $\varepsilon_0, \beta_0, \eta_0, \mathbf{C}, T$ but there is a point $Z_k = (\chi_k,\xi_k,\zeta_k) \in \op{spt} T_k \cap \mathbf{B}_{1/2}(0)$ such that $\Theta(T_k,Z_k) \geq q$
\begin{equation}\label{branchpt close2plane1 eqn1}
	|\chi_k|^2 + \widehat{E}_k^2 |\xi_k|^2 \geq \delta^2 \widehat{E}_k^2 , 
\end{equation}
where $\widehat{E}_k = E(T_k,\mathbf{P}_0,\mathbf{B}_1(0))$.  Let $\mathbf{C}_k = \sum_{i=1}^p q_i \llbracket P^{(k)}_i \rrbracket$ where $q_i$ are integers with $\sum_{i=1}^p q_i = q$ and $P^{(k)}_i = \{ (z,x,y) \in \mathbb{R}^{n+m} : z = A^{(k)}_i x \}$ for some distinct $m \times 2$ matrices $A^{(k)}_i$ with $\|A^{(k)}_i\| \leq C(n,m,q) \,\widehat{E}_k$.  After passing to a subsequence assume that $q_i$ are independent of $k$.  Arguing as in~\cite[Section~6]{KrumWica}, we can blow-up $T_k$ relative to $\mathbf{P}_0$ by using Almgren's Strong Lipschitz Approximation Theorem (Theorem~\ref{lip approx thm}) to find a Lispchitz approximating function $f_k : B_{3/4}(0) \rightarrow \mathcal{A}_q(\mathbb{R}^m)$ of $T_k$ in $\mathbf{C}_{3/4}(0),$ followed by standard arguments giving, after passing to a subsequence, $f_k/\widehat{E}_k \rightarrow w$ in $L^2(B_{3/4}(0), \mathcal{A}_q(\mathbb{R}^m))$.  After passing to a subsequence for each $i$ let $A^{(k)}_i/\widehat{E}_k \rightarrow \Lambda_i$ as $m \times 2$ matrices.  By Theorem~\ref{graphrep close2plane thm}(c) and \eqref{main hyp eqn3}, $w$ is also a blow-up of $\mathbf{C}_k$ relative to $\mathbf{P}_0$ and in particular 
\begin{equation}\label{branchpt close2plane1 eqn2}
	w(x,y) = \sum_{i=1}^p q_i \llbracket \Lambda_i x \rrbracket . 
\end{equation}
Using the Hardt-Simon inequality (\cite[Lemma~6.4]{KrumWica}), after passing to a subsequence $(\xi_k,\zeta_k) \rightarrow (\xi,\zeta)$ in $\overline{B_{1/2}(0)}$ and $\chi_k/\widehat{E}_k \rightarrow \lambda$ in $\mathbb{R}^m$ such that 
\begin{equation*}
	\int_{B_{1/4}(\xi,\zeta)} R^{2-n} \left| \frac{\partial}{\partial R} \left(\frac{w-\lambda}{R}\right)\right|^2 \leq C 
\end{equation*}
for some constant $C = C(n,m,q) \in (0,\infty)$ and thus $w(\xi,\zeta) = q \llbracket \lambda \rrbracket$.  But by \eqref{branchpt close2plane1 eqn2} and \eqref{close2plane eqn15}, we have $w(\xi,\zeta) = q \llbracket \lambda \rrbracket$ only if $\xi = 0$ and $\lambda = 0$.  Therefore, $\xi_k \rightarrow \xi = 0$ and $\chi_k/\widehat{E}_k \rightarrow \lambda = 0$, contradicting \eqref{branchpt close2plane1 eqn1}.
\end{proof}

\begin{lemma} \label{branchpt close2plane lemma2}
Let $3 \leq p_0 \leq q$ be integers such that (H4) holds true.  Given $M \in [1,\infty)$ and $\delta > 0$ there exist $\varepsilon,\beta,\gamma,\eta \in (0,1)$ depending only on $n,m,q,p_0,M,\delta$ such that the following holds true.  Let $\mathbf{C}$ and $T$ satisfy Hypothesis~$(\star)$, Hypothesis~$(\star\star)$ (of Section~\ref{notation-and-graphical}) and Hypothesis~$(\dagger)$ with $\varepsilon,\beta,\eta$ in place of $\varepsilon_0,\beta_0,\eta_0$.  Let $Z = (\chi,\xi,\zeta) \in \op{spt} T \cap \mathbf{B}_{1/2}(0)$ be such that $\Theta(T,Z) \geq q$.  Suppose that for some $s \in \{2,3,\ldots,p_0-1\}$
\begin{equation} \label{branchpt close2plane2 hyp}
	\inf_{\mathbf{C}' \in \bigcup_{p'=1}^s \mathcal{C}_{q,p'}} Q(T,\mathbf{C}',\mathbf{B}_1(0))
	\leq \gamma \inf_{\mathbf{C}' \in \bigcup_{p'=1}^{s-1} \mathcal{C}_{q,p'}} Q(T,\mathbf{C}',\mathbf{B}_1(0)) . 
\end{equation}
Then 
\begin{equation} \label{branchpt close2plane2 concl}
	|\chi|^2 + \widehat{E}^2 |\xi|^2 \leq \delta^2 \inf_{\mathbf{C}' \in \bigcup_{p'=1}^s \mathcal{C}_{q,p'}} Q(T,\mathbf{C}',\mathbf{B}_1(0))^2 , 
\end{equation}
where $\widehat{E} = E(T,\mathbf{P}_0,\mathbf{B}_1(0))$. 
\end{lemma}

\begin{remark}{\rm
By \eqref{close2plane eqn2} and \eqref{branchpt close2plane2 hyp}, $\inf_{\mathbf{C}' \in \bigcup_{p'=1}^s \mathcal{C}_{q,p'}} Q(T,\mathbf{C}',\mathbf{B}_1(0)) \leq \gamma_0 E(T, \mathbf{P}_0, \mathbf{B}_1(0))$. 
}\end{remark}

\begin{proof}[Proof of Lemma~\ref{branchpt close2plane lemma2}]
The proof is similar to that of~\cite[Proposition~10.7]{Wic14} with some minor changes and we sketch the proof for completion.  Without loss of generality, fix $M \in [1,\infty)$, $\delta > 0$ and $2 \leq s < p_0 \leq q$.  Suppose to the contrary that for $k = 1,2,3,\ldots$ there exists $\varepsilon_k \rightarrow 0^+$, $\beta_k \rightarrow 0^+$, $\gamma_k \rightarrow 0^+$, $\eta_k \rightarrow 0^+$, $\mathbf{C}_k \in \mathcal{C}_{q,p_0}$, and a locally area minimizing rectifiable current $T_k$ of $\mathbf{B}_1(0)$ such that Hypothesis~$(\star)$, Hypothesis~$(\star\star)$, Hypothesis~$(\dagger)$, \eqref{branchpt close2plane2 hyp} hold true with $\varepsilon_k,\beta_k,\gamma_k,\eta_k,\mathbf{C}_k, T_k$ in place of $\varepsilon_0,\beta_0,\gamma_0,\eta_0,\mathbf{C},T$ but for some $Z_k = (\chi_k,\xi_k,\zeta_k) \in \op{spt} T_k \cap \mathbf{B}_{1/2}(0)$ such that $\Theta(T_k,Z_k) \geq q$
\begin{equation}\label{branchpt close2plane2 eqn1}
	|\chi_k|^2 + \widehat{E}_k^2 |\xi_k|^2 \geq \delta^2 \widetilde{Q}_k^2 , 
\end{equation}
where $\widehat{E}_k = E(T_k,\mathbf{P}_0,\mathbf{B}_1(0))$ and $\widetilde{Q}_k = \inf_{\mathbf{C}' \in \bigcup_{p'=1}^s \mathcal{C}_{q,p'}} Q(T,\mathbf{C}',\mathbf{B}_1(0))$.  Let $\widetilde{\mathbf{C}}_k \in \mathcal{C}_{q,s}$ such that $Q(T_k, \widetilde{\mathbf{C}}_k, \mathbf{B}_1(0)) \leq 2 \widetilde{Q}_k$.  Express $\widetilde{\mathbf{C}}_k$ as $\widetilde{\mathbf{C}}_k = \sum_{i=1}^p \widetilde{q}_i \llbracket \widetilde{P}^{(k)}_i \rrbracket$ where $\widetilde{q}_i$ are integers with $\sum_{i=1}^p \widetilde{q}_i = q$ and $\widetilde{P}^{(k)}_i = \{ (z,x,y) \in \mathbb{R}^{n+m} : z = A^{(k)}_i x \}$ for some $m \times 2$ matrices $A^{(k)}_i$ with $\|A^{(k)}_i\| \leq C(n,m,q) \,\widehat{E}_k$.  After passing to a subsequence assume that $q_i$ are independent of $k$.  Let $\tau_k \rightarrow 0^+$ and $u^{(k)}_i : B_{3/4}(0) \cap \{r > \tau_k\} \rightarrow \mathcal{A}_{q_i}(\mathbb{R}^m)$ be as in Theorem~\ref{graphrep close2plane thm}(c) with $\widetilde{\mathbf{C}}_k, T_k, u^{(k)}_i$ in place of $\widetilde{\mathbf{C}}, T, u_i$.  After passing to a subsequence, blow-up $T_k$ relative to $\widetilde{\mathbf{C}}_k$ by letting $u^{(k)}_i/\widetilde{Q}_k \rightarrow w_i$ in $L^2(B_{3/4}(0) \cap \{r > \sigma\}, \mathcal{A}_{q_i}(\mathbb{R}^m))$ for each $\sigma > 0$.  By Theorem~\ref{graphrep close2plane thm}(c) and \eqref{main hyp eqn3}, $w_i$ is also a blow-up of $\mathbf{C}_k$ relative to $\widetilde{\mathbf{C}}_k$ and in particular 
\begin{equation}\label{branchpt close2plane2 eqn2}
	w_i(x,y) = \sum_{j=1}^{q_i} \llbracket M_{i,j} x \rrbracket . 
\end{equation}
for some $m \times 2$ matrices $M_{i,j}$ (not necessarily distinct).  By (H4) we can apply Theorem~\ref{branchpt close2plane thm}(a) with $T_k$ and $\widetilde{\mathbf{C}}_k$ in place of $T$ and $\mathbf{C}$ to obtain $|\chi_k|^2 + \widehat{E}_k^2 |\xi_k|^2 \leq C \widehat{E}_k^2$ for some constant $C = C(n,m,q,p_0,M) \in (0,\infty)$.  After passing to a subsequence let $\chi_k/\widetilde{Q}_k \rightarrow \lambda$ in $\mathbb{R}^m$, $\widehat{E}_k \xi_k/\widetilde{Q}_k \rightarrow \kappa$ in $\mathbb{R}^2$, $\zeta_k \rightarrow \zeta$ in $\mathbb{R}^{n-2}$, and $A^{(k)}_i/\widehat{E}_k \rightarrow \Lambda_i$ as $m \times 2$ matrices.  By (H4), we can apply Theorem~\ref{branchpt close2plane thm}(d) with $T_k$ and $\widetilde{\mathbf{C}}_k$ in place of $T$ and $\mathbf{C}$ to obtain 
\begin{equation*}
	\int_{B_{3/4}(0)} \frac{|w_i(x) - \lambda + \Lambda_i \kappa|^2}{|(x,y) - (0,\zeta)|^{n+3/2}} \leq C 
\end{equation*}
for some constant $C = C(n,m,q,p_0,M) \in (0,\infty)$ and thus by \eqref{branchpt close2plane2 eqn2} we must have that $\lambda - \Lambda_i \kappa = 0$ for all $1 \leq i \leq p$.  Hence $\Lambda_i \kappa = \Lambda_j \kappa$ for all $1 \leq i < j \leq p$, which by \eqref{close2plane eqn15} implies that $\kappa = 0$, and thus $\lambda = 0$.  Therefore, $\chi_k/\widetilde{Q}_k \rightarrow \lambda = 0$ and $\widehat{E}_k \xi_k/\widetilde{Q}_k \rightarrow \kappa = 0$, contradicting \eqref{branchpt close2plane2 eqn1}.
\end{proof}

\begin{corollary} \label{branchpt close2plane cor3}
Let $2 \leq p_0 \leq q$ be integers such that either $p_0 = 2$ or $p_0 > 2$ and (H4) holds true.  Given $M \in [1,\infty)$ and $\delta > 0$ there exists $\varepsilon_0,\beta_0,\eta_0 \in (0,1)$ depending only on $n,m,q,p_0,M,\delta$ such that if $\mathbf{C}$ and $T$ satisfy Hypothesis~$(\star)$, Hypothesis~$(\star\star)$, and Hypothesis~$(\dagger)$ and $Z = (\chi,\xi,\zeta) \in \op{spt} T \cap \mathbf{B}_{1/2}(0)$ such that $\Theta(T,Z) \geq q$, then 
\begin{equation} \label{branchpt close2plane3 concl}
	|\chi|^2 + \widehat{E}^2 |\xi|^2 \leq \delta^2 \inf_{\mathbf{C}' \in \bigcup_{p'=1}^{p_0-1} \mathcal{C}_{q,p'}} Q(T,\mathbf{C}',\mathbf{B}_1(0))^2 , 
\end{equation}
where $\widehat{E} = E(T,\mathbf{P}_0,\mathbf{B}_1(0))$. 
\end{corollary}

\begin{proof}
Corollary~\ref{branchpt close2plane cor3} follows from Lemma~\ref{branchpt close2plane lemma1} and Lemma~\ref{branchpt close2plane lemma2} in the same way that~\cite[Lemma~10.4]{Wic14} follows from~\cite[Lemma~10.6]{Wic14} and~\cite[Proposition~10.7]{Wic14}. 
\end{proof}

\begin{proof}[Proof of Theorem~\ref{branchpt close2plane thm}]
We proceed by induction.  Suppose that $p_0 \in \{2,3,\ldots,q\}$ such that either $p_0 = 2$ or $p_0 > 2$ and (H4) holds true.  We want to show that the conclusion of Theorem~\ref{branchpt close2plane thm} holds true when $p = p_0$.  Arguing as in~\cite[Proposition~10.5]{Wic14} using Theorem~\ref{graphrep close2plane thm}(c), Corollary~\ref{branchpt close2plane cor3}, 
\eqref{close2plane eqn11}, and \eqref{close2plane eqn12}, for every $\rho \in (0,1/4]$ and every $\varepsilon, \beta, \eta \in (0,1)$ there exists $\varepsilon_0, \beta_0, \eta_0 \in (0,1)$ depending only on $n,m,q,M,\rho, \varepsilon, \beta, \eta$ such that if $\mathbf{C} \in \mathcal{C}_{q,p_0}$ and $T$ satisfy Hypotheses~$(\star)$, $(\star\star)$, and $(\dagger)$ and if $Z \in \op{spt} T \cap \mathbf{B}_{1/2}(0)$ with $\Theta(T,Z) \geq q$, then $\mathbf{C}$ and $\eta_{Z,\rho\#} T$ satisfy Hypotheses~$(\star)$, $(\star\star)$, and $(\dagger)$ with $\varepsilon,\beta,\eta, C(n,m,q)\,M$ in place of $\varepsilon_0,\beta_0,\eta_0,M$.  In particular, by Corollary~\ref{keyest cor1}  
\begin{equation}\label{branchpt close2plane eqn1}
	\int_{\mathbf{B}_{3\rho/4}(Z)} \frac{\op{dist}^2(X-Z, \op{spt} \mathbf{C}) \,d\|T\|(X)}{|X-Z|^{n+2-\sigma}}
		\leq C \rho^{-n-2+ \sigma} \int_{\mathbf{B}_{\rho}(Z)} \op{dist}^2(X-Z, \op{spt} \mathbf{C}) \,d\|T\|(X)
\end{equation}
for some constant $C = C(n,m,q,p_0,\sigma) \in (0,\infty)$. 

Now the conclusion of Theorem~\ref{branchpt close2plane thm} follows by arguing as in~\cite[Lemma~6.21]{Wic04}.  Let $\mathbf{C} = \sum_{i=1}^p q_i \llbracket P_i \rrbracket$ where $q_i \geq 1$ are integers such that $\sum_{i=1}^p q_i = q$ and $P_i$ are $n$-dimensional oriented planes given by \eqref{close2plane eqn4} for some $m \times 2$ matrix $A_i$ satisfying \eqref{close2plane eqn5}.  Let $0 < \tau < \gamma < 1$ and, assuming $\varepsilon_0,\beta_0,\eta_0$ are sufficiently small, let $T_i$ be as in Theorem~\ref{graphrep close2plane thm}(b).  Note that given $\delta > 0$ and assuming $\varepsilon_0,\beta_0,\eta_0$ are sufficiently small depending on $n,m,q,p_0,M,\delta$, by Lemma~\ref{branchpt close2plane cor3}, \eqref{close2plane eqn4}, and \eqref{close2plane eqn5}, 
\begin{equation*}
	\op{dist}_{\mathcal{H}}(\op{spt} \mathbf{C} \cap \mathbf{B}_1(0), (Z+\op{spt} \mathbf{C}) \cap \mathbf{B}_1(0)) 
	\leq \max_{1 \leq i \leq p} |A_i \xi - \chi| \leq C \delta \inf_{\mathbf{C}' \in \bigcup_{p'=1}^{p_0-1} \mathcal{C}_{q,p'}} Q(T,\mathbf{C}',\mathbf{B}_1(0)) . 
\end{equation*}
Thus by Theorem~\ref{graphrep close2plane thm}(a)(b), provided we take $\delta = \delta(n,m,q,p_0,M,\tau)$ to be sufficiently small, for each $X \in \op{spt} T_i$ the closest point $X'$ to $X$ on $\op{spt} \mathbf{C}$ lies on $P_i$ and the closest point to $X$ on $Z+\op{spt} \mathbf{C}$ lies on $Z + P_i$.  Hence 
\begin{equation*}
	\op{dist}(X, Z+\op{spt} \mathbf{C}) = |X' - X - \pi_{P_i^{\perp}} (\chi,\xi,0)|
\end{equation*}
(as in (6.31) of~\cite{Wic04}).  The main change from~\cite[Lemma~6.21]{Wic04} is that we obtain (6.34) of~\cite{Wic04} as follows: we want to show that for some $i_0 \in \{1,2,\ldots,p\}$ 
\begin{equation}\label{branchpt close2plane eqn5}
	|\pi_{P_{i_0}^{\perp}} (\chi,\xi,0)| \geq c \,\widehat{E} \,|(\chi,\xi)| 
\end{equation}
for some constant $c = c(n,m,q,M) > 0$.  Suppose to the contrary that $|\pi_{P_i^{\perp}} (\chi,\xi,0)| < \widehat{E} \,|(\chi,\xi)|$ for all $i \in \{1,2,\ldots,p\}$.  By taking $P = P_i$ to be one of the planes of $\mathbf{C}$ in \eqref{close2plane eqn9} and using \eqref{graphrep concl a2}, 
\begin{equation*}
	\widehat{E} \leq CM \max_{1 \leq i < j \leq p} \op{dist}_{\mathcal{H}}(P_i \cap \mathbf{B}_1(0), P_j \cap \mathbf{B}_1(0))
		\leq CM \max_{1 \leq i < j \leq p} \inf_{X \in P_i \cap (\mathbf{S}^{m+1} \times \mathbb{R}^{n-2})} \op{dist}(X,P_j) ,
\end{equation*}
where $\widehat{E} = E(T, \mathbf{P}_0, \mathbf{B}_1(0))$ and $C = C(n,m,q) \in (0,\infty)$ are constants.  Thus there exists $i,j \in \{1,2,\ldots,p\}$ such that 
\begin{equation*}
	|\pi_{P_i^{\perp}} (\chi,\xi,0) - \pi_{P_j^{\perp}} (\chi,\xi,0)| \geq |\pi_{P_i} (\chi,\xi,0) - \pi_{P_j} (\chi,\xi,0)| \geq 2c \,\widehat{E} \,|(\chi,\xi)|
\end{equation*}
for some constant $c = c(n,m,q,M) > 0$.  Therefore \eqref{branchpt close2plane eqn5} holds true for either $i_0 = i$ or $i_0 = j$.  Now let $\rho_0 = \rho_0(n,m,q,p_0,M) \in (0,1/8]$ and $\tau = \tau(n,m,q,p_0,M,\rho_0) \in (0,\rho_0/2]$ be constants to be later determined.  Provided $\tau$ is sufficiently small, there exists a $\|T_{i_0}\|$-measurable set $S \subset \op{spt} T_{i_0} \cap \mathbf{B}_{\rho_0}(Z) \cap \{|x| > \tau\}$ such that $\|T_{i_0}\|(S) \geq \tfrac{1}{2} \,\omega_n \rho_0^n$.  
Integrating \eqref{branchpt close2plane eqn5} over $S$ gives us 
\begin{equation}\label{branchpt close2plane eqn6}
	\widehat{E}^2 |(\chi,\xi)|^2 \leq \frac{C}{\rho_0^n} \sum_{i=1}^p \int_{\mathbf{B}_{\rho_0}(Z) \cap \{r > \tau\}} |\pi_{P_i^{\perp}} (\chi,\xi,0)|^2 \,d\|T_i\|(X) 
\end{equation}
for some constant $C = C(n,m,q,p_0,M) \in (0,\infty)$.  Inequality (6.34) of~\cite{Wic04} follows from \eqref{branchpt close2plane eqn6} and (6.35) of~\cite{Wic04}.
\end{proof}

\begin{corollary} \label{nonconcentration close2plane cor}
For all $\delta \in (0,1/4)$, $\sigma \in (0,1)$, and $M \in [1,\infty)$ there exists $\varepsilon_0,\beta_0,\eta_0 \in (0,1)$ depending only on $n,m,q,M,\delta,\sigma$ such that if $\mathbf{C}$ and $T$ satisfy Hypothesis~$(\star)$, Hypothesis~$(\star\star)$, and Hypothesis~$(\dagger)$ and if 
\begin{equation*}
	\mathbf{B}_{\delta}(0,y_0) \cap \{ X \in \mathbf{B}_{3/4}(0) : \Theta(T,X) \geq q \} \neq \emptyset 
\end{equation*}
for all $y_0 \in B^{n-2}_{1/2}(0)$, then 
\begin{equation} \label{noncencentration close2plane concl} 
	\int_{\mathbf{B}_{1/2}(0)} \frac{\op{dist}^2(X, \op{spt} \mathbf{C})}{r_{\delta}^{2-\sigma}} \,d\|T\|(X) 
		\leq C \int_{\mathbf{B}_1(0)} \op{dist}^2(X, \op{spt} \mathbf{C}) \,d\|T\|(X) , 
\end{equation}
where $r_{\delta} = \max\{|x|,\delta\}$ and $C = C(n,m,q,M,\sigma) \in (0,\infty)$ is a constant (independent of $\delta$). 
\end{corollary}

\begin{proof}
This follows from Theorem~\ref{branchpt close2plane thm} by arguing as in the proof of~\cite[Corollary~3.2]{Sim93}.
\end{proof}

\subsection{A priori estimates for area-minimizers not close to a plane} \label{sec:apriori nonplanar}

Here we will consider the case where, in contrast with Hypothesis~$(\dagger)$, $T$ is not close to a $n$-dimensional plane.  Thus, given a constant $\eta_0 \in (0,1),$ we shall assume that: 

\noindent\textbf{Hypothesis~$(\dagger\dagger)$.}  For $q \geq 2$ an integer and $T$ an $n$-dimensional locally area-minimizing rectifiable current in $\mathbf{B}_1(0)$ with $T \llcorner {\mathbf B}_{1}(0) = 0,$ suppose that 
\begin{equation}\label{main hyp eqn6} 
	\inf_{\mathbf{P} \in \mathcal{C}_{q,1}} E(T, \mathbf{P}, \mathbf{B}_1(0)) \geq \eta_0 . 
\end{equation}

We shall represent points $X \in \mathbb{R}^{n+m}$ by $X = (x,y)$ where $x \in \mathbb{R}^{m+2}$ and $y \in \mathbb{R}^{n-2}$.  
Let $p$ be an integer with $2 \leq p \leq q$, and let $\mathbf{C} = \sum_{i=1}^p q_i \llbracket P_i \rrbracket \in {\mathcal C}_{q, p}$ and suppose that $\mathbf{C}$ and $T$ satisfy Hypothesis~$(\star)$, Hypothesis~$(\star\star)$, and Hypothesis~$(\dagger\dagger)$.    By the triangle inequality, the assumption $\|T\|(\mathbf{B}_1(0)) \leq (q+1/2)\,\omega_n$, and \eqref{main hyp eqn3}, we have that for each $i = 1,2,\ldots,p,$ 
\begin{align*}
	&\int_{\mathbf{B}_1(0)} \op{dist}^2(X, P_i) \,d\|T\|(X) 
	\\ &\hspace{1in}\leq 2 \int_{\mathbf{B}_1(0)} \op{dist}^2(X, \op{spt} \mathbf{C}) \,d\|T\|(X) 
		+ 2 (q+1/2) \omega_n \sup_{X \in \op{spt} \mathbf{C} \cap \mathbf{B}_1(0)} \op{dist}^2(X, P_i) 
	\\ & \hspace{1.5in}\leq 2 \beta_0^2 \int_{\mathbf{B}_1(0)} \op{dist}^2(X, P_i) \,d\|T\|(X) 
		+ 2 (q+1/2) \omega_n \sup_{X \in \op{spt} \mathbf{C} \cap \mathbf{B}_1(0)} \op{dist}^2(X, P_i) .
\end{align*}
Hence assuming $\beta_0 < 1/2$ and using \eqref{main hyp eqn4}, 
\begin{align}\label{nonplanar eqn1}
	\eta_0^2 \leq \frac{1}{\omega_n} \int_{\mathbf{B}_1(0)} \op{dist}^2(X, P_i) \,d\|T\|(X) 
	&\leq 4 (q+1/2) \sup_{X \in \op{spt} \mathbf{C} \cap \mathbf{B}_1(0)} \op{dist}^2(X, P_i) 
	\\&\leq 4 (q+1/2) \max_{1 \leq j \leq p} \op{dist}_{\mathcal{H}}(P_i \cap \mathbf{B}_1(0), P_j \cap \mathbf{B}_1(0)) . \nonumber 
\end{align}

We wish to obtain estimates similar to those of Corollary~\ref{keyest cor2}, Theorem~\ref{branchpt close2plane thm}, and Corollary~\ref{nonconcentration close2plane cor} in this setting.  

First, as a consequence of Theorem~\ref{keyest thm} and Theorem~\ref{graphrep thm}, we have: 

\begin{corollary} \label{keyest cor3}
For each $0 < \tau < \gamma < 1$  there exists $\varepsilon_0, \beta_0 \in (0,1)$ depending only on $n,m,q,\gamma,\tau$ such that if $\mathbf{C}$ and $T$ satisfy Hypothesis~$(\star)$ and Hypothesis~$(\star\star)$, then: 
\begin{equation} \label{keyest cor3 concl}
	\sum_{i=1}^p \int_{\mathbf{B}_{\gamma}(0,P_i) \cap \{r > \tau\}} R^{2-n} \left| \frac{\partial (u_i/R)}{\partial R} \right|^2 
		\leq C \int_{\mathbf{B}_1(0)} \op{dist}^2(X, \op{spt} \mathbf{C}) \,d\|T\|(X) , 
\end{equation}
where $R(x,y) = |(x,y)|$, $u_i : B_{\gamma}(0,P_i) \cap \{r > \tau\} \rightarrow \mathcal{A}_{q_i}(P_i^{\perp})$ are as in Theorem~\ref{graphrep thm}(c), and $C = C(n,m,q,\gamma) \in (0,\infty)$ is a constant. 
\end{corollary}

\begin{proof}
Let $(x',y)$ denote points on $P_i$, where $x' \in P_i \cap (\mathbb{R}^{m+2} \times \{0\})$ and $y \in \mathbb{R}^{n-2}$.  For each $(x',y) \in B_{\gamma}(0) \cap \{r > \tau\}$ let $u_i(x',y) = \sum_{j=1}^{q_i} \llbracket u_{i,j}(x',y) \rrbracket$ where $u_{i,j}(x',y) \in \mathbb{R}^m$.  Arguing as in the proof of Corollary~\ref{keyest cor2}, for $\mathcal{H}^n$-a.e.~$(x',y) \in K_i$ at the point $(x,y) = (x',y) + u_{i,j}(x',y)$ 
\begin{equation*}
	|(x,y)^{\perp}| \geq \frac{1}{2} \,R^2 \,\Bigg|\frac{\partial}{\partial R}\Bigg(\frac{u_{i,j}(x',y)}{R}\Bigg)\Bigg| . 
\end{equation*}
Now \eqref{keyest cor3 concl} follows from Theorem~\ref{keyest thm} and Theorem~\ref{graphrep thm}(c). 
\end{proof}

Next were derive a priori estimates based on~\cite[Theorem~3.1]{Sim93} and~\cite[Corollary~10.2]{Wic14}. 

\begin{theorem} \label{branchpt nonplanar thm}
For all integers $p$, $q$ with $2 \leq p \leq q$, each $\eta_0 \in (0,1)$, each $\tau, \gamma$ with $0 < \tau < \gamma < 1$, and each $\sigma \in (0,1)$ there exists $\varepsilon_0,\beta_0 \in (0,1)$ depending only on $n,m,q,p,\eta_0,\gamma,\tau,\sigma$ such that the following holds true.  Suppose that $\mathbf{C}$ and $T$ satisfy Hypothesis~$(\star)$, Hypothesis~$(\star\star)$, and Hypothesis~$(\dagger\dagger)$.  Let $Z = (\xi,\zeta) \in \op{spt} T \cap \mathbf{B}_{1/2}(0)$ be such that $\Theta(T,Z) \geq q$.  Then 
\begin{align*} 
	&(a) \quad |\xi|^2 \leq C \int_{\mathbf{B}_1(0)} \op{dist}^2(X, \op{spt} \mathbf{C}) \,d\|T\|(X) , \\
	&(b) \quad \int_{\mathbf{B}_{1/4}(Z)} \op{dist}(X-Z,\op{spt} \mathbf{C})^2 \,d\|T\|(X) 
		\leq C \int_{\mathbf{B}_1(0)} \op{dist}^2(X, \op{spt} \mathbf{C}) \,d\|T\|(X) ,
\end{align*}
where $C = C(n,m,q,p) \in (0,\infty)$ is a constant (independent of $\tau$).  Moreover, 
\begin{align*} 
	&(c) \quad \int_{\mathbf{B}_{\gamma}(0)} \frac{\op{dist}^2(X,\op{spt} \mathbf{C}) \,d\|T\|(X)}{|X-Z|^{n+2-\sigma}}
		\leq C \int_{\mathbf{B}_1(0)} \op{dist}^2(X, \op{spt} \mathbf{C}) \,d\|T\|(X) , \\
	&(d) \quad \int_{B_{\gamma}(0) \cap \{r > \tau\}} \frac{|u_i(x,y) - \pi_{P_i^{\perp}} \xi|^2}{|(x,y) - (\xi,\zeta)|^{n+2-\sigma}} 
		\leq C \int_{\mathbf{B}_1(0)} \op{dist}^2(X, \op{spt} \mathbf{C}) \,d\|T\|(X) \text{ for all $i$}
\end{align*}
where $u_i : B_{\gamma}(0,P_i) \cap \{r > \tau\} \rightarrow \mathcal{A}_{q_i}(P_i^{\perp})$ are as in Theorem~\ref{graphrep close2plane thm}(c) and $C = C(n,m,q,p,\gamma,\sigma) \in (0,\infty)$ is a constant (independent of $\tau$). 
\end{theorem}

The proof of Theorem~\ref{branchpt nonplanar thm} is similar to that of Theorem~\ref{branchpt close2plane thm} above.  We shall assume the inductive hypothesis that $p_0 \in \{2,3,\ldots,q\}$ such that either $p_0 = 2$ or $p_0 > 2$ and 
\begin{enumerate}
	\item[(H5)]  Theorem~\ref{branchpt nonplanar thm} holds true for all $p \in \{2,3,\ldots,p_0-1\}$ . 
\end{enumerate}
Theorem~\ref{branchpt nonplanar thm} in the case $p = p_0$ will follow from two preliminary results Lemma~\ref{branchpt nonplanar lemma2} and Corollary~\ref{branchpt nonplanar cor3}.  Note that by Theorem~\ref{graphrep thm}(a), for each $\tau > 0$ there exists $\varepsilon_0 = \varepsilon_0(\tau) \in (0,1)$ and $\beta_0 = \beta_0(\tau) \in (0,1)$ with $\varepsilon_0(\tau) \rightarrow 0^+$ and $\beta_0(\tau) \rightarrow 0^+$ as $\tau \rightarrow 0^+$ such that if $\mathbf{C}$ and $T$ satisfy Hypothesis~$(\star)$ and Hypothesis~$(\star\star)$ and $Z = (\xi,\zeta) \in \op{spt} T \cap \mathbf{B}_{1/2}(0)$ such that $\Theta(T,Z) \geq q$ then 
\begin{equation} \label{branchpt nonplanar2 eqn1}
	|\xi| \leq \tau . 
\end{equation}

\begin{lemma} \label{branchpt nonplanar lemma2}
Let $2 \leq p_0 \leq q$ be integers such that either $p_0 = 2$ or $p_0 > 2$ and (H5) holds true.  Given $\eta \in (0,1)$ and $\delta > 0$ there exists $\varepsilon,\beta,\gamma \in (0,1)$ depending only on $n,m,q,p_0,\eta,\delta$ such that the following holds true.  Let $\mathbf{C}$ and $T$ satisfy Hypothesis~$(\star)$, Hypothesis~$(\star\star)$, and Hypothesis~$(\dagger\dagger)$ with $\varepsilon,\beta,\eta$ in place of $\varepsilon_0,\beta_0,\eta_0$.  Let $Z = (\xi,\zeta) \in \op{spt} T \cap \mathbf{B}_{1/2}(0)$ be such that $\Theta(T,Z) \geq q$.  Suppose that for some $s \in \{1,2,\ldots,p_0-1\}$ either $s = 1$ or $s > 1$ and 
\begin{equation} \label{branchpt nonplanar2 hyp}
	\inf_{\mathbf{C}' \in \bigcup_{p'=1}^s \mathcal{C}_{q,p'}} Q(T,\mathbf{C}',\mathbf{B}_1(0))
	\leq \gamma \inf_{\mathbf{C}' \in \bigcup_{p'=1}^{s-1} \mathcal{C}_{q,p'}} Q(T,\mathbf{C}',\mathbf{B}_1(0)) . 
\end{equation}
Then 
\begin{equation} \label{branchpt nonplanar2 concl}
	|\xi| \leq \delta \inf_{\mathbf{C}' \in \bigcup_{p'=1}^s \mathcal{C}_{q,p'}} Q(T,\mathbf{C}',\mathbf{B}_1(0)) . 
\end{equation}
\end{lemma}

\begin{proof}[Proof of Lemma~\ref{branchpt nonplanar lemma2}]
The proof is similar to that of~\cite[Proposition~10.7]{Wic14} with some minor changes and we sketch the proof for completion.  Without loss of generality, fix $\eta,\delta > 0$ and $1 \leq s < p_0 \leq q$.  Suppose to the contrary that for $k = 1,2,3,\ldots$ there exists $\varepsilon_k \rightarrow 0^+$, $\beta_k \rightarrow 0^+$, $\gamma_k \rightarrow 0^+$, $\mathbf{C}_k \in \mathcal{C}_{q,p_0}$, and a locally area minimizing rectifiable current $T_k$ of $\mathbf{B}_1(0)$ such that Hypothesis~$(\star)$, Hypothesis~$(\star\star)$, and Hypothesis~$(\dagger\dagger)$ hold true with $\varepsilon_k,\beta_k,\gamma_k,\eta,\mathbf{C}_k, T_k$ in place of $\varepsilon_0,\beta_0,\gamma_0,\eta_0,\mathbf{C},T$ but for some $Z_k = (\chi_k,\xi_k,\zeta_k) \in \op{spt} T_k \cap \mathbf{B}_{1/2}(0)$ such that $\Theta(T_k,Z_k) \geq q$
\begin{equation}\label{branchpt nonplanar2 eqn2}
	|\xi_k|^2 \geq \delta^2 \widetilde{Q}_k^2 , 
\end{equation}
where $\widetilde{Q}_k = \inf_{\mathbf{C}' \in \bigcup_{p'=1}^s \mathcal{C}_{q,p'}} Q(T,\mathbf{C}',\mathbf{B}_1(0))$.  Note that by \eqref{branchpt nonplanar2 eqn1}, we know that $\xi_k \rightarrow 0$ as $k \rightarrow \infty$.  Thus by \eqref{branchpt nonplanar2 eqn2} we must have that $\widetilde{Q}_k \rightarrow 0^+$.  Moreover, by Hypothesis~$(\dagger\dagger)$ we must have that $s \geq 2$.  Let $\widetilde{\mathbf{C}}_k \in \mathcal{C}_{q,s}$ such that $Q(T_k, \widetilde{\mathbf{C}}_k, \mathbf{B}_1(0)) \leq 2 \widetilde{Q}_k$.  Express $\widetilde{\mathbf{C}}_k$ as $\widetilde{\mathbf{C}}_k = \sum_{i=1}^p \widetilde{q}_i \llbracket \widetilde{P}^{(k)}_i \rrbracket$ where $\widetilde{q}_i$ are integers with $\sum_{i=1}^p \widetilde{q}_i = q$ and $\widetilde{P}^{(k)}_i$ are $n$-dimensional oriented planes of $\mathbb{R}^{n+m}$ with $\{0\} \times \mathbb{R}^{n-2} \subset \widetilde{P}^{(k)}_i$.  After passing to a subsequence, assume $\widetilde{q}_i$ are independent of $k$ and find $n$-dimensional planes $\widetilde{P}^{(\infty)}_i$ such that 
\begin{equation*}
	\lim_{k \rightarrow \infty} \op{dist}_{\mathcal{H}}(\widetilde{P}^{(k)}_i \cap \mathbf{B}_1(0), \widetilde{P}^{(\infty)}_i \cap \mathbf{B}_1(0)) = 0 . 
\end{equation*}
Thus we can regard $\widetilde{\mathbf{C}}_k$ as converging weakly to $\widetilde{\mathbf{C}}_{\infty} = \sum_{i=1}^p \widetilde{q}_i \llbracket \widetilde{P}^{(\infty)}_i \rrbracket$ in $\mathbb{R}^{n+m}$ for an appropriate choice of orientation of $\widetilde{P}^{(\infty)}_i$.  Let $\tau_k \rightarrow 0^+$ and $u^{(k)}_i : B_{3/4}(0,\widetilde{P}^{(k)}_i) \cap \{r > \tau_k\} \rightarrow \mathcal{A}_{q_i}((\widetilde{P}^{(k)}_i)^{\perp})$ be as in Theorem~\ref{graphrep thm}(c) with $\widetilde{\mathbf{C}}_k, T_k, u^{(k)}_i$ in place of $\widetilde{\mathbf{C}}, T, u_i$.  Let $\mathfrak{q}^{(k)}_i : \mathbb{R}^{n+m} \rightarrow \mathbb{R}^{n+m}$ be an orthogonal linear transformation such that 
\begin{equation*}
	\mathfrak{q}^{(k)}_i(\{0\} \times \mathbb{R}^{n-2}) = \{0\} \times \mathbb{R}^{n-2}, \quad 
	\mathfrak{q}^{(k)}_i(\widetilde{P}^{(k)}_i) = \widetilde{P}^{(\infty)}_i, \quad 
	\mathfrak{q}^{(k)}_i((\widetilde{P}^{(k)}_i)^{\perp}) = (\widetilde{P}^{(\infty)}_i)^{\perp} . 
\end{equation*}
After passing to a subsequence, blow-up $T_k$ relative to $\widetilde{\mathbf{C}}_k$ by letting $(\mathfrak{q}^{(k)}_i \circ u^{(k)}_i \circ (\mathfrak{q}^{(k)}_i)^{-1})/\widetilde{Q}_k \rightarrow w_i$ in $L^2(B_{3/4}(0,\widetilde{P}^{(\infty)}_i) \cap \{r > \sigma\}, \mathcal{A}_{q_i}((\widetilde{P}^{(\infty)}_i)^{\perp})$ for each $\sigma > 0$.  By Theorem~\ref{graphrep thm}(c) and \eqref{main hyp eqn3}, $w_i$ is also a blow-up of $\mathbf{C}_k$ relative to $\widetilde{\mathbf{C}}_k$ and in particular 
\begin{equation}\label{branchpt nonplanar2 eqn3}
	w_i(x,y) = \sum_{j=1}^{q_i} \llbracket \psi_{i,j}(x) \rrbracket . 
\end{equation}
for each $(x,y) \in \widetilde{P}^{(\infty)}_i$, where $\psi_{i,j} : \widetilde{P}^{(\infty)}_i \cap (\mathbb{R}^{m+2} \times \{0\}) \rightarrow (\widetilde{P}^{(\infty)}_i)^{\perp}$ are linear maps (not necessarily distinct).  By (H5) we can apply Theorem~\ref{branchpt nonplanar thm}(a) with $T_k$ and $\widetilde{\mathbf{C}}_k$ in place of $T$ and $\mathbf{C}$ to obtain $|\xi_k| \leq C \widetilde{Q}_k$ for some constant $C = C(n,m,q,p_0) \in (0,\infty)$.  After passing to a subsequence let $\xi_k/\widetilde{Q}_k \rightarrow \kappa$ in $\mathbb{R}^m$ and $\zeta_k \rightarrow \zeta$ in $\mathbb{R}^{n-2}$.  By (H5), we can apply Theorem~\ref{branchpt close2plane thm}(d) with $T_k$ and $\widetilde{\mathbf{C}}_k$ in place of $T$ and $\mathbf{C}$ to obtain 
\begin{equation*}
	\int_{B_{3/4}(0)} \frac{|w_i(x) - \pi_{(\widetilde{P}^{(\infty)}_i)^{\perp}} (\kappa,0)|^2}{|(x,y) - (0,\zeta)|^{n+3/2}} \leq C 
\end{equation*}
for some constant $C = C(n,m,q,p_0) \in (0,\infty)$ and thus by \eqref{branchpt nonplanar2 eqn3} we must have that $\pi_{(\widetilde{P}^{(\infty)}_i)^{\perp}} \kappa = 0$ for all $1 \leq i \leq p$.  In other words, $(\kappa,0) \in \widetilde{P}^{(\infty)}_i$ for all $i$.  However, by \eqref{nonplanar eqn1} and \eqref{graphrep concl a2} with $\widetilde{\mathbf{C}}_k$ in place of $\mathbf{C}$, after passing to a subsequence there exists $i \neq j$ such that 
\begin{equation*}
	\inf_{X \in \widetilde{P}^{(k)}_i \cap (\mathbb{S}^{m+1} \times \mathbb{R}^{n-2})} \op{dist}(X,\widetilde{P}^{(k)}_j) \geq c \eta
\end{equation*}
for some constant $c = c(n,m,q) > 0$.  Thus letting $k \rightarrow \infty$, 
\begin{equation*}
	\inf_{X \in \widetilde{P}^{(\infty)}_i \cap (\mathbb{S}^{m+1} \times \mathbb{R}^{n-2})} \op{dist}(X,\widetilde{P}^{(\infty)}_j) \geq c \eta
\end{equation*}
for some constant $c = c(n,m,q) > 0$.  That is, $\widetilde{P}^{(\infty)}_i \cap \widetilde{P}^{(\infty)}_j = \{0\} \times \mathbb{R}^{n-2}$.  Since $(\kappa,0) \in \widetilde{P}^{(\infty)}_i \cap \widetilde{P}^{(\infty)}_j$, we must have that $\kappa = 0$.  Therefore, $\xi_k/\widetilde{Q}_k \rightarrow \kappa = 0$, contradicting \eqref{branchpt nonplanar2 eqn2}.
\end{proof}

\begin{corollary} \label{branchpt nonplanar cor3}
Let $2 \leq p_0 \leq q$ be integers such that either $p_0 = 2$ or $p_0 > 2$ and (H5) holds true.  Given $\eta_{0} \in (0,1)$ and $\delta > 0$ there exists $\varepsilon_0,\beta_0 \in (0,1)$ depending only on $n,m,q,p_0,\eta_{0},\delta$ such that if $\mathbf{C}$ and $T$ satisfy Hypothesis~$(\star)$, Hypothesis~$(\star\star)$, and Hypothesis~$(\dagger\dagger)$ and if $Z = (\xi,\zeta) \in \op{spt} T \cap \mathbf{B}_{1/2}(0)$ is such that $\Theta(T,Z) \geq q$, then 
\begin{equation} \label{branchpt close2plane3 concl}
	|\xi| \leq \delta \inf_{\mathbf{C}' \in \bigcup_{p'=1}^{p_0-1} \mathcal{C}_{q,p'}} Q(T,\mathbf{C}',\mathbf{B}_1(0)) .
\end{equation}
\end{corollary}

\begin{proof}
Corollary~\ref{branchpt nonplanar cor3} follows from Lemma~\ref{branchpt nonplanar lemma2} in much the same way that~\cite[Lemma~10.4]{Wic14} follows from~\cite[Lemma~10.6]{Wic14} and~\cite[Proposition~10.7]{Wic14}. 
\end{proof}

\begin{proof}[Proof of Theorem~\ref{branchpt nonplanar thm}]
We proceed by induction.  Suppose that $p_0 \in \{2,3,\ldots,q\}$ such that either $p_0 = 2$ or $p_0 > 2$ and (H5) holds true.  We want to show that the conclusion of Theorem~\ref{branchpt nonplanar thm} holds true when $p = p_0$.  We claim that for every $\rho \in (0,1/4]$ and every $\eta_0, \varepsilon, \beta \in (0,1)$ there exists $\varepsilon_0, \beta_0 \in (0,1)$ depending only on $n,m,q,\rho$ such that if $\mathbf{C} \in \mathcal{C}_{q,p_0}$ and $T$ satisfy Hypotheses~$(\star)$, $(\star\star)$, and $(\dagger\dagger)$ and if $Z \in \op{spt} T \cap \mathbf{B}_{1/2}(0)$ with $\Theta(T,Z) \geq q$, then $\mathbf{C}$ and $\eta_{Z,\rho\#} T$ satisfy Hypotheses~$(\star)$ and $(\star\star)$ with $\varepsilon,\beta$ in place of $\varepsilon_0,\beta_0$.  Hypotheses~$(\star)$ is readily verified using \eqref{branchpt nonplanar2 eqn1}.  Hypotheses~$(\star\star)$ is verified by arguing as in~\cite[Proposition~10.5]{Wic14} using Theorem~\ref{graphrep close2plane thm}(c) and Corollary~\ref{branchpt nonplanar cor3}.  Hence by Corollary~\ref{keyest cor1}, \eqref{branchpt close2plane eqn1} holds true.  

Now let $\mathbf{C} = \sum_{i=1}^p q_i \llbracket P_i \rrbracket$ where $q_i \geq 1$ are integers such that $\sum_{i=1}^p q_i = q$ and $P_i$ are $n$-dimensional oriented planes with $\{0\} \times \mathbb{R}^{n-2} \subset P_i$.  Let $0 < \tau < \gamma < 1$ and, assuming $\varepsilon_0,\beta_0$ are sufficiently small, let $T_i$ be as in Theorem~\ref{graphrep thm}(b).  Arguing as we did for Theorem~\ref{branchpt close2plane thm}, it follows using Theorem~\ref{graphrep thm} and Lemma~\ref{branchpt nonplanar cor3} that for each $X \in \op{spt} T_i$ the closest point $X'$ to $X$ on $\op{spt} \mathbf{C}$ lies on $P_i$ and the closest point to $X$ on $Z+\op{spt} \mathbf{C}$ lies on $Z + P_i$.  Hence 
\begin{equation}\label{branchpt nonplanar eqn1}
	\op{dist}(X, Z+\op{spt} \mathbf{C}) = |X' - X - \pi_{P_i^{\perp}} (\xi,0)|
\end{equation}
In particular, since $\op{dist}(X,\op{spt}\mathbf{C}) = |X'-X|$, 
\begin{gather}
	\label{branchpt nonplanar eqn2} |\pi_{P_i^{\perp}} (\xi,0)| \leq \op{dist}(X, \op{spt} \mathbf{C}) + \op{dist}(X, Z+\op{spt} \mathbf{C}) , \\
	\label{branchpt nonplanar eqn3} \big| \op{dist}(X, \op{spt} \mathbf{C}) - \op{dist}(X, Z+\op{spt} \mathbf{C}) \big| \leq |\xi| . 
\end{gather}
By \eqref{graphrep concl a2} and \eqref{nonplanar eqn1} there exists $i,j \in \{1,2,\ldots,p\}$ such that 
\begin{equation*}
	|\pi_{P_i^{\perp}} (\xi,0) - \pi_{P_j^{\perp}} (\xi,0)| \geq |\pi_{P_i} (\xi,0) - \pi_{P_j} (\chi,\xi,0)| \geq 2c \eta \,|\xi|
\end{equation*}
for some constant $c = c(n,m,q) > 0$.  Thus for $i_0 = i$ or $i_0 = j$ 
\begin{equation}\label{branchpt nonplanar eqn4}
	|\pi_{P_{i_0}^{\perp}} (\xi,0)| \geq c \eta \,|\xi| .
\end{equation}
Let $\rho_0 = \rho_0(n,m,q,p_0) \in (0,1/8]$ and $\tau = \tau(n,m,q,p_0,\rho_0) \in (0,\rho_0/2]$ be constants to be later determined.  Provided $\tau$ is sufficiently small, there exists a $\|T_{i_0}\|$-measurable set $S \subset \op{spt} T_{i_0} \cap \mathbf{B}_{\rho_0}(Z) \cap \{|x| > \tau\}$ such that $\|T_{i_0}\|(S) \geq \tfrac{1}{2} \,\omega_n \rho_0^n$.  Integrating \eqref{branchpt nonplanar eqn4} over $S$ and arguing as in~\cite[Lemma~3.9]{Sim93} using \eqref{branchpt close2plane eqn1}, \eqref{branchpt nonplanar eqn2}, and \eqref{branchpt nonplanar eqn3} prove (a) and (b).  Arguing as in ~\cite[Theorem~3.1]{Sim93} again using \eqref{branchpt close2plane eqn1}, \eqref{branchpt nonplanar eqn1}, and \eqref{branchpt nonplanar eqn3} proves (c) and (d).
\end{proof}

\begin{corollary} \label{noncencentration nonplanar cor}
For all $\eta_0 \in (0,1)$, $\delta \in (0,1/4)$, and $\sigma \in (0,1)$ there exists $\varepsilon_0,\beta_0 \in (0,1)$ depending only on $n,m,q,\eta_0,\delta,\sigma$ such that if $\mathbf{C}$ and $T$ satisfy Hypothesis~$(\star)$, Hypothesis~$(\star\star)$, and Hypothesis~$(\dagger\dagger)$ and if 
\begin{equation*}
	\mathbf{B}_{\delta}(0,y_0) \cap \{ X \in \mathbf{B}_{3/4}(0) : \Theta(T,X) \geq q \} \neq \emptyset 
\end{equation*}
for all $y_0 \in B^{n-2}_{1/2}(0)$, then 
\begin{equation} \label{noncencentration nonplanar concl} 
	\int_{\mathbf{B}_{1/2}(0)} \frac{\op{dist}^2(X, \op{spt} \mathbf{C})}{r_{\delta}^{2-\sigma}} \,d\|T\|(X) 
		\leq C \int_{\mathbf{B}_1(0)} \op{dist}^2(X, \op{spt} \mathbf{C}) \,d\|T\|(X) , 
\end{equation}
where $r_{\delta} = \max\{|x|,\delta\}$ and $C = C(n,m,q,\sigma) \in (0,\infty)$ is a constant (independent of $\delta$). 
\end{corollary}

\begin{proof}
The estimate \eqref{noncencentration nonplanar concl} follows from Theorem~\ref{branchpt nonplanar thm} by arguing as in the proof of~\cite[Corollary~3.2]{Sim93}.
\end{proof}

\section{Fine blow-ups and decay of fine excess}\label{blow-up-analysis}
\setcounter{equation}{0}

\subsection{Preliminaries, fine blow up class and notation}\label{fine-blowup-prelim}

Let $q$ be an integer  $\geq 2$, $M \geq 1$ and let $(\e_{k})$, $(\b_{k}),$ $(\eta_{k})$ be sequences of positive numbers converging to $0$. For each $k=1, 2, \ldots$, let $T_{k}$ be an $n$-dimensional locally area minimizing rectifiable current   in ${\mathbf B}_{1}(0)$ with $\left(\partial \, T_{k} \right)\llcorner {\mathbf B}_{1}(0) = 0$ and let ${\mathbf C}_{k} \in {\mathcal C}_{q, p_{k}}$, where $p_{k}$ is an integer with $2 \leq p_{k} \leq q$, such that Hypothesis~($\star$), Hypothesis~($\star\star$), Hypothesis~($\dag$) and conditions (\ref{close2plane convention1}) and (\ref{close2plane convention2}) hold true with $T_{k}$, ${\mathbf C}_{k}$, $p_{k}$, $\e_{k}$, $\b_{k}, \eta_{k}$  in place of $T$, ${\mathbf C}$, $p$, $\e$, $\b$, $\eta$ respectively. Thus, for each $k=1, 2, \ldots,$ we suppose: 
\begin{itemize}
\item[$(1)$] ${\mathbf C}_{k} = \sum_{i=1}^{p_{k}} q_{i}^{(k)}\llbracket P_{i}^{(k)}\rrbracket$ where for each $k$ and $1 \leq i \leq p_{k},$ $P_{i}^{(k)}$  are distinct $n$-dimensional oriented planes with $\{0\} \times {\mathbb R}^{n-2} \subset P_{i}^{(k)}$ and orienting $n$-vector 
$\vec{P}_{i}^{(k)}$, and $q_{i}^{(k)},$ $i \in \{1, 2, \ldots, p_{k}\}$, are positive integers with $\sum_{i=1}^{p_{k}} q_{i}^{(k)} = q$; 
\item[$(2)$]  $(\partial T_{k}) \res {\mathbf B}_{1}(0) = 0$, $\Theta \, (T_{k}, 0) \geq q$, \; $\|T_{k}\|({\mathbf B}_{1}(0)) < (q + 1/2)\omega_{n};$  
\item[$(3)$] $0 < E_{k}= E(T_{k}, {\mathbf C}_{k}, {\mathbf B}_{1}(0)) < \e_{k}$; 
\item[$(4)$] $Q(T_{k}, {\mathbf C}_{k}, {\mathbf B}_{1}(0)) < \b_{k} \, \inf_{{\mathbf C}^{\prime} \in \cup_{p^{\prime} = 1}^{p_{k} - 1} {\mathcal C}_{q, p^{\prime}}} \, Q(T_{k}, {\mathbf C}^{\prime}, {\mathbf B}_{1}(0));$
\item[$(5)$] $E(T_{k}, {\mathbf P}_{0}, {\mathbf B}_{1}(0)) < \eta_{k};$
\item[$(6)$] $E(T_{k}, {\mathbf P}_{0}, {\mathbf B}_{1}(0)) \leq M \inf_{{\mathbf P} \in {\mathcal C}_{q, 1}} \, E(T_{k}, {\mathbf P}, {\mathbf B}_{1}(0));$
\item[$(7)$] $\pi_{P_0\#} (T_{k} \llcorner \mathbf{B}_{7/8}) \llcorner \mathbf{B}_{3/4}(0) = q \llbracket P_0 \rrbracket \llcorner \mathbf{B}_{3/4}(0)$;
\item[$(8)$] $|\vec P_i^{(k)} - \vec P_0| < 1/2$. 
\end{itemize}

Suppose further that  for each $k=1, 2, \ldots$, there is $\delta_{k}>0$ with $\delta_{k} \to 0$ such that 
\begin{equation}\label{excess-1}
{\mathbf B}_{\d_{k}}(0,z) \cap \{Z \, : \, \Theta \, (T_{k}, Z) \geq q\} \neq \emptyset
\end{equation}
for each point $(0, z) \in \{0\} \times {\mathbb R}^{n-2} \cap B_{1/2}(0).$ 

Let $(\t_{k})$ be a sequences of decreasing positive numbers converging to 0, and let $\gamma \in (0, 1)$. By passing to appropriate 
subsequences of $(T_{k})$, $({\mathbf C}_{k})$, and possibly replacing ${\mathbf C}_{k}$ with a cone 
${\mathbf C}_{k}^{\prime}$ with ${\rm spt} \, \|{\mathbf C}_{k}^{\prime}\| = {\rm spt} \, \|{\mathbf C}_{k}\|$ without changing notation (see Remark~\ref{cone-multiplicity}), we find an integer $p \in \{2, \ldots, q\}$, and for each $i \in \{1, 2, \ldots, p\}$, integers $q_{i} \geq 1$ with $\sum_{i=1}^{p} q_{i} = q$ such that (1)-(8) above hold with 
$p_{k} = p$ and $q_{i}^{(k)} = q_{i}$ for each $k =1, 2, 3, \ldots$; furthermore, we find that the assertions $(A)$-$(E)$ below hold for each $k=1, 2, \ldots$:
 \begin{itemize}
 \item[(A)] Writing $\widehat{E}_{k}  = E(T_{k}, {\mathbf P}_{0}, {\mathbf B}_{1}(0))$ and $P_{i}^{(k)} = \{(z, x, y) \in {\mathbb R}^{n+m}\, : \, z = A_{i}^{(k)}x\}$ for some $m \times 2$ matrices $A_{i}^{(k)}$, we have by (\ref{close2plane eqn5}), (\ref{close2plane eqn6}), (\ref{close2plane eqn7}) and (\ref{close2plane eqn12}) 
 that  
  
 \begin{gather}
 \label{close2plane eqn1k} \|A_{i}^{(k)}\| \leq C \widehat{E}_{k}, \\
	\label{close2plane eqn2k} \min_{1 \leq i < j \leq p} \,\inf_{x \in \mathbb{S}^1} |A_i^{(k)} x - A_j^{(k)} x| 
		\geq c \inf_{\mathbf{C}' \in \bigcup_{p'=1}^{p-1} \mathcal{C}_{p^{\prime}, }} Q(T_{k}, \mathbf{C}', {\mathbf B}_{1}(0)) , \\
		\label{close2plane eqn3k} 
	\max_{1 \leq i < j \leq p} \inf_{x \in \mathbb{S}^1} |A_i ^{(k)}x - A_j^{(k)} x| \geq \frac{c}{M} \widehat{E}_{k}, \\
	\label{close2plane eqn4k} \|A_i^{(k)} - A_j^{(k)}\| \leq C \inf_{x \in \mathbb{S}^1} |A_i^{(k)} x - A_j^{(k)} x| \text{ for all $1 \leq i < j \leq p$,} 
\end{gather}
where $c = c(n,m,q)$ and $C = C(n,m,q)$.\\ 
 \item[(B)] By Corollary~\ref{nonconcentration close2plane cor},
\begin{equation}\label{excess-2}
\int_{{\mathbf B} _{1/2}(0) \cap \{|(x^{1}, x^{2})| < \delta\}} {\rm dist}^{2} \, (X, {\rm spt} \, {\mathbf C}_{k}) \, d\|T_{k}\|(X) 
\leq C\delta^{3/2} E_{k}^{2}
\end{equation}
for each $\delta \in [\d_{k}, 1/4),$ where $C= C(n,m, q, M) \in (0, \infty).$\\
\item[(C)] By Theorem~\ref{graphrep close2plane thm},  for $1 \leq i \leq p$, there exists an $n$-dimensional locally area minimizing rectifiable current $T_i^{(k)}$ of $\mathbf{C}_{(3+\gamma)/4}(0) \cap \{r > \tau_{k}/4\}$  such that: 
\begin{enumerate} \setlength{\itemsep}{5pt}
	\item[(a)] 
	\begin{gather}\label{graphrep close2plane concl b} 
		T_{k} \llcorner \mathbf{B}_{(7+\gamma)/8}(0) \cap \mathbf{C}_{(3+\gamma)/4}(0) \cap \{r > \tau_{k}/4\} = \sum_{i=1}^{p} T_i^{(k)} , \\
		(\partial T_i^{(k)}) \llcorner \mathbf{C}_{(3+\gamma)/4}(0) \cap \{r > \tau_{k}/4\} = 0, \nonumber \\ 
		(\pi_{P_0 \#} T_i^{(k)}) \llcorner \mathbf{C}_{(1+\gamma)/2}(0) \cap \{r > \tau_{k}/2\} 
			= q_i \llbracket P_0 \rrbracket \llcorner \mathbf{C}_{(1+\gamma)/2}(0) \cap \{r > \tau_{k}/2\} , \nonumber \\ 
		\sup_{X \in \op{spt} T_i^{(k)} \cap \{r > \sigma\}} \op{dist}(X, P_i^{(k)}) \leq C_{\sigma} E_{k} \text{ for all } \sigma \in [\tau_{k}/2,1/2] , \nonumber 
	\end{gather}
	where $r(X) = \op{dist}(X, \{0\} \times \mathbb{R}^{n-2})$ and $C_{\sigma} = C_{\sigma}(n,m,q,\gamma,\sigma) \in (0,\infty)$ is a constant; 
	
	\item[(b)] for each $i \in \{1,2,\ldots,p\}$ there exists a Lipschitz $q_i$-valued function $u_i^{(k)} : B_{\gamma}(0) \cap \{r > \tau_{k}\} \rightarrow \mathcal{A}_{q_i}({\mathbb R}^{m})$ and a
	measurable set $K_{i}^{(k)} \subseteq B_{\gamma}(0) \cap \{r > \tau_{k}\}$ such that 
	\begin{gather}\label{graphrep close2plane concl c} 
		T_i^{(k)} \llcorner ({\mathbb R}^{m} \times K_{i}^{(k)}) = \op{graph} \, (A_{i}^{(k)}x + u_i^{(k)}) \cap ({\mathbb R}^{m} \times K_{i}^{(k)}) , \\
		\mathcal{L}^n(B_{\gamma}(0)  \cap \{r > \sigma\} \setminus K_{i}^{(k)})
			+ \|T_i^{(k)}\|({\mathbb R}^{m} \times (B_{\gamma}(0) \cap \{r > \sigma\} \setminus K_{i}^{(k)})) \leq C_{\sigma} E_{k}^{2+\alpha} , \nonumber \\
		\sup_{B_{\gamma}(0) \cap \{r > \sigma\}} |u_i^{(k)}| \leq C_{\sigma} E_{k} , \quad 
			\sup_{B_{\gamma}(0) \cap \{r > \sigma\}} |\nabla \,  u_i^{(k)}| \leq C_{\sigma} E_{k}^{\alpha} \nonumber
	\end{gather}
	for all $\sigma \in [\tau_{k},1/2]$, where $\alpha = \alpha(n,m,q) \in (0,1)$, $C_{\sigma} = C_{\sigma}(n,m,q,\gamma,\sigma) \in (0,\infty)$ are constants.\\
\end{enumerate}
\item[(D)]  By Theorem~\ref{branchpt close2plane thm},  for each point $Z = (\chi,\xi,\zeta) \in \op{spt} T_{k} \cap \mathbf{B}_{1/2}(0)$ with $\Theta(T_{k},Z) \geq q,$
\begin{equation}\label{goodpt-est} 
\quad |\chi|^2 + \widehat{E}_{k}^2 |\xi|^2 \leq C E_{k}^{2}
\end{equation}
\end{itemize} 
where $C = C(n, m, q, M) \in (0, \infty)$. Furthermore, by Theorem~\ref{keyest thm}, Corollary~\ref{keyest cor2} and Theorem~\ref{branchpt close2plane thm}, for each $\rho \in (0, 1/2)$ and sufficiently large $k$ depending on $\rho$, and for each point $Z = (\chi, \xi, \zeta) \in 
\op{spt} T_{k} \cap {\mathbf B}_{1/2}(0)$ with $\Theta (T_{k}, Z) \geq q$, 
\begin{itemize}
\item [(E)] 
\begin{eqnarray}\label{estimate-E1}
&&\int_{B_{\gamma\rho}(0) \cap \{r > \t_{k}\}} R^{2-n}\left(\frac{\partial(u_{i}^{(k)}/R)}{\partial R}\right)^{2}  + \rho^{-n}\int_{B_{\gamma\rho}(0) \cap \{r > \t_{k}\}} |D_{y}u_{i}^{(k)}|^{2} \\
&&\hspace{2in}		\leq C \rho^{-n-2} \int_{{\mathbf B}_{\rho}(0)} {\rm dist}^{2} \, (X, {\rm spt} \, \|{\mathbf C}_{k}\|) \, d\|T_{k}\| \nonumber
\end{eqnarray}
\begin{eqnarray}\label{estimate-E2}
&&\int_{B_{\gamma\rho}(\xi, \z) \cap \{r > \tau_{k}\}} \frac{|u_i^{(k)}(x,y) - (\chi - A_i^{(k)} \xi)|^{2}}{|(u_{i}^{(k)}(x, y) - \chi, x - \xi,y - \zeta)|^{n+2-\sigma}} \\
&&\hspace{1.5in}		\leq C_{1}  \rho^{-n-2 + \sigma} \int_{{\mathbb R}^{m} \times B_{\rho}(\xi, \z)} {\rm dist}^{2} \, (X, {\rm spt} \, \|\nu_{Z \, \#} {\mathbf C}_{k}\|) \, d\|T_{k}\| \nonumber
\end{eqnarray}
for $ 1 \leq i \leq p$, where $R(X) = |X|$; $\nu_{Z}$ is the translation $X \mapsto X + Z$; $C = C(n, m, q, \gamma) \in (0, \infty)$, $C_{1} = C_{1}(n, m, q, M, \gamma, \sigma) \in (0, \infty)$ are constants; 
and, writing $$u_{i}^{(k)}(x, y) = \sum_{\kappa=1}^{q_{i}} \llbracket u_{i, \kappa}^{(k)}(x, y) \rrbracket$$ with $u^{(k)}_{i, \kappa}(x, y)\in {\mathbb R}^{m}$,  the integrand on the left hand side of \eqref{estimate-E2} is the single-valued function 
$(x, y) \mapsto \sum_{\kappa=1}^{q_{i}} \frac{|u_{i, \kappa}^{(k)}(x,y) - (\chi - A_i^{(k)} \xi)|^{2}}{|(u_{i, \kappa}^{(k)}(x, y) - \chi, x - \xi,y - \zeta)|^{n+2-\sigma}}$.  
\end{itemize}

Extend $u_{i}^{(k)}$ to all of $B_{\gamma}(0)$ by setting $u_{i}^{(k)}(x,y) = q_{i}\llbracket 0 \rrbracket$ for $(x,y) \in B_{\gamma}(0) \cap \{r \leq \t_{k}\}.$ 

By (\ref{close2plane eqn1k}) and (\ref{close2plane eqn3k}), for each $i \in \{1, 2, \ldots p\}$ there exists an $m \times 2$ matrix $A_{i}$ such  that,  
after passing to a subsequence of $(k)$ without changing notation, 
\begin{equation}\label{excess-7}
\widehat{E}_{k}^{-1} A_{i}^{(k)}  \to A_{i} 
\end{equation}
for each $i=1, 2, \ldots, p,$ and 
\begin{equation}\label{matrix-upperbound}
\|A_{i}\| \leq C; 
\end{equation}
\begin{equation}\label{rank2}
\max_{1 \leq i < j \leq p} \, \inf_{x \in {\mathbb S}^{1}} \, |A_{i}x - A_{j}x| \geq  c/M
\end{equation} 
where $C = C(n, m ,q) \in (0, \infty)$, $c = c(n, m, q) \in (0, \infty);$
hence in particular there exist $i, j \in \{1, 2, \ldots, p\}$ with $i \neq j$ such that $A_{i} - A_{j}$ has rank 2.  

By (C) above, Theorem~\ref{harmonic thm} and Lemma~\ref{energy est lemma}, there exist locally Dirichlet energy minimizing functions $w_{i} \, : \, B_{\gamma}(0) \setminus \{0\} \times {\mathbb R}^{n-2}  
\to {\mathcal A}_{q_{i}}({\mathbb R}^{m})$  such that, writing $v^{(k)}_{i} = E_{k}^{-1}u^{(k)}_{i}$ and passing to a subsequence of $(k)$ without changing notation, we have 
\begin{equation}\label{excess-8}
\int_{K} {\mathcal G}(v^{(k)}_{i}, w_{i})^{2} + (|Dv^{(k)}_{i}| - |Dw_{i}|)^{2} \to 0 \;\; \mbox{as $k \to \infty$}
\end{equation}
for each compact set $K \subset B_{\gamma}(0) \setminus \{0\} \times {\mathbb R}^{n- 2}.$ From (\ref{excess-2}) it follows that 
$\int_{{\mathbf B}_{1/2}(0) \cap \{0 < r < \s\}} |u_{i}^{(k)}|^{2} \leq C \s^{1/2}E_{k}^{2}$ and consequently that $\int_{B_{1/2}(0) \cap \{0 < r < \s\}}|w_{i}|^{2} \leq C\s^{1/2}$
for each $\s \in (0, 1/4)$,  where $C = C(n, m, q) \in (0, \infty)$, and hence that $w_{i} \in L^{2}(B_{1/2}(0); {\mathcal A}_{q_{i}}({\mathbb R}^{m}))$ and 
that $$v^{(k)}_{i} \to w_{i} \;\;\; {\rm in}  \;\;L^{2}\,(B_{1/2}(0); {\mathcal A}_{q_{i}}({\mathbb R}^{m})).$$ Taking a sequence $\gamma_{\ell} \to 1^{-}$ and applying this procedure with $\gamma_{\ell}$ in place of $\gamma$, in conjunction with a diagonal subsequence argument to select a further subsequence of $(k)$ without relabeling, we obtain, for each $i=1, 2, \ldots, p$, a locally Dirichlet energy minimizing function $w_{i} \, : \, B_{1}(0) \setminus \{0\} \times {\mathbb R}^{n-2} \to {\mathcal A}_{q_{i}}({\mathbb R}^{m})$ such that 
\begin{equation}\label{excess-8-1}
\int_{B_{1/2}(0)} {\mathcal G}(v^{(k)}_{i}, w_{i})^{2} + \int_{K} (|Dv^{(k)}_{i}| - |Dw_{i}|)^{2} \to 0 \;\; \mbox{as $k \to \infty$}
\end{equation}
for each compact set $K \subset B_{1}(0) \setminus \{0\} \times {\mathbb R}^{n-2}$. Furthermore, we have by (C) and \eqref{excess-2} that for each $\rho \in (0, 1/2)$ and each $\sigma \in (0, \rho/8)$, 
\begin{align*}
&\int_{B_{\rho - \sigma}(0)} \sum_{i=1}^{p}|w_{i}|^{2}  = \lim_{k \to \infty} E_{k}^{-2} \int_{{\mathbf C}_{\rho - \sigma}(0)} {\rm dist}^{2}(X, {\rm spt} \, {\mathbf C}_{k}) \, d\|T_{k}\| \\
&\leq \liminf_{k \to \infty} E_{k}^{-2} \int_{{\mathbf B}_{\rho}(0)} {\rm dist}^{2}(X, {\rm spt} \, {\mathbf C}_{k})  \, d\|T_{k}\| \leq \limsup_{k \to \infty} E_{k}^{-2} \int_{{\mathbf B}_{\rho}(0)} {\rm dist}^{2}(X, {\rm spt} \, {\mathbf C}_{k})  \, d\|T_{k}\|\\ 
&\leq \lim_{k \to \infty} E_{k}^{-2} \int_{{\mathbf C}_{\rho}(0)} {\rm dist}^{2}(X, {\rm spt} \, {\mathbf C}_{k})  \, d\|T_{k}\| = \int_{B_{\rho}(0)} \sum_{i=1}^{p} |w_{i}|^{2},
\end{align*} 
and thus 
\begin{equation}\label{blowup_norm_conv_eqn2_0}
\lim_{k \to \infty} E_{k}^{-2} \int_{{\mathbf B}_{\rho}(0)} {\rm dist}^{2}(X, {\rm spt} \, {\mathbf C}_{k})  \, d\|T_{k}\| = \int_{B_{\rho}}\sum_{i=1}^{p} |w_{i}|^{2};
\end{equation}
hence 
in particular $\omega_{n}^{-1}\int_{B_{1/2}(0)} |w_{i}|^{2} \leq 1$, and we have that 
\begin{eqnarray} \label{blowup_norm_conv_eqn2-1}
&&w_{i}  \in L^{2}(B_{1/2}(0); {\mathcal A}_{q_{i}}({\mathbb R}^{m})) \cap W^{1, 2}_{\rm loc}(B_{1}(0) \setminus \{0\} \times {\mathbb R}^{n-2}; {\mathcal A}_{q_{i}}({\mathbb R}^{m})), \;\; i=1, 2, \ldots, p, \;\; \mbox{and}\\
&&\omega_{n}^{-1}\int_{B_{1/2}(0)} \sum_{i=1}^{p} |w_{i}|^{2} \leq 1.\nonumber
\end{eqnarray}

\noindent
\begin{definition}[fine blow-up]\label{blowupclass_defn} 
Fix integers $q \geq 2$, $p \in \{2, \ldots, q\}$ and $q_{1}, q_{2}, \ldots, q_{p} \geq 1$ such that $\sum_{j=1}^{p} q_{j} = q,$ and write ${\mathbf q} = (q_{1}, q_{2}, \ldots, q_{p}).$ Let $w= (w_{1}, w_{2}, \ldots, w_{p}),$ where for each $i \in \{1, 2, \ldots, p\}$: 
\begin{itemize}
\item[(a)] $w_{i}  \in L^{2}(B_{1/2}(0); {\mathcal A}_{q_{i}}({\mathbb R}^{m})) \cap W^{1, 2}_{\rm loc}(B_{1}(0) \setminus \{0\} \times {\mathbb R}^{n-2}; {\mathcal A}_{q_{i}}({\mathbb R}^{m}));$ 
\item[(b)] $w_{i}$ is as in (\ref{excess-8-1}), so that $w_{i}$ arises in the manner described above, corresponding to (a subsequence of)  a  sequence $(T_{k})$ of $n$-dimensional locally area minimizing rectifiable currents  in ${\mathbf B}_{1}(0)$, a sequence of cones $({\mathbf C}_{k})$ with 
${\mathbf C}_{k} \in {\mathcal C}_{p, {\mathbf q}},$ 
sequences of positive numbers $(\e_{k}),$ $(\b_{k}),$ $(\eta_{k}),$ $(\delta_{k})$ converging to 0, satisfying,  for some $M \geq 1$ (independent of $k$), every $\gamma \in (0, 1)$ and sufficiently large $k$ (depending on $\gamma$):  
\begin{itemize} 
\item[(i)] Hypothesis~($\star$), Hypothesis~($\star\star$), Hypothesis~($\dag$) and conditions (\ref{close2plane convention1}) and (\ref{close2plane convention2}) with $T_{k}$, ${\mathbf C}_{k}$, $\e_{k}$, $\b_{k}, \eta_{k}$  in place of $T$, ${\mathbf C}$, $\e$, $\b$, $\eta$ respectively (or equivalently, after possibly reversing orientation of $T_{k}$ and of ${\mathbf C}_{k}$, conditions (1)-(8) above with $p_{k} = p$ and $q_{i}^{(k)} = q_{i}$ for all $k=1, 2, \ldots$ and $i = 1, \ldots, p$), and
\item[(ii)] condition (\ref{excess-1}). 
\end{itemize}
\end{itemize}
{\rm We call $w$ a  {\em fine blow-up} of the sequence $(T_{k})$ relative to $({\mathbf C}_{k})$.} 
\end{definition}

\begin{definition}[fine blow-up classes]\label{blowupclass_defn_1} With $q$, $p$, ${\mathbf q}$ as in Definition~\ref{blowupclass_defn} above, and for $M \geq 1$, we let ${\mathfrak B}_{p, {\mathbf q}}(M)$ denote the set of all functions $w  = (w_{1}, \ldots, w_{p})$ with $w_{i} \in  L^{2}(B_{1/2}(0); {\mathcal A}_{q_{i}}({\mathbb R}^{m})) \cap W^{1, 2}_{\rm loc}(B_{1}(0) \setminus \{0\} \times {\mathbb R}^{n-2}; {\mathcal A}_{q_{i}}({\mathbb R}^{m}))$ such that $w$ is a fine blow-up, per Definition~\ref{blowupclass_defn},  of a sequences of currents $(T_{k})$ relative to a sequences of cones $({\mathbf C}_{k})$, where $M$ corresponds to $(T_{k})$ as in condition (6) above.

Given $m \times 2$ matrices $A_{1}, \ldots, A_{p}$ and $M \geq 1,$ we let ${\mathfrak B}_{p, {\mathbf q}}(M, A_{1}, \ldots, A_{p})$ denote the set of all 
$w \in {\mathfrak B}_{p, {\mathbf q}}(M)$ such that $w$ is a fine blow-up, per Definition~\ref{blowupclass_defn},  of  a sequences of currents $(T_{k})$ relative to a sequence of cones 
$({\mathbf C}_{k})$ with ${\mathbf C}_{k} \in {\mathcal C}_{p, {\mathbf q}},$  where: 
\begin{itemize}
\item $M$ corresponds to $(T_{k})$ as in condition (6) above, and 
\item the $m \times 2$ matrices $A_{1}^{(k)}, \ldots, A_{p}^{(k)}$ corresponding to ${\mathbf C}_{k}$ (as in (A) above)
satisfy $${\hat E}_{k}^{-1} A_{i}^{(k)} \to A_{i}$$ as $k \to \infty$, for each $i \in \{1, 2, \ldots, p\}$, where ${\hat E}_{k} = E(T_{k}, {\mathbf P}_{0}, {\mathbf B}_{1}(0))$. 
\end{itemize} 
\end{definition} 

\begin{remark} \label{conditions-on-Aj}
By the discussion above, it follows that ${\mathfrak B}_{p, {\mathbf q}}(M, A_{1}, \ldots, A_{p})$ is defined only for matrices $A_{1}, \ldots, A_{p}$ satisfying 
\eqref{matrix-upperbound} and \eqref{rank2}.   
\end{remark}

\begin{remark} \label{rescaling-fine-blowup}
Let $q$, $p$, ${\mathbf q}$ be as in Definition~\ref{blowupclass_defn}, let $M \geq 1$  and let $w = (w_{1}, \ldots, w_{p}) \in {\mathfrak B}_{p, {\mathbf q}}(M, A_{1}, \ldots, A_{p})$ where  $A_{1}, \ldots, A_{p}$ are constant $m \times 2$ matrices satisfying the requirements of Definition~\ref{blowupclass_defn_1} for some sequences  
$(T_{k})$, $({\mathbf C}_{k})$ of currents and cones giving rise to $w$ per Definition~\ref{blowupclass_defn_1}. If $\rho \in (0, 1)$ and 
$w_{i} \not\equiv q_{i}\llbracket 0 \rrbracket$ in $B_{\rho}(0)$ for some $i \in \{1, \ldots, p\}$, then $$\widetilde{w} \equiv 
\|w(\rho(\cdot))\|_{L^{2}(B_{1}(0))}^{-1}w(\rho(\cdot)) \in {\mathfrak B}_{p, {\mathbf q}}(CM, A_{1}, \ldots, A_{p}),$$ where $C= C(n, m, q)$, $\|w(\rho(\cdot))\|^{2}_{L^{2}(B_{1}(0))} = 
\rho^{-n}\int_{B_{\rho}(0)} \sum_{j=1}^{p} |w_{j}|^{2}$ and $\widetilde{w} = (\widetilde{w}_{1}, \widetilde{w}_{2}, \ldots, \widetilde{w}_{p})$ with 
$\widetilde{w}_{j} = \|w(\rho(\cdot))\|^{-1}_{L^{2}(B_{1}(0))} w_{j}(\rho(\cdot))$; indeed, $E(T_{k}, {\mathbf C}_{k}, {\mathbf B}_{\rho}(0)) \neq 0$ for all sufficiently large $k$, and $\widetilde{w}$ is a fine blow up of the sequence 
$\widetilde{T}_{k} = \eta_{0, \rho \, \#} \, T_{k}$ relative to ${\mathbf C}_{k}.$ {\rm To see this, note that by \eqref{blowup_norm_conv_eqn2_0} we have that 
$E(T_{k}, {\mathbf C}_{k}, {\mathbf B}_{\rho}(0)) \neq 0$ for all sufficiently large $k$. By assumption conditions (1)-(8) above (with $p_{k}=p$) and condition ~\eqref{excess-1} hold, and hence it is readily checked that conditions (1)-(3), (5), (7), (8) above and condition ~\eqref{excess-1} are satisfied with $\widetilde{T}_{k}$ in place of $T_{k}$ and with $\rho^{-(n+2)/2}\varepsilon_{k}$, 
$\rho^{-(n+2)/2}\eta_{k}$, $\rho^{-1}\delta_{k}$ in place of $\varepsilon_{k}$, $\eta_{k}$, $\delta_{k}$; by arguing as in \cite[pp.\ 910-914]{Wic14}, using Theorem~\ref{graphrep close2plane thm} in places where the argument depended on \cite[Theorem~10.1]{Wic14}, we can also verify that condition (4) holds with $\widetilde{T}_{k}$ in place of $T_{k}$ and $C\rho^{-(n+2)/2}\beta_{k}$ in place of $\beta_{k}$ where $C = C(n, m, q) \in (0, \infty)$, and that condition (6) holds with $\widetilde{T}_{k}$ in place of $T_{k}$ and 
$CM$ in place of $M$, where again $C = C(n, m, q) \in (0, \infty)$. Hence we can construct a fine blow up of $(\widetilde{T}_{k})$ relative to ${\mathbf C}_{k}$, 
which can readily be checked to be equal to $\|w(\rho(\cdot))\|_{L^{2}(B_{1}(0))}^{-1}w(\rho(\cdot))$.} 
\end{remark}

In subsequent sections, we shall use the following notation:

\begin{itemize}
\item ${\mathcal S}$ denotes the subspace of $(m+n) \times (m+n)$ skew symmetric matrices spanned by the set of skew symmetric matrices corresponding to the transformations $\Gamma_{ij} \, : \, {\mathbb R}^{m+n} \to {\mathbb R}^{m+n}$ ($1 \leq i \leq 2+m$,  $1 \leq j \leq n-2$) given by 
$\Gamma_{ij}(x, y) = x^{i}e_{2+m+j} - y^{j}e_{i},$ $(x,y) \in {\mathbb R}^{2+m}  \times {\mathbb R}^{n-2}$. 
\end{itemize}

In the following, we fix an integer $q \geq 2.$ 
\begin{itemize}
\item For $p \in \{2, \ldots, q\}$, let $$\mathfrak{M}_{p} = \{(q_{1}, q_{2}, \ldots, q_{p}) \, : \, q_{1}, q_{2}, \ldots, q_{p} \;\; \mbox{are integers $\geq 1,$ with} \;\sum_{j=1}^{p} q_{j} = q\}.$$ 

\item For $p \in \{2, \ldots, q\}$ and ${\mathbf q} \in {\mathfrak M}_{p}$, let $$\widetilde{\mathcal C}_{p, {\mathbf q}} = \{(e^{A})_{\#} \, {\mathbf C} \, : \, {\mathbf C} \in {\mathcal C}_{p, {\mathbf q}}, \;A \in {\mathcal S}\}$$ and let 
$$\widetilde{\mathcal C} = \cup_{p=2}^{q} \cup_{{\mathbf q} \in {\mathfrak M}_{p}} \widetilde{\mathcal C}_{p, {\mathbf q}}.$$

\item For $p \in \{2, \ldots, q\}$ and ${\mathbf q}  = (q_{1}, q_{2}, \ldots, q_{p}) \in {\mathfrak M}_{p}$, let ${\mathfrak L}_{p, {\mathbf q}}$  be the space of functions $\psi = (\psi_{1}, \psi_{2}, \ldots, \psi_{p})$ such that 
for each $j \in \{1, 2, \ldots, p\}$, $\psi_{j} \, : \, {\mathbb R}^{n} \to {\mathcal A}_{q_{j}}({\mathbb R}^{m})$ is a (linear) function of the form 
$$\psi_{j}(x,y) = \sum_{\ell=1}^{q_{j}} \llbracket \psi_{j, \ell}(x) \rrbracket \hspace{.4in} \forall (x, y) \in {\mathbb R}^{2} \times {\mathbb R}^{n-2},$$ 
where for $\ell=1, \ldots, q_{j}$, 
$\psi_{j, \ell} \, : \, {\mathbb R}^{2} \to {\mathbb R}^{m}$  are linear functions such that for $k \neq \ell$, either $\psi_{j, k}(x) \equiv \psi_{j, \ell}(x)$ for all $x \in {\mathbb R}^{2}$ or $\psi_{j, k}(x) \neq \psi_{j, \ell}(x)$ for all $x \in {\mathbb R}^{2} \setminus \{0\}$.  

\item Let ${\mathfrak L} = \cup_{p=2}^{q} \cup_{{\mathbf q} \in {\mathfrak M}_{p}} {\mathfrak L}_{p, {\mathbf q}}.$

\item For $p \in \{2, \ldots, q\}$ and ${\mathbf q}  = (q_{1}, q_{2}, \ldots, q_{p}) \in {\mathfrak M}_{p}$,  let $\widetilde{\mathfrak L}_{p, {\mathbf q}}$ be the space of functions 
$\widetilde{\psi} = (\widetilde{\psi}_{1}, \widetilde{\psi}_{2}, \ldots, \widetilde{\psi}_{p})$ such that: 
\begin{itemize}
\item[(a)] $\widetilde{\psi}_{j} \, : \, {\mathbb R}^{n} \to {\mathcal A}_{q_{j}}({\mathbb R}^{m})$ for each $j \in \{1, 2, \ldots, p\};$
\item[(b)] there are $\psi = (\psi_{1}, \psi_{2}, \ldots, \psi_{p}) \in {\mathfrak L}_{p, {\mathbf q}}$, two (single-valued) linear functions $L_{1} \,: \, {\mathbb R}^{n-2} \to {\mathbb R}^{m}$ and $L_{2} \, : \, {\mathbb R}^{n-2} \to {\mathbb R}^{2},$ and constant $m \times 2$ matrices $A_{1}, A_{2}, \ldots, A_{p}$ such that for each $j \in \{1,2, \ldots, p\}$,  
$\widetilde{\psi}_{j}(x, y) = \sum_{\ell=1}^{q_{j}} \llbracket \psi_{j, \ell}(x) +  L_{1}(y) + A_{j} L_{2}(y)\rrbracket$ 
$\forall (x, y) \in {\mathbb R}^{2} \times {\mathbb R}^{n-2}.$ 
\end{itemize}
\item For $p \in \{2, \ldots, q\}$, ${\mathbf q}  = (q_{1}, q_{2}, \ldots, q_{p}) \in {\mathfrak M}_{p}$ and $m \times 2$ matrices $\overline{A}_{1}, \ldots, \overline{A}_{p}$, let 
$$\widetilde{\mathfrak L}_{p, {\mathbf q}}(\overline{A}_{1}, \ldots, \overline{A}_{p}) = \{\widetilde{\psi} \in \widetilde{\mathfrak L}_{p, {\mathbf q}} \, : \, \mbox{(b) above holds with $A_{j} = \overline{A}_{j}$ for $j=1, \ldots, p$}\}.$$ 
\item Let $\widetilde{\mathfrak L} =  \cup_{p=2}^{q} \cup_{{\mathbf q} \in {\mathfrak M}_{p}} \widetilde{\mathfrak L}_{p, {\mathbf q}}.$
\end{itemize}

In the following lemma and subsequently, we let $r(X) = |(x_{1}, x_{2})|$ for $X  = (x_{1}, x_{2}, y) \in {\mathbb R}^{n}$ where $(x_{1}, x_{2}) \in {\mathbb R}^{2}, y \in {\mathbb R}^{n-2}$. The lemma establishes the elementary fact that corresponding to each $\psi \in \widetilde{{\mathfrak L}}_{p, {\mathbf q}}$, each sequence of cones ${\mathbf C}_{k} \in {\mathcal C}_{p, {\mathbf q}}$  converging to $\mathbf{P}_0 = q \llbracket \{0\} \times {\mathbb R}^{n} \rrbracket$ and 
each sequence of numbers $E_{k} \to 0^{+}$,  
there is  a sequence of cones $\widetilde{\mathbf C}_{k} \in \widetilde{{\mathcal C}}_{p^{\prime}, {\mathbf q}^{\prime}}$ (for some $p^{\prime} \in \{p, \ldots, q\}$ and 
${\mathbf q}^{\prime} \in {\mathfrak M}_{p^{\prime}}$) such that $\psi$ is the blow-up of $(\widetilde{\mathbf C}_{k})$ relative to $({\mathbf C}_{k})$. 

\begin{lemma} \label{blowup2L_lemma} 
Let $p \in \{2, \ldots, q\}$, ${\mathbf q} \in {\mathfrak M}_{p}$  and  let $({\mathbf C}_{k})$  be a sequence in ${\mathcal C}_{p, {\mathbf q}}$ with  
$$\e_{k}^{2} = \int_{{\mathbf B}_{1}(0)} {\rm dist}^{2}(X, \{0\} \times {\mathbb R}^{n}) \, d\|{\mathbf C}_{k}\|(X) \to 0.$$ Write ${\mathbf C}_{k} = \sum_{j=1}^{p} q_{j}\llbracket P_{j}^{(k)} \rrbracket$ and $P_{j}^{(k)} = \{(A_{j}^{(k)}x, x, y) \, : \, (x, y) \in {\mathbb R}^{2} \times {\mathbb R}^{n-2}\}$ where $A_{j}^{(k)}$ ($j=1, 2, \ldots, p$) are distinct constant $m \times 2$ matrices, and suppose that  for each $j \in \{1, 2, \ldots, p\}$, 
$\e_{k}^{-1}A_{j}^{(k)} \to A_{j}$ as $k \to \infty$ for some $m \times 2$ matrices $A_{1}, \ldots, A_{p}$.   
Let $\widetilde{\psi}  = (\widetilde{\psi}_{j})_{j=1}^{p} \in \widetilde{\mathfrak{L}}_{p, {\mathbf q}}(A_{1}, \ldots, A_{p})$ and let $\psi = (\psi_{j})_{j=1}^{p} \in \mathfrak{L}_{p, {\mathbf q}}$, $L_{1} \, : \, {\mathbb R}^{n-2} \to {\mathbb R}^{m}$ and $L_{2} \, : \, {\mathbb R}^{n-2} \to {\mathbb R}^{2}$ correspond to $\widetilde{\psi}$ in accordance with the definition of  
$\widetilde{\mathfrak{L}}_{p, {\mathbf q}}(A_{1}, \ldots, A_{p})$ above.  
Let $(E_{k})$ be a sequence of positive numbers with 
$E_{k} \leq \b_{k} \e_{k}$ where $\b_{k} \to 0^{+}$. 
For  some $p^{\prime} \in \{p, \ldots, q\}$, some ${\mathbf q}^{\prime} \in {\mathfrak M}_{p^{\prime}}$ and each $k = 1, 2, 3, \ldots,$ there is a cone $\widetilde{\mathbf C}_{k} \in \widetilde{\mathcal C}_{p^{\prime}, {\mathbf q}^{\prime}}$ with 
$\int_{{\mathbf B}_{1}(0)} {\rm dist}^{2}(X, \{0\} \times {\mathbb R}^{n}) \, d \|\widetilde{\mathbf C}_{k}\|(X)  \leq (\e_{k} + CE_{k})^{2}$ for sufficiently large $k$, where   $C = C(n, q, \widetilde{\psi})$ such that: 
\begin{itemize}
\item[(i)] $\widetilde{\mathbf C}_k = \sum_{j=1}^p \sum_{l=1}^{q_j} \llbracket \widetilde{P}^{(k)}_{j,l} \rrbracket$ where  
 \begin{equation}\label{blowup2L_eqn5}
 \widetilde{P}^{(k)}_{j,l}= \{ (A_j^{(k)} x + E_k \widetilde{\psi}_{j,l}(x) + \mathcal{R}^{(k)}_{j,l}(x,y), x, y) : (x,y) \in \mathbb{R}^n \}
 \end{equation} 
 with 
 \begin{equation}\label{homog_error}
 |\mathcal{R}^{(k)}_{j,l}(x,y)| \leq C (\epsilon_k + \beta_k + |\e_{k}^{-1} A_{j}^{(k)} - A_{j}|) \,E_k \|L_2\| + C E_k^2 \|L_1\| (|\psi_{j,l}| + \|L_1\| + \|L_2\|)
  \end{equation}
 for $(x, y) \in B_{1}(0)$, where $C = C(n, q)$.
\item[(ii)] For $\tau \in (0, 1/4)$ and $k$ sufficiently large, 
\begin{equation}\label{blowup2L_eqn6}
	\int_{({\mathbb R}^{m} \times B_{1/2}(0)) \cap \{r < \tau\}} {\rm dist}^{2}(X, {\rm spt} \, {\mathbf C}_{k}) \, d\|\widetilde{\mathbf C}_{k}\| 
	\leq C \|\psi\|^{2}_{L^{2}(B_{1}(0))}\tau^{2} E_{k}^2
\end{equation}
where $C = C(n, q)$; 
\item[(iii)] \begin{equation} \label{blowup2L_eqn1}
	\int_{B_{1/2}(0)} \sum_{j=1}^p |\widetilde{\psi}_{j}(X)|^2 \,dX 
		= \lim_{k \rightarrow \infty} \frac{1}{E_{k}^2} \int_{{\mathbb R}^{m} \times B_{1/2}(0)} {\rm dist}^{2} \, (X, {\rm spt} \, {\mathbf C}_{k}) \, d\|{\widetilde{\mathbf C}_{k}}\|.
\end{equation}
\end{itemize}
\end{lemma}

\begin{proof} Note that $\e_{k} >0$ since $p \geq 2$ and $A_{1}^{(k)}, \ldots, A_{p}^{(k)}$ are distinct. Let $\widetilde{\mathbf C}_k = \sum_{j=1}^p \sum_{l=1}^{q_j} \llbracket \widetilde{P}^{(k)}_{j,l} \rrbracket$ where  
	\begin{equation*}
		\widetilde{P}^{(k)}_{j,l} = \{ e^{E_k M_1 + E_k M_2/\e_k} (A_j^{(k)} x + E_k \psi_{j,l}(x), x,y) : (x,y) \in \mathbb{R}^n \} 
	\end{equation*}
	and 
	\begin{equation*}
		M_1 = \left[\begin{array}{ccc} 0 & 0 & L_1 \\ 0 & 0 & 0 \\ -L_1^T & 0 & 0 \end{array}\right] , \quad 
		M_2 = \left[\begin{array}{ccc} 0 & 0 & 0 \\ 0 & 0 & -L_2 \\ 0 & L_2^T & 0 \end{array}\right] 
	\end{equation*}
	representing $L_1$ as an $m \times (n-2)$ matrix and $L_2$ as a $2 \times (n-2)$ matrix.  It is clear that $\widetilde{P}^{(k)}_{j,l}$ is close to $\{0\} \times \mathbb{R}^n$ and hence is a graph over $\{0\} \times \mathbb{R}^n$.  $e^{E_k M_2/\e_k}$ is a rotation of the $(x,y)$-coordinates.  Since $E_{k} \leq \b_{k} \e_{k}$,
	we have that  
	\begin{equation*}
		\widetilde{P}^{(k)}_{j,l} = \Big\{ e^{E_k M_1} \Big( (A_j^{(k)} x + E_k \psi_{j,l}(x) + \e_{k}^{-1}E_k A_j^{(k)} L_2(y), x,y) 
			+ O(\beta_k E_k \|L_2\||(x, y)|) \Big) : (x,y) \in \mathbb{R}^n \Big\}.
	\end{equation*}
	$e^{E_k M_1}$ is a rotation of the $(y,z)$-coordinates.  Hence 
	\begin{align*}
		\widetilde{P}^{(k)}_{j,l} =\,& \{ (A_j^{(k)} x + E_k \psi_{j,l}(x) + E_k L_1(y) + \e_{k}^{-1}E_k A_j^{(k)} L_2(y), x, y - E_k L_1^T A_{j}^{(k)}x) 
		 	\\&\hspace{10mm} + O(\beta_k E_k \|L_2\| + E_k^2 \|L_1\| (|\psi_{j,l}| + \|L_1\| + \|L_2\|)|(x, y)|) : (x,y) \in \mathbb{R}^n \} 
		 \\=\,& \{ (A_j^{(k)} x + E_k \psi_{j,l}(x) + E_k L_1(y) + \e_{k}^{-1}E_k A_j^{(k)} L_2(y), x, y) 
		 	\\&\hspace{10mm} + O(\beta_k E_k \|L_2\| + E_k^2 \|L_1\| (|\psi_{j,l}| + \|L_1\| + \|L_2\|)|(x, y)|): (x,y) \in \mathbb{R}^n \} 
		 \\=\,& \{ (A_j^{(k)} x + E_k \widetilde{\psi}_{j,l}(x) + \mathcal{R}^{(k)}_{j,l}(x,y), x, y) : (x,y) \in \mathbb{R}^n \} 
	\end{align*}
	where 
	\begin{align*}
		|\mathcal{R}^{(k)}_{j,l}(x,y)| \leq\,& C (\beta_k + |\e_{k}^{-1} A_{j}^{(k)} - A_{j}|) \,E_k \|L_2\| + C E_k^2 \|L_1\| (|\psi_{j,l}| + \|L_1\| + \|L_2\|)
	\end{align*}
with $C = C(n, q)$. The rest of the conclusions are now immediate.
\end{proof}

\begin{lemma} \label{blowup_norm_conv_lemma}
Let $p \in \{2, \ldots, q\}$, ${\mathbf q} \in {\mathfrak M}_{p}$, $M \geq 1$, $w = (w_{j})_{j=1}^{p} \in \mathfrak{B}_{p, {\mathbf q}}(M, A_{1}, \ldots, A_{p})$ with associated sequences of locally area minimizing rectifiable currents $T_{k}$ in ${\mathbf B}_{1}(0)$, cones ${\mathbf C}_{k} \in 
{\mathcal C}_{p, {\mathbf q}}$ and $m \times 2$ matrices $A^{(k)}_{j}$ ($1 \leq j \leq p$) corresponding to ${\mathbf C}_{k}$ (as in condition (A) above) such that ${\hat E}_{k}^{-1} A_{j}^{(k)} \to A_{j}$ as $k \to \infty$ for each $j \in \{1, 2, \ldots, p\}$, where ${\hat E}_{k} = E(T_{k}, {\mathbf P}_{0}, {\mathbf B}_{1}(0))$. Let $E_{k} = E(T_{k}, {\mathbf C}_{k}, {\mathbf B}_{1}(0))$ and let  $\widetilde{\psi}  = (\widetilde{\psi}_{j})_{j=1}^{p} \in \widetilde{{\mathfrak L}}_{p, {\mathbf q}}(A_{1}, \ldots, A_{p})$.   For each $k = 1, 2, 3, \ldots,$ let  $\widetilde{\mathbf C}_{k}$ be the cone 
corresponding to ${\mathbf C}_{k}$, $\widetilde{\psi}$, $E_{k}$ given by Lemma~\ref{blowup2L_lemma}.  Then for every $\rho \in (0,1)$, 
\begin{equation} \label{blowup_norm_conv_eqn1}
	\lim_{k \rightarrow \infty} E_{k}^{-2} E(T_{k}, \widetilde{\mathbf  C}_{k}, {\mathbf B}_{\rho}(0))^{2} 
	\leq \frac{1}{\omega_n \rho^{n+2}}\int_{B_{\rho}(0)} \sum_{j=1}^p \mathcal{G}(w_{j},\widetilde{\psi}_{j})^2; 
\end{equation} 
\begin{equation}\label{blowup_norm_conv_eqn2}
\lim_{k \to \infty} E_{k}^{-2}\int_{{\mathbf B}_{\rho}(0) \setminus \{r \leq \rho/8\}} {\rm dist}^{2} \, (X, {\rm spt} \, T_{k}) d\|\widetilde{\mathbf C}_{k}\| \leq \int_{B_{\rho}(0) \setminus \{r \leq \rho/8\}} \sum_{j=1}^{p}\mathcal{G}(w_{j}, \widetilde{\psi}_{j})^{2}.
\end{equation}
\end{lemma}

\begin{proof}
By (\ref{excess-2}),  for every $\delta \in (0,1/4),$ $\rho \in (0, 1)$ and $k$ sufficiently large we have that 
$$\int_{{\mathbf B}_{\rho}(0) \cap \{r < \delta\}} {\rm dist}^{2} \, (X, {\rm spt}\, {\mathbf C}_{k}) \, d\|T_{k}\| \leq C \delta^{1/2} E_{k}^2,$$ and consequently, by the triangle inequality,  
\begin{equation*} 
\int_{{\mathbf B}_{\rho}(0) \cap \{r < \delta\}} {\rm dist}^{2} \, (X, {\rm spt} \, \widetilde{\mathbf C}_{k}) \, d\|T_{k}\| \leq C \delta^{1/2} E_{k}^2,  
\end{equation*}
for some constant $C = C(n,m,q, \rho) \in (0,\infty)$.  In view of this, the conclusions \eqref{blowup_norm_conv_eqn1} and \eqref{blowup_norm_conv_eqn2} follow from the definition of fine blow-up, \eqref{blowup2L_eqn5} and the estimates \eqref{homog_error} and  
\eqref{blowup2L_eqn6}. 
\end{proof}

\subsection{Main estimates for fine blow-ups}
We shall now derive estimates for $w \in \mathfrak{B}_{p, {\mathbf q}}$ that correspond to the estimates in Section~\ref{estimates}.  These estimates will form the basis of our asymptotic analysis of the fine blow-ups (carried out in Sections~\ref{sec:homogblowup_sec} and \ref{sec:asymptotics_sec} below) which in turn will play a key role in the proof of the main excess decay result, Lemma~\ref{excess-improvement}. 

\begin{lemma} \label{blowup_est_lemma} 
Let $p \in \{2, \ldots, q\}$, ${\mathbf q} \in {\mathfrak M}_{p}$, $M \geq 1$ and $w = (w_{i})_{i=1}^{p} \in {\mathfrak B}_{p, {\mathbf q}}(M, A_{1}, \ldots, A_{p})$ for some $m \times 2$ matrices $A_{1}, \ldots, A_{p}.$ The following estimates hold: 
\begin{enumerate}
\item[(a)] for each $\psi  = (\psi_{i})_{i=1}^{p}\in \widetilde{\mathfrak{L}}_{p, {\mathbf q}}(A_{1}, \ldots, A_{p})$ and each $\rho \in (0, 1/2),$
\begin{align} 
	\label{blowup_est1} &\int_{B_{\rho/2}(0)} \sum_{i=1}^p 
		R^{2-n} \left| \frac{\partial (w_{i}/R)}{\partial R} \right|^2 
		\leq C \rho^{-n-2} \int_{B_{\rho}(0)} \sum_{i=1}^p \mathcal{G}(w_{i},\psi_{i})^2, 
	\end{align}
where  $R = |X|$ and $C = C(n,m,q) \in (0,\infty)$ is a constant;
\item[(b)] for each $\rho \in (0, 1/2)$, 
\begin{equation}\label{blowup_est3} 
	\rho^{-n}\int_{B_{\rho/2}(0)} \sum_{i=1}^p  |D_y w_{i}|^2 \leq C \rho^{-n-2}\int_{B_{\rho}(0)} \sum_{i=1}^{p}|w_{i}|^{2}, 
\end{equation}
where $X = (x,y)$ for $x = (x_1,x_2)$ and $y = (x_3,\ldots,x_n)$ and $C = C(n,m,q) \in (0,\infty)$ is a constant; 
\item[(c)] there exist functions $\lambda_{1} : B^{n-2}_{1/4}(0) \rightarrow \mathbb{R}^m$, $\lambda_{2} : B^{n-2}_{1/4}(0) \rightarrow \mathbb{R}^2$ with 
\begin{equation}\label{lambda_bound}
\sup_{B^{n-2}_{1/4}(0)} \left(|\lambda_{1}(z)| + |\lambda_{2}(z)| \right)\leq C
\end{equation}
such that for  each $\sigma \in (0,1/2)$, $\rho \in (0, 1/2)$ and $z \in B^{n-2}_{1/4}(0)$, 
\begin{equation} \label{blowup_est4}
	\int_{B_{\rho/2}(0,z)} \sum_{i=1}^p \frac{|w_{i}(X) - \left(\lambda_{1}(z) -  A_{i} \lambda_{2}(z)\right)|^2}{
		|X-(0,z)|^{n+2-\sigma}} dX \leq C \rho^{-n-2+\sigma}\int_{B_{\rho}(0,z)} \sum_{i=1}^{p} |w_{i} - (\lambda_{1}(z) - A_{i}\lambda_{2}(z))|^{2}
\end{equation}
where $C = C(n,m,q,M, \sigma) \in (0,\infty)$.  Moreover, $\lambda_{1}(z), \lambda_{2}(z)$ are uniquely determined (by $z$, $w$ and $A_{1}, A_{2}, \ldots, A_{p}$) 
subject only to the condition  
$$\int_{B_{\rho}(0,z)} \sum_{i=1}^p \frac{|w_{i}(X) - (\lambda_{1}(z) - A_{i} \lambda_{2}(z))|^2}{
		|X-(0,z)|^{n+2-\sigma}} dX < \infty \;\;\; \mbox{for some $\sigma \in (0, 1/2)$ and some $\rho \in (|z|, 1/2).$}$$ 
\end{enumerate}
\end{lemma}

\begin{proof}
Let $w = (w_{i})_{i=1}^{p} \in \mathfrak{B}_{p, {\mathbf q}}(M, A_{1}, \ldots, A_{p})$. Choose  a sequence $(T_{k})$ of locally area minimizing rectifiable currents in ${\mathbf B}_{1}(0)$ and a sequence $({\mathbf C}_{k})$ of cones in ${\mathcal C}_{p, {\mathbf q}}$ such that 
the $m \times 2$ matrices $A^{(k)}_{i}$ ($1 \leq i \leq p$) corresponding to ${\mathbf C}_{k}$ (as in condition (A) above) satisfy 
${\hat E}_{k}^{-1} A_{i}^{(k)} \to A_{i}$ as $k \to \infty$ for each $i \in \{1, 2, \ldots, p\},$ and such that 
$w$ is the fine blow-up of of $(T_{k})$ relative to 
$({\mathbf C}_{k})$ in accordance with Definition~\ref{blowupclass_defn}. 

To see part (a), let $\psi  = (\psi_{i})_{i=1}^{p}\in \widetilde{\mathfrak{L}}_{p, {\mathbf q}}(A_{1}, \ldots, A_{p})$. For each $k=1, 2, \ldots,$ let $\widetilde{\mathbf C}_{k} \in \widetilde{\mathcal C}_{p^{\prime}, {\mathbf q}^{\prime}}$  be the cone 
given by Lemma~\ref{blowup2L_lemma} corresponding to ${\mathbf C}_{k}$, $\psi$ and $E_{k} = E(T_{k}, {\mathbf C}_{k}, {\mathbf B}_{1}(0)),$ and let $R_{k} \in {\mathcal S}$ be such that 
$\widehat{\mathbf C}_{k} \equiv (e^{-R_{k}})_{\#} \, \widetilde{\mathbf C}_{k} \in {\mathcal C}_{p^{\prime}, {\mathbf q}^{\prime}}.$ 
Let $\epsilon_{0} = \epsilon_{0}(n, m, q, 1/2)$, $\beta_{0} = \beta_{0}(n,m,q, 1/2)$ be as in Theorem~\ref{keyest thm} taken with $\gamma = 1/2$. 
Fix $\rho \in (0, 1/2).$ Note that for sufficiently large $k$, Hypothesis~($\star$) is satisfied with $\eta_{0, \rho \, \#} \, T_{k}$, ${\mathbf C}_{k}$ in place of $T$, ${\mathbf C}$ (since $E(T_{k}, {\mathbf P}_{0}, {\mathbf B}_{1}(0)) \to 0$ and hence, by condition (4) of Section~\ref{fine-blowup-prelim}, $Q(T_{k}, {\mathbf C}_{k}, {\mathbf B}_{1}(0)) \to 0$), and hence Hypothesis~($\star$) also holds with $\widehat{T}_{k} \equiv \eta_{0, \rho \, \#} \, (e^{-R_{k}})_{\#} \, T_{k},$ $\widehat{\mathbf C}_{k}$ in place of $T$, ${\mathbf C}$. Suppose that for sufficiently large $k$, 
Hypothesis~($\star\star$) is satisfied with $\widehat{T}_{k}$, $\widehat{\mathbf C}_{k}$ in place of 
$T$, ${\mathbf C}$, i.e.\ that 
\begin{equation}\label{blowup_est_lemma-1}
Q(\widehat{T}_{k},\widehat{\mathbf{C}}_{k},\mathbf{B}_1(0)) \leq \beta_0 \inf_{\mathbf{C}' \in \bigcup_{p'=1}^{p-1} \mathcal{C}_{q,p'}} Q(\widehat{T}_{k},\mathbf{C}',\mathbf{B}_1(0)) 
\end{equation} 
for sufficiently large $k$. 
Then we can apply Theorem~\ref{keyest thm}(a)  (with $\widehat{T}_{k},$ $\widehat{\mathbf C}_{k}$ in place of $T$, ${\mathbf C}$) to deduce that
$$\int_{{\mathbf B}_{\rho/2}(0)} \frac{|X^{\perp}|^{2}}{|X|^{n+2}} d\|T_{k}\|(X) \leq C\rho^{-n-2}\int_{{\mathbf B}_{\rho}(0)} {\rm dist}^{2}(X, {\rm spt} \, \widetilde{\mathbf C}_{k}) \, d\|T_{k}\|(X)$$
for all sufficiently large $k$, where $C = C(n,m, q) \in (0, \infty)$. Reasoning exactly as in the proof of Corollary~\ref{keyest cor2}, this implies, for any $\tau \in (0, \rho/4)$ and all sufficiently large $k$, 
$$\int_{B_{\rho/2}(0) \setminus \{r < \tau\}} \sum_{i=1}^{p}R^{2-n}\left|\frac{\partial (u_{i}^{(k)}/R)}{\partial R}\right|^{2} \leq C\rho^{-n-2}\int_{{\mathbf B}_{\rho}(0)} {\rm dist}^{2}(X, {\rm spt} \, \widetilde{\mathbf C}_{k}) \, d\|T_{k}\|(X)$$ where $u_{i}^{(k)}$ are as in (\ref{graphrep close2plane concl c}). Dividing this by $E_{k}$, letting $k \to \infty$ and applying Lemma~\ref{blowup_norm_conv_lemma}, and then letting $\tau \to 0$, this yields \eqref{blowup_est1}. If on the other hand \eqref{blowup_est_lemma-1} fails 
for infinitely many $k$, then (see Remark~\ref{tildeC rmk}) we can choose $p^{\prime} \in \{1, \ldots, p-1\}$, a subsequence of $(k)$ without relabeling, and cones ${\mathbf C}_{k}^{\prime} \in 
{\mathcal C}_{q, p^{\prime}}$ such that: 
\begin{itemize}
\item[(i)] $Q(\widehat{T}_{k}, {\mathbf C}^{\prime}_{k}, {\mathbf B}_{1}(0)) < (2\beta_{0}^{-1})^{p-1} Q(\widehat{T}_{k}, \widehat{\mathbf C}_{k}, {\mathbf B}_{1}(0))$ and 
\item[(ii)] either $p^{\prime} = 1$ or Hypothesis~($\star\star$) holds with $\widehat{T}_{k},$ ${\mathbf C}^{\prime}_{k}$ in place of $T$, ${\mathbf C}$. 
\end{itemize}
We clearly also have that Hypothesis~($\star$) is satisfied with $\widehat{T}_{k}$, ${\mathbf C}_{k}^{\prime}$ in place of $T$, ${\mathbf C}$
(by (i), and the fact that $Q(T_{k}, {\mathbf C}_{k}, {\mathbf B}_{1}(0)) \to 0$ whence $Q(\widehat{T}_{k}, \widehat{\mathbf C}_{k}, {\mathbf B}_{1}(0)) \to 0$). 
We can thus apply Theorem~\ref{keyest thm}(a) again, with $\widehat{T}_{k}$, ${\mathbf C}_{k}^{\prime}$ in place of $T$, ${\mathbf C}$ (noting that if we have in (ii) above that $p^{\prime} = 1,$ then Theorem~\ref{keyest thm}(a) holds by a standard argument based on the monotonicity formula, Theorem~\ref{lip approx thm} and Lemma~\ref{energy est lemma}), and combine the resulting inequality with (i) above to deduce that 
\begin{align*}
&\int_{{\mathbf B}_{\rho/2}(0)} \frac{|X^{\perp}|^{2}}{|X|^{n+2}} d\|T_{k}\|(X) \leq C\rho^{-n-2}\int_{{\mathbf B}_{\rho}(0)} {\rm dist}^{2}(X, {\rm spt} \, \widetilde{\mathbf C}_{k}) \, d\|T_{k}\|(X)\\
& \hspace{2in}+ C\rho^{-n-2}\int_{e^{-R_k}({\mathbf B}_{\rho/2}(0) \setminus \{r(X) < \rho/16\})} {\rm dist}^{2}(X, {\rm spt} \, T_{k}) \, d\|\widetilde{\mathbf C}_{k}\|(X)
\end{align*}
for infinitely many $k$. Again reasoning as in the proof of Corollary~\ref{keyest cor2}, this implies, for any $\tau \in (0, \rho/4)$,  
 \begin{align*}
&\int_{B_{\rho/2}(0) \setminus \{r < \tau\}} R^{2-n}\left|\frac{\partial (u_{i}^{(k)}/R)}{\partial R}\right|^{2} \leq C\rho^{-n-2}\int_{{\mathbf B}_{\rho}(0)} {\rm dist}^{2}(X, {\rm spt} \, \widetilde{\mathbf C}_{k}) \, d\|T_{k}\|(X)\\
&\hspace{2.2in} + C\rho^{-n-2}\int_{e^{-R_k}({\mathbf B}_{\rho/2}(0) \setminus \{r(X) < \rho/16\})} {\rm dist}^{2}(X, {\rm spt} \, T_{k}) \, d\|\widetilde{\mathbf C}_{k}\|(X)
\end{align*}
for infinitely many $k$. Dividing this by $E_{k}$, letting $k \to \infty$ and then letting $\tau \to 0$, in view of Lemma~\ref{blowup_norm_conv_lemma} and the fact that $e^{-R_k} \to I$, this yields \eqref{blowup_est1}.

Part (b)  follows directly from (\ref{estimate-E1}) and \eqref{blowup_norm_conv_eqn2_0}.

To see Part (c), let $\rho  \in (0, 1/2)$ and $z \in B^{n-2}_{1/4}(0)$.  By \eqref{excess-1}, for each $k= 1,2,\ldots$ there exists $Z_{k}(z)  = (\chi_{k}(z), \xi_{k}(z), \zeta_{k}(z)) \in {\mathbf B}_{\delta_{k}}(0,z)$ with $\Theta(\|T_{k}\|, Z_{k}(z)) \geq q$.  By (\ref{estimate-E2}), 
\begin{eqnarray}
&&\int_{B_{\rho/2}(\xi_{k}(z), \z_{k}(z)) \cap \{r > \tau_{k}\}} \frac{|u_i^{(k)}(x,y) - (\chi_{k}(z) -  A_i^{(k)} \xi_{k}(z))|^{2}}{|(u_{i}^{(k)}(x, y) - \chi_{k}(z), x - \xi_{k}(z),y - \z_{k}(z))|^{n+2-\sigma}} 
\nonumber\\
&&\hspace{2in}\leq C_{1}  \rho^{-n-2 + \sigma} \int_{{\mathbb R}^{m} \times B_{\rho}(\xi_{k}(z), \z_{k}(z))} {\rm dist}^{2} \, (X, {\rm spt} \, \nu_{Z_{k}(z) \, \#} {\mathbf C}_{k}) \, d\|T_{k}\|
\end{eqnarray}
for all sufficiently large $k$, where $C_{1} = C_{1}(n, m, q, M, 1/2, \sigma)$. 
By (\ref{goodpt-est}), there are 
$\lambda_{1}(z) \in {\mathbb R}^{m}$ and $\lambda_{2}(z) \in {\mathbb R}^{2}$ with $|\lambda_{1}(z)|, |\lambda_{2}(z)| \leq C = C(n, m, q, M)$ such that, passing to a further subsequence,  $E_{k}^{-1}\chi_{k}(z) \to \lambda_{1}(z)$ and 
$E_{k}^{-1} {\widehat E}_{k} \xi_{k}(z) \to \lambda_{2}(z).$ Hence dividing the above by $E_{k}$ and letting $k \to \infty$, we obtain the existence of $\lambda_{1}(z)$, $\lambda_{2}(z)$ such that (\ref{blowup_est4}) holds with $C = C_{1}$.  	
Moreover $\lambda_{1}(z)$, $\lambda_{2}(z)$ are uniquely determined by $z$, $w$ and $A_{1}, A_{2}, \ldots, A_{p}$ (independently of $\rho$, $\sigma$ and the chosen sequence $(Z_{k}(z))$); to see this, let $z \in B^{n-2}_{1/4}(0)$, $|z| < \rho_{1}, \rho_{2} < 1/2$ and suppose that for $\ell=1, 2$, there are 
$\lambda_{1}^{(\ell)} \in {\mathbb R}^{m}$, $\lambda_{2}^{(\ell)} \in {\mathbb R}^{2}$ such that 
$$\int_{B_{\rho_{\ell}/2}(0,z)} \sum_{i=1}^p \frac{|w_{i}(X) - (\lambda_{1}^{(\ell)} -  A_{i} \lambda_{2}^{(\ell)})|^2}{
		|X-(0,z)|^{n+2-\sigma_{\ell}}} dX < \infty$$ 
for some $\sigma_{1}, \sigma_{2}$ with $0 < \sigma_{1} \leq \sigma_{2} <1/2.$ 
By replacing $\rho_{1}, \rho_{2}$ with $\min \, \{\rho_{1}, \rho_{2}\}$, we may assume $\rho_{1} = \rho_{2}$. By the triangle inequality this implies that 
\begin{equation*}
	\int_{B_{\rho_{1}/2}(0,z)} \sum_{i=1}^p  \frac{|(\lambda_{1}^{(1)} - \lambda_{1}^{(2)}) - A_{i}(\lambda_2^{(1)} - \lambda_2^{(2)})|^2}{|X-(0,z)|^{n+2-\sigma_{2}}} dX < \infty 
\end{equation*}
whence (since $\int_{B_{\rho_{1}/2}(0,z)} |X - (0, z)|^{-n-2+ \sigma_{2}} \, dX = \infty$) we must have $(\lambda_{1}^{(1)} - \lambda_{1}^{(2)}) -  A_{i}(\lambda_2^{(1)} - \lambda_2^{(2)}) = 0$ for each 
$i=1, 2, \ldots, p.$  By (\ref{rank2}) (which holds automatically, see Remark~\ref{conditions-on-Aj}) this gives first that $\lambda_{2}^{(1)} = \lambda_{2}^{(2)}$ and consequently also that $\lambda_{1}^{(1)} = \lambda_{1}^{(2)}.$
\end{proof}

\subsection{Classification of homogeneous degree 1 fine blow-ups}\label{sec:homogblowup_sec} 

In this section we establish the following classification theorem for homogeneous degree 1 fine blow-ups:

\begin{theorem} \label{classification_thm} 
Let $p \in \{2, \ldots, q\}$, ${\mathbf q}  = (q_{1}, \ldots , q_{p}) \in {\mathfrak M}_{p}$ and $M \geq 1$. Let $w  = (w_{j})_{j=1}^{p} \in \mathfrak{B}_{p, {\mathbf q}}(M)$ be homogeneous of degree 1 in the sense that $X \mapsto w_{j}(X)$ is an ${\mathcal A}_{q_{j}}({\mathbb R}^{m})$-valued homogeneous degree 1 function of $X \in B_{1}(0) \setminus \{0\} \times {\mathbb R}^{n-2}$ for every $j \in \{1,  \ldots, p\}.$  Then $w \in \widetilde{\mathfrak{L}}_{p, {\mathbf q}},$ and moreover 
there are a function 
$\psi = (\psi_{1}, \ldots, \psi_{p}) \in \mathfrak{L}_{p, {\mathbf q}}$, constant $(m \times 2)$ matrices $A_{1}, \ldots, A_{p}$, and  two linear functions 
$L_{1} \, : \, {\mathbb R}^{n-2} \to {\mathbb R}^{m}$, $L_{2} \, : \, {\mathbb R}^{n-2} \to {\mathbb R}^{2}$ with 
\begin{equation}\label{classification-bound} 
\|A_{j}\| + \|L_{1}\| + \|L_{2}\| \leq C, \;\;\; \max_{1 \leq i < j \leq p} \, \inf_{x \in {\mathbb S}^{1}} \, |A_{i}x - A_{j}x| \geq  c/M 
\end{equation}
where $C = C(n, m, q)$, $c = c(n, m, q),$ such that 
$w_{j}(x, y) = \sum_{k =1}^{q_{j}}\llbracket \psi_{j, k}(x) + L_{1}(y) + A_{j}L_{2}(y)\rrbracket$ for each $j \in \{1, 2, \ldots, p\}$ and $(x, y) \in {\mathbb R}^{2} \times {\mathbb R}^{n-2}.$
\end{theorem}

\begin{proof} Note that $w$ satisfies \eqref{blowup_est4} for some functions $\lambda_{1} \, : \, B_{1/4}(0) \cap \{0\} \times {\mathbb R}^{n-2} \to {\mathbb R}^{m}$, $\lambda_{2} \, : \, B_{1/4}(0) \cap \{0\} \times {\mathbb R}^{n-2} \to {\mathbb R}^{2}$ and some constant $(m \times 2)$ matrices $A_{i}$, $1 \leq i \leq p$. First consider the case $\lambda_{1} = 0,$ $\lambda_{2} = 0$, i.e.\  the case where 
for each $\sigma \in (0, 1/2)$, $\rho \in (0, 1/2)$ and $z \in B_{1/4}^{n-2}(0)$, 
\begin{equation} \label{homogrep2_eqn1}
	\int_{B_{\rho/2}((0,z))} \sum_{j=1}^p\frac{|w_{j}(X)|^2}{|X-(0,z)|^{n+2-\sigma}} dX \leq C \rho^{-n-2+\sigma}\int_{B_{\rho}((0,z))} \sum_{j=1}^{p} |w_{j}|^{2}
 \end{equation}
for some constant $C = C(n,m,q,M, \sigma) \in (0,\infty)$.  Since each $w_{j}$ is locally Dirichlet energy minimizing in $B_{1/2}(0) \setminus \{0\} \times {\mathbb R}^{n-2}$, 
it follows from standard sup and energy estimates for such functions (\cite{Almgren}) and  \eqref{homogrep2_eqn1}  that for every $\sigma \in (0,1/2),$ $\rho \in (0, 1/2)$ and every $X = (x,y) \in B^{2}_{\rho/2}(0) \times B_{1/4}^{n-2}(0) \setminus \{0\} \times {\mathbb R}^{n-2},$
\begin{equation} \label{homogrep2_eqn2}
	\sum_{j=1}^p  |w_{j}(X)|^{2} + |x|^{2-n} \int_{B_{|x|/2}(X)} \sum_{j=1}^p |Dw_{j}|^2  \leq C \left(\frac{|x|}{\rho}\right)^{2-\sigma} \rho^{-n} \int_{B_{\rho}(0, y)} \sum_{j=1}^{p} |w_{j}|^{2}
\end{equation}
where $C = C(n,m,q,\sigma) \in (0,\infty)$. From this, with the help of a standard covering argument, it follows that 
for each $\rho \in (0, 1/2)$, $\delta \in (0, 1/16)$ and $z \in B^{n-2}_{1/16}(0)$, 
\begin{equation}\label{homogrep2_eqn2_1_0}
\rho^{2-n}\int_{B_{\rho/16}(0, z) \cap \{\delta\rho/2 \leq r \leq \delta\rho\}} \sum_{j=1}^{p}|Dw_{j}|^{2} < C\delta^{2 - \sigma} \rho^{-n} \int_{B_{\rho}(0, z)} \sum_{j=1}^{p} |w_{j}|^{2} 
\end{equation}
where $C = C(n, m, q, \sigma),$ whence for each $\rho \in (0, 1/2)$, $\delta \in (0, 1/16)$ and $z \in B^{n-2}_{1/16}(0)$, 
\begin{equation}\label{homogrep2_eqn2_1_1}
\rho^{2-n}\int_{B_{\rho/16}(0, z) \cap \{r \leq \delta\rho\}} \sum_{j=1}^{p}|Dw_{j}|^{2} < C\delta^{2 - \sigma} \rho^{-n} \int_{B_{\rho}(0, z)} \sum_{j=1}^{p} |w_{j}|^{2};
\end{equation}
in particular, for each $\rho \in (0, 1/2)$ and $z \in B^{n-2}_{1/16}(0)$, 
\begin{equation}\label{homogrep2_eqn2_1}
\rho^{2-n}\int_{B_{\rho/16}(0,z)} \sum_{j=1}^{p}|Dw_{j}|^{2} < C\rho^{-n} \int_{B_{\rho}(0, z)} \sum_{j=1}^{p} |w_{j}|^{2}
\end{equation}
where $C = C(n, m, q).$ 

Now fix $j \in \{1, 2, \ldots, p\}$, and set $\overline{w} = w_{j}$ and $\overline{q} = q_{j}$. Since $\overline{w}$ is ${\mathcal A}_{\overline{q}}({\mathbb R}^{m})$-valued, we may write $$\overline{w}(X) = \sum_{\ell=1}^{\overline{q}} \llbracket \overline{w}_{\ell}(X) \rrbracket,$$ where 
$\overline{w}_{\ell}(X) = (\overline{w}_{\ell}^{1}(X), \overline{w}_{\ell}^{2}(X), \ldots, \overline{w}_{\ell}^{m}(X)) \in {\mathbb R}^{m}.$ With this notation, we next verify  the two identities: 
\begin{gather}
	\int_{\mathbb{R}^n} |D\overline{w}|^2 \zeta = -\int_{\mathbb{R}^n} \overline{w}_{\ell}^{\kappa} D_i \overline{w}_{\ell}^{\kappa} D_i \zeta, 
		\label{homogrep2_freqid1} \\
	\int_{\mathbb{R}^n}  \left(\tfrac{1}{2} |D\overline{w}|^2 \delta_{ik} - D_i \overline{w}_{\ell}^{\kappa} D_{k}\overline{w}_{\ell}^{\kappa}\right)D_{i}\zeta_{k} = 0, \label{homogrep2_freqid2}
\end{gather}
for all $\zeta, \zeta_{1}, \ldots, \zeta_{n} \in C^1_c(B_{1/16}(0)),$ where  we use the convention of summing over repeated indices. 
Indeed, since $\overline{w}$ is locally energy minimizing in $B_{1/16}(0) \setminus \{0\} \times {\mathbb R}^{n-2}$, \eqref{homogrep2_freqid1}, \eqref{homogrep2_freqid2} hold 
whenever $\zeta, \zeta_{1}, \ldots, \zeta_{n}  \in C^1_c(B_{1/16}(0) \setminus \{0\} \times \mathbb{R}^{n-2})$.  For $\delta \in (0,1/32)$, let $\widetilde{\chi}_{\delta} : [0,\infty) \rightarrow \mathbb{R}$ be a smooth function such that $0 \leq \widetilde{\chi}_{\delta}(r) \leq 1$, $\widetilde{\chi}_{\delta}(r) = 0$ for all $r \in [0,\delta/2]$, $\widetilde{\chi}_{\delta}(r) = 1$ for all $r \geq \delta$, and $|\widetilde{\chi}'_{\delta}(r)| \leq 3/\delta$. Define $\chi_{\delta} : \mathbb{R}^n \rightarrow \mathbb{R}$ by $\chi_{\delta}(x,y) = \widetilde{\chi}_{\delta}(|x|)$.  Let $\zeta \in C^1_c(B_{1/16}(0))$ be arbitrary and replace $\zeta$ in \eqref{homogrep2_freqid1} by $\chi_{\delta} \zeta$ to get, using \eqref{homogrep2_eqn2}, \eqref{homogrep2_eqn2_1_0} and \eqref{blowup_norm_conv_eqn2-1}, that
\begin{align} \label{homogrep2_eqn3}
	\left| \int_{B_{1/16}(0)} (|D\overline{w}|^2 \zeta + \overline{w}^{\kappa}_{\ell} D_i \overline{w}^{\kappa}_{\ell} D_i \zeta) \chi_{\delta} \right|
	&\leq \int_{B_{1/16}(0)} |\overline{w}| |D\overline{w}| |D\chi_{\delta}| |\zeta| 
	\\&\leq C \delta^{2-\sigma} \sup_{B_{1/16}(0)} |\zeta| \nonumber 
\end{align}
for every $\sigma \in (0,1/2)$, where $C = C(n,m,q,\sigma) \in (0,\infty)$.  In view of \eqref{homogrep2_eqn2_1} and \eqref{homogrep2_eqn2}, we may let $\delta \to 0^{+}$ in this to conclude \eqref{homogrep2_freqid1}. A similar argument (using both \eqref{homogrep2_eqn2_1_0} and \eqref{homogrep2_eqn2_1}) yields  \eqref{homogrep2_freqid2}. 

It is standard (see e.g.\ \cite[Section~4.4]{KrumWic2}) that \eqref{homogrep2_freqid1}, \eqref{homogrep2_freqid2} imply that either $\overline{w} \equiv \overline{q}\llbracket 0 \rrbracket$ in $B_{1/16}(0)$ or for each $Z \in B_{1/16}(0)$ and $\rho \in (0, 1/16)$, 
the frequency function 
\begin{equation} \label{homogrep2_freq}
	N_{\overline{w}, \, Z}(\rho) = \frac{\rho^{2-n} \int_{B_{\rho}(Z)} |D\overline{w}|^2}{\rho^{1-n} \int_{\partial B_{\rho}(Z)} |\overline{w}|^2}
\end{equation}
is well-defined and monotone nondecreasing as a function of $\rho \in (0,1/16)$.  Suppose that $\overline{w}$ is not identically zero, and write ${\mathcal N}_{\overline{w}}(Z) = \lim_{\rho \to 0^{+}} \, N_{\overline{w}, Z}(\rho)$. Since $\overline{w}$ is homogeneous of degree 1 from the origin, 
$1 = {\mathcal N}_{\overline{w}}(0) \geq {\mathcal N}_{\overline{w}}(Z)$ for each $Z \in B_{1/16}(0)$. 

Let  $Z = (0, z) \in B_{1/16}(0) \cap \{0\} \times {\mathbb R}^{n-2}$. By \eqref{homogrep2_eqn1} we have that for each $\sigma \in (0, 1/2)$, 
$\sup_{\rho \in (0, 1/4)} \, \rho^{-n-2+\sigma} \int_{B_{\sigma}(Z)} |\overline{w}|^{2} < \infty,$ while by the monotonicity of the frequency 
function \eqref{homogrep2_freq} we have  that 
$\inf_{\rho \in (0, \rho_{1})} \rho^{-n-2N_{\overline{w}, Z}(\rho_{1})} \int_{B_{\rho}(Z)} |\overline{w}|^{2} >0$ for each $\rho_{1} \in (0, 1/16)$. It follows that  $2N_{\overline{w}, Z}(\rho_{1}) \geq 2 - \sigma$ for each $\rho_{1} \in (0, 1/16)$ and $\sigma \in (0, /12)$, whence, letting $\sigma, \rho_{1} \to 0$, we 
deduce that ${\mathcal N}_{\overline{w}}(Z) \geq 1$.   
Thus ${\mathcal N}_{\overline{w}}(Z)= {\mathcal N}_{\overline{w}}(0) = 1$ for every $Z \in \{0\} \times {\mathbb R}^{n-2}$, and hence by standard arguments again, $\overline{w}$ is invariant under translations along the subspace $\{0\} \times {\mathbb R}^{n-2}$. Thus there  is a function $\overline{w}_{1} \in W^{1, 2}_{\rm loc}({\mathbb R}^{2}; {\mathcal A}_{\overline{q}}({\mathbb R}^{m}))$ such that 
$\overline{w}(x, y) = \overline{w}_{1}(x)$ for all $(x, y) \in {\mathbb R}^{2} \times {\mathbb R}^{n-2}.$ Since $\overline{w}_{1}$ is homogeneous of degree 1 and locally energy minimizing in ${\mathbb R}^{2} \setminus \{0\},$ $\overline{w}_{1}$ is given by $\overline{q}$ linear functions on ${\mathbb R}^{2},$ of which any two distinct functions take distinct values at every point ${\mathbb R}^{2} \setminus \{0\}$.  We have thus shown that a homogeneous degree 1 element $w \in {\mathfrak B}_{p, {\mathbf q}}$ belongs to ${\mathfrak L}_{p, {\mathbf q}}$ subject to the additional assumption that the two functions  $\lambda_{1}, \lambda_{2}$ corresponding to $w$ as in \eqref{blowup_est4} are both zero. 

To complete the proof of the theorem, let now $w  = (w_{1}, w_{2}, \ldots, w_{p})$ be an arbitrary homogeneous degree 1 element of  
${\mathfrak B}_{p, {\mathbf q}}(M).$ Then $w \in {\mathfrak B}_{p, {\mathbf q}}(M, A_{1}, \ldots, A_{p})$ for some $m \times 2$ matrices $A_{1}, \ldots, A_{p}$ satisfying 
\eqref{matrix-upperbound} and \eqref{rank2}. By Lemma~\ref{blowup_est_lemma}(c), there are functions $\lambda_{1} \, : \, B_{1/4}^{n-2}(0) \to {\mathbb R}^{m}$, $\lambda_{2} \, : \, B_{1/4}^{n-2}(0) \to {\mathbb R}^{2}$  satisfying \eqref{lambda_bound} and \eqref{blowup_est4}.  Fix $j \in \{1, 2, \ldots, p\}$ and write 
$w_{j}(X) = \sum_{\ell=1}^{q_{j}} \llbracket w_{j, \ell}(X) \rrbracket$ with $w_{j, \ell}(X) \in {\mathbb R}^{m}.$ Let $w_{j, a}$ denote the single-valued average given by 
$w_{j, a}(X) = q_{j}^{-1}\sum_{\ell=1}^{q_{j}} w_{j, \ell}(X)$ for $X \in B_{1/2}(0) \setminus \{0\} \times {\mathbb R}^{n-2}$. Since $w_{j}$ is locally Dirichlet energy minimizing on $B_{1/2}(0) \setminus \{0\} \times {\mathbb R}^{n-2}$, we have that $w_{j, a}$ is harmonic on  $B_{1/2}(0) \setminus \{0\} \times {\mathbb R}^{n-2}$. By \eqref{blowup_est4} (with $\rho = 1/4$ and $\sigma = 1/4$), \eqref{blowup_norm_conv_eqn2-1}, \eqref{lambda_bound} and \eqref{matrix-upperbound}, we have that  for  each  $z \in B^{n-2}_{1/4}(0)$, 
\begin{equation} \label{blowup_est4_avg}
	\int_{B_{1/4}(0,z)} \frac{|w_{j, a}(X) - \left(\lambda_{1}(z) -  A_{j} \lambda_{2}(z)\right)|^2}{
		|X-(0,z)|^{n+2-1/4}} dX \leq C
\end{equation}
where $C = C(n,m,q, M) \in (0,\infty)$. Since $|\lambda_{1}(z) - A_{j}\lambda_{2}(z)| \leq C$, this implies that $$\rho^{-n}\int_{B_{\rho}(0, z)} |w_{j, a}|^{2} \leq C$$ for any $z \in B_{1/4}^{n-2}(0)$ and any $\rho \in (0, 1/4)$, where $C = C(n,m, q, M).$ Consequently, by the mean value property for harmonic functions, we have that $w_{j, a}$ is bounded on 
$B_{1/4}(0) \setminus \{0\} \times {\mathbb R}^{n-2}.$ It is then standard to see that $w_{j, a}$ extends to $B_{1/4}(0)$ as a harmonic function (e.g.\ by using the fact that $\int_{B_{1/4}(0)} |Dw_{j, a}|^{2} \z^{2} \leq 4\int_{B_{1/4}(0)} |w_{j, a}|^{2}|D\z|^{2}$ for all $\z \in C^{1}_{c}(B_{1/4}(0) \setminus (\{0\} \times {\mathbb R}^{n-2}))$ and the fact that $\{0\} \times {\mathbb R}^{n-2}$ has vanishing 2-capacity to verify first that $w_{j, a} \in W^{1, 2}_{\rm loc}(B_{1/4}(0))$, followed by another use of the same vanishing 2-capcity property to verify that $w_{j, a}$ is weakly harmonic in $B_{1/4}(0)$). Since $w_{j}$ is homogeneous of degree 1 on $B_{1/2}(0) \setminus \{0\} \times {\mathbb R}^{n-2}$ by assumption, this extended function (which we shall continue to denote by $w_{j, a}$) is homogeneous of degree 1, so it is a linear function on $B_{1/4}(0)$. Since \eqref{blowup_est4_avg} implies that 
$\left.w_{j, a}\right|_{\{0\} \times B_{1/4}^{n-2}(0)} = \lambda_{1} - A_{j} \lambda_{2},$ we conclude that for each $j \in \{1, 2, \ldots, p\}$, the function $\lambda^{(j)}  =  \lambda_{1} - A_{j} \lambda_{2} \, : \, B_{1/4}^{n-2}(0) \to {\mathbb R}^{m}$ is linear. Since by  \eqref{rank2} there are $i, j \in \{1, 2, \ldots, p\}$ such that $A_{i} - A_{j}$ has full rank $(=2)$, it follows that $\lambda_{2}$ is linear and consequently so is $\lambda_{1}.$  Thus, writing points $z \in {\mathbb R}^{n-2}$ as column vectors, there is a constant $m \times (n-2)$ matrix $D_{1}$ and a constant $2 \times (n-2)$ matrix $D_{2}$ such that $\lambda_{1}(z) = D_{1}z$ and $\lambda_{2}(z) = D_{2} z.$ Now, there is $(T_{k})$ a sequence of locally area minimizing rectifiable currents in ${\mathbf B}_{1}(0)$ with $\partial \, T_{k} \llcorner {\mathbf B}_{1}(0) = 0$, and $({\mathbf C}_{k})$ a sequence of cones in ${\mathcal C}_{p, {\mathbf q}}$ such that $w$ is the fine blow-up of 
$(T_{k})$ relative to $({\mathbf C}_{k})$. Let $E_{k} = E(T_{k}, {\mathbf C}_{k}, {\mathbf B}_{1}(0))$. 
Define rotations $R_{k}^{(1)},  R_{k}^{(2)}\, : \, {\mathbb R}^{m+n} \to {\mathbb R}^{m+n}$  by 
$R_{k}^{(1)} = e^{E_{k}M_{1}}$ and $R_{k}^{(2)} = e^{\widehat{E}_{k}^{-1}E_{k}M_{2}}$ where $M_{1}$, $M_{2}$ are as in the proof of Lemma~\ref{blowup2L_lemma}
with $\lambda_{1}$, $\lambda_{2}$ in place of $L_{1}$, $L_{2}$ respectively and $\widehat{E}_{k} = E(T_{k}, {\mathbf P}_{0}, {\mathbf B}_{1}(0))$. 
By the argument of \cite[Theorem~10.1]{KrumWic2} (see the end of \cite[Section 10.2]{KrumWic2}) and Lemma~\ref{blowup2L_lemma}, 
we see that the fine blow-up $\widetilde{w}  = (\widetilde{w}_{1}, \widetilde{w}_{2}, \ldots, \widetilde{w}_{p}) \in {\mathfrak B}_{p, {\mathbf q}}$ of the rotated sequence $\widetilde{T}_{k}  = R^{(1)}_{k \, \#} \, R^{(2)}_{k \, \#} \, T_{k}$ relative to the (same) sequence of cones $({\mathbf C}_{k})$ (by the excess $\widetilde{E}_{k} = E(\widetilde{T}_{k}, {\mathbf C}_{k}, {\mathbf B}_{1}(0))$)  satisfies 
\begin{equation}\label{rotated-blowup}
\widetilde{c} \,\widetilde{w}_{j}(x, y) = w_{j}(x, y) - (\lambda_{1}(y) - A_{j} \lambda_{2}(y))
\end{equation}
for some constant $\widetilde{c} \in [0, C]$ (with $C = C(n, m, q)$) and all $j =1, 2, \ldots, p;$ in fact 
$$\widetilde{c} = \limsup_{k \to \infty} \, E_{k}^{-1}\widetilde{E}_{k}.$$ 
If $\widetilde{c} = 0$ then $w_{j}(x, y) = q_{j} \llbracket \lambda_{1}(y) - A_{j}\lambda_{2}y\rrbracket$ for each $j=1, 2, \ldots, p$, so in this case $w \in \widetilde{\mathfrak L}_{p, {\mathbf q}}$ trivially. If on the other hand $\widetilde{c} >0$, then by \eqref{blowup_est4} and \eqref{rotated-blowup}, we see that $\widetilde{w}$ is a homogeneous degree 1 element of ${\mathfrak B}_{p, {\mathbf q}}$ satisfying, 
for  each $\sigma \in (0,1/2)$, $\rho \in (0, 1/2)$ and $z \in B^{n-2}_{1/4}(0)$, 
\begin{equation*} 
	\int_{B_{\rho/2}(0,z)} \sum_{j=1}^p \frac{|\widetilde{w}_{j}(X)|^{2}}{
		|X-(0,z)|^{n+2-\sigma}} dX \leq C \rho^{-n-2+\sigma}\int_{B_{\rho}(0,z)} \sum_{j=1}^{p} |\widetilde{w}_{j}|^{2}
\end{equation*}
where $C = C(n,m,q,M,\sigma) \in (0,\infty)$. So by applying the special case (i.e.\ the case $\lambda_{1} = 0, \lambda_{2} = 0$) of the theorem just proved to $\widetilde{w}$, we see that $\widetilde{w} \in {\mathfrak L}_{p, {\mathbf q}}.$ By \eqref{rotated-blowup}, this means that 
$w \in \widetilde{\mathfrak L}_{p, {\mathbf q}}$, i.e.\ 
$w_{j}(x, y) = \sum_{k =1}^{q_{j}}\llbracket \psi_{j, k}(x) + L_{1}(y) + A_{j}L_{2}(y)\rrbracket$ for each 
$j \in \{1, 2, \ldots, p\}$ and $(x, y) \in {\mathbb R}^{2} \times {\mathbb R}^{n-2}$ where $L_{i} = \lambda_{1}$ and $L_{2} = \lambda_{2}$. 
The asserted bounds \eqref{classification-bound} follow from  \eqref{lambda_bound}, \eqref{matrix-upperbound} and \eqref{rank2}. \end{proof}

\subsection{Asymptotic decay of fine blow-ups} \label{sec:asymptotics_sec}

The main result of this section is a decay estimate for the fine blow-ups, Theorem~\ref{blowupdecay_lemma}. 
Broadly speaking, our proof of Theorem~\ref{blowupdecay_lemma} will employ the well-known hole-filling technique, in a manner similar to its use in \cite{Sim93} for the multiplicity 1 counterpart of the result, and will be based on the following lemma.
\begin{lemma} \label{holefilling_lemma}
Let $p \in \{2, \ldots, q\}$ and ${\mathbf q} \in {\mathfrak M}_{p}$. There is a constant $C = C(n, m, q) \in (0, \infty)$ such that if $w \in \mathfrak{B}_{p, {\mathbf q}}$ and if $\psi \in \widetilde{\mathfrak L}_{p, {\mathbf q}}$ is such that 
\begin{equation}\label{holefilling_est-0}
\int_{B_{1/2}(0)} \sum_{j=1}^p  \mathcal{G}(w_{j},\psi_{j})^2  \leq 2 \inf_{\psi^{\prime} \in \widetilde{\mathfrak L}_{p, {\mathbf q}}} \, \int_{B_{1/2}(0)} \sum_{j=1}^p  \mathcal{G}(w_{j},\psi^{\prime}_{j})^2
\end{equation} 
then 
\begin{equation} \label{holefilling_est}
	\int_{B_{1/2}(0)} \sum_{j=1}^p  \mathcal{G}(w_{j},\psi_{j})^2  
	\leq C \int_{B_{1/2}(0) \setminus B_{1/8}(0)} \sum_{j=1}^p \left| \frac{\partial (w_{j}/R)}{\partial R} \right|^2. 
\end{equation}
\end{lemma}
Because ours is a higher multiplicity setting and the functions involved are multi-valued, the proof of Lemma~\ref{holefilling_lemma} will have to be different from that of the analogous result in \cite{Sim93}.  We proceed in two steps, where the first step is Lemma~\ref{holefilling0_lemma} below giving the same conclusions as Lemma~\ref{holefilling_lemma} subject to a weaker condition 
on $\psi$ (namely, inequality \eqref{holefilling0_hyp1}) than \eqref{holefilling_est-0}, together with an additional hypothesis on $w$ (namely, condition \eqref{holefilling0_hyp2}) that is the analogue, for fine blow-ups, of Hypothesis~($\star\star$) on area minimizing currents associated with fine blow-ups. Note also that in Lemma~\ref{holefilling0_lemma}, we work under the assumption that $\psi \in {\mathfrak L}_{p, {\mathbf q}}$ (rather than $\psi \in \widetilde{\mathfrak L}_{p, {\mathbf q}}$). 

Let $p \in \{2, \ldots, q\}$, ${\mathbf q} = (q_{1}, q_{2}, \ldots, q_{p})  \in {\mathfrak M}_{p}$ and $\psi  = (\psi_{1}, \psi_{2}, \ldots, \psi_{p}) \in {\mathfrak L}_{p, {\mathbf q}}.$ Recall that then $\sum_{j=1}^{p} q_{j} = q$ and for $j \in \{1, 2, \ldots, p\}$, we have that $\psi_{j}(x, y) = \sum_{k=1}^{q_{j}}\llbracket \psi_{j, k}(x) \rrbracket$ for all $(x, y) \in {\mathbb R}^{n} \approx {\mathbb R}^{2} \times {\mathbb R}^{n-2}$, 
where $\psi_{j, k} \, : \, {\mathbb R}^{2} \to {\mathbb R}^{m}$ are linear functions. In Lemma~\ref{holefilling0_lemma}, we shall use the following additional notation: 

For $j \in \{1, 2, \ldots, p\}$, let $d_{j}(\psi)$ be the number of \emph{distinct} functions in the collection $\{\psi_{j, k} \, : \,  k = 1, \ldots, q_{j}\}.$ 

Let $d(\psi) = \sum_{j=1}^{p} d_{j}(\psi).$ Note that then $1 \leq d_{j}(\psi) \leq q_{j}$ and $p \leq d(\psi) \leq q$. 

For $p \leq s \leq q$, let ${\mathfrak L}_{p, {\mathbf q}}(s) = \{\psi \in {\mathfrak L}_{p, {\mathbf q}} \, : \, d(\psi) = s\}$ and note that 
${\mathfrak L}_{p, {\mathbf q}} = \cup_{s=p}^{q} {\mathfrak L}_{p, {\mathbf q}}(s)$.

\begin{lemma} \label{holefilling0_lemma}
For every $\overline{M} \geq 1$ and $M \geq 1,$ there exist constants $\overline{\beta}  = \overline{\b}(n, m, q, M, \overline{M}) \in (0,1)$ and $\overline{C}  = \overline{C}(n, m, q, M, \overline{M}) \in (0,\infty)$ such that the following holds: if 
$p \in \{2, \ldots, q\}$, ${\mathbf q} \in {\mathfrak M}_{p}$,  $w \in \mathfrak{B}_{p, {\mathbf q}}$ and $\psi \in \mathfrak{L}_{p,{\mathbf q}}$ 
are such that  
\begin{equation} \label{holefilling0_hyp1}
	\int_{B_{1/2}(0)} \sum_{j=1}^p  \mathcal{G}(w_{j},\psi_{j})^2 
		\leq 2\overline{M}^{2} \inf_{\psi' \in \widetilde{\mathfrak{L}}_{p, {\mathbf q}}} \int_{B_{1/2}(0)} \sum_{j=1}^p \mathcal{G}(w_{j},\psi'_{j})^2,
\end{equation}
and either (i) $d(\psi)  = p$, or (ii) $d(\psi) > p$ and 
\begin{equation} \label{holefilling0_hyp2}
	\int_{B_{1/2}(0)} \sum_{j=1}^p \mathcal{G}(w_{j},\psi_{j})^2 
		\leq \overline{\beta}^{2} \inf_{\psi' \in \bigcup_{s=p}^{d(\psi) -1} \mathfrak{L}_{p, {\mathbf q}}(s)} 
			\int_{B_{1/2}(0)} \sum_{j=1}^p\mathcal{G}(w_{j},\psi'_{j})^2, 
\end{equation}
then \eqref{holefilling_est} holds with $C = \overline{C}$. 
\end{lemma}

\begin{proof}
Note that it suffices to fix $p \in \{2, \ldots, q\}$, ${\mathbf q}  = (q_{1}, \ldots, q_{p}) \in {\mathfrak M}_{p}$, $s \in \{p,\ldots,q\}$ and prove the lemma for $\psi \in \mathfrak{L}_{p,{\mathbf q}}(s),$ with $\overline{\beta}$ and $\overline{C}$ depending on $n$, $m$, $q$, $M$, $\overline{M}$, $p,$ ${\mathbf q}$ and $s$.
We argue by contradiction, so suppose that for some fixed $\overline{M} \geq 1$, $M \geq 1$, $p \in \{2, \ldots, q\},$ ${\mathbf q} \in {\mathfrak M}_{p}$  and $s \in \{p,\ldots,q\},$ and for each $\nu =1, 2, 3, \ldots,$ there exist $\overline{\beta}_{\nu}>0$, $w^{(\nu)} = (w^{(\nu)}_{1}, \ldots, w^{(\nu)}_{p}) \in 
\mathfrak{B}_{p, {\mathbf q}}(M)$, $\psi^{(\nu)} = (\psi^{(\nu)}_{1}, \ldots, \psi^{(\nu)}_{p}) \in \mathfrak{L}_{p,{\mathbf q}}(s)$  such that $\overline{\beta}_{\nu} \downarrow 0$ and for each $\nu$, 
\begin{equation} \label{holefilling0_hyp1-0}
	\int_{B_{1/2}(0)} \sum_{j=1}^p  \mathcal{G}(w^{(\nu)}_{j},\psi_{j}^{(\nu)})^2 
		\leq 2\overline{M}^{2} \inf_{\psi' \in \widetilde{\mathfrak{L}}_{p, {\mathbf q}}} \int_{B_{1/2}(0)} \sum_{j=1}^p \mathcal{G}(w^{(\nu)}_{j},\psi'_{j})^2,
\end{equation}
and either (i) $d(\psi^{(\nu)})  = p$ (i.e.\ $s=p$), or (ii) $d(\psi^{(\nu)}) > p$ (i.e.\ $s > p$) and 
\begin{equation} \label{holefilling0_hyp2-0}
	\int_{B_{1/2}(0)} \sum_{j=1}^p \mathcal{G}(w_{j}^{(\nu)},\psi_{j}^{(\nu)})^2 
		\leq \overline{\beta}_{\nu}^{2} \inf_{\psi' \in \bigcup_{s=p}^{d(\psi) -1} \mathfrak{L}_{p, {\mathbf q}}(s)} 
			\int_{B_{1/2}(0)} \sum_{j=1}^p\mathcal{G}(w_{j}^{(\nu)},\psi'_{j})^2, 
\end{equation}
and yet  
\begin{equation} \label{holefilling0_eqn1}
	\int_{B_{1/2}(0) \setminus B_{1/8}(0)} \sum_{j=1}^p 
		\left| \frac{\partial (w^{(\nu)}_{j}/R)}{\partial R} \right|^2 
	< \frac{1}{\nu} \int_{B_{1/2}(0)} \sum_{j=1}^p \mathcal{G}\left(w^{(\nu)}_{j},\psi^{(\nu)}_{j}\right)^2. 
\end{equation} 

Write
\begin{equation*}
	F_{\nu} = \left( \int_{B_{1/2}(0)} \sum_{j=1}^p  \mathcal{G}\left(w^{(\nu)}_{j},\psi^{(\nu)}_{j}\right)^2 \right)^{1/2} 
\end{equation*}
and note that $F_{\nu}>0$ by \eqref{holefilling0_eqn1}. To obtain a contradiction, we shall proceed in 5 steps. 

\noindent
{\bf Step 1:} \emph{Selection of currents associated with $w^{(\nu)},$ construction of cones associated with $\psi^{(\nu)}$ and some preliminary bounds}. Since $w^{(\nu)} \in \mathfrak{B}_{p, {\mathbf q}}$, there exist a sequence of locally area minimizing rectifiable currents $(T^{(\nu,k)})_{k = 1}^{\infty}$ in ${\mathbf B}_{1}(0)$ with 
$\partial \, T^{(\nu,k)} \llcorner {\mathbf B}_{1}(0)= 0$, a sequence of cones $({\mathbf C}^{(\nu,k)})_{k=1}^{\infty}$ in ${\mathcal C}_{p, {\mathbf q}}$  and sequences of positive numbers $(\epsilon^{(\nu, k)})_{k=1}^{\infty}$, $(\beta^{(\nu, k)})_{k=1}^{\infty},$ $(\eta^{(\nu, k)})_{k=1}^{\infty}$ and $(\delta^{(\nu, k)})_{k=1}^{\infty}$  tending to 0, all such that conditions (1)-(8) of Section~\ref{fine-blowup-prelim} and condition~\eqref{excess-1} hold with $T_{k} = T^{(\nu, k)}$, ${\mathbf C}_{k} = {\mathbf C}^{(\nu, k)}$, $p_{k}=p$, $q_{j}^{(k)} = q_{j}$, $\epsilon_{k} = \epsilon^{(\nu, k)}$, $\beta_{k} = \beta^{(\nu, k)}$, $\eta_{k} = \eta^{(\nu, k)}$, $\delta_{k} = \delta^{(\nu, k)}$, and such that $w^{(\nu)}$ is the fine blow-up (as in Definition~\ref{blowupclass_defn}) of $T^{(\nu,k)}$ relative to ${\mathbf C}^{(\nu,k)};$ denote the relevant excess by $E_{\nu,k} = E(T^{(\nu,k)},\mathbf{C}^{(\nu,k)},\mathbf{B}_1(0))$.  
 By the definition of ${\mathcal C}_{p, {\mathbf q}}$, we have that 
\begin{equation}\label{holefilling0_eqn2-0}
{\mathbf C}^{(\nu, k)} = \sum_{i=1}^{p} q_{i} \llbracket P_{i}^{(\nu, k)}\rrbracket
\end{equation}
where for each $\nu$ and each $k$,  
$P_{i}^{(\nu, k)}$ ($1 \leq i \leq p$) are distinct $n$-dimensional oriented planes with $\{0\} \times {\mathbb R}^{n-2}  = P_{i}^{(\nu, k)} \cap P_{j}^{(\nu, k)}$ whenever $i \neq j,$ and with orienting $n$-vector 
$\vec{P}_{i}^{(\nu, k)}$,  and moreover, $P_{i}^{(\nu, k)} = \{(z, x, y) \in {\mathbb R}^{n+m}\, : \, z = A_{i}^{(\nu, k)}x\}$ for some constant $m \times 2$ matrices $A_{i}^{(\nu, k)}.$  
By the definition of fine blow-up, for each $\nu$, there is a sequence of positive numbers $(\tau^{(\nu,k)})_{k=1}^{\infty}$ with $\tau^{(\nu, k)} \to 0$ and functions $u^{(\nu,k)}_{i} : \, B_{1/2}(0) \setminus \{r < \t_{\nu, k}\} \to {\mathcal A}_{q_{i}}({\mathbb R}^{m})$ for $i =1, 2, \ldots, p$ corresponding to the sequences $(\tau_{k})$, $(u^{(k)}_{i})$ respectively as in the discussion of Section~\ref{fine-blowup-prelim} 
 taken with $\gamma = 1/2$, $T = T^{(\nu,k)}$, and ${\mathbf C} = {\mathbf C}^{(\nu,k)}$; and moreover, for each $\nu =1, 2, \dots$ and each $j =1, 2, \ldots, p$, 
we have (by \eqref{excess-8-1}) that 
\begin{equation} \label{holefilling0_eqn2}
			 \frac{u^{(\nu,k)}_{i}}{E_{\nu,k}} \to w^{(\nu)}_{i}
\end{equation}
in $L^{2}(B_{1/2}; {\mathcal A}_{q_{i}}({\mathbb R}^{m}))$ and also that $E_{\nu,k}^{-1} |Du^{(\nu, k)}_{i}| \to |Dw^{(\nu)}_{i}|$ locally in $L^{2}(K)$ for each 
compact set $K \subset B_{1/2}(0) \setminus \{0\} \times {\mathbb R}^{n-2}$.

By the definition of ${\mathfrak L}_{p, {\mathbf q}}(s)$, we have that $d(\psi^{(\nu)}) = s$ and 
\begin{equation*}
    \psi^{(\nu)}_i(x,y) = \sum_{\ell=1}^{m_i^{(\nu)}} \widetilde{q}_{\ell}^{(\nu)} \llbracket \widetilde{\psi}^{(\nu)}_{i,\ell}(x) \rrbracket
\end{equation*}
where $m_i^{(\nu)}$ are positive integers with $\sum_{i=1}^{p} m^{(\nu)}_{i} = s$, $\widetilde{q}_{i, \ell}^{(\nu)}$ are positive integers with $\sum_{\ell=1}^{m^{(\nu)}_{i}} \widetilde{q}_{i, \ell}^{(\nu)} = q_{i}$ and $\widetilde{\psi}^{(\nu)}_{i,\ell} : \mathbb{R}^n \rightarrow \mathbb{R}^m$ are distinct linear functions.  By passing to a subsequence without changing notation, we may and shall assume that for each $i \in \{1, \ldots, p\}$, $m_{i}^{(\nu)} = m_{i}$  for some fixed integer $m_{i} \geq 1$ and all $\nu = 1, 2, \ldots$, and also for each $i \in \{1, \ldots, p\}$ and $\ell \in \{1, \ldots, m_{i}\}$ that $\widetilde{q}^{(\nu)}_{i, \ell} = \widetilde{q}_{i, \ell}$ some fixed integer $\widetilde{q}_{i, \ell} \geq 1$ and all $\nu =1, 2, \ldots.$  For each $\nu$ and $k$, let 
\begin{equation}\label{holefilling0_eqn2-1}
	\widetilde{\mathbf C}^{(\nu,k)} = \sum_{i=1}^{p} \sum_{\ell=1}^{m_i} \widetilde{q}_{\ell} \llbracket \widetilde{P}^{(\nu,k)}_{i,\ell} \rrbracket 
\end{equation}
where 
\begin{equation*}
	\widetilde{P}^{(\nu, k)}_{i, \ell} = \{(z, x, y) \in {\mathbb R}^{m} \times {\mathbb R}^{2} \times {\mathbb R}^{n-2}\, : \, z = A_{i}^{(\nu, k)} x + E_{\nu, k} \widetilde{\psi}_{i, \ell}^{(\nu)}(x)\}.\end{equation*}
(Thus $\widetilde{\mathbf C}^{(\nu, k)}$ is as constructed in Lemma~\ref{blowup2L_lemma} with ${\mathbf C}^{(\nu,k)}$, $\psi^{(\nu)}$, and 
$E_{\nu,k}$ in place of ${\mathbf C}_{k}$, $\psi$ and $E_{k}$).  Note that then, since $\psi^{(\nu)} \in {\mathfrak L}_{p, {\mathbf q}}(s)$, we have that $\widetilde{\mathbf C}^{(\nu, k)} \in {\mathcal C}_{s, {\mathbf q}^{\prime}}$ for some ${\mathbf q}^{\prime} \in {\mathfrak M}_{s},$ and by Lemma~\ref{blowup_norm_conv_lemma} we have
\begin{equation} \label{holefilling0_eqn3}
	F_{\nu}^2 \geq \lim_{k \rightarrow \infty} E_{\nu,k}^{-2} \int_{{\mathbf B}_{1/2}(0)} {\rm dist}^{2} \, (X, {\rm spt} \, \widetilde{\mathbf C}^{(\nu,k)}) d\|T^{(\nu, k)}\|\;\; \mbox{and}
\end{equation}

\begin{equation}\label{holefilling0_eqn3-1}
F_{\nu}^{2} \geq \lim_{k \to \infty} E_{\nu, k}^{-2}\int_{{\mathbf B}_{1/4}(0) \setminus \{r < 1/32\}} {\rm dist}^{2} \, (X, {\rm spt} \, T^{(\nu, k)}) d\|\widetilde{\mathbf C}^{(\nu, k)}\|.
\end{equation}

Observe that by taking $\psi' = 0$ in \eqref{holefilling0_hyp1-0} and applying  \eqref{blowup_norm_conv_eqn2-1} we have that  
\begin{equation}\label{holefilling0_eqn5-1}
	F_{\nu}^2 \leq 2\overline{M}^{2} \int_{B_{1/2}(0)} \sum_{j=1}^p  |w_{j}^{(\nu)}|^2 \leq 2\overline{M}^{2} 
\end{equation}
whence by \eqref{holefilling0_eqn3}, \eqref{holefilling0_eqn3-1}, for each $\nu$ and sufficiently large $k$ (depending on $\nu$), 
\begin{equation} \label{holefilling0_eqn6}
	Q^{2}(T^{(\nu, k)}, \widetilde{\mathbf C}^{(\nu, k)}, {\mathbf B}_{1/2}(0)) 
    \leq \omega_n^{-1} 2^{n+4} E_{\nu, k}^2 F_{\nu}^2 
    \leq \omega_n^{-1} 2^{n+5} \overline{M}^{2} E_{\nu, k}^2.
\end{equation}
We also have from \eqref{holefilling0_eqn5-1}, \eqref{blowup_norm_conv_eqn2-1} and the triangle inequality that 
\begin{equation}\label{holefilling0_eqn5-2}
\int_{B_{1/2}(0)} \sum_{j=1}^{p} |\psi_{j}^{(\nu)}|^{2} \leq 4\overline{M}^{2} +2.
\end{equation}

\noindent\textbf{Step 2:}  We claim the following: \emph{if $s >p$ (possible only if $2 \leq p < q$) then for each sufficiently large $\nu$ and sufficiently large $k$ (depending on $\nu$),}
\begin{equation}\label{eqn1}
	Q(T^{(\nu,k)}, \widetilde{\mathbf C}^{(\nu,k)}, \mathbf{B}_{1/2}(0)) \leq C_1\overline{\beta}_{\nu} \inf_{\mathbf{C}' \in \bigcup_{s'=p}^{s-1} \mathcal{C}_{q,s'}} 
		Q(\widetilde{\mathbf C}^{(\nu,k)}, \mathbf{C}', \mathbf{B}_{1/2}(0)) ,
\end{equation}
\emph{where $C_1 = C_1(n,m,q,\overline{M}) \in (0,\infty).$} To see this, with $C_{1} = C_{1}(n, m, q, \overline{M})$ to be determined suppose to the contrary that 
\begin{equation}\label{eqn2}
	Q(T^{(\nu,k)}, \widetilde{\mathbf C}^{(\nu,k)}, \mathbf{B}_{1/2}(0)) > C_1 \overline{\beta}_{\nu} \inf_{\mathbf{C}' \in \bigcup_{s'=p}^{s-1} \mathcal{C}_{q,s'}} 
		Q(\widetilde{\mathbf C}^{(\nu,k)}, \mathbf{C}', \mathbf{B}_{1/2}(0)) .
\end{equation}
For each $\nu$ and $k$ choose $\widehat{s}^{(\nu,k)} \in \{p,\ldots,s-1\}$ and a cone $\widehat{\mathbf C}^{(\nu,k)} \in \mathcal{C}_{q,\widehat{s}^{(\nu,k)}}$ such that 
\begin{equation}\label{eqn3}
	Q(\widetilde{\mathbf C}^{(\nu,k)}, \widehat{\mathbf C}^{(\nu,k)}, \mathbf{B}_{1/2}(0)) < 2 \inf_{\mathbf{C}' \in \bigcup_{s'=p}^{s-1} \mathcal{C}_{q,s'}}
		Q(\widetilde{\mathbf C}^{(\nu,k)}, \mathbf{C}', \mathbf{B}_{1/2}(0)) . 
\end{equation}
Taking $\mathbf{C}' = \mathbf{C}^{(\nu,k)}$ in \eqref{eqn3} gives us, in view of the definition \eqref{holefilling0_eqn2-1} of $\widetilde{\mathbf C}^{(\nu,k)},$ and \eqref{holefilling0_eqn6} that 
\begin{equation*}
	Q(\widetilde{\mathbf C}^{(\nu,k)}, \widehat{\mathbf C}^{(\nu,k)}, \mathbf{B}_{1/2}(0)) 
	< 2 Q(\widetilde{\mathbf C}^{(\nu,k)}, \mathbf{C}^{(\nu,k)}, \mathbf{B}_{1/2}(0)) 
	\leq C E_{\nu,k}, 
\end{equation*}
where $C = C(n,q,\overline{M}) \in (0,\infty)$ is a constant.  Since $\widetilde{\mathbf C}^{(\nu,k)}, \widehat{\mathbf C}^{(\nu,k)}, \mathbf{C}^{(\nu,k)}$ are supported on unions of distinct planes intersecting along $\{0\} \times \mathbb{R}^{n-2}$, it follows that 
\begin{equation*}
	\op{dist}_{\mathcal H}(\op{spt} \widehat{\mathbf C}^{(\nu,k)} \cap \mathbf{B}_1(0), \op{spt} \mathbf{C}^{(\nu,k)} \cap \mathbf{B}_1(0)) \leq C E_{\nu,k} 
\end{equation*}
for some constant $C = C(n,m,q,\overline{M}) \in (0,\infty)$.  Recalling that by Theorem~\ref{graphrep thm}(a) and condition (4) of Section~\ref{fine-blowup-prelim} we have $\text{minsep} \, \mathbf{C}^{(\nu,k)} \geq c(n,m,q) \,(\beta^{(\nu,k)})^{-1} E_{\nu,k}$, we can express $\widehat{\mathbf C}^{(\nu,k)}$ as 
\begin{equation*}
	\widehat{\mathbf C}^{(\nu,k)} = \sum_{i=1}^p \sum_{\ell=1}^{\widehat{s}^{(\nu,k)}_i} \widehat{q}^{(\nu,k)}_{\ell} \llbracket \widehat{P}_{i,\ell}^{(\nu,k)} \rrbracket 
\end{equation*}
where $\widehat{s}^{(\nu,k)}_i$ and $\widehat{q}^{(\nu,k)}_i$ are positive integers such that $\sum_{i=1}^p \widehat{s}^{(\nu,k)}_i = \widehat{s}^{(\nu,k)}$ and $\sum_{i=1}^p \sum_{\ell=1}^{\widehat{s}^{(\nu,k)}_i} \widehat{q}_{\ell}^{(\nu,k)} = q$ and where $\widehat{P}_{i,\ell}^{(\nu,k)}$ are distinct oriented planes with 
\begin{equation*}
	\op{dist}_{\mathcal H}(\widehat{P}_{i,\ell}^{(\nu,k)} \cap \mathbf{B}_1(0), P_i^{(\nu,k)} \cap \mathbf{B}_1(0)) \leq C E_{\nu,k} 
\end{equation*}
for each $i \in \{1,\ldots,p\}$ and $\ell \in \{1,\ldots,\widehat{s}^{(\nu,k)}_i\}$.  (Note that we do not claim that $\sum_{\ell=1}^{\widehat{s}^{(\nu,k)}_i} \widehat{q}_{\ell}^{(\nu,k)} = q_i.$)  Moreover, for each $i \in \{1,\ldots,p\}$ and $\ell \in \{1,\ldots,\widehat{s}^{(\nu,k)}_i\}$ there exists a linear function $\widehat{\psi}_{i,\ell}^{(\nu,k)} : \mathbb{R}^2 \rightarrow \mathbb{R}^n$ such that 
\begin{equation*}
	\widehat{P}_{i,\ell}^{(\nu,k)} = \{(z,x,y) : z = A_i^{(\nu,k)} x + E_{\nu,k} \widehat{\psi}_{i,\ell}^{(\nu,k)}(x)\} . 
\end{equation*}
Since $\widehat{\mathbf C}^{(\nu,k)}$ has fewer planes than $\widetilde{\mathbf C}^{(\nu,k)}$, for each $x \in \mathbb{R}^2$ there exists $i \in \{1,\ldots,p\}$, $\ell_1,\ell_2 \in \{1,\ldots,m_i\}$, and $\widehat{\ell} \in \{1,\ldots,\widehat{s}^{(\nu,k)}_i\}$ such that $\ell_1 \neq \ell_2$ and 
\begin{equation}\label{eqn4}
	E_{\nu,k} |\widetilde{\psi}_{i,\ell_j}^{(\nu)}(x) - \widehat{\psi}_{i,\widehat{\ell}}^{(\nu,k)}(x)| \leq 2 \op{dist}((A_i^{(\nu,k)} x + E_{\nu,k} \widetilde{\psi}_{i,\ell_j}^{(\nu)}(x),x,0), \op{spt} \widehat{\mathbf C}^{(\nu,k)}) 
    \text{ for } j = 1,2.
\end{equation}
Thus for some $i \in \{1,\ldots,p\}$, $\ell_1,\ell_2 \in \{1,\ldots,m_i\}$, and $\widehat{\ell} \in \{1,\ldots,\widehat{s}^{(\nu,k)}_i\}$ we have that $\ell_1 \neq \ell_2$ and the set $S = S_{i,\ell_1,\ell_2,\widehat{\ell}}$ of all $x \in B^2_{1/4}(0)$ for which \eqref{eqn4} holds true satisfies $\mathcal{L}^2(S) \geq \pi/(32q^3)$.  By \eqref{eqn4} for each $x \in S$ and $y \in \mathbb{R}^{n-2}$ 
\begin{align*}
	E_{\nu,k} |\widetilde{\psi}_{i,\ell_1}^{(\nu)}(x) - \widetilde{\psi}_{i,\ell_2}^{(\nu)}(x)| \leq\,& E_{\nu,k} |\widetilde{\psi}_{i,\ell_1}^{(\nu)}(x) - \widehat{\psi}_{i,\widehat{\ell}}^{(\nu,k)}(x)| 
		+ E_{\nu,k} |\widetilde{\psi}_{i,\ell_2}^{(\nu)}(x) - \widehat{\psi}_{i,\widehat{\ell}}^{(\nu,k)}(x)| 
	\\ \leq\,& 2 \sum_{\ell=1}^{q_i} \op{dist}((A_i^{(\nu,k)} x + E_{\nu,k} \widetilde{\psi}_{i,\ell}^{(\nu)}(x),x,y), \op{spt} \widehat{\mathbf C}^{(\nu,k)}) .
\end{align*}
Hence squaring both sides and integrating over $(x,y) \in S \times B^{n-2}_{1/4}(0)$, 
\begin{equation*}
	E_{\nu,k}^2 \int_{B^{n-2}_{1/4}(0)} \int_S |\widetilde{\psi}^{(\nu)}_{i,\ell_1}(x) - \widetilde{\psi}^{(\nu)}_{i,\ell_2}(x)|^2 \,dx \,dy 
	\leq 4q \int_{\mathbf{B}_1(0)} \op{dist}^2(X, \op{spt} \widehat{\mathbf C}^{(\nu,k)}) \,d\|\widetilde{\mathbf C}^{(\nu,k)}\|(X) .
\end{equation*}
Since $\widetilde{\psi}^{(\nu)}_{i,\ell_1} - \widetilde{\psi}^{(\nu)}_{i,\ell_2}$ is a linear function, there exists unit vectors $v_1,v_2 \in \mathbb{R}^2$ such that \eqref{eqn4} holds true for $x = v_1,v_2$ and 
\begin{equation*}
	\frac{\pi \omega_{n-2}}{4^{n+1} q^3} \,E_{\nu,k}^2 |\widetilde{\psi}^{(\nu)}_{i,\ell_1}(x) - \widetilde{\psi}^{(\nu)}_{i,\ell_2}(x)|^2 
	\leq 4q \int_{\mathbf{B}_1(0)} \op{dist}(X, \op{spt} \widehat{\mathbf C}^{(\nu,k)}) \,d\|\widetilde{\mathbf C}^{(\nu,k)}\|(X) 
\end{equation*}
for $x = v_1,v_2$ and the angle between $v_1,v_2$ is in $[\pi/(4q^3), \pi-\pi/(4q^3)]$.  Hence 
\begin{align}\label{eqn5}
	E_{\nu,k}^2 \int_{B_{1/2}(0)} |\widetilde{\psi}^{(\nu)}_{i,\ell_1} - \widetilde{\psi}^{(\nu)}_{i,\ell_2}|^2 
	\leq\,& C E_{\nu,k}^2 \sup_{B_{1/2}(0)} |\widetilde{\psi}^{(\nu)}_{i,\ell_1} - \widetilde{\psi}^{(\nu)}_{i,\ell_2}|^2 
	\\ \leq\,& C E_{\nu,k}^2 \int_{\mathbf{B}_1(0)} \op{dist}(X, \op{spt} \widehat{\mathbf C}^{(\nu,k)}) \,d\|\widetilde{\mathbf C}^{(\nu,k)}\|(X) \nonumber 
\end{align}
where $C = C(n,q) \in (0,\infty)$ are constants.  On the other hand, setting $\psi'_i(x) = \widetilde{q}_{i,\ell_2} \llbracket \widetilde{\psi}_{i,\ell_1}(x) \rrbracket + \sum_{\ell\neq\ell_2} \widetilde{q}_{i,\ell} \llbracket \widetilde{\psi}_{i,\ell}(x) \rrbracket$ and $\psi'_j(x) = \psi_j(x)$ if $j \neq i$, 
\begin{equation}\label{eqn6}
	\inf_{\psi' \in \bigcup_{s'=p}^{s-1} \mathfrak{L}_{p,\mathbf{q}}(s')} \int_{B_{1/2}(0)} \sum_{i=1}^p \mathcal{G}(\psi^{(\nu)}_i,\psi'_i)^2 
	\leq q \int_{B_{1/2}(0)} |\widetilde{\psi}^{(\nu)}_{i,\ell_1} - \widetilde{\psi}^{(\nu)}_{i,\ell_2}|^2 . 
\end{equation}
Combining \eqref{eqn2}, \eqref{eqn3}, \eqref{eqn5} and \eqref{eqn6} gives
\begin{equation}\label{eqn7}
	Q^2(T^{(\nu,k)}, \widetilde{\mathbf C}^{(\nu,k)}, \mathbf{B}_{1/2}(0)) 
	> c \,C_1^2 \overline{\beta}_{\nu}^2 E_{\nu,k}^2 \inf_{\psi' \in \bigcup_{s'=p}^{s-1} \mathfrak{L}_{p,\mathbf{q}}(s')} 
		\int_{B_{1/2}(0)} \sum_{i=1}^p \mathcal{G}(\psi^{(\nu)}_i,\psi'_i)^2 
\end{equation}
for some constant $c = c(n,m,q,\overline{M}) > 0$.  Dividing both sides of \eqref{eqn7} by $E_{\nu,k}^2$ and letting $k \rightarrow \infty$ using \eqref{holefilling0_eqn3} and \eqref{holefilling0_eqn3-1} gives
\begin{equation*}
	\int_{B_{1/2}(0)} \sum_{i=1}^p \mathcal{G}(w^{(\nu)}_i,\psi^{(\nu)}_i)^2
	\geq c \,C_1^2 \overline{\beta}_{\nu}^2 \inf_{\psi' \in \bigcup_{s'=p}^{s-1} \mathfrak{L}_{p,\mathbf{q}}(s')} 
		\int_{B_{1/2}(0)} \sum_{i=1}^p \mathcal{G}(\psi^{(\nu)}_i,\psi'_i)^2 
\end{equation*}
for some constant $c = c(n,m,q,\overline{M}) > 0$.  Hence by the triangle inequality  
\begin{equation}\label{eqn8}
	\int_{B_{1/2}(0)} \sum_{i=1}^p \mathcal{G}(w^{(\nu)}_i,\psi^{(\nu)}_i)^2
	\geq \frac{c \,C_1^2 \overline{\beta}_{\nu}^2}{4} \inf_{\psi' \in \bigcup_{s'=p}^{s-1} \mathfrak{L}_{p,\mathbf{q}}(s')} 
		\int_{B_{1/2}(0)} \sum_{i=1}^p \mathcal{G}(w^{(\nu)}_i,\psi'_i)^2 
\end{equation}
for all sufficiently large $\nu$.  Choosing $C_1 > 2c^{-1/2}$, \eqref{eqn8} contradicts \eqref{holefilling0_hyp2-0}.  Therefore, \eqref{eqn1} must hold true.

\noindent\textbf{Step 3:}  Next we claim the following: \emph{for each sufficiently large $\nu$ and sufficiently large $k$ (depending on $\nu$)}
\begin{equation}\label{eqn9}
	Q(T^{(\nu,k)}, \widetilde{\mathbf C}^{(\nu,k)}, \mathbf{B}_{1/2}(0)) \leq C_2\overline{\beta}_{\nu} \inf_{\mathbf{C}' \in \bigcup_{s'=1}^{s-1} \mathcal{C}_{q,s'}} 
		Q(T^{(\nu, k)}, \mathbf{C}', \mathbf{B}_{1/2}(0)) ,
\end{equation}
where $C_{2} = C_{2}(n,m,q,\overline{M}) \in (0,\infty)$ is a constant.  To see this, note first that by $F_{\nu} > 0$ and \eqref{holefilling0_hyp1-0} (with $\psi' = 0$), $w^{(\nu)}$ is not identically zero in $B_{1/2}(0)$.  Hence by Remark~\ref{rescaling-fine-blowup} for each fixed $\nu$ 
 $$Q(T^{(\nu,k)}, {\mathbf C}^{(\nu,k)}, \mathbf{B}_{1/2}(0)) \leq C\beta^{(\nu, k)} \inf_{\mathbf{C}' \in \bigcup_{s'=1}^{p-1} \mathcal{C}_{q,s'}} Q(T^{(\nu, k)}, \mathbf{C}', \mathbf{B}_{1/2}(0))$$
where $\b^{(\nu, k)}$ is as in Step~1 (so that $\b^{(\nu, k)} \to 0$ as $k \to \infty$) and $C = C(n) \in (0,\infty)$ is a constant.  Hence in case $s = p$, \eqref{eqn9} (with $C_{2} = 1$ and for sufficiently large $k$ depending on $\nu$) follows from this and the definition \eqref{holefilling0_eqn2-1} of $\widetilde{\mathbf C}^{(\nu, k)}.$ So suppose that $s >p$. Let $\beta_{\star\star} = \beta_{\star\star}(n,m,q) = \min\{\beta_0(n,m,$ $q,p,3/4,1/8) : p = 2,\ldots,q\}$ where $\beta_0$ is as in Theorem~\ref{graphrep close2plane thm}.  Arguing as in Remark~\ref{tildeC rmk}, we can find $\overline{s}^{(\nu,k)} \in \{p,\ldots,s\}$ and a cone $\overline{\mathbf C}^{(\nu,k)} \in \mathcal{C}_{q,\overline{s}^{(\nu,k)}}$ such that 
\begin{gather}
	\label{eqn10} Q(T^{(\nu,k)}, \overline{\mathbf C}^{(\nu,k)}, \mathbf{B}_{1/2}(0)) \leq 2^{q-1} \beta_{\star\star}^{2-q} 
		\inf_{\mathbf{C}' \in \bigcup_{s'=1}^{s-1} \mathcal{C}_{q,s'}} Q(T^{(\nu,k)}, \mathbf{C}', \mathbf{B}_{1/2}(0)) , \\
	\label{eqn11} Q(T^{(\nu,k)}, \overline{\mathbf C}^{(\nu,k)}, \mathbf{B}_{1/2}(0)) \leq 
		\beta_{\star\star} \inf_{\mathbf{C}' \in \bigcup_{s'=1}^{\overline{s}^{(\nu,k)}-1} \mathcal{C}_{q,s'}} Q(T^{(\nu,k)}, \mathbf{C}', \mathbf{B}_{1/2}(0)) .
\end{gather}
(Note that by setting $\mathbf{C}' = \mathbf{C}^{(\nu,k)}$ in \eqref{eqn10}, $Q(T^{(\nu,k)}, \overline{\mathbf C}^{(\nu,k)}, \mathbf{B}_{1/2}(0)) \leq C E_{\nu,k}$.  Thus by condition~(4) of Section~\ref{fine-blowup-prelim}, we must have that $\overline{s}^{(\nu,k)} \geq p$.)  In light of condition~(3) of 
Section~\ref{fine-blowup-prelim} and \eqref{eqn11}, we can apply Theorem~\ref{graphrep close2plane thm}  to deduce that 
\begin{align}\label{eqn12}
	&\op{dist}_{\mathcal H}(\op{spt} T^{(\nu,k)} \cap \mathbf{B}_{1/4}(0) \cap \{r \geq 1/16\}, 
		\op{spt} \overline{\mathbf C}^{(\nu,k)} \cap \mathbf{B}_{1/4}(0) \cap \{r \geq 1/16\}) 
	\\ \leq\,& C E(T^{(\nu,k)}, \overline{\mathbf C}^{(\nu,k)}, \mathbf{B}_{1/2}(0)) \nonumber 
\end{align}
for some constant $C = C(n,m,q) \in (0,\infty)$.  Since $\widetilde{\mathbf C}^{(\nu,k)}, \overline{\mathbf C}^{(\nu,k)}$ are supported on unions of distinct planes intersecting along $\{0\} \times \mathbb{R}^{n-2}$, 
\begin{equation*}
	Q^2(\widetilde{\mathbf C}^{(\nu,k)}, \overline{\mathbf C}^{(\nu,k)}, \mathbf{B}_{1/2}(0)) 
	\leq C \int_{\mathbf{B}_{1/4}(0) \cap \{r > 1/16\}} \op{dist}^2(X, \op{spt} \overline{\mathbf C}^{(\nu,k)}) \,d\|\widetilde{\mathbf C}^{(\nu,k)}\|(X) 
\end{equation*}
for some constants $C = C(n,m,q) \in (0,\infty)$.  Noting that by the triangle inequality 
\begin{align*}
	&\op{dist}(X, \op{spt} \overline{\mathbf C}^{(\nu,k)}) 
	\leq \op{dist}(X, \op{spt} T^{(\nu,k)}) \\&\hspace{15mm} + \op{dist}_{\mathcal H}(\op{spt} T^{(\nu,k)} \cap \mathbf{B}_{1/4}(0) \cap \{r \geq 1/16\}, 
		\op{spt} \overline{\mathbf C}^{(\nu,k)} \cap \mathbf{B}_{1/4}(0) \cap \{r \geq 1/16\}) 
\end{align*}
for each $X \in \op{spt} \widetilde{\mathbf C}^{(\nu,k)} \cap \mathbf{B}_{1/4}(0) \cap \{r \geq 1/16\}$ and hence using \eqref{eqn12}
\begin{align}\label{eqn13}
	&Q^2(\widetilde{\mathbf C}^{(\nu,k)}, \overline{\mathbf C}^{(\nu,k)}, \mathbf{B}_{1/2}(0)) 
	\\ \leq\,& C \int_{\mathbf{B}_{1/4}(0) \cap \{r > 1/16\}} \op{dist}^2(X, \op{spt} T^{(\nu,k)}) \,d\|\widetilde{\mathbf C}^{(\nu,k)}\|(X) 
		+ C E^2(T^{(\nu,k)}, \overline{\mathbf C}^{(\nu,k)}, \mathbf{B}_{1/2}(0)) \nonumber
	\\ \leq\,& C Q^2(T^{(\nu,k)}, \widetilde{\mathbf C}^{(\nu,k)}, \mathbf{B}_{1/2}(0)) + C Q^2(T^{(\nu,k)}, \overline{\mathbf C}^{(\nu,k)}, \mathbf{B}_{1/2}(0))  \nonumber
\end{align}
for some constants $C = C(n,m,q) \in (0,\infty)$.  Thus combining \eqref{eqn1}, \eqref{eqn13} and \eqref{eqn10}
\begin{align*}
	Q(T^{(\nu,k)}, \widetilde{\mathbf C}^{(\nu,k)}, \mathbf{B}_{1/2}(0)) 
	\leq\,& C \overline{\beta}_{\nu} Q(T^{(\nu,k)}, \widetilde{\mathbf C}^{(\nu,k)}, \mathbf{B}_{1/2}(0))
	\\&+ C \overline{\beta}_{\nu} \inf_{\mathbf{C}' \in \bigcup_{s'=1}^{s-1} \mathcal{C}_{q,s'}} Q(T^{(\nu,k)}, \mathbf{C}', \mathbf{B}_{1/2}(0))
\end{align*}
for some constant $C = C(n,m,q,\overline{M}) \in (0,\infty)$, which, provided $C\overline{\beta}_{\nu} < 1/2,$ proves \eqref{eqn9}.

\noindent
{\bf Step 4:} \emph{Fine blow up relative to the cones associated with $\psi^{(\nu)}$.} Using \eqref{eqn9} and the fact that condition (6) in Section~\ref{fine-blowup-prelim} holds with $T^{(\nu, k)}$ in place $T,$ we can now verify (by arguing as in \cite[pp.\ 910-914]{Wic14}) the following: for sufficiently large $\nu$ and $k$, \eqref{eqn9} holds with $\eta_{0, 1/2 \, \#} \, T^{(\nu, k)}$ in place of $T^{(\nu, k)}$ and with a larger but fixed constant (depending only on $n$, $m$, $q$, $\overline{M}$) in place of $C_{2}$, and condition (6) in Section~\ref{fine-blowup-prelim} holds with $\eta_{0, 1/2 \, \#} \, T^{(\nu, k)}$ in place of $T$ and $CM$ in place of $M$, where $C = C(n, m, q)$.   
This allows us to apply Corollary~\ref{nonconcentration close2plane cor} with $\eta_{0, 1/2 \#} \, T^{(\nu, k)}$ in place of $T$ and 
$\widetilde{\mathbf C}^{(\nu, k)}$ in place of ${\mathbf C}$ to deduce that for any $\delta \in (0, 1/2)$, 
\begin{equation} \label{holefilling0_eqn19-1} 
	\int_{\mathbf{B}_{1/4}(0) \cap \{r < \delta\}} \op{dist}^2(X, \widetilde{\mathbf C}^{(\nu, k)}) \,d\|T^{(\nu, k)}\|(X) 
		\leq C_{1} \delta \int_{\mathbf{B}_{1/2}(0)} \op{dist}^2(X, \widetilde{\mathbf C}^{(\nu, k)}) \,d\|T^{(\nu, k)}\|(X)  
\end{equation}
for all sufficiently large $\nu$ and $k$ (depending on $\delta$), where $C_{1} = C_{1}(n, m, q, M)$. Using Theorem~\ref{graphrep close2plane thm} (applied to the left hand side of this), dividing both sides by $E_{\nu, k}^{2}$ and letting $k \to \infty$ (for fixed $\nu$), we obtain with the help of Lemma~\ref{blowup_norm_conv_lemma}
(applied to the right hand side) that for all sufficiently large $\nu$, 
\begin{equation} \label{holefilling0_eqn19-2} 
	\int_{B_{1/4}(0) \cap \{r < \delta\}} \sum_{i=1}^{p}{\mathcal G}(w^{(\nu)}_{i}, \psi^{(\nu)}_{i})^{2} 
	\leq C_{1} \delta \int_{B_{1/2}(0)} \sum_{i=1}^{p}{\mathcal G}(w^{(\nu)}_{i}, \psi^{(\nu)}_{i})^{2}  =C_{1}\delta F_{\nu}^{2}
\end{equation}
where $C_{1} = C_{1}(n, m, q, M).$

In light of \eqref{eqn9}, there exist sequences $(\widetilde{\tau}_{\nu})_{\nu=1}^{\infty}$ and $(\widetilde{\gamma}_{\nu})_{\nu=1}^{\infty}$ such that $0 < \widetilde{\tau}_{\nu} < \widetilde{\gamma}_{\nu} < 1$ and $\widetilde{\tau}_{\nu} \rightarrow 0$ and $\widetilde{\gamma}_{\nu} \rightarrow 1$ as $\nu\rightarrow\infty$ and for each sufficiently large $\nu$ there exists a Lipschitz $\widetilde{q}_{i,l}$-valued functions $\widetilde{u}^{(\nu,k)}_{i, \ell} : \, B_{\gamma}(0) \setminus \{r < \widetilde{\tau}_{\nu}\} \to {\mathcal A}_{\widetilde{q}_{i, \ell}}({\mathbb R}^{m})$ and closed sets $\widetilde{K}^{(\nu,k)} \subset B_{\widetilde{\gamma}_{\nu}}(0) \setminus \{r < \widetilde{\tau}_{\nu}\}$ ($1 \leq i \leq p, \;\; 1 \leq \ell \leq m_{i}$) such that for each sufficiently large $k$ conclusion~(C) of Theorem~\ref{graphrep close2plane thm} holds true with $\widetilde{\gamma}_{\nu},\widetilde{\tau}_{\nu},\eta_{0,1/2\#} T^{(\nu,k)}, \widetilde{\mathbf C}^{(\nu,k)}, \widetilde{u}^{(\nu,k)}_{i,\ell},\widetilde{K}^{(\nu,k)}$ in place of $\gamma,\tau,T,\mathbf{C},u_i,K$.  By Theorem~\ref{graphrep close2plane thm} and the definition of $\widetilde{\mathbf C}^{(\nu)}$ (as in \eqref{holefilling0_eqn2-1}), we have that for each $i \in \{1, \ldots, p\}$, $\ell \in \{1, \ldots, m_{i}\}$ and $(x,y) \in B_{{\gamma}}(0) \cap \widetilde{K}^{(\nu,k)} \cap (2K^{(\nu,k)})$, 
\begin{equation}\label{holefilling0_eqn28}
    2u_{i}^{(\nu,k)}(x/2,y/2) = \sum_{\ell=1}^{m_{i}} \sum_{j=1}^{\widetilde{q}_{i, \ell}} \llbracket E_{\nu,k} \widetilde{\psi}_{i, \ell}^{(\nu)}(x) + \widetilde{u}^{(\nu,k)}_{i, \ell, j}(x,y) \rrbracket
\end{equation}
where $K^{(\nu,k)}$ is the set corresponding to $K$  in Theorem~\ref{graphrep close2plane thm} with $1/2,\tau_{\nu,k},T^{(\nu)},{\mathbf C}^{(\nu)}$ in place of $\gamma,\tau,T,\mathbf{C}$.  By Theorem~\ref{graphrep close2plane thm}, Lemma~\ref{energy est lemma} and \eqref{holefilling0_eqn3}, for each $0 < \tau < \gamma < 1$
\begin{equation}\label{holefilling0_eqn34}
    \sup_{B_{\gamma}(0) \setminus \{r \leq \tau\}} |\widetilde{u}^{(\nu,k)}_{i,\ell}|^2 + \int_{B_{\gamma}(0) \setminus \{r \leq \tau\}} |D\widetilde{u}^{(\nu,k)}_{i,\ell}|^2 \leq C E(T^{(\nu,k)},\widetilde{\mathbf C}^{(\nu,k)},\mathbf{B}_{1/2}(0))^2 \leq C E_{\nu,k}^2 F_{\nu}^2 
\end{equation}
whenever $\widetilde{\tau}_{\nu} < \tau < \gamma < \widetilde{\gamma}_{\nu}$ and $k$ is sufficiently large, where $C = C(n,m,q,\gamma,\tau,\overline{M}) \in (0,\infty)$ are constants.  Hence by the Rellich compactness lemma for multi-valued Sobolev functions~\cite[Proposition~2.11]{DeLSpaDirMin} and Theorem~\ref{harmonic thm}, after passing to a subsequence there exist locally Dirichlet energy minimizing functions $\overline{w}^{(\nu)}_{i, \ell} : \, B_{\widetilde{\gamma}_{\nu}}(0) \setminus \{r \leq \widetilde{\tau}_{\nu}\} \to {\mathcal A}_{\widetilde{q}_{i, \ell}}({\mathbb R}^{m})$ such that 
\begin{equation}\label{holefilling0_eqn35}
    \frac{\widetilde{u}^{(\nu,k)}_{i,\ell}}{E_{\nu,k}} \rightarrow \overline{w}^{(\nu)}_{i,\ell} 
\end{equation}
in $L^2(B_{\widetilde{\gamma}_{\nu}}(0) \setminus \{r \leq \widetilde{\tau}_{\nu}\};\mathcal{A}_{\widetilde{q}_{i,\ell}}(\mathbb{R}^m))$ and also $E_{\nu,k}^{-1} |D\widetilde{u}^{(\nu,k)}_{i,\ell}| \rightarrow |D\overline{w}^{(\nu)}_{i,\ell}|$ locally in $L^2(K)$ for each compact set $K \subset B_{\widetilde{\gamma}_{\nu}}(0) \setminus \{r \leq \widetilde{\tau}_{\nu}\}$ as $k\rightarrow\infty$.  By dividing both sides of \eqref{holefilling0_eqn28} by $E_{\nu,k}$ and letting $k \to \infty$ 
\begin{equation}\label{holefilling0_eqn36}
    2w_{i}^{(\nu)}(x/2,y/2) = \sum_{\ell=1}^{m_{i}} \sum_{j=1}^{\widetilde{q}_{i, \ell}} \llbracket \widetilde{\psi}_{i, \ell}^{(\nu)}(x) + \overline{w}^{(\nu)}_{i, \ell, j}(x,y) \rrbracket
\end{equation}
for all $(x,y) \in B_{1/2}(0) \setminus \{r \leq \widetilde{\tau}_{\nu}\}$ and $i \in \{1,2,\ldots,p\}$.  By dividing both sides of \eqref{holefilling0_eqn34} by $E_{\nu,k}$ and letting $k \to \infty$, for each $0 < \tau < \gamma < 1$
\begin{equation*}
    \sup_{B_{\gamma}(0) \cap \{r > \tau\}} |\overline{w}^{(\nu)}_{i,\ell}| \leq C F_{\nu} 
\end{equation*}
whenever $\widetilde{\tau}_{\nu} < \tau < \gamma < \widetilde{\gamma}_{\nu}$, where $C = C(n,m,q,\gamma,\tau,\overline{M}) \in (0,\infty)$ is a constant.  Hence by applying the compactness of Dirichlet energy minimizing multi-valued functions, after passing to a subsequence there exist locally Dirichlet energy minimizing functions $\overline{w}_{i, \ell} : \, B_1(0) \setminus \{0\} \times \mathbb{R}^{n-2} \to {\mathcal A}_{\widetilde{q}_{i, \ell}}({\mathbb R}^{m})$ such that 
\begin{equation}\label{holefilling0_eqn37}
    \frac{\overline{w}^{(\nu)}_{i,\ell}}{F_{\nu}}  \rightarrow \overline{w}_{i,\ell} 
\end{equation}
uniformly on $K$ for each compact set $K \subset B_1(0) \setminus \{0\} \times \mathbb{R}^{n-2}$ as $\nu\to\infty$.  We let $\overline{w}_{i,\ell}(x,y) = \sum_{j=1}^{\widetilde{q}_{i,\ell}} \llbracket \overline{w}_{i,\ell,j}(x,y)\rrbracket$ for each $(x,y) \in B_1(0) \setminus \{0\} \times \mathbb{R}^{n-2}$.

Select diagonal sequences $T^{(\nu,k(\nu))}$, ${\mathbf C}^{(\nu,k(\nu))}$   by choosing $k = k(\nu)$ large enough such that conditions (1)-(8) of Section~\ref{fine-blowup-prelim} and condition~\eqref{excess-1} are satisfied with $T^{(\nu, k(\nu))}$, ${\mathbf C}^{(\nu, k(\nu))}$ in place of 
$T_{k}$, ${\mathbf C}_{k}$; with $p_{k}=p$, $q_{j}^{(k)} = q_{j}$ and with $\epsilon^{(\nu, k(\nu))}$, $\beta^{(\nu, k(\nu))}$, $\eta^{(\nu, k(\nu))}$, $\delta^{(\nu, k(\nu))}$ $\downarrow 0$ (as $\nu 
\to \infty$) in place of $\epsilon_{k}$, $\beta_{k}$, $\eta_{k}$, $\delta_{k}$ respectively. Furthermore, by taking $\nu$ large enough, in view of \eqref{holefilling0_eqn2}, \eqref{holefilling0_eqn3}, \eqref{holefilling0_eqn3-1} and the claim in Step~3, we may, and shall, require also that
\begin{gather}
	\int_{B_{1/2}(0) \setminus \{r < \nu^{-1}\}} \sum_{i=1}^p \mathcal{G}\left( \frac{u^{(\nu,k(\nu))}_{i}}{E_{\nu,k(\nu)}}, w^{(\nu)}_{i} \right)^{2}
    < \frac{1}{\nu} F_{\nu}^{2}; \label{holefilling0_eqn20}\\
	\lim_{\nu \rightarrow \infty} \frac{1}{E_{\nu,k(\nu)}^2 F_{\nu}^2} \int_{{\mathbf B}_{1/2}(0)} 
		{\rm dist}^{2} \, (X, {\rm spt} \, \widetilde{\mathbf C}^{(\nu,k(\nu))}) \, d\|T^{(\nu, k(\nu))}\| \leq 1; \label{holefilling0_eqn21}\\
		\lim_{\nu \to \infty} \frac{1}{E_{\nu, k(\nu)}^{2}F_{\nu}^{2}}\int_{{\mathbf B}_{1/4}(0) \setminus \{r < 1/32\}} {\rm dist}^{2} \, (X, {\rm spt} \, T^{(\nu, k(\nu))}) \, d\|\widetilde{\mathbf C}^{(\nu, k(\nu))}\| \leq 1 \;\; \mbox{and}\label{holefilling0_eqn22}
\end{gather}
\begin{eqnarray} \label{holefilling0_eqn23}
&&\mbox{either (i) $s = p$ or (ii) $s > p$ and}\\
&&	Q(T^{(\nu, k(\nu))}, \widetilde{\mathbf C}^{(\nu, k(\nu))}, {\mathbf B}_{1/2}(0)) 
		\leq C\overline{\beta}_{\nu} \inf_{{\mathbf C}^{\prime} \in \bigcup_{s^{\prime}=1}^{s-1} \bigcup_{{\mathbf q}^{\prime} \in {\mathfrak M}_{s^{\prime}}} {\mathcal C}_{s^{\prime}, {\mathbf q}^{\prime}}} Q(T^{(\nu, k(\nu))}, {\mathbf C}^{\prime}, {\mathbf B}_{1/2}(0)), \nonumber
\end{eqnarray}
where $C = C(n, m, q, \overline{M}).$  Here $E_{\nu,k(\nu)} = E(T^{(\nu,k(\nu))},\mathbf{C}^{(\nu,k(\nu))},\mathbf{B}_1(0))$ and $\widetilde{\mathbf C}^{(\nu,k(\nu))}$ and $u_{i}^{(\nu, k(\nu))}$ are as in  \eqref{holefilling0_eqn2-1} and \eqref{holefilling0_eqn2} with $k = k(\nu)$. 
Set $T^{(\nu)} = T^{(\nu,k(\nu))}$, ${\mathbf C}^{(\nu)} = {\mathbf C}^{(\nu,k(\nu))}$, $u_{i}^{(\nu)} = u_{i}^{(\nu, k(\nu))}$, $\widetilde{\mathbf C}^{(\nu)} = \widetilde{\mathbf C}^{(\nu,k(\nu))}$, $E_{\nu} = E_{\nu,k(\nu)}$, 
$\epsilon_{\nu} = \epsilon^{(\nu, k(\nu))}$, $\beta_{\nu} = \beta^{(\nu, k(\nu))}$, $\eta_{\nu} = \eta^{(\nu, k(\nu))}$ and $\delta_{\nu} = \delta^{(\nu, k(\nu))}$. 
After relabelling, assume the above hold for $\nu=1, 2, 3, \ldots$. 
Write 
\begin{equation*}
{\mathbf C}^{(\nu)} = \sum_{i=1}^{p} q_{i} \llbracket P_{i}^{(\nu)}\rrbracket
\end{equation*}
where $P_{i}^{(\nu)} = P_{i}^{(\nu, k(\nu))},$ with $P_{i}^{(\nu, k)} = \{(z, x, y) \in {\mathbb R}^{m} \times {\mathbb R}^{2} \times {\mathbb R}^{n-2}\, : \, 
z = A_{i}^{(\nu, k)} x\}$ as in \eqref{holefilling0_eqn2-0}, where $A_{i}^{(\nu, k)}$ is a constant $(m \times 2)$ matrix. By the definition of $\widetilde{\mathbf C}^{(\nu)}$ (see \eqref{holefilling0_eqn2-1}), we can write 
\begin{equation*}
\widetilde{\mathbf C}^{(\nu)} = \sum_{i=1}^{p} \sum_{\ell=1}^{m_{i}} \widetilde{q}_{i, \ell} \llbracket \widetilde{P}_{i, \ell}^{(\nu)} \rrbracket 
\end{equation*}
where $\widetilde{P}_{i, \ell}^{(\nu)} = \widetilde{P}_{i, \ell}^{(\nu,k(\nu))}$ are distinct planes such that $\widetilde{P}^{(\nu)}_{i, \ell} = \{(z, x, y) \in {\mathbb R}^{m} \times {\mathbb R}^{2} \times {\mathbb R}^{n-2}\, : \, z = A_{i}^{(\nu)} x + E_{\nu} \widetilde{\psi}_{i, \ell}^{(\nu)}(x)\}$ and 
$${\rm dist} \, (X, {\rm spt} \, {\mathbf C}^{(\nu)}) = {\rm dist} \, (X, P_{i}^{(\nu)}) \;\;\;\; \forall X \in \widetilde{P}_{i, \ell}^{(\nu)}.$$

In view of \eqref{holefilling0_eqn23}, all conditions necessary to produce a fine blow-up of (a subsequence of) $(\eta_{0, 1/2 \, \#} \, T^{(\nu)})$ relative to the (corresponding subsequence of) cones $(\widetilde{\mathbf C}^{(\nu)})$ are met.  Thus, we obtain $\widetilde{w} = (\widetilde{w}_{i,l}) \in {\mathfrak B}_{s, \widetilde{\mathbf q}},$ where $\widetilde{\mathbf q} = (\widetilde{q}_{1, 1},  \ldots, \widetilde{q}_{1, m_{1}}, \widetilde{q}_{2, 1} \ldots, \widetilde{q}_{2, m_{2}}, \ldots, \widetilde{q}_{p, 1}, \ldots, \widetilde{q}_{p, m_{p}}),$ 
which can be written as $$\widetilde{w}_{i,\ell}(x) = \sum_{j=1}^{\widetilde{q}_{i, \ell}}\llbracket \widetilde{w}_{i, \ell, j}(x)\rrbracket$$ 
with $\widetilde{w}_{i, \ell, j}(x) \in {\mathbb R}^{m}$ for $x \in B_{1}(0),$ such that for each $i=1, 2, \dots, p$ and each $\ell =1, 2, \ldots,$ 
\begin{equation} \label{holefilling0_eqn26}
			\widetilde{E}_{\nu}^{-1} \widetilde{u}^{(\nu)}_{i, \ell} \to \widetilde{w}_{i, \ell}
\end{equation}
in $L^{2}_{\rm loc}(B_1; {\mathcal A}_{\widetilde{q}_{i, \ell}}({\mathbb R}^{m}))$ and $\widetilde{E}_{\nu}^{-1} |D\widetilde{u}^{(\nu)}_{i, \ell}| \to |D\widetilde{w}_{i, \ell}|$ locally in $L^2(K)$ for each compact set $K \subset B_1(0) \setminus \{0\}\times\mathbb{R}^{n-2}$ as $\nu \rightarrow \infty$, where $\widetilde{E}_{\nu} = E(T^{(\nu)},\widetilde{\mathbf C}^{(\nu)},\mathbf{B}_{1/2}(0)).$   
By  \eqref{holefilling0_eqn28}, \eqref{holefilling0_eqn36} and \eqref{holefilling0_eqn20} we see that for each $0 < \tau < \gamma < 1$, 
\begin{equation*}
\int_{B_{\gamma}(0) \setminus \{r < \tau\}} \sum_{i, \ell} {\mathcal G} \left(\frac{\widetilde{u}^{(\nu)}_{i, \ell}}{E_{\nu}F_{\nu}}, \frac{\overline{w}^{(\nu)}_{i,\ell}}{F_{\nu}}\right)^{2}
= \frac{2^{n+2}}{F_{\nu}^2}\int_{B_{\gamma/2}(0) \setminus \{r < \tau/2\}} 
\sum_{i} {\mathcal G} \left( \frac{u^{(\nu)}_{i}}{E_{\nu}},  w^{(\nu)}_{i}\right)^{2}\to 0 
\end{equation*}
as $\nu \to \infty$, and hence by \eqref{holefilling0_eqn37} and \eqref{holefilling0_eqn26} we have \eqref{holefilling0_eqn26} 
$c \, \widetilde{w}_{i, \ell} = \overline{w}_{i, \ell}$ on $B_{1} \setminus \{0\} \times{\mathbb R}^{n-2}$ where the constant 
$c  = \lim_{\nu \to \infty} \, (E_{\nu}F_{\nu})^{-1} \widetilde{E}_{\nu}\in [0, 2^{n+2}]$ (where the limit exists, after passing to a subsequence, by \eqref{holefilling0_eqn21}).

\noindent
{\bf Step 5:} \emph{Non-triviality of the fine blow up relative to the cones associated with $\psi^{(\nu)}$, its homogeneity, and the contradiction that completes the proof.}  By \eqref{holefilling0_eqn1} and the fact that $\widetilde{\psi}_{i, \ell}$ is homogeneous of degree~1, it follows that for each $\delta \in (0, 1/4)$, 
\begin{eqnarray*}
    &&\int_{(B_1(0) \setminus B_{1/4}(0)) \setminus \{r < \delta\}} \sum_{i=1}^p \sum_{\ell = 1}^{m_{i}} \left|\frac{\partial (\overline{w}_{i, \ell}/R)}{\partial R}\right|^{2} \\ 
    &&\hspace{1in}\leq \lim_{\nu \to \infty} \frac{1}{F_{\nu}^2}\int_{(B_1(0) \setminus B_{1/4}(0)) \setminus \{r < \delta\}} \sum_{i=1}^p \sum_{\ell=1}^{m_{i}}
	\left| \frac{\partial (\overline{w}^{(\nu)}_{i, \ell} /R)}{\partial R}\right|^{2}\\ 
    &&\hspace{1in}= \lim_{\nu \to 0} \frac{1}{2^{n+2} F_{\nu}^2} \int_{(B_{1/2}(0) \setminus B_{1/8}(0)) \setminus \{r < \delta\}} \sum_{i=1}^p \left| \frac{\partial (w^{(\nu)}_{i}/R)}{\partial R}\right|^{2} = 0.
\end{eqnarray*}	
Thus for each $i$ and $\ell$, $\frac{\partial (\overline{w}_{i, \ell}/R)}{\partial R} = 0$ at ${\mathcal H}^{n}$-a.e.\ $X \in B_{1/2}(0) \setminus B_{1/8}(0)$ and hence 
$\overline{w}_{i, \ell}$ is homogeneous of degree 1 in $(B_1(0) \setminus B_{1/4}(0)) \setminus \{0\} \times {\mathbb R}^{n-2}$.  Recall that $c \, \widetilde{w}_{i, \ell} = \overline{w}_{i, \ell}$ on $B_1 \setminus \{0\} \times{\mathbb R}^{n-2}$ where $c  = \lim_{\nu \to \infty} \, (E_{\nu}F_{\nu})^{-1} \widetilde{E}_{\nu}\in [0, 2^{n+2}]$. Since $\sum_{i, \ell} |\overline{w}_{i, \ell}|^{2} \not\equiv 0$ , it follows that $c >0$, and hence, since the homogeneous degree 1 extension of $\widetilde{w} \llcorner ((B_1 \setminus B_{1/2}) \setminus \{0\} \times {\mathbb R}^{n-2})$ is locally Dirichlet energy minimizing in ${\mathbb R}^{n} \setminus \{0\}\times {\mathbb R}^{n-2}$, it follows from unique continuation for locally Dirichlet energy minimizing functions that $\widetilde{w}$ is (non-zero and) homogeneous of degree 1 in $B_{1} \setminus \{0\} \times {\mathbb R}^{n-2}$. By Theorem~\ref{classification_thm}, $\widetilde{w} \in \widetilde{\mathfrak L}_{s, \widetilde{\mathbf q}}$ so that $\widetilde{w}_{i, \ell}(x,y)  = \sum_{k=1}^{\widetilde{q}_{i, \ell}} \llbracket \widetilde{\varphi}_{i, \ell, k}(x,y) \rrbracket$ where $\widetilde{\varphi}_{i, \ell, k} \, : \, {\mathbb R}^{n} \to{\mathbb R}^{m}$ are single valued linear functions of the form 
$\widetilde{\varphi}_{i, \ell, k}(x, y) =  \varphi_{i, \ell, k}(x) +  L_{1}(y) + A_{i, \ell} L_{2}(y)$ for each $(x, y) \in {\mathbb R}^{2} \times {\mathbb R}^{n-2}$, where $\varphi_{i, \ell, k} \, : \, {\mathbb R}^{2} \to {\mathbb R}^{m}$ are linear functions such that: for $k_{1} \neq k_{2}$, either 
$\varphi_{i, \ell, k_{1}}(x) \equiv \varphi_{i, \ell, k_{2}}(x)$ for all $x \in {\mathbb R}^{2}$ or $\varphi_{i, \ell, k_{1}}(x) \neq \varphi_{i, \ell, k_{2}}(x)$ for all $x \in {\mathbb R}^{2} \setminus \{0\};$ $L_{1} \, : \, {\mathbb R}^{n-2} \to {\mathbb R}^{m}$, $L_{2} \, : \, {\mathbb R}^{n-2} \to {\mathbb R}^{2}$ are linear functions, and $A_{i, \ell}$ are constant $m \times 2$ matrices. By the triangle inequality, we then have $\int_{B_1(0) \setminus \{r < \nu^{-1}\}} \sum_{i=1}^p \sum_{\ell=1}^{m_i} {\mathcal G} \left(\frac{\overline{w}^{(\nu)}_{i,\ell}}{F_{\nu}}, c \, \widetilde{\varphi}_{i, \ell}\right)^{2} \to 0$ as $\nu \to \infty$, 
which by \eqref{holefilling0_eqn36} says that  
$$\frac{1}{F_{\nu}^2} \int_{B_{1/2}(0) \setminus \{r < \nu^{-1}\}} \sum_{i=1}^p \sum_{\ell=1}^{m_i} {\mathcal G} \left(\sum_{k=1}^{\widetilde{q}_{i, \ell}} \llbracket w^{(\nu)}_{i, \ell, k}(x, y)  - \widetilde{\psi}^{(\nu)}_{i, \ell}(x)\rrbracket, \sum_{k=1}^{\widetilde{q}_{i, \ell}}\llbracket c \, F_{\nu} \widetilde{\varphi}_{i, \ell, k}(x, y)\rrbracket\right)^{2} \to 0$$ 
as $\nu \to \infty.$ If we define $\overline{\psi}^{(\nu)}_{i} \, : \, {\mathbb R}^{n} \to {\mathcal A}_{q_{i}}({\mathbb R}^{m})$ by 
$\overline{\psi}^{(\nu)}_{i}(x, y) = \sum_{\ell=1}^{m_{i}}\sum_{j=1}^{\widetilde{q}_{i, \ell}} \llbracket \widetilde{\psi}_{i, \ell}^{(\nu)}(x)+ c \, F_{\nu} (\varphi_{i, \ell, j}(x) +  L_{1}(y) + A_{i, \ell} L_{2}(y)) \rrbracket$ and set $\overline{\psi}^{(\nu)} = (\overline{\psi}^{(\nu)}_{1}, \ldots, \overline{\psi}^{(\nu)}_{p}),$ then 
$\overline{\psi}^{(\nu)} \in \widetilde{\mathfrak L}_{p, {\mathbf q}}$ and the above says 
\begin{equation}\label{holefilling0_eqn32-1} 
    \frac{1}{F_{\nu}^2}\int_{B_{1/2}(0) \setminus \{r < \nu^{-1}\}} \sum_{i=1}^{p} {\mathcal G} \left(w^{(\nu)}_{i}, \overline{\psi}^{(\nu)}_{i}\right)^{2} \to 0.
\end{equation}

By the regularity theory for locally Dirichlet energy minimizing functions (\cite{Almgren}), there is a 
closed set $\Sigma_{i}^{(\nu)} \subset B_{1/2}(0)$ (consisting of the union of the singular set of $w^{(\nu)}_{i} \llcorner (B_{1/2} \setminus (\{0\} \times {\mathbb R}^{n-2}))$ and $B_{1/2} \cap (\{0\} \times {\mathbb R}^{n-2})$) of Hausdorff dimension $\leq n-2$ such that locally about any point $y \in B_{1/2} \setminus \Sigma^{(\nu)}_{i}$, the function $w^{(\nu)}_{i}$ is given by $q_{i}$ 
smooth ${\mathbb R}^{m}$-valued harmonic functions. Letting $\sigma^{(\nu)}_{i} \subset {\mathbb S}^{n-1}$ be the radial projection of $\Sigma_{i}^{(\nu)} \cap (B_{1/2} \setminus B_{1/8})$ into ${\mathbb S}^{n-1}$, we then have that ${\rm dim}_{\mathcal H} \, (\sigma^{(\nu)}_{i}) \leq n-2$ and that for each $\omega \in {\mathbb S}^{n-1} \setminus \sigma^{(\nu)}_{i}$, there is a simply connected neighborhood $K_{\omega}$ of $\{t\omega \, : \, t \in [1/8, 1/2]\}$ such that 
$w^{(\nu)}_{i} \llcorner K_{\omega}$ decomposes as $q_{i}$ smooth ${\mathbb R}^{m}$-valued functions on $K_{\omega}$. Using this decomposition, we see  
that for each $\omega \in {\mathbb S}^{n-2} \setminus \sigma^{(\nu)}_{i}$ and $1/8 \leq \rho_{1} \leq \rho_{2} \leq 1/2,$  
${\mathcal G}(w^{(\nu)}_{i}(\rho_{2}\omega)/\rho_{2}, w^{(\nu)}_{i}(\rho_{1}\omega)/\rho_{1}) \leq \int_{1/8}^{1/2}  \left| \frac{\partial (w^{(\nu)}_{i}/R)(t\omega)}{\partial R}\right| dt.$ Hence we have by homogeneity of $\overline{\psi}^{(\nu)}_{i}$, triangle inequality and the Cauchy--Schwarz inequality, 
$$\frac{{\mathcal G}(w^{(\nu)}_{i}(\rho_{2}\omega), \overline{\psi}^{(\nu)}_{i}(\rho_{2}\omega))^{2}}{\rho_{2}^{2}} \leq  \frac{2\,{\mathcal G}(w^{(\nu)}_{i}(\rho_{1}\omega), \overline{\psi}^{(\nu)}_{i}(\rho_{1}\omega))^{2}}{\rho_{1}^{2}}  + \frac{9 \cdot 8^{n-1}}{32}\int_{1/8}^{1/2}  t^{n-1}\left|\frac{\partial (w^{(\nu)}_{i}/R)(t\omega)}{\partial R}\right|^{2} dt.$$ Multiplying this by 
$\rho_{2}^{n+1}$ and integrating with respect to $\rho_{2} \in [1/8, 1/2]$, and then multiplying the resulting inequality by  $\rho_{1}^{n+1}$ and integrating with respect to $\rho_{1} \in [1/8, 1/4]$, followed by integration with respect to $\omega \in {\mathbb S}^{n-1} \setminus \sigma^{(\nu)}_{i},$ adding 
$\int_{B_{1/8}} {\mathcal G}(w^{(\nu)}_{i}(x), \overline{\psi}^{(\nu)}_{i}(x))^{2} \, dx$ to both sides and summing over $i$, we obtain that 
\begin{align}\label{holefilling0_eqn32}
    \int_{B_{1/2}} \sum_{i=1}^{p}{\mathcal G}(w^{(\nu)}_{i}(x), \overline{\psi}^{(\nu)}_{i}(x))^{2} \, dx 
    \leq\,& C \int_{B_{1/4}} \sum_{i=1}^{p}{\mathcal G}(w^{(\nu)}_{i}(x), \overline{\psi}^{(\nu)}_{i}(x))^{2} \, dx 
    \\&+ \, C\int_{B_{1/2}(0) \setminus B_{1/8}(0)} \sum_{i=1}^p \left| \frac{\partial (w^{(\nu)}_{i}/R)}{\partial R}\right|^{2} \nonumber
\end{align} 
for each $\nu$, where $C = C(n).$  
Since $\overline{\psi}^{(\nu)} \in \widetilde{\mathfrak L}_{p, {\mathbf q}}$  satisfies \eqref{holefilling0_eqn32-1} 
\begin{equation*}
    \int_{B_{1/4}\cap \{r < \delta\}} \sum_{i=1}^{p}{\mathcal G}({\psi}^{(\nu)}_{i}, \overline{\psi}_{i}^{(\nu)})^{2} \, dx \leq C\delta^2 \int_{B_{1/2}} \sum_{i=1}^{p}{\mathcal G}({\psi}^{(\nu)}_{i}, \overline{\psi}_{i}^{(\nu)})^{2} \, dx \leq C \delta^{2} F_{\nu}^2
\end{equation*}
for some $C = C(n).$  Hence by \eqref{holefilling0_eqn19-2}, for any $\delta \in (0, 1/4)$ and sufficiently large $\nu$, 
\begin{align}\label{holefilling0_eqn38}
     \int_{B_{1/2} \cap \{r < \delta\}} \sum_{i=1}^{p}{\mathcal G}(w^{(\nu)}_{i}, \overline{\psi}^{(\nu)}_{i})^{2} \, dx  
     \leq\,& 2 \int_{B_{1/4}\cap \{r < \delta\}} \sum_{i=1}^{p}{\mathcal G}(w^{(\nu)}_{i}, {\psi}^{(\nu)}_{i})^{2} \, dx \\&+ 2 \int_{B_{1/4}\cap \{r < \delta\}} \sum_{i=1}^{p}{\mathcal G}({\psi}^{(\nu)}_{i}, \overline{\psi}_{i}^{(\nu)})^{2} \, dx \nonumber\\ 
     \leq\,& C\delta F_{\nu}^{2} , \nonumber
\end{align} 
where $C = C(n, m, q,M).$  Taking $\psi^{\prime} = \overline{\psi}^{(\nu)}$ in \eqref{holefilling0_hyp1-0} and using \eqref{holefilling0_eqn38}, \eqref{holefilling0_eqn32} gives us 
\begin{align*}
    F_{\nu}^{2} \leq\,& C \delta \overline{M}^2 F_{\nu}^{2} +  2C \overline{M}^2 \int_{B_{1/4}(0) \setminus \{r < \delta\}} \sum_{i=1}^{p} {\mathcal G} \left(w^{(\nu)}_{i}, \overline{\psi}^{(\nu)}_{i}\right)^{2} \, dx \\&+ C \overline{M}^2 \int_{B_{1/2}(0) \setminus B_{1/8}(0)} \sum_{i=1}^p \left| \frac{\partial (w^{(\nu)}_{i}/R)}{\partial R}\right|^{2} ,
\end{align*}
where $C = C(n, m, q,M).$  In view of \eqref{holefilling0_eqn1} and \eqref{holefilling0_eqn32-1}, dividing this by $F_{\nu}^{2}$ and letting $\nu \to \infty$ leads to a contradiction if we choose $\delta = \delta(n, m, q,\overline{M})$ such that $C \delta \overline{M}^2 = 1/2$.
This completes the proof. \end{proof}

\begin{proof}[Proof of Lemma~\ref{holefilling_lemma}] 
 Let $w \in \mathfrak{B}_{p, {\mathbf q}}(M)$ and $\psi \in \widetilde{\mathfrak{L}}_{p, {\mathbf q}}$ be as in the lemma, and let $(T_{k}),$ 
$({\mathbf C}_{k})$ be sequences corresponding to $w$ as in Defintion~\ref{blowupclass_defn}. 
By considering a fine blow up of appropriately rotated $T_{k}$ (as in the proof of Theorem~\ref{classification_thm}), we may assume without loss of generality that 
$\psi \in \mathfrak{L}_{p, {\mathbf q}}$.  
  
Let $\overline{\b} = \overline{\b}(n, m, q, M, \overline{M})$, $\overline{C} = \overline{C}(n, m, q, M, \overline{M})$ be the constants as in Lemma~\ref{holefilling0_lemma} (which depend only on $n$, $m$, $q,$ $M$ and the parameter $\overline{M} \geq 1$). For $i \in \{1,\ldots,q-p+1\}$, inductively define $\beta_i$ and $C_i$ by 
setting $\beta_1 = \overline{\beta}(n, m, q, M, 1)$ and $C_1 = \overline{C}(n, m, q, M, 1),$ 
and for each $i \geq 2$, $\beta_i = \overline{\beta}(n, m, q, M, M_{i})$ and $C_i = 2\overline{C}(n, m, q, M, M_{i})$ 
where  $M_{i} = 2^{i-1} (\beta_1 \beta_2 \cdots \beta_{i-1})^{-2}$.  

Observe that \eqref{holefilling0_hyp1} with $\overline{M} = 1$ holds true by hypothesis of the present lemma.  If either $d(\psi) = p,$ or if $d(\psi) > p$ and $w$ and $\psi$ satisfy \eqref{holefilling0_hyp2} with $\overline{\beta} = \beta_1$, then by Lemma~\ref{holefilling0_lemma}, $w$ and $\psi$ satisfy \eqref{holefilling_est} with $C = C_1$.  If instead 
$d(\psi)  > p$ and $w$ and $\psi$ do not satisfy 
\eqref{holefilling0_hyp2} with $\overline{\beta} = \beta_1$, then for $i = 1,2,\ldots,i_0$ inductively select 
$\psi^{(i)} \in \mathfrak{L}_{p, {\mathbf q}}(s_{i})$ such that when $i = 1$ we have $p \leq s_1 < d(\psi)$ and 
\begin{equation*}
	\int_{B_{1/2}(0)} \sum_{j=1}^p \mathcal{G}(w_{j},\psi^{(1)}_{j})^2 
		< 2 \inf_{\psi' \in \bigcup_{s'=p}^{d(\psi)-1} \mathfrak{L}_{p,{\mathbf q}}(s^{\prime})} 
			\int_{B_{1/2}(0)} \sum_{j=1}^p \mathcal{G}(w_{j},\psi'_{j})^2 
\end{equation*}
and for each $i \geq 2$ we have $p \leq s_i < s_{i-1}$ and 
\begin{equation*}
	\int_{B_{1/2}(0)} \sum_{j=1}^p \mathcal{G}(w_{j},\psi^{(i)}_{j})^2 
		< 2 \inf_{\psi' \in \bigcup_{s'=p}^{s_{i-1}-1} \mathfrak{L}_{p,{\mathbf q}}(s^{\prime})} 
			\int_{B_{1/2}(0)} \sum_{j=1}^p  \mathcal{G}(w_{j},\psi'_{j})^2 
\end{equation*}
and terminate either when $i$ equals the smallest $i_0$ for which $s_{i_0} > p$ and 
\begin{equation} \label{holefilling_eqn1}
	\int_{B_{1/2}(0)} \sum_{j=1}^p\mathcal{G}(w_{j},\psi^{(i_0)}_{j})^2 \leq \beta_{i_0+1}^{2} \inf_{\psi' \in \bigcup_{s'=p}^{s_{i_0}-1} \mathfrak{L}_{p, {\mathbf q}}(s^{\prime})} 
			\int_{B_{1/2}(0)} \sum_{j=1}^p \mathcal{G}(w_{j},\psi'_{j})^2  
\end{equation}
or (if no such $i_{0}$ exists) when $i$ equals the value $i_{0}$ for which $s_{i_{0}} = p.$   
By choice of $\psi^{(i)}$, one readily checks that $\psi^{(i_0)}$ satisfies 
\begin{align} \label{holefilling_eqn2}
	&\int_{B_{1/2}(0)} \sum_{j=1}^p \mathcal{G}(w_{j},\psi_{j})^2 
	\leq 2\int_{B_{1/2}(0)} \sum_{j=1}^p \mathcal{G}(w_{j},\psi^{(i_0)}_{j})^2 
	\\&\hspace{10mm} \leq \frac{2^{i_0}}{\beta_1^{2} \beta_2^{2} \cdots \beta_{i_0}^{2}} \int_{B_{1/2}(0)} 
		\sum_{j=1}^p \mathcal{G}(w_{j},\psi_{j})^2, \nonumber 
\end{align}
where the first inequality follows from \eqref{holefilling_est-0} and the second inequality follows by the construction of $\psi^{(i)}$, including the failure, when $i < i_{0}$, of \eqref{holefilling_eqn1} with $i$ in place of $i_{0}$. Using Lemma~\ref{holefilling0_lemma} and \eqref{holefilling_eqn2} gives us 
\begin{equation*}
\int_{B_{1/2}(0)} \sum_{j=1}^p \mathcal{G}(w_{j},\psi_{j})^2 
	\leq 2\int_{B_{1/2}(0)} \sum_{j=1}^p \mathcal{G}(w_{j},\psi^{(i_0)}_{j})^2 \leq 2C_{i_{0}} \int_{B_{1/2}(0) \setminus B_{1/8}(0)} \sum_{i=1}^p \left| \frac{\partial (w_{i}/R)}{\partial R}\right|^{2}.
\end{equation*}
\end{proof}

\begin{lemma} \label{blowupdecay_lemma}
Let $p \in \{2, \ldots, q\}$, ${\mathbf q} \in {\mathfrak M}_{p}$ and $M \geq 1$. Let $w \in \mathfrak{B}_{p, {\mathbf q}}(M)$. For each $\vartheta\in (0, 1/8]$, let 
$\psi^{(\vartheta)}$ be any element in $\widetilde{\mathfrak{L}}_{p, {\mathbf q}}$ such that 
$\int_{B_{\vartheta}(0)} \sum_{j=1}^{p} \mathcal{G}(w_{j},\psi^{(\vartheta)}_{j})^2 \leq 2 \inf_{\psi^{\prime} \in \widetilde{\mathfrak{L}}_{p, {\mathbf q}}} \, \int_{B_{\vartheta}(0)} \sum_{j=1}^{p} \mathcal{G}(w_{j},\psi^{\prime}_{j})^2.$
Then for each $\vartheta_{1}$, $\vartheta_{2}$ with $0 < \vartheta_{1} \leq \vartheta_{2} \leq 1/8,$
\begin{align} \label{blowupdecay_eqn1}
	\vartheta_{1}^{-n-2} \int_{B_{\vartheta_{1}}(0)} \sum_{j=1}^{p} \mathcal{G}(w_{j},\psi^{(\vartheta_{1})}_{j})^2 
	\leq C \left(\frac{\vartheta_{1}}{\vartheta_{2}}\right)^{2\mu} \vartheta_{2}^{-n-2}\int_{B_{\vartheta_{2}}(0)} \sum_{j=1}^{p} \mathcal{G}(w_{j},\psi^{(\vartheta_{2})}_{j})^2, 
\end{align}
where $\mu \in (0,1)$ and $C \in (0,\infty)$ depend only on $n$, $m$, $q$ and $M$ (in particular independent of $\vartheta_{1}$ and $\vartheta_{2}$).
\end{lemma}
\begin{proof}
Let $\rho \in [\vartheta_{1}, 1/2]$. By \eqref{blowup_est1} with $\psi^{(\rho)}$ in place of $\psi$, 
\begin{equation} \label{blowupdecay_eqn2}
	\int_{B_{\rho/4}(0)} \sum_{j=1}^p R^{2-n} \left| \frac{\partial (w_{j}/R)}{\partial R} \right|^2 
	\leq C \rho^{-n-2} \int_{B_{\rho}(0)} \sum_{j=1}^p \mathcal{G}(w_{j},\psi^{(\rho)}_{j})^2 
\end{equation}
where $C = C(n,m,q) \in (0,\infty)$. Let $T_{k}$, ${\mathbf C}_{k}$ be as in Definition~\ref{blowupclass_defn} such that $w$ is the blow-up of $T_{k}$ relative to 
${\mathbf C}_{k}$ (by the excess $E_{k} = \sqrt{\int_{B_{1}(0)} {\rm dist}^{2} \, (X, {\rm spt} \, \|{\mathbf C}_{k}\|) \, d\|T_{k}\|}$). We may 
 assume that $\left.w\right|_{B_{{\vartheta}_{1}}} \neq 0$ (else there is nothing to prove) and hence $\left.w\right|_{B_{\rho}} \neq 0$. By 
 Remark~\ref{rescaling-fine-blowup}, we have that $w_{1} \equiv \|w(\rho(\cdot))\|_{L^{2}(B_{1}(0))}^{-1}w(\rho(\cdot)) \in {\mathfrak B}_{p, {\mathbf q}}(CM)$ for 
 some $C = C(n, m, q)$.
Hence by applying Lemma~\ref{holefilling_lemma} to $w_{1}$, we see that 
\begin{equation} \label{blowupdecay_eqn3}
\rho^{-n-2} \int_{B_{\rho}(0)} \sum_{j=1}^p \mathcal{G}(w_{j},\psi^{(\rho)}_{j})^2 \leq C \int_{B_{\rho}(0) \setminus B_{\rho/4}(0)} \sum_{j=1}^p R^{2-n} \left| \frac{\partial (w_{j}/R)}{\partial R} \right|^2  
\end{equation}
where $C = C(n,m,q, M) \in (0,\infty)$.  Thus by \eqref{blowupdecay_eqn2} and \eqref{blowupdecay_eqn3}, 
\begin{equation} \label{blowupdecay_eqn4}
\int_{B_{\rho/4}(0)} \sum_{j=1}^p R^{2-n} \left| \frac{\partial (w_{j}/R)}{\partial R} \right|^2	\leq C_0 \int_{B_{\rho}(0) \setminus B_{\rho/4}(0)} \sum_{j=1}^p R^{2-n} \left| \frac{\partial (w_{j}/R)}{\partial R} \right|^2 \end{equation}
for all $\rho \in [\vartheta_{1},1/2]$ and some constant $C_0 = C_0(n,m,q, M) \in (0,\infty)$.  By adding $C_0$ times the left-hand side of \eqref{blowupdecay_eqn4} to both sides of \eqref{blowupdecay_eqn4}, 
\begin{equation} \label{blowupdecay_eqn5}
	\int_{B_{\rho/4}(0)} \sum_{j=1}^p R^{2-n} \left| \frac{\partial (w_{j}/R)}{\partial R} \right|^2
	\leq \gamma \int_{B_{\rho}(0)} \sum_{j=1}^p R^{2-n} \left| \frac{\partial (w_{j}/R)}{\partial R} \right|^2
\end{equation}
for all $\rho \in [\vartheta_{1},1/2]$, where $\gamma = C_0/(1 + C_0) \in (0,1)$.  Noting that the assertion of the lemma easily follows if $\vartheta_{1} \geq \vartheta_{2}/32$, we may assume that $\vartheta_{1} <\vartheta_{2}/32$ and iteratively apply \eqref{blowupdecay_eqn5} with $\rho = 2^{-2i-1}\vartheta_{2}$ for $i = 1,2,\ldots,N-1$, where $N$ is the positive integer such that $2^{-2N-3}\vartheta_{2} < \vartheta_{1} \leq 2^{-2N-1}\vartheta_{2}$, to obtain 
\begin{equation} \label{blowupdecay_eqn6}
	\int_{B_{\vartheta_{1}}(0)} \sum_{j=1}^p R^{2-n} \left| \frac{\partial (w_{j}/R)}{\partial R} \right|^2
	\leq C \left(\frac{\vartheta_{1}}{\vartheta_{2}}\right)^{2\mu} \int_{B_{\vartheta_{2}/8}(0)} \sum_{j=1}^p R^{2-n} \left| \frac{\partial (w_{j}/R)}{\partial R} \right|^2, 
\end{equation}
where $\mu = -\log \gamma/\log 16$ and $C = C(n, m, q, M)$.  By combining \eqref{blowupdecay_eqn2} with $\rho = \vartheta_{2}/2$, \eqref{blowupdecay_eqn3} with $\rho = \vartheta_{1}$, and \eqref{blowupdecay_eqn6}, we obtain \eqref{blowupdecay_eqn1}.
\end{proof}

\begin{theorem}\label{blowup-decay-unique}
Let $p \in \{2, \ldots, q\}$, ${\mathbf q}  = (q_{1}, q_{2}, \ldots, q_{p}) \in {\mathfrak M}_{p},$ $M \geq 1$ and $w \in \mathfrak{B}_{p, {\mathbf q}}(M).$ There exists unique 
$\widetilde{\psi} \ = (\widetilde{\psi}_{1}, \ldots, \widetilde{\psi}_{p}) \in  \widetilde{\mathfrak{L}}_{p, {\mathbf q}}$ such that 
\begin{align} \label{blowup-decay-unique-1}
	\rho^{-n-2} \int_{B_{\rho}(0)} \sum_{j=1}^{p} \mathcal{G}(w_{j},\widetilde{\psi}_{j})^2 
	\leq C \rho^{2\mu} \int_{B_{1/2}(0)} \sum_{j=1}^{p} |w_{j}|^2, 
\end{align}
for all $\rho \in (0, 1/8]$, where $\mu \in (0,1)$ and $C \in (0,\infty)$ depend only on $n$, $m$,  $q$ and $M$. Additionally, 
\begin{itemize}
\item[(i)] we have that   
\begin{equation}\label{blowup-decay-unique-2}
\int_{B_{1}(0)} \sum_{j=1}^{p}|\widetilde{\psi}_{j}|^{2} \leq \overline{C}
\end{equation}
for some constant $\overline{C} = \overline{C}(n, m, q) \in (0, \infty)$, and 
\item[(ii)] there are $\psi = (\psi_{1}, \ldots, \psi_{p}) \in \mathfrak{L}_{p, {\mathbf q}}$, 
constant $(m \times 2)$ matrices $A_{1}, \ldots, A_{p}$ and  two linear functions 
$L_{1} \, : \, {\mathbb R}^{n-2} \to {\mathbb R}^{m}$, $L_{2} \, : \, {\mathbb R}^{n-2} \to {\mathbb R}^{2}$ with  
\begin{equation}\label{blowup-decay-unique-3} 
	\sum_{j=1}^p \|A_{j}\| + \|L_{1}\| + \|L_{2}\| \leq C
\end{equation}
for some $C = C(n, m, q)$ such that $w \in {\mathfrak B}_{p, {\mathbf q}}(M, A_{1}, \ldots, A_{p})$ and $\widetilde{\psi} \in \widetilde{\mathfrak{L}}_{p, {\mathbf q}}(A_{1}, \ldots, A_{p})$ with $\widetilde{\psi}_{j}(x, y) = \sum_{k =1}^{q_{j}}\llbracket \psi_{j, k}(x) + L_{1}(y) + A_{j}L_{2}(y)\rrbracket$ for $(x, y) \in {\mathbb R}^{2} \times {\mathbb R}^{n-2}.$
\end{itemize}  
\end{theorem}

\begin{proof}
For $k=1, 2, \ldots,$, let $\psi^{k} = \psi^{(8^{-k})}$ be an element in $\widetilde{\mathfrak{L}}_{p, {\mathbf q}}$ such that 
$\int_{B_{8^{-k}}(0)} \sum_{j=1}^{p} \mathcal{G}(w_{j},\psi^{k}_{j})^2 \leq 2 \inf_{\psi^{\prime} \in \widetilde{\mathfrak{L}}_{p, {\mathbf q}}} \, \int_{B_{8^{-k}}(0)} \sum_{j=1}^{p} \mathcal{G}(w_{j},\psi^{\prime}_{j})^2.$
Applying Lemma~\ref{blowupdecay_lemma} with $\vartheta_{2} = 1/2$, $\vartheta_{1}  \in \{8^{-(k+1)}, 8^{-k}\}$, and 
using the triangle inequality and homogeneity of $\psi_{k}$, we see that   
$\int_{B_{1}(0)} \sum_{j=1}^{p}{\mathcal G}(\psi_{j}^{k+1}, \psi_{j}^{k})^{2} \leq C8^{-2k\mu}\int_{B_{1/2}(0)} \sum_{j=1}^{p} |w_{j}|^2$ where $C = C(n, m, q) \in (0, \infty)$ and $\mu = \mu(n, m, q) \in (0, 1)$. This implies that for each $j=1, 2, \ldots, p$, $(\psi_{j}^{k})$ is a sequence of ${\mathcal A}_{q_{j}}({\mathbb R}^{m})$-valued linear functions on $B_{1}(0)$ which is Cauchy with respect to the uniform metric, and hence there is $\widetilde{\psi} = (\widetilde{\psi}_{1}, \ldots, \widetilde{\psi}_{p}) \in \widetilde{\mathfrak L}_{p, {\mathbf q}}$ such that $\psi^{k} \to \widetilde{\psi}$ as $k \to \infty$ uniformly 
on $B_{1}(0)$.  It can now be readily checked that $(8^{-k})^{-n-2}\int_{B_{8^{-k}}(0)} \sum_{j=1}^{p}{\mathcal G}(w_{j}, \widetilde{\psi}_{j})^{2} \leq C8^{-2k\mu}\int_{B_{1/2}(0)} \sum_{j=1}^{p} |w_{j}|^2$ for all $k=1, 2, \ldots$, and hence that $\rho^{-n-2}\int_{B_{\rho}(0)} \sum_{j=1}^{p}{\mathcal G}(w_{j}, \widetilde{\psi}_{j})^{2} \leq C_{1}\rho^{2\mu_{1}}\int_{B_{1/2}(0)} \sum_{j=1}^{p} |w_{j}|^2$ for all $\rho \in (0, 1/8]$, where $C_{1} = C_{1}(n, m, q) \in (0, \infty)$ and $\mu_{1} = \mu_{1}(n, m, q) \in (0, 1)$. Uniqueness of $\widetilde{\psi}$ is clear from this estimate. Taking $\rho  = 8^{-1}$ in the estimate and using the triangle inequality and the fact that $\omega_{n}^{-1}\int_{B_{1/2}(0)} |w_{j}|^{2} \leq 1$ imply that $\int_{B_{1}(0)} \sum_{j=1}^{p}|\widetilde{\psi}_{j}|^{2} \leq \overline{C}$, where $\overline{C} = \overline{C}(n, m, q)$. To see that conclusion (ii) holds, we may assume that $\widetilde{\psi} \not\equiv 0$ whence by 
(\ref{blowup-decay-unique-1}), $0 < c \equiv \int_{B_{1}(0)} \sum_{j=1}^{p}|\widetilde{\psi}_{j}|^{2}  = \rho^{-n-2}\int_{B_{\rho}(0)} \sum_{j=1}^{p} |\widetilde{\psi}_{j}|^{2} \leq 
2C\rho^{2\mu} + 2\rho^{-n-2}\int_{B_{\rho}(0)} \sum_{j=1}^{p}|w_{j}|^{2}$ and hence $\rho^{-n-2}\int_{B_{\rho}(0)} \sum_{j=1}^{p}|w_{j}|^{2} >c/2$ for all sufficiently small $\rho >0$. We also have, again by (\ref{blowup-decay-unique-1}) and the triangle inequality that 
$\rho^{-n-2}\int_{B_{\rho}(0)} \sum_{j=1}^{p}|w_{j}|^{2} \leq  2C\rho^{2\mu} \int_{B_{1/2}(0)} \sum_{j=1}^{p} |w_{j}|^{2} + \overline{C} \leq 2C + \overline{C}$. Hence if $\rho_{\ell}$ is a sequence with $\rho_{\ell} \to 0^{+}$, then after passing to a subsequence $c_{\ell} \equiv \rho_{\ell}^{-n-2}\int_{B_{\rho_{\ell}}(0)}\sum_{j=1}^{p}|w_{j}|^{2}  \to c_{\star} \in 
[c/2, 2C + \overline{C}].$  Now note that $w \in {\mathfrak B}_{p, {\mathbf q}}(M, A_{1}, \ldots, A_{p})$ for some $m \times 2$ matrices $A_{1}, \ldots, A_{p}$, so by Remark~\ref{rescaling-fine-blowup}, we have that  
$w^{(\ell)} \equiv \|w(\rho_{\ell}(\cdot))\|_{L^{2}(B_{1}(0))}^{-1}w(\rho_{\ell}(\cdot)) \in {\mathfrak B}_{p, {\mathbf q}}(CM, A_{1}, \ldots, A_{p})$ for each $\ell$ and some fixed 
$C = C(n, m, q)$, and hence by a diagonal sequence argument, 
there exists $\overline{w} \in {\mathfrak B}_{p, {\mathbf q}}(CM, A_{1}, \ldots, A_{p})$ such that $w^{(\ell)} \to \overline{w}$ locally in $L^{2}$. 
Taking $\rho = \rho_{\ell}$ in (\ref{blowup-decay-unique-1}), dividing both sides by $c_{\ell}$ and letting $\ell \to \infty$, we deduce that $\widetilde{\psi} = \sqrt{c_{\star}} \overline{w}$. In particular, $\overline{w}$ is a homogeneous degree 1 element of ${\mathfrak B}_{p, {\mathbf q}}(CM, A_{1}, \ldots, A_{p})$. Conclusion (ii), and in particular the estimate (\ref{blowup-decay-unique-3}), follows from Theorem~\ref{classification_thm} since $c_{\star} \leq 2C + \overline{C}$ (with $C$, $\overline{C}$ depending only on $n$, $m$ and $q$). 
\end{proof}

\subsection{Decay of fine excess of area minimizing currents}\label{finedecay}
We start with the following preliminary excess decay result. 
\begin{lemma}\label{excess-improvement} 
Let $q$ be an integer such that $q \geq 2$, $\th \in (0, 1/4)$ and $M \in [1, \infty)$. 
There exist numbers  $\overline\b = \overline\b(n, m, q, M, \th) \in (0, 1/2)$, $\overline\eta = \overline\eta(n, m, q, M, \th) \in (0, 1/2)$  and $\overline{\delta} = \overline{\delta}(n, m, q, M, \th) \in (0, 1/2)$ such that the following holds true: if ${\mathbf C} \in {\mathcal C}_{q, p}$ for some $p \in \{2, \dots, q\}$ and $T$ is an $n$-dimensional locally area minimizing rectifiable current in ${\bf B}_{1}(0)$ such that 
$$\partial \, T \llcorner {\mathbf B}_{1}(0) = 0, \;\; \omega_{n}^{-1}\|T\|({\mathbf B}_{1}(0)) < q + 1/2, \;\; \Theta(T, 0) \geq q,$$ Hypothesis~($\star\star$) (of Section~\ref{notation-and-graphical}) and Hypothesis~($\dag$)  (of Section~\ref{sec:apriori close2plane}) hold with $\overline\b$, $\overline\eta$  in place of 
$\beta_{0}$, $\eta_{0}$ respectively,   
then either 
\begin{itemize}
\item[(A)] ${\mathbf B}_{\overline{\d}}(0,z) \cap \{Z \, : \, \Theta \, (T, Z) \geq q\}  = \emptyset$
for some point $(0, z) \in \{0\} \times {\mathbb R}^{n-2} \cap {\mathbf B}_{1/2}(0)$ or, 
 
\item[(B)] there exist an orthogonal rotation $\G$ of ${\mathbb R}^{n+m}$ and a cone ${\mathbf C}^{\prime}$ with
${\mathbf C}^{\prime} \in {\mathcal C}_{q, p^{\prime}}$ for some $p^{\prime} \in \{2, \ldots, q\}$ such that, writing 
$${\hat E}_{T}^{2} = \int_{{\mathbf B}_{1}(0)} {\rm dist}^{2}(X,  P_{0}) \, d\|T\|(X) \;\;\; {and} \;\;\; 
E_{T}^{2} =  \int_{{\mathbf B}_{1}(0)} {\rm dist}^{2} \, (X, {\rm spt} \, {\mathbf C}) \, d\|T\|(X),$$ 
the following hold:
\begin{eqnarray*}
&&({\rm a}) \;\;\;\; |e_{i} - \G(e_{i})| \leq \overline{\k}E_{T}  \;\; {for} \;\; i=1,  \ldots, m\;\; {and}  \;\; |e_{m+j} - \G(e_{m+j})| \leq \overline{\k}{\hat E}_{T}^{-1}E_{T}\nonumber\\ 
&&\;\;\;\;\; \mbox{for} \;\; j=1,  \ldots, n;\nonumber\\   
&&({\rm b})\;\;\;\; {\rm dist}_{\mathcal H} \, ({\rm spt} \, {\mathbf C}^{\prime} \cap {\mathbf B}_{1}, {\rm spt} \, {\mathbf C} \cap {\mathbf B}_{1}) \leq \overline{C}_{0}E_{T};\nonumber\\
&&({\rm c})\;\;\;\;\th^{-n-2}\int_{\G\left({\mathbf B}_{\th/2} \setminus \{r(X) \leq \th/16\}\right)} {\rm dist}^{2} \, (X, {\rm spt} \, T) \, d\|\G_{\#} \, {\mathbf C}^{\prime}\|(X)\nonumber\\
&&\hspace{1.5in}+ \;\th^{-n-2}\int_{{\mathbf B}_{\th}} {\rm dist}^{2} \, (X, {\rm spt} \, \G_{\#} \, {\mathbf C}^{\prime}) \, d\|T\|(X)  \leq \overline{\nu}\th^{2\mu}E_{T}^{2};\nonumber\\
&&({\rm d}) \;\;\;\; \left(\th^{-n-2}\int_{{\mathbf B}_{\th}} {\rm dist}^{2}\, (X, P)\, d\|\G^{-1}_{\#} \,T\|(X)\right)^{1/2} \geq \overline{C}_{1} \,{\rm dist}_{\mathcal H}\, ({\rm spt} \, {\mathbf C} \cap {\mathbf B}_{1}, P \cap {\mathbf B}_{1}) - \overline{C}_{2}E_{T} \nonumber\\
&&\;\;\;\;\mbox{for any $n$-dimensional plane $P \subset {\mathbb R}^{n+m}$ containing $\{0\} \times {\mathbb R}^{n-2}$};\nonumber\\ 
\end{eqnarray*}
\end{itemize}
Here the constant $\overline{C}_{1} \in (0, \infty)$ depends only on $n$ and $m$, and the constants $\overline{\k}, \overline{\nu}, \overline{C}_{0}, \overline{C}_{2} \in (0, \infty)$  and $\mu \in (0, 1)$ each depends only on 
$n$, $m$ and $q$. 
\end{lemma}

\begin{proof}  Fix $q \geq 2$, $\theta \in (0, 1/2)$ and $M \in [1, \infty).$ For $k=1, 2, \ldots$, let $\eta_{k}, \b_{k}, \delta_{k} \in (0, 1)$ be such that $\eta_{k} \to 0$, $\b_{k} \to 0$, $\delta_{k} \to 0$ as $k \to \infty$; $T_{k}$ be a locally area minimizing $n$-dimensional rectifiable current in ${\mathbf B}_{1}(0)$;  
${\mathbf C}_{k} \in {\mathcal C}_{p, {\mathbf q}}$ where $p \in \{2, \ldots, q\}$ and ${\mathbf q}  = (q_{1}, \ldots, q_{p}) \in {\mathfrak M}_{p}$ are independent of 
$k,$ such that conditions (1)-(8) of Section~\ref{fine-blowup-prelim} hold for each $k$, with $q_{i}$ in place of $q_{i}^{(k)}$ and with $\e_{k} \to 0^{+}$, and also such that conclusion (A) of the present lemma with $T_{k}$ in place of $T$ and $\d_{k}$ in place of $\overline{\d}$ fails and hence condition (\ref{excess-1}) holds for each $k$. To prove the lemma, it suffices to show that for each $k$, there are an orthogonal rotation $\G_{k}$ of ${\mathbb R}^{n+m}$ and a cone ${\mathbf C}_{k}^{\prime} \in {\mathcal C}_{p^{\prime}, {\mathbf q}^{\prime}}$ for some $p^{\prime} \in \{2, \ldots, q\}$ and some ${\mathbf q}^{\prime} \in {\mathfrak M}_{p^{\prime}},$ such that after passing to an appropriate subsequences of $(k)$, conclusion (B) holds with $T_{k},$ ${\mathbf C}_{k}$, ${\mathbf C}^{\prime}_{k}$, $\G_{k}$ in place of $T$, ${\mathbf C}$, ${\mathbf C}^{\prime}$, $\G$ respectively, and with fixed constants $\overline{\k}, \overline{\nu}, \overline{C}_{0}, \overline{C}_{2} \in (0, \infty)$  and $\mu \in (0, 1)$ depending only on 
$n$, $m$ and $q$ and a fixed constant $\overline{C}_{1} \in (0, \infty)$ depending only on $n$, $m$. 

By the definition of ${\mathcal C}_{p, {\mathbf q}}$, we have that 
$${\mathbf C}_{k} = \sum_{j=1}^{p} q_{j}\llbracket P^{(k)}_{j} \rrbracket$$
where for each $k$, $P_{j}^{(k)} = \{A_{j}^{(k)}x, x, y) \, : \, (x, y) \in {\mathbb R}^{2} \times {\mathbb R}^{n-2}\}$ for $j \in \{1, 2, \ldots, p\}$ and 
for some constant $m \times 2$ distinct matrices $A_{1}^{(k)}, A_{2}^{(k)}, \ldots, A^{(k)}_{p}$. Choosing an appropriate sequence of numbers $\t_{k}$ with $\t_{k} \to 0^{+}$, we 
obtain a fine blow-up $w \in {\mathfrak B}_{p, {\mathbf q}}(M)$ of $(T_{k})$ relative to $({\mathbf C}_{k})$ as described in  Section~\ref{fine-blowup-prelim}. By Theorem~\ref{blowup-decay-unique}, there is $\widetilde{\psi}  = (\widetilde{\psi}_{1}, \widetilde{\psi}_{2}, \ldots, \widetilde{\psi}_{p}) \in \widetilde{\mathfrak L}_{p, {\mathbf q}},$ with $\widetilde{\psi}_{j}(x, y) = \sum_{\ell=1}^{q_{j}} \llbracket \psi_{j, \ell}(x) +  L_{1}(y) + A_{j} L_{2}(y)\rrbracket$ for $(x, y) \in {\mathbb R}^{2} \times {\mathbb R}^{n-2}$ and $j \in \{1,2, \ldots, p\}$, where $\psi = (\psi_{1}, \psi_{2}, \ldots, \psi_{p}) \in {\mathfrak L}_{p, {\mathbf q}}$, $L_{1} \,: \, {\mathbb R}^{n-2} \to {\mathbb R}^{m}$ and $L_{2} \, : \, {\mathbb R}^{n-2} \to {\mathbb R}^{2}$  are (single-valued) linear functions and $A_{1}, A_{2}, \ldots, A_{p}$ are  constant $m \times 2$ matrices,  
such that (\ref{blowup-decay-unique-1}), (\ref{blowup-decay-unique-2}) and (\ref{blowup-decay-unique-3}) hold for some constants $C, \overline{C} \in (0, \infty)$ and $\mu \in (0, 1)$ depending only on $n$, $m$ and $q$. Set 
$$\G_{k} = e^{E_k M_1 + E_k M_2/\widehat{E}_k}$$
where 
$${\widehat E}_{k}^{2} = {\widehat E}_{T_{k}}^{2}= \int_{{\mathbf B}_{1}(0)} {\rm dist}^{2}(X,  P_{0}) \,d\|T_{k}\|(X), \;\;\; 
E_{k}^{2} =  E_{T_{k}}^{2} = \int_{{\mathbf B}_{1}(0)} {\rm dist}^{2} (X, {\rm spt} \, {\mathbf C}_{k}) \, d\|T_{k}\|(X),$$  
	\begin{equation*}
		M_1 = \left[\begin{array}{ccc} 0 & 0 & L_1 \\ 0 & 0 & 0 \\ -L_1^T & 0 & 0 \end{array}\right] \quad 
		{\rm and } \quad M_2 = \left[\begin{array}{ccc} 0 & 0 & 0 \\ 0 & 0 & -L_2 \\ 0 & L_2^T & 0 \end{array}\right], 
	\end{equation*}
	representing $L_1$ as an $m \times (n-2)$ matrix and $L_2$ as a $2 \times (n-2)$ matrix. Set 
$${\mathbf C}_{k}^{\prime} = \sum_{j=1}^{p}\sum_{\ell =1}^{q_{j}} \llbracket P^{(k)}_{j, \ell} \rrbracket$$
where $P^{(k)}_{j, \ell}  = \{(A_j^{(k)} x + E_k \psi_{j,\ell}(x), x,y) \,: \, (x,y) \in {\mathbb R}^{2} \times {\mathbb R}^{n-2}\}.$ In view of the bounds (\ref{blowup-decay-unique-2}), (\ref{blowup-decay-unique-3}), it readily follows that conclusions (B)(a) and B(b) hold with $\G_{k}$ in place of $\G$, $T_{k}$ in place of $T$, 
${\mathbf C}_{k}$ in place of ${\mathbf C}$ and ${\mathbf C}_{k}^{\prime}$ in place of ${\mathbf C}^{\prime},$ and  for some constants $\overline{\k} = \overline{\k}(n, m ,q)$, $\overline{C}_{0} = \overline{C}_{0}(n, m, q)$. It follows from assertions (B), (C) of Section~\ref{fine-blowup-prelim}  and the decay estimate (\ref{blowup-decay-unique-1}) that for sufficiently large $k$, conclusion (c) holds with $T_{k}$, ${\mathbf C}^{\prime}_{k}$, $\G_{k}$ in place of $T$, ${\mathbf C}$, $\G$ and 
with appropriate constants $\overline{\nu}$, $\mu$ depending only on $n$, $m$, $q$. Finally, arguing exactly as in \cite[p.~941]{Wic14} (giving estimate (13.14) therein), utilising again assertion (C) of Section~\ref{fine-blowup-prelim}, we see that conclusion (d) holds with $T_{k}$, ${\mathbf C}_{k}$, $\G_{k}$ in place of 
$T$, ${\mathbf C}$, $\G$, and with appropriate constants $\overline{C}_{1} = \overline{C}_{1}(n, m)$ and $\overline{C}_{2} = \overline{C}_{2}(n, m, q)$. 
\end{proof}

For our purposes in the next section, where we establish for an area minimizing current $T$, ${\mathcal H}^{n-2}$ a.e.\ uniqueness of tangent cones and countable $(n-2)$ rectifiability of the set of all singularities where $T$ does not rapidly decay to a plane, we need a version of Theorem~\ref{excess-improvement} in which Hypothesis~($\star\star$) is relaxed to the following weaker assumption: the fine excess of $T$ relative to a cone ${\mathbf C} \in \cup_{p=2}^{q}{\mathcal C}_{q, p}$ is significantly smaller than the coarse excess of $T$ relative to any plane ${\mathbf P}$, i.e.\ the condition  
$Q(T, {\mathbf C}, {\mathbf B}_{1}(0)) \leq \beta \inf_{{\mathbf P} \in {\mathcal C}_{q, 1}} \, E(T, {\mathbf P}, {\mathbf B}_{1}(0))$ for a fixed, appropriately small constant $\beta$. Relaxing Hypothesis~($\star\star$) in this manner can readily be achieved by employing Theorem~\ref{excess-improvement} itself, provided we are content to require that conclusion~(c) (the improvement of the fine excess) and conclusion~(d) of Theorem~\ref{excess-improvement} hold at one of a fixed number of (in fact $q-1$) smaller scales 
$\theta_{1}, \theta_{2}, \ldots, \theta_{q-1}.$ For the purpose of deducing a uniform decay estimate for $T$ by iteratively applying the lemma (as we do in the next section), allowing a fixed number of scales to choose from at each step of the iteration is just as good as single scale improvement at each stage.

\begin{lemma}\label{excess-improvement-final} 
Let $q$ be an integer such that $q \geq 2,$ $M \in [1, \infty)$ and let $\th_{1}, \th_{2}, \ldots, \th_{q-1} \in (0, 1/4)$ be distinct numbers. 
There exist numbers 
$$\eta = \eta(n, m, q, M, \th_{1}, \ldots, \th_{q-1}) \in (0, 1/2), \;\; \b = \b(n, m, q, M, \th_{1}, \ldots, \th_{q-1}) \in (0, 1/2)$$ and $\delta = \delta(n, m, q, M, \th_{1}, \ldots, \th_{q-1}) \in (0, 1/2)$ such that the following holds true: let ${\mathbf C} \in 
\cup_{p=2}^{q} \, {\mathcal C}_{q, p}$ and let $T$ be an $n$-dimensional locally area minimizing rectifiable current in ${\bf B}_{1}(0)$ with 
$$\partial \, T \llcorner {\mathbf B}_{1}(0) = 0,  \omega_{n}^{-1}\|T\|({\mathbf B}_{1}(0)) < q + 1/2 \;\; \mbox{and} \;\; \Theta(T, 0) \geq q.$$ Suppose that the following two conditions hold:  
\begin{itemize}
\item[(i)] Hypothesis~($\dag$) holds with $\eta$ in place of $\eta_{0}$, i.e.\ we have that $$E(T, {\mathbf P}_{0}, {\mathbf B}_{1}(0)) < \eta \;\; \mbox{and} \;\; E(T, {\mathbf P}_{0}, {\mathbf B}_{1}(0)) \leq  M \inf_{{\mathbf P} \in {\mathcal C}_{q, 1}} E(T, {\mathbf P}, {\mathbf B}_{1}(0));$$ 
\item [(ii)] $$Q(T, {\mathbf C}, {\mathbf B}_{1}(0)) < \beta E(T, {\mathbf P}_{0}, {\mathbf B}_{1}(0)).$$ 
\end{itemize}
Then either 
\begin{itemize}
\item[(A)] ${\mathbf B}_{\d}(0,z) \cap \{Z \, : \, \Theta \, (T, Z) \geq q\}  = \emptyset$
for some point $(0, z) \in \{0\} \times {\mathbb R}^{n-2} \cap {\mathbf B}_{1/2}(0)$ or, 
 
\item[(B)] there exist an orthogonal rotation $\G$ of ${\mathbb R}^{n+m}$ and a cone ${\mathbf C}^{\prime} \in \cup_{p=2}^{q} \, {\mathcal C}_{q, p}$ such that, writing 
${\hat E}_{T} =E(T, {\mathbf P}_{0}, {\mathbf B}_{1}(0))$,
the following hold:
\begin{eqnarray*}
&&\hspace{-.5in}({\rm a})\;\;\;\; |e_{i} - \G(e_{i})| \leq \k \, Q(T, {\mathbf C}, {\mathbf B}_{1}(0))  \;\; {for} \;\; i=1,  \ldots, m \;\; \mbox{and}\nonumber\\ 
&&\;\; |e_{m+j} - \G(e_{m+j})| \leq \k \, {\hat E}_{T}^{-1}Q(T, {\mathbf C}, {\mathbf B}_{1}(0)) \;\; \mbox{for  $j=1,  \ldots, n$};\nonumber\\   
&&\hspace{-.5in}({\rm b})\;\;\;\; {\rm dist}_{\mathcal H} \, ({\rm spt} \, {\mathbf C}^{\prime} \cap {\mathbf B}_{1}, {\rm spt} \, {\mathbf C} \cap {\mathbf B}_{1}) 
\leq C_{0} \, Q(T, {\mathbf C}, {\mathbf B}_{1}(0));  
\end{eqnarray*}
and for some $j \in \{1, 2, \ldots, q-1\},$
\begin{eqnarray*}
&&\hspace{.23in}({\rm c})\;\;\;\;\th_{j}^{-n-2}\int_{\G\left({\mathbf B}_{\th_{j}/2} \setminus \{r(X) \leq \th_{j}/16\}\right)} {\rm dist}^{2} \, (X, {\rm spt} \, T) \, d\|\G_{\#} \, {\mathbf C}^{\prime}\|(X)\nonumber\\
&&\hspace{1.5in}+ \;\th_{j}^{-n-2}\int_{{\mathbf B}_{\th_{j}}} {\rm dist}^{2} \, (X, {\rm spt} \, \G_{\#} \, {\mathbf C}^{\prime}) \, d\|T\|(X)  \leq \nu_{j}\th_{j}^{2\mu} \, 
Q^{2}(T, {\mathbf C}, {\mathbf B}_{1}(0));\nonumber\\
&&\hspace{.23in}({\rm d})\;\;\left(\th_{j}^{-n-2}\int_{{\mathbf B}_{\th_{j}}} {\rm dist}^{2}\, (X, P)\, d\|\G^{-1}_{\#} \,T\|(X)\right)^{1/2}\nonumber\\ 
&&\hspace{2.5in}\geq C_{1} \,{\rm dist}_{\mathcal H}\, ({\rm spt} \, {\mathbf C} \cap {\mathbf B}_{1}, P \cap {\mathbf B}_{1}) - C_{2} \, Q(T, {\mathbf C}, {\mathbf B}_{1}(0))\\
&&\hspace{.5in}\mbox{for any $n$-dimensional plane $P \subset {\mathbb R}^{n+m}$ containing $\{0\} \times {\mathbb R}^{n-2}$}.
\end{eqnarray*}
\end{itemize}
Here $\mu = \mu(n, m, q) \in (0, 1)$; $C_{1} = C_{1}(n, m) \in (0, \infty)$; the constants $\k, C_{0}, C_{2} \in (0, \infty)$ depend only on $n$ and $m$ in case $q=2$ and only on $n$, $m$, $q$ and $\th_{1}, \th_{2}, \ldots, \th_{q-2}$ in case $q \geq 3$; $\nu_{1} = \nu_{1}(n, m, q) \in (0, \infty)$ and, in case $q\geq 3$,  $\nu_{j} = \nu_{j}(n, m, q, \th_{1}, \ldots, \th_{j-1}) \in (0, \infty)$ for each $j = 2, \ldots, q-1$. (In particular, for $q \geq 3$, $\nu_{j}$ is independent of $\th_{j}, \th_{j+1}, \ldots, \th_{q-1}$ for each $j=2, \ldots, q-1$.)  
\end{lemma}

\begin{proof} 
We assert the following slightly more refined version of the lemma:

\noindent 
\emph{Claim: Let $p \in \{2, \ldots, q\}$ and let $\th_{1}, \th_{2}, \ldots, \th_{p-1} \in (0, 1/4).$ There exist $$\eta^{(p)} = \eta^{(p)}(n, m, q, M, \th_{1}, \ldots, \th_{p-1}) \in (0, 1/2), \;\; \b^{(p)} = \b^{(p)}(n, m, q, M, \th_{1}, \ldots, \th_{p-1}) \in (0, 1/2) \;\; \mbox{and}$$  $$\delta^{(p)} = \delta^{(p)}(n, m, q, M, \th_{1}, \ldots, \th_{p-1}) \in (0, 1/2)$$ 
such that if ${\mathbf C}$, $T$ with ${\mathbf C} \in \cup_{p^{\prime}=2}^{p} \, {\mathcal C}_{q, p^{\prime}}$ satisfy the hypotheses of the lemma with $\eta^{(p)}$, $\beta^{(p)}$ in place of $\eta$, $\beta$, then conclusion~(A) or conclusion (B) of the lemma holds, with the following choices: in conclusion (A), $\delta^{(p)}$ is taken in place of $\delta$; in conclusion (B), the cone ${\mathbf C}^{\prime} \in \cup_{p=2}^{q} {\mathcal C}_{q, p}$; in conclusions~(B)(a), (B)(b), constants $\k^{(p)}$, $C_{0}^{(p)}$ are taken in place of $\k$, $C_{0}$ and in conclusions~(B)(c), (B)(d), the index $j$ is such that $j \in \{1, \ldots, p-1\}$, 
the constants $\mu$, $C_{1}$ are such that $\mu = \mu(n, m, q) \in (0, 1),$ $C_{1} = C_{1}(n, m) \in (0, \infty)$ and constants $\nu_{j}^{(p)}$, $C_{2}^{(p)} \in (0, \infty)$ are taken in place of $\nu_{j}$, $C_{2}$, where:  $\k^{(p)}, C_{0}^{(p)}, C_{2}^{(p)}$ depend only on $n$, $m$ and $q$ in case $p=2$ and only on $n$, $m$, $q$ and $\th_{1}, \th_{2}, \ldots, \th_{p-2}$ in case $3\leq p \leq q$; $\nu_{1} = \nu_{1}(n, m, q) \in (0, \infty)$ and, in case $q\geq 3$, $\nu_{j}^{(p)} = \nu_{j}^{(p)}(n, m, q, \th_{1}, \ldots, \th_{j-1})$ for each $j = 2, 3, \ldots, p-1$.} 

It is clear that the lemma as stated follows from the claim, by simply setting $\eta = \min \{\eta^{(p)}: p=2, \ldots, q\}$, $\beta = \min \{\beta^{(p)}: p=2, \ldots, q\}$, $\delta = \min \{\delta^{(p)}: p=2, \ldots, q\},$ $\k = \max \{\k^{(p)}: p=2, \ldots, q\}$, $C_{0} = \max \{C_{0}^{(p)}: p=2, \ldots, q\}$, $C_{2} = \max  \{C_{2}^{(p)}: p=2, \ldots, q\}$ and, for $q \geq 3$ and $j \in \{2, \ldots, q-1\}$, $\nu_{j} = \max \{\nu_{j}^{(p)}: p=2, \ldots, q\}.$

To see the claim, we argue by induction on $p$ (keeping $q \geq 2$ fixed). First set $\mu  = \mu(n, m, q)$ to be the constant as in Lemma~\ref{excess-improvement}, and also set 
$C_{1} = \overline{C}_{1}(n,m)$ and $\nu_{1} = \overline{\nu}(n, m, q)$ where $\overline{C}_{1}$, $\overline{\nu}$ are as in Lemma~\ref{excess-improvement}. If $p=2$, the claim follows directly from Lemma~\ref{excess-improvement} taken with 
$\th = \th_{1}$, provided we take $\delta^{(2)} = \overline{\delta}(n, m, q, M, \th_{1})$, $\eta^{(2)}  = \overline{\eta}(n, m, q, M, \th_{1})$ and 
$\b^{(2)} = \overline{\b}(n, m, q, M, \th_{1})$, $\k^{(2)} = \overline{\k}(n, m, q)$, $C_{0}^{(2)} = \overline{C}_{0}(n, m, q)$, 
$C_{2}^{(2)} = \overline{C}_{2}(n, m , q)$ where $\overline{\delta}$, $\overline{\eta}$, $\overline{\b},$ $\overline{\k}$,  $\overline{C}_{0},$ $\overline{C}_{2}$ are as in Lemma~\ref{excess-improvement}. 

Let $p_{1} \in \{3, \ldots, q\}$ and assume (the induction hypothesis) that the claim holds with $p_{1}-1$ in place of $p$. We wish to show that the claim holds with $p=p_{1}$, so let $\th_{1}, \th_{2}, \ldots, \th_{p_{1}-1} \in (0, 1/4).$ 
By applying the induction hypothesis (with $\th_{2}, \ldots, \th_{p_{1}-1}$ in place of $\th_{1}, \ldots, \th_{p_{1}-2}$), we obtain constants $\eta^{(p_{1}-1)} = \eta^{(p_{1}-1)}(n, m, q, M, \th_{2}, \ldots, \th_{p_{1}-1})$, $\b^{(p_{1}-1)} = \b^{(p_{1}-1)}(n, m, q, M, \th_{2}, \ldots, \th_{p_{1}-1})$, 
$\delta^{(p_{1}-1)} = \delta^{(p_{1}-1)}(n, m, q, M, \th_{2}, \ldots, \th_{p_{1}-1}),$ $\mu = \mu(n, m, q),$ $C_{1} = C_{1}(n, m)$, 
$\k^{(p_{1}-1)} = \k^{(p_{1}-1)}(n, m, q, \th_{2}, \ldots, \th_{p_{1}-2})$, $C_{0}^{(p_{1}-1)} =  C_{0}^{(p_{1}-1)}(n, m, q, \th_{2}, \ldots, \th_{p_{1}-2})$, 
$C_{2}^{(p_{1}-1)} =  C_{2}^{(p_{1}-1)}(n, m, q, \th_{2}, \ldots, \th_{p_{1}-2})$ and $\nu_{j}^{(p_{1}-1)} = \nu_{j}^{(p_{1}-1)}(n, m, q, \th_{2}, \ldots, \th_{j})$ for $j = 2, 3, \ldots, p_{1}-2$ so that the claim is true with $p = p_{1}-1$ and with $\th_{2}, \ldots, \th_{p_{1}-1}$ in place of $\th_{1}, \ldots, \th_{p_{1}-2}$.  We assert that the claim is true with $p = p_{1}$ and with  constants  
$$\eta^{(p_{1})} = \min\{\eta^{(p_{1} - 1)}, \overline{\eta}(n, m, q, M, \th_{1})\}, \;\;  \b^{(p_{1})} = \frac{1}{2^{q}}(\overline{\b}(n, m, q, M, \th_{1}))^{q-1}\b^{(p_{1} - 1)},$$ 
$$\delta^{(p_{1})} = \min\{\delta^{(p_{1} - 1)}, \overline{\delta}(n, m, q, M, \th_{1})\},$$ 
$$\k^{(p_{1})} = \max\left\{\frac{2\k^{(p_{1}-1)}}{\overline{\b}(n, m, q, M, \th_{1})}, \overline{\k}(n, m, q)\right\},  \;\;  
C_{0}^{(p_{1})} = \max\left\{\overline{C} + \frac{2(C_{0}^{(p_{1}-1)} +\overline{C})}{\overline{\b}(n, m, q, M, \th_{1})}, \overline{C_{0}}(n, m, q)\right\},$$ 
$$C_{2}^{(p_{1})} = \max\left\{C_{1}\overline{C} + \frac{2(C_{2}^{(p_{1}-1)} +C_{1}\overline{C})}{\overline{\b}(n, m, q, M, \th_{1})}, \overline{C_{2}}(n, m, q)\right\} \;\; \mbox{and}$$  
$$\nu_{j}^{(p_{1})} = \max\left\{\frac{4\nu_{j}^{(p_{1}-1)}}{\overline{\b}^{2}(n, m, q, M, \th_{1})}, \overline{\nu}(n, m, q)\right\} \;\; \mbox{for $j=1, 2, \ldots, p_{1}-1$},$$ where $\overline{\eta}, \overline{\b}$, $\overline{\d},$ $\overline{\k}$, $\overline{C_{0}}$, $\overline{C_{2}}$ and $\overline{\nu}$ are as in Lemma~\ref{excess-improvement} taken with 
$\th = \th_{1}$, and $\overline{C} = \overline{C}(n, q)$ is to be specified momentarily. To see this,  let  
${\mathbf C} \in \cup_{p^{\prime} = 1}^{p_{1}} {\mathcal C}_{q, p^{\prime}}$ and suppose that ${\mathbf C}$, $T$ are such that 
$T$ is an $n$-dimensional locally area minimizing rectifiable current in ${\bf B}_{1}(0)$, $\partial \, T \llcorner {\mathbf B}_{1}(0) = 0,$  
$\omega_{n}^{-1}\|T\|({\mathbf B}_{1}(0)) < q + 1/2,$  
$E(T, {\mathbf P}_{0}, {\mathbf B}_{1}(0)) < \eta^{(p_{1})}$,  $E(T, {\mathbf P}_{0}, {\mathbf B}_{1}(0)) < M \inf_{{\mathbf P} \in {\mathcal C}_{q, 1}} E(T, {\mathbf P}, {\mathbf B}_{1}(0))$ and $$Q(T, {\mathbf C}, {\mathbf B}_{1}(0)) < \beta^{(p_{1})} E(T, {\mathbf P}_{0}, {\mathbf B}_{1}(0)).$$ 
If ${\mathbf C} \in \cup_{p^{\prime} = 1}^{p_{1}-1} {\mathcal C}_{q, p^{\prime}},$ 
then the conclusions of the claim are immediate by the induction hypothesis, so assume that ${\mathbf C} \in {\mathcal C}_{q, p_{1}}.$ In this case, 
if additionally we have Hypothesis~($\star\star$) i.e.\ that  
\begin{equation*} \label{excess-improvement-final-1}
Q(T, {\mathbf C}, {\mathbf B}_{1}(0)) < \overline{\b}(n, m, q, M, \th_{1}) 
\inf_{{\mathbf C}^{\prime} \in \cup_{p^{\prime}=1}^{p_{1} - 1} {\mathcal C}_{q, p^{\prime}}} \, Q(T, {\mathbf C}^{\prime}, {\mathbf B}_{1}(0)),  
\end{equation*}
then we can apply Lemma~\ref{excess-improvement} (with $\th = \th_{1}$) to see that the conclusions of the claim hold again. If on the other hand 
\begin{equation} \label{excess-improvement-final-2}
Q(T, {\mathbf C}, {\mathbf B}_{1}(0)) \geq \overline{\b}(n, m, q, M, \th_{1}) 
\inf_{{\mathbf C}^{\prime} \in \cup_{p^{\prime}=1}^{p_{1} - 1} {\mathcal C}_{q, p^{\prime}}} \, Q(T, {\mathbf C}^{\prime}, {\mathbf B}_{1}(0)),  
\end{equation}
then arguing as in Remark~\ref{tildeC rmk}, we can find $\overline{p} \in \{1,\ldots,p_{1}-1\}$ and a cone $\overline{\mathbf C} \in \mathcal{C}_{q,\overline{p}}$ such that 
\begin{equation}
	\label{excess-improvement-final-2-1} Q(T, \overline{\mathbf C}, \mathbf{B}_{1}(0)) \leq 2^{q-1} (\overline{\beta}(n, m, q, M, \th_{1}))^{2-q} 
		\inf_{\mathbf{C}' \in \bigcup_{p'=1}^{p_{1}-1} \mathcal{C}_{q,p'}} Q(T, \mathbf{C}', \mathbf{B}_{1}(0))
\end{equation}
and either $\overline{p} = 1$ or 
\begin{equation}
\label{excess-improvement-final-2-2} 
\mbox{$\overline{p} \geq 2$ and} \;\; 	Q(T, \overline{\mathbf C}, \mathbf{B}_{1}(0)) \leq 
		\overline{\beta}(n, m, q, M, \th_{1})\inf_{\mathbf{C}' \in \bigcup_{p'=1}^{\overline{p}-1} \mathcal{C}_{q,p'}} Q(T, \mathbf{C}', \mathbf{B}_{1}(0)) .
\end{equation}
Note that then we must have  $$Q(T, \overline{\mathbf C}, {\mathbf B}_{1}(0)) < \beta^{(p_{1}-1)} E(T, {\mathbf P}_{0}, {\mathbf B}_{1}(0)),$$ 
for otherwise, by \eqref{excess-improvement-final-2} and \eqref{excess-improvement-final-2-1},  we would have 
\begin{eqnarray*}
&\hspace{-1in}\beta^{(p_{1}-1)}E(T, {\mathbf P}_{0}, {\mathbf B}_{1}(0)) \leq Q(T, \overline{\mathbf C}, {\mathbf B}_{1}(0))\\ 
&\hspace{2.1in}\leq  2^{q-1}(\overline{\b}(n, m, q, \th_{1}))^{1-q} Q(T, {\mathbf C}, {\mathbf B}_{1}(0))\\ 
&\hspace{3in}\leq  2^{q-1}(\overline{\beta}(n, m, q, M, \th_{1}))^{1-q}\beta^{(p_{1})} E(T, {\mathbf P}_{0}, {\mathbf B}_{1}(0)),
\end{eqnarray*}
 which is impossible in view of the definition of $\b^{(p_{1})}$. Thus we can apply the induction hypothesis again to deduce first that the conclusions of the claim hold (for some cone ${\mathbf C}^{\prime} \in \cup_{p=2}^{q} {\mathcal C}_{q, p}$) with $\overline{\mathbf C}$ in place of ${\mathbf C}$ and with  the constants $\eta^{(p_{1} - 1)}$, $\b^{(p_{1}-1)}$, $\d^{(p_{1}-1)}$, $\k^{(p_{1}-1)}$, $C_{0}^{(p_{1}-1)}$, $C_{2}^{(p_{1}-1)}$, $\nu_{j}^{(p_{1}-1)}$ in place of 
$\eta^{(p)}$, $\b^{(p)}$, $\d^{(p)}$, $\k^{(p)}$, $C_{0}^{(p)}$, $C_{2}^{(p)}$, $\nu_{j}^{(p)}$; consequently, in view of \eqref{excess-improvement-final-2} and \eqref{excess-improvement-final-2-1}, together with the fact that 
\begin{equation}\label{excess-improvement-final-dist}
{\rm dist}_{\mathcal H} \, ({\rm spt} \, {\mathbf C} \cap {\mathbf B}_{1}(0), {\rm spt} \, \overline{\mathbf C} \cap {\mathbf B}_{1}(0)) \leq \overline{C} (Q(T, {\mathbf C}, {\mathbf B}_{1}(0)) + Q(T, \overline{\mathbf C}, {\mathbf B}_{1}(0)))
\end{equation}
where $\overline{C} = \overline{C}(n, q)$ (which we shall justify momentarily), 
the claim also holds with $p= p_{1}$, and with the same cone ${\mathbf C}^{\prime}$ and with the choice of constants $\eta^{(p_{1})}$, $\b^{(p_{1})}$, $\d^{(p_{1})}$, $\k^{(p_{1})}$, $C_{0}^{(p_{1})}$, $C_{2}^{(p_{1})}$, $\nu_{j}^{(p_{1})}$ as defined above (taking $\overline{C}$ in the definitions of $C_{0}^{(p_{1})},$ $C_{2}^{(p_{1})}$ to be the constant in \eqref{excess-improvement-final-dist}). 

To see \eqref{excess-improvement-final-dist}, write $\overline{\mathbf C} = \sum_{i=1}^{\overline{p}} {q}_i \llbracket {P}_i \rrbracket$ for 
${P}_i= \{z = {A}_i \,x\}$.  Note that in view of \eqref{excess-improvement-final-2-2}, the current $T$ has a ``graphical representation relative to 
$\overline{\mathbf C}$'' in the sense of Theorem~\ref{graphrep close2plane thm}.  In particular, using the notation of Theorem~\ref{graphrep close2plane thm}, there is a good set $K$ and Lipschitz approximation $u_i: B_{1/2}(0) \cap \{r \geq 1/8\} \rightarrow \mathcal{A}_{q_i}(\mathbb{R}^m)$ ($1 \leq i \leq \overline{p}$) relative to $\overline{\mathbf C}$.  Hence by Theorem~\ref{graphrep close2plane thm}(b), 
	\begin{align*}
		&\int_{\mathbf{B}_{1/2}(0) \cap \{r > 1/8\}} \op{dist}^2(X,\op{spt} \overline{\mathbf C}) \,d\|\mathbf{C}\|(X) 
		\\ \leq\,& 2 \int_{\mathbf{B}_{1/2}(0) \cap \{r > 1/8\}} \op{dist}^2(X,\op{spt} T) \,d\|\mathbf{C}\|(X) 
			\\&+ C \sup_{X \in \op{spt} T \cap \mathbf{B}_{3/4}(0) \cap \{r > 1/16\}} \op{dist}^2(X,\op{spt} \overline{\mathbf C})
		\\ \leq\,& 2 \int_{\mathbf{B}_{1/2}(0) \cap \{r > 1/8\}} \op{dist}^2(X,\op{spt} T) \,d\|\mathbf{C}\|(X) 
			+ C E^2(T,\overline{\mathbf C},\mathbf{B}_1(0)) 
		\\ \leq\,& 2 Q^2(T,\mathbf{C},\mathbf{B}_1(0)) + C Q^2(T,\overline{\mathbf C},\mathbf{B}_1(0)) 
	\end{align*}
and by Theorem~\ref{graphrep close2plane thm}(c), 
	\begin{align*}
		&\int_{\mathbf{B}_{1/2}(0) \cap \{r > 1/8\}} \op{dist}^2(X,\op{spt} \mathbf{C}) \,d\|\overline{\mathbf C}\|(X) 
		\\ \leq\,& 2 \int_{B_{1/2}(0) \cap \{r > 1/8\}} \sum_{i=1}^{\overline{p}} \sum_{\ell=1}^{q_i} 
			\op{dist}^2( (A_i \, x + u_{i,\ell}(x,y),x,y), \op{spt} \mathbf{C}) \,d\mathcal{L}^n(x,y) 
			\\&+ C \sum_{i=1}^{\overline{p}} \sup_{B_{1/2}(0) \cap \{r \geq 1/8\}} |u_i|^2 
		\\ \leq\,& 2 \int_{\mathbf{B}_1(0)} \op{dist}^{2}(X,\op{spt} \mathbf{C}) \,d\|T\|(X) + C E^2(T,\overline{\mathbf C},\mathbf{B}_1(0)) 
		\\ \leq\,& 2 Q^2(T,\mathbf{C},\mathbf{B}_1(0)) + C Q^2(T,\overline{\mathbf C},\mathbf{B}_1(0)) ,
	\end{align*}
	where $C = C(n,m,q) \in (0,\infty)$ are constants.  Since the cones ${\mathbf C}$, $\overline{\mathbf C}$ are made up of planes, the Hausdorff distance bound as in \eqref{excess-improvement-final-dist} immediately follows.

This completes the inductive proof of the claim. The claim readily implies the lemma as already indicated.\end{proof}

Lemma~\ref{excess-improvement-final}  is the main excess decay result that handles the ``degenerate'' case, i.e.\ the case when the current $T$ is close to a plane (and much closer to a cone ${\mathbf C} \in \cup_{p=2}^{q} {\mathcal C}_{q, p}$).  We also have the following 
version of excess decay that is applicable to the non-degenerate case, i.e.\ when $T$ lies far from any plane. 
  
\begin{lemma}\label{excess-improvement-final-1} 
Let $q$ be an integer such that $q \geq 2,$ $\eta \in (0, 1)$, and let $\th_{1}, \th_{2}, \ldots, \th_{q-1} \in (0, 1/4)$ be distinct numbers.  
There exist numbers  $\b_{1} = \b_{1}(n, m, q, \eta, \th_{1}, \ldots, \th_{q-1}) \in (0, 1/2)$,   and $\delta_{1} = \delta_{1}(n, m, q, \eta, \th_{1}, \ldots, \th_{q-1}) \in (0, 1/2)$ such that the following holds true: if ${\mathbf C} \in 
\cup_{p=2}^{q} \, {\mathcal C}_{q, p}$, $T$ is an $n$-dimensional locally area minimizing rectifiable current in ${\bf B}_{1}(0)$ with 
$$\partial \, T \llcorner {\mathbf B}_{1}(0) = 0, \;\; \omega_{n}^{-1}\|T\|({\mathbf B}_{1}(0)) < q + 1/2 \;\; \mbox{and} \;\; \Theta(T, 0) \geq q,$$ and if 
\begin{itemize}
\item[{\rm(i)}] Hypothesis~($\dag\dag$) holds with $\eta$ in place of $\eta_{0}$, i.e.\ 
$$\inf_{{\mathbf P} \in {\mathcal C}_{q, 1}} E(T, {\mathbf P}, {\mathbf B}_{1}(0)) \geq \eta, \;\; \mbox{and}$$ 
\item [{\rm(ii)}] $$Q(T, {\mathbf C}, {\mathbf B}_{1}(0)) < \beta_{1},$$ 
\end{itemize}
then either 
\begin{itemize}
\item[(A)] ${\mathbf B}_{\d_{1}}(0,z) \cap \{Z \, : \, \Theta \, (T, Z) \geq q\}  = \emptyset$
for some point $(0, z) \in \{0\} \times {\mathbb R}^{n-2} \cap {\mathbf B}_{1/2}(0)$ or, 
 
\item[(B)] there exist an orthogonal rotation $\G$ of ${\mathbb R}^{n+m}$ and a cone ${\mathbf C}^{\prime} \in \cup_{p=2}^{q} \, {\mathcal C}_{q, p}$ such that,
the following hold:
\begin{eqnarray*}
&&\hspace{-1.8in}({\rm a})\;\;\;\; |I - \G| \leq \k_{1} \, Q(T, {\mathbf C}, {\mathbf B}_{1}(0));\\  
&&\hspace{-1.8in}({\rm b})\;\;\;\; {\rm dist}_{\mathcal H} \, ({\rm spt} \, {\mathbf C}^{\prime} \cap {\mathbf B}_{1}, {\rm spt} \, {\mathbf C} \cap {\mathbf B}_{1}) 
\leq C_{0}^{(1)} \, Q(T, {\mathbf C}, {\mathbf B}_{1}(0));  
\end{eqnarray*}
and for some $j \in \{1, 2, \ldots, q-1\},$
\begin{eqnarray*}
&&\hspace{.2in}({\rm c})\;\;\;\;\th_{j}^{-n-2}\int_{\G\left({\mathbf B}_{\th_{j}/2} \setminus \{r(X) \leq \th_{j}/16\}\right)} {\rm dist}^{2} \, (X, {\rm spt} \, T) \, d\|\G_{\#} \, {\mathbf C}^{\prime}\|(X)\nonumber\\
&&\hspace{1.5in}+ \;\th_{j}^{-n-2}\int_{{\mathbf B}_{\th_{j}}} {\rm dist}^{2} \, (X, {\rm spt} \, \G_{\#} \, {\mathbf C}^{\prime}) \, d\|T\|(X)  \leq \nu_{j}^{(1)}\th_{j}^{2\mu_{1}} \, 
Q^{2}(T, {\mathbf C}, {\mathbf B}_{1}(0)).\nonumber\\
\end{eqnarray*}
\end{itemize}
Here $\mu_{1} = \mu(n, m, q, \eta) \in (0, 1)$; the constants $\k_{1}, C_{0}^{(1)} \in (0, \infty)$ depend only on $n$, $m$ and $\eta$ in case $q=2$ and only on $n$, $m$, $q$, $\eta$ and $\th_{1}, \th_{2}, \ldots, \th_{q-2}$ in case $q \geq 3$; $\nu_{1}^{(1)} = \nu_{1}(n, m, q, \eta) \in (0, \infty)$ and, in case $q\geq 3$,  $\nu_{j}^{(1)} = \nu_{j}^{(1)}(n, m, q, \eta, \th_{1}, \ldots, \th_{j-1})$ for each $j = 2, \ldots, q-1$. (In particular, for $q \geq 3$, $\nu_{j}^{(1)}$ is independent of $\th_{j}, \th_{j+1}, \ldots, \th_{q-1}$ for each $j=2, \ldots, q-1$.)  
\end{lemma}

\begin{proof} It is clear that given any $\eta \in (0, 1)$, there is a fixed cone ${\mathbf C}_{0}  = {\mathbf C}_{0}(\eta) \in \cup_{p=2}^{q} {\mathcal C}_{q, p}$ with 
$\Theta({\mathbf C}_{0}, 0) = q$ such that the following holds: for any $\epsilon  \in (0, 1)$, we can choose $\b_{1} = \b_{1}(n, m, q, \eta, \epsilon) \in (0, 1)$ such that if $T$, ${\mathbf C}$ satisfy the hypotheses of the lemma then 
${\rm dist}_{\mathcal H} \, ({\rm spt} \, {\mathbf C} \cap {\mathbf B}_{1}(0), {\rm spt} \, {\mathbf C}_{0} \cap {\mathbf B}_{1}(0)) < \epsilon$. 
To prove Lemma~\ref{excess-improvement-final-1}, repeat the  entire argument leading to Lemma~\ref{excess-improvement-final} with obvious (minor) modifications; in particular,  we utilise in places where that argument depended on results of Section~\ref{sec:apriori close2plane} the corresponding results from Section~\ref{sec:apriori nonplanar}; moreover, we use ${\mathbf C}_{0}$ (in place of ${\mathbf P}_{0}$) as the parameter space for the blow-ups analysis corresponding to the blow-up analysis of Section~\ref{blow-up-analysis}, for sequences of currents $T_{k}$ subject to the hypotheses of the present lemma with $T_{k}$ in place of $T$, and with $Q(T_{k}, {\mathbf C}_{k}, {\mathbf B}_{1}(0)) \to 0$ for a sequence of cones 
${\mathbf C}_{k} \in \cup_{p=2}^{q} {\mathcal C}_{q, p}$ with ${\mathbf C}_{k} \to {\mathbf C}_{0}$.

More specifically, the argument involves incorporating the following changes:
\begin{itemize}
		\item[(i)]  We replace Theorem~\ref{graphrep close2plane thm} (graphical representation), Corollary~\ref{keyest cor2} (Hardt-Simon inequality), and Corollary~\ref{nonconcentration close2plane cor} (excess non-concentration) with Theorem~\ref{graphrep thm}, Corollary~\ref{keyest cor3}, and Corollary~\ref{noncencentration nonplanar cor}.
		
		\item[(ii)]  We replace Theorem~\ref{branchpt close2plane thm} with Theorem~\ref{branchpt nonplanar thm}.  In particular, we use coordinates $X = (x,y) \in \mathbb{R}^{m+2} \times \mathbb{R}^{n-2}$ and in place of \eqref{goodpt-est} we get that $Z = (\xi,\zeta)$ with 
		$$|\xi| \leq C E_k.$$
		In place of \eqref{estimate-E2} we get 
		$$\int_{B_{\gamma}(0) \cap \{r > \tau_k\}} \frac{|u_i(x,y) - \pi_{P_i^{(k)}}^{\perp} \xi|^2}{|(x,y) - (\xi,\zeta)|^{n+2-\sigma}} \leq C E_k^2.$$
		In place of \eqref{blowup_est4}, we get that there exists $\lambda : B^{n-2}_{1/4}(0) \rightarrow \mathbb{R}^{m+2}$ with $\sup_{B^{n-2}_{1/4}(0)} |\lambda| \leq C$ and 
		\begin{align*}
			\hspace{20mm} &\int_{\mathbf{B}_{\rho/2}(0,z)} \sum_{i=1}^p \frac{|w_i(X) - \pi_{P_i^{(0)}}^{\perp} \lambda(z)|^2}{|X - (0,z)|^{n+2-\sigma}} 
				\,d\|\mathbf{C}_0\|(X)
			\\&\hspace{15mm} \leq C \rho^{-n-2+\sigma} \int_{\mathbf{B}_{\rho}(0,z)} \sum_{i=1}^p |w_i(X) - \pi_{P_i^{(0)}}^{\perp} \lambda(z)|^2 
				\,d\|\mathbf{C}_0\|(X),
		\end{align*}
		where $P_i^{(0)}$ are the planes making up the reference cone $\mathbf{C}_0$. 
		
		\item[(iii)]  The blow-up class $\widetilde{\mathfrak L}$ should now be defined as the set of $\widetilde{\psi}_j(x,y) = \sum_{\ell=1}^{q_j} \llbracket \psi_{j,\ell}(x) + \pi_{P_i^{(0)}}^{\perp} (L(y),0) \rrbracket$, where $L : \mathbb{R}^{n-2} \rightarrow \mathbb{R}^{m+2}$ is a linear function.

		\item[(iv)]  Lemma~\ref{blowup2L_lemma} changes in a corresponding way given item~(iii), with the rotation given by $e^{E_k M}$ where 
		$$M = \left[\begin{array}{cc} 0 & L \\ -L^T & 0 \end{array}\right].$$
		Hence $\widetilde{P}_{j,\ell}^{(k)}$ is parameterized by 
		\begin{align*}
			\hspace{25mm} &e^{E_k M} ((x,y) + E_k \psi_{j,\ell}(x)) 
			\\=\,& (x',y') + \pi_{P_i^{(k)}}^{\perp} e^{E_k M} (x,y) + E_k \psi_{j,\ell}(x) + O(E_k^2 \|\psi_{j,\ell}\| |(x,y)|)
			\\ \approx\,& (x',y') + E_k \pi_{P_i^{(k)}}^{\perp} (L(y'),0) + E_k \psi_{j,\ell}(x') + O(E_k^2 (\|\psi_{j,\ell}\| + \|L\|^2) \,|(x,y)|)
			\\ =\,& (x',y') + E_k \widetilde{\psi}_{j,\ell}(x') + O(E_k^2 (\|\psi_{j,\ell}\| + \|L\|^2) \,|(x,y)|)
		\end{align*}
		where $(x',y') = \pi_{P_i^{(k)}} e^{E_k M} (x,y)$.  A similar rotation is used in the proof of Theorem~\ref{classification_thm} and Lemma~\ref{blowup_est_lemma}.  This rotation gives us conclusion~(B)(a) of the present lemma.
\end{itemize}
\end{proof}

\section{${\mathcal H}^{n-2}$ a.e.\ uniqueness of tangent cones and $(n-2)$-rectifiability of the set of singularities where rapid decay to a plane fails}\label{rectifiability} 

\begin{theorem} \label{decomposition}
Let $q$ be an integer $\geq 2$.  There exist $\b_{\star}, \gamma_{\star} \in (0, 1)$ depending only on $n$, $m$ and $q$ such that the following holds: if ${\mathbf C} \in \cup_{p=2}^{q} {\mathcal C}_{q, p}$ and if $T$ is an $n$-dimensional locally area minimizing rectifiable current in ${\mathbf B}_{1}(0)$ with $\partial \, T \llcorner {\mathbf B}_{1}(0) = 0$ such that  $\omega_{n}^{-1} \|T\|({\mathbf B}_{1}(0)) < q + 1/2$, $\Theta (T, 0) \geq q$  and  
$$Q(T, {\mathbf C}, {\mathbf B}_{1}(0)) < \beta_{\star} \inf_{{\mathbf P} \in {\mathcal C}_{q, 1}} \, E(T, {\mathbf P}, {\mathbf B}_{1}(0))$$ then 
\begin{equation*}
	\{ Z \in {\mathbf B}_{1/2}(0) \, : \, \Theta(T, Z) \geq q \}  = \Sigma \cup \Gamma
\end{equation*}
where $\Sigma \subset L$ with $L$ a properly embedded $(n-2)$-dimensional $C^{1,\mu_{\star}}$-submanifold of ${\mathbf B}_{1/2}(0)$ with $\mathcal{H}^{n-2}(L) \leq 2\omega_{n-2}\left(\frac{1}{2}\right)^{n-2}$ and $\Gamma \subseteq \bigcup_{j=1}^{\infty} {\mathbf B}_{\rho_j}(Y_{j})$ for a countable family of balls $\{{\mathbf B}_{\rho_j}(Y_{j})\}$ with $\rho_{j} < 1/2$ and $\sum_{j=1}^{\infty} \rho_j^{n-2} \leq 1-\gamma_{\star}$.  Here $\mu_{\star} = \mu_{\star}(n, m , q) \in (0, 1)$. Moreover, for each $Z \in \Sigma$, the current $T$ at $Z$ has a unique tangent cone ${\mathbf C}_{Z}$ which has the form $\mathbf{C}_Z = \sum_{j=1}^{p_Z} q_j^{(Z)} \llbracket P_j^{(Z)} \rrbracket$, where $p_Z \geq 2$ and $q^{(Z)}_{1}, \ldots, q^{(Z)}_{p_Z} \geq 1$ are integers and $P^{(Z)}_{1}, \ldots, P^{(Z)}_{p_Z}$ are distinct $n$-dimensional oriented planes such that there is an $(n-2)$-dimensional subspace $L_Z$ with $P^{(Z)}_{i} \cap P^{(Z)}_{j} = L_Z$ for all $i \neq j$, and which satisfies the estimates 
\begin{equation} \label{thm1 concl2}
	\rho^{-n-2} \int_{B_{\rho}(Z)} {\rm dist}^{2}(X,  Z + {\rm spt} \, {\mathbf C}_{Z})\, d\|T\|(X) \leq C \rho^{2\mu_{\star}} Q^{2}(T, {\mathbf C}, {\mathbf B}_{1}(0))
\end{equation}
holds for all $\rho \in (0,1/4]$ and for some $\rho_Z \in (0,1/4]$ (depending on $Z$) 
\begin{align}\label{thm1 concl3}
	&\rho^{-n-2} \int_{B_{\rho}(Z)} {\rm dist}^{2}(X,  Z + {\rm spt} \, {\mathbf C}_{Z})\, d\|T\|(X) \\&\hspace{1.2in} \leq C \Big(\frac{\rho}{\sigma}\Big)^{2\mu_{\star}} \sigma^{-n-2} \int_{B_{\sigma}(Z)} {\rm dist}^{2}(X,  Z + {\rm spt} \, {\mathbf C}_{Z})\, d\|T\|(X) \nonumber 
\end{align}
for all $0 < \rho \leq \sigma \leq \rho_Z$, where $C = C(n,m,q) \in (0,\infty)$ is a constant.
\end{theorem}

\begin{proof} 
Write 
$${\rm sing}_{q}^{\star} \, T = \{Z \in {\mathbf B}_{1/2}(0) \, : \, \Theta \, (T, Z) \geq q\}.$$
Note first that  for any given $\b \in (0, 1)$, $\varepsilon \in (0, 1)$, we can choose $\overline{\b}_{\star} = \overline{\b}_{\star}(n, m, q, \b, \varepsilon) \in (0, 1)$
such that if ${\mathbf C}$, $T$ satisfy the hypotheses of the lemma with $\overline{\b}_{\star}$ in place of $\b_{\star}$, then for any 
$Z \in {\rm sing}^{\star}_{q} \, T$,  we have that

\begin{equation}\label{decomposition-0-0}
{\rm dist} \, (Z, {\rm spine} \, {\mathbf C}) \leq \varepsilon,
\end{equation}

\begin{equation}\label{decomposition-0}
Q(\eta_{Z, 1/2 \, \#} \, T, {\mathbf C}, {\mathbf B}_{1}(0)) < \b \inf_{{\mathbf P} \in {\mathcal C}_{q, 1}} \, E(\eta_{Z, 1/2 \, \#} \, T, {\mathbf P}, {\mathbf B}_{1}(0))   \;\; \mbox{and}
\end{equation}

\begin{equation}\label{decomposition-1-1}
\omega_{n}^{-1}\|\eta_{Z, 1/2 \, \#} \, T\|({\mathbf B}_{1}(0)) < q + 1/2.
\end{equation}

To check this, we argue by contradiction. If the claim is false, there are sequences ${\mathbf C}_{k} \in \cup_{p=1}^{q} {\mathcal C}_{q, p}$ and $T_{k}$ such that the hypotheses of the lemma hold with ${\mathbf C}_{k}$, $T_{k},$ $k^{-1}$ in place of ${\mathbf C}$, $T,$ $\b_{\star}$ and so in particular  

\begin{equation}\label{decomposition-1}
Q(T_{k}, {\mathbf C}_{k}, {\mathbf B}_{1}(0)) < k^{-1} \inf_{{\mathbf P} \in {\mathcal C}_{q, 1}} \, E(T_{k}, {\mathbf P}, {\mathbf B}_{1}(0)),
\end{equation}
 and yet for each $k$ and some $Z_{k} \in {\rm sing}^{\star}_{q} \, T_{k}$, either 
 \begin{equation}\label{decomposition-1-0}
 {\rm dist} \, (Z_{k}, {\rm spine} \, {\mathbf C}_{k}) > \epsilon \;\; \mbox{or}
 \end{equation}
 
\begin{equation}\label{decomposition-2}
Q(\eta_{Z_{k}, 1/2 \, \#} \, T_{k}, {\mathbf C}_{k}, {\mathbf B}_{1}(0)) \geq  \b \inf_{{\mathbf P} \in {\mathcal C}_{q, 1}} \, E(\eta_{Z_{k}, 1/2 \, \#} \, T_{k}, {\mathbf P}, {\mathbf B}_{1}(0)) \;\; \mbox{or}
\end{equation}
 
\begin{equation}\label{decomposition-2-1}
\omega_{n}^{-1}\|\eta_{Z_{k}, 1/2 \, \#} \, T_{k}\|({\mathbf B}_{1}(0)) \geq  q + 1/2.
\end{equation}

It follows from \eqref{decomposition-1} that for each $k \geq 2$, 
\begin{equation}\label{decomposition-3}
\inf_{{\mathbf P} \in {\mathcal C}_{q, 1}} \, {\rm dist}_{\mathcal H} \, ({\rm spt} \, {\mathbf C}_{k} \cap {\mathbf B}_{1}(0), {\rm spt} \, {\mathbf P} \cap {\mathbf B}_{1}(0)) \geq c \inf_{{\mathbf P} \in {\mathcal C}_{q, 1}} \, E(T_{k}, {\mathbf P}, {\mathbf B}_{1}(0)) 
\end{equation}
where $c = c(n, q) >0$. Choose a plane ${\mathbf P}_{k} \in {\mathcal C}_{q, 1}$ so that $E(T_{k}, {\mathbf P}_{k}, {\mathbf B}_{1}(0))  = \inf_{{\mathbf P} \in {\mathcal C}_{q, 1}} \, E(T_{k}, {\mathbf P}, {\mathbf B}_{1}(0)).$ By passing to a subsequence of $(k)$, we obtain $Z \in \overline{\mathbf B}_{1/2}(0)$ with $Z_{k} \to Z,$ ${\mathbf Q} \in {\mathcal C}_{q, 1}$ and ${\mathbf C} \in \cup_{p=1}^{q} {\mathcal C}_{q, p}$ such that, after possibly changing multiplicities and orientations of planes constituting ${\mathbf C}_{k}$, we have that ${\mathbf P}_{k} \to {\mathbf Q},$ 
${\mathbf C}_{k} \to {\mathbf C}$  as currents, and by the Federer--Fleming compactness theorem, \eqref{decomposition-1} and the fact that $\Theta \, (T_{k} , 0) \geq q$, also $T_{k} \to {\mathbf C}$ as currents, so in particular ${\mathbf C}$ is locally area minimizing. 

Consider the cases:
\begin{itemize}
\item[(i)] $\liminf_{k \to \infty} \, \inf_{{\mathbf P} \in {\mathcal C}_{q, 1}} \, E(T_{k}, {\mathbf P}, {\mathbf B}_{1}(0)) >0,$ or  
\item[(ii)] $\liminf_{k \to \infty} \, \inf_{{\mathbf P} \in {\mathcal C}_{q, 1}} \, E(T_{k}, {\mathbf P}, {\mathbf B}_{1}(0)) =0.$ 
\end{itemize}
If case (i) occurs, then by \eqref{decomposition-3} we see that ${\mathbf C} \in \cup_{p=2}^{q} \, {\mathcal C}_{q, p}$ (so ${\mathbf C}$ is not a plane), and by upper semi-continuity of density, that $Z \in {\rm spine} \, {\mathbf C}$. Thus $\eta_{Z_{k}, 1/2 \, \#} \, T_{k} \to \eta_{Z, 1/2 \, \#} \, {\mathbf C} = {\mathbf C}$, which, together with mass convergence, leads to a contradiction by letting $k \to \infty$ in \eqref{decomposition-1-0}, \eqref{decomposition-2} or \eqref{decomposition-2-1}, whichever hold for infinitely many $k$.

If case (ii) occurs, assume without loss of generality (by rotating) that ${\mathbf P}_{k}  = {\mathbf P}_{0} \equiv  q\llbracket \{0\} \times {\mathbf R}^{n} \rrbracket$ for each $k$,  
and let $w \in W^{1, 2}(B_{1}^{n}(0); {\mathcal A}_{q}({\mathbb R}^{m}))$ be a (coarse) blow-up of a subsequence of $(T_{k})$ relative to  ${\mathbf P}_{0}$ (produced using Theorem~\ref{lip approx thm}, as described in \cite[Section~6]{KrumWica}).  Let $p_{k} \in \{2, \ldots, q\}$ be such that ${\mathbf C}_{k} \in {\mathcal C}_{q, p_{k}}$. Let $\t \in (0, 1/16)$ and $\g \in (3/4, 1)$ be arbitrary, and let $\b_{0} = \b_{0}(n, m, q, \g, \t) \in (0, 1)$ be as in Theorem~\ref{graphrep close2plane thm}. In view of \eqref{decomposition-1}, for each $k$ we can find a cone ${\mathbf C}_{k}^{\prime}$ 
such that ${\mathbf C}^{\prime}_{k}  \in {\mathcal C}_{q, s_{k}}$ for some $s_{k} \in \{2, \ldots, p_{k}\}$,  
\begin{equation}\label{decomposition-4}
Q(T_{k}, {\mathbf C}_{k}^{\prime}, {\mathbf B}_{1}(0)) < \b_{0} \inf_{{\mathbf C}^{\prime} \cup_{p=1}^{s_{k} - 1} \, {\mathcal C}_{q, p}} \, 
Q(T_{k}, {\mathbf C}^{\prime}, {\mathbf B}_{1}(0))
\end{equation}
and $Q(T_{k}, {\mathbf C}_{k}^{\prime}, {\mathbf B}_{1}(0)) \leq \b_{0}^{-(q-2)} 
Q(T_{k}, {\mathbf C}_{k}, {\mathbf B}_{1}(0)).$ Writing $\widehat{E}_{k} = E(T_{k}, {\mathbf P}_{0}, {\mathbf B}_{1}(0))$, by \eqref{decomposition-1} again, we then have that 
\begin{equation} \label{decomposition-5}
Q(T_{k}, {\mathbf C}_{k}^{\prime}, {\mathbf B}_{1}(0)) \leq \b_{0}^{-(q-2)} k^{-1} \widehat{E}_{k} 
\end{equation}
which in particular implies that 
\begin{equation}\label{decomposition-6}
{\rm dist}_{\mathcal H} \, ({\rm spt} \, {\mathbf C}_{k}^{\prime} \cap ({\mathbb R}^{m} \times B_{1}(0)), {\rm spt} \, {\mathbf P} \cap ({\mathbb R}^{m} \times B_{1}(0))) \geq c \, \widehat{E}_{k}
\end{equation}
for any ${\mathbf P} \in {\mathcal C}_{q, 1}$, where $c = c(n, q) > 0.$ The condition \eqref{decomposition-4} allows us to apply Theorem~\ref{graphrep close2plane thm}, whence, in view of \eqref{decomposition-5} and the fact that the blow-up $w$ is locally Dirichlet energy minimizing, we see that ${\rm graph} \, w  = {\rm spt} \, {\mathbf C}_{w} \cap ({\mathbb R}^{m} \times B_{1}(0)))$ for some cone ${\mathbf C}_{w} \in \cup_{p=1}^{q} {\mathcal C}_{q, p}.$ (In particular, the possibility that ${\rm graph} \, w$ consists of distinct planes intersecting along an $(n-1)$-dimensional subspace is ruled out by the energy minimizing property of $w$, and we have that ${\rm spine} \, ({\mathbf C}_{w}) = \{0\} \times {\mathbf R}^{n-2}$, although these facts are not needed for the rest of the argument).  Moreover, given any linear function $L \, : \, {\mathbb R}^{n} \to {\mathbb R}^{m}$, we can take 
${\mathbf P} = q\llbracket{\rm graph} \, (\widehat{E}_{k}L)\rrbracket$ in \eqref{decomposition-6}, divide both sides by 
$\widehat{E}_{k}$ and pass to the limit as $k \to \infty$ to see that  
${\rm dist}_{\mathcal H} \, ({\rm spt} \, ({\mathbf C}_{w}) \cap {\mathbf B}_{1}(0), {\rm graph} \, L \cap {\mathbf B}_{1}(0)) \geq c$ whence 
${\rm spt} \, {\mathbf C}_{w}$ is not a single plane, so that ${\mathbf C}_{w} \in {\mathcal C}_{q, p}$ for some $p \geq 2$. Also, by \cite[Lemma~6.2]{KrumWica}, we have that $w_{a} \equiv 0$ (where $w_{a}(x)$ is the average of the $q$ values of $w(x)$). On the other hand, since $\Theta(T_{k}, Z_{k}) \geq q$,  writing $Z_{k} = (\chi_{k}, \xi_{k}, \z_{k})$ with $\chi_{k} \in {\mathbb R}^{m}$ and $(\xi_{k}, \z_{k}) \in {\mathbb R}^{2} \times {\mathbb R}^{n-2},$ we see using \cite[Lemma~6.4]{KrumWica} that there is a point $\z \in {\rm spine} \, ({\mathbf C}_{w})$ such that passing to a subsequence without changing notation, 
$(\widehat{E}_{k}^{-1}\chi_{k}, \xi_{k}, \z_{k}) \to (w_{a}(0, \z), 0, \z) = (0, 0, \z)$. 
In particular, this means that ${\rm dist} \, (Z_{k}, {\rm spine} \, {\mathbf C}_{k}) \to 0,$ so \eqref{decomposition-1-0} must fail for sufficiently large $k$. It is also clear by mass convergence that \eqref{decomposition-2-1} must fail for sufficiently large $k$. So we must have \eqref{decomposition-2} for infinitely many $k$.

Now, in view of \eqref{decomposition-4}, we can argue exactly as for \eqref{excess-improvement-final-dist} to to see that 
\begin{equation*} 
{\rm dist}_{\mathcal H} \, ({\rm spt} \, {\mathbf C}_{k}^{\prime} \cap {\mathbf B}_{1}(0), {\rm spt} \, {\mathbf C}_{k} \cap {\mathbf B}_{1}(0)) \leq C(Q(T_{k}, {\mathbf C}_{k}^{\prime}, {\mathbf B}_{1}(0)) + Q(T_{k}, {\mathbf C}_{k}, {\mathbf B}_{1}(0))) \leq Ck^{-1}\widehat{E}_{k}.
\end{equation*} 
This implies that ${\rm dist}_{\mathcal H} \, ({\rm spt} \, {\mathbf C}_{k} \cap {\mathbf B}_{1}(0), {\mathbf P}_{0} \cap {\mathbf B}_{1}(0)) \leq C\widehat{E}_{k}$, and hence 
$${\rm dist}^{2}_{\mathcal H} \, ({\rm spt} \, \eta_{Z_{k}, 1 \, \#} \, {\mathbf C}_{k} \cap {\mathbf B}_{1}(0), {\rm spt} \, {\mathbf C}_{k} \cap {\mathbf B}_{1}(0)) \leq C(|\chi_{k}|^{2} + \widehat{E}_{k}^{2}|\xi_{k}|^{2}).$$ Therefore 
\begin{eqnarray*}
&&\int_{{\mathbf B}_{1}(0)} {\rm dist}^{2}\, (X, {\rm spt} \, {\mathbf C}_{k}) \, d\|\eta_{Z_{k}, 1/2 \, \#} \, T_{k}\| = 
\left(\frac{1}{2}\right)^{-n-2} \int_{{\mathbf B}_{1/2}(Z_{k})} {\rm dist}^{2} \, (X, Z_{k} + {\rm spt} \, {\mathbf C}_{k}) \, d\|T_{k}\|\\ 
&&\leq C\left(\int_{{\mathbf B}_{1}(0)} {\rm dist}^{2} \, (X, {\rm spt} \, {\mathbf C}_{k}) \, d\|T_{k}\| +  (|\chi_{k}|^{2} + \widehat{E}_{k}^{2}|\xi_{k}|^{2})\right). \nonumber
\end{eqnarray*} 
Similarly, 
\begin{eqnarray*}
&&\hspace{-.2in}\int_{{\mathbf B}_{1/2}(0) \setminus \{r(X) < 1/16\}} {\rm dist}^{2}\, (X, {\rm spt} \, \eta_{Z_{k}, 1/2 \, \#} \, T_{k}) \, d\|{\mathbf C}_{k}\|\\ 
&&\hspace{-.2in}= \left(\frac{1}{2}\right)^{-n-2} \int_{{\mathbf B}_{1/2}(Z_{k}) \setminus \{r(X - Z_{k}) < 1/32\}} {\rm dist}^{2} \, (X, {\rm spt} \, T_{k}) \, d\| \eta_{Z_{k}, 1 \, \#} \, {\mathbf C}_{k}\|\nonumber\\ 
&&\hspace{-.2in}\leq \left(\frac{1}{2}\right)^{-n-2}\int_{{\mathbf B}_{1/2}(0) \setminus \{r(X) < 1/64\}} {\rm dist}^{2} \, (X, {\rm spt} \, {T}_{k}) \, d\|{\mathbf C}_{k}\| +  C(|\chi_{k}|^{2} + \widehat{E}_{k}^{2}|\xi_{k}|^{2})\nonumber\\
&&\hspace{-.2in}\leq C\int_{{\mathbf B}_{1}(0)} {\rm dist}^{2} \, (X, {\rm spt} \, {\mathbf C}_{k}^{\prime}) \, d\|{\mathbf C}_{k}\| + C(|\chi_{k}|^{2} + \widehat{E}_{k}^{2}|\xi_{k}|^{2}) +\nonumber\\ 
&& C\,{\rm dist}^{2}_{\mathcal H} \, ({\rm spt} \, {\mathbf C}_{k}^{\prime} \cap ({\mathbf B}_{1/2}(0) \setminus \{r(X) < 1/64\}), {\rm spt} \, T_{k} \cap ({\mathbf B}_{1/2}(0) \setminus \{r(X) < 1/64\}))\nonumber\\
&&\hspace{-.2in}\leq C\left({\rm dist}^{2}_{\mathcal H} \, ({\rm spt} \, {\mathbf C}_{k}^{\prime} \cap {\mathbf B}_{1}(0), {\rm spt} \, {\mathbf C}_{k} \cap {\mathbf B}_{1}(0)) + (|\chi_{k}|^{2} + \widehat{E}_{k}^{2}|\xi_{k}|^{2}) + Q(T_{k}, {\mathbf C}_{k}^{\prime}, {\mathbf B}_{1}(0))\right)
\end{eqnarray*} 
where $C = C(n, m, q),$ and where in the last inequality we have used Theorem~\ref{graphrep close2plane thm}(b) (with $T_k$ and $\mathbf{C}'_k$ in place of $T$ and $\mathbf{C}$). Combining the above estimates with Theorem~\ref{branchpt close2plane thm}(a) 
we see, in view of \eqref{decomposition-1} and \eqref{decomposition-5}, that $${\hat E}_{k}^{-1} Q(\eta_{Z_{k}, 1/2 \, \#} \, T_{k}, {\mathbf C}_{k}, {\mathbf B}_{1}(0)) \to 0$$
as $k \to \infty$. Hence by \eqref{decomposition-2}, 
\begin{equation}\label{decomposition-7}
\widehat{E}_{k}^{-1}\inf_{{\mathbf P} \in {\mathcal C}_{q, 1}} \, E(\eta_{Z_{k}, 1/2 \, \#} \, T_{k}, {\mathbf P}, {\mathbf B}_{1}(0)) \to 0
\end{equation}
as $k \to \infty.$ Now choosing for each $k$ a plane ${\mathbf Q}_{k} \in {\mathcal C}_{q, 1}$ such that 
$$E(\eta_{Z_{k}, 1/2 \, \#} \, T_{k}, {\mathbf Q}_{k}, {\mathbf B}_{1}(0)) = \inf_{{\mathbf P} \in {\mathcal C}_{q, 1}} \, E(\eta_{Z_{k}, 1/2 \, \#} \, T_{k}, {\mathbf P}, {\mathbf B}_{1}(0)),$$ we see that ${\rm dist}_{\mathcal H} \, ({\rm spt} \, {\mathbf Q}_{k} \cap {\mathbf B}_{1}(0), B_{1}(0)) \leq C \widehat{E}_{k}$ where 
$C = C(n)$, and hence ${\rm spt} \, {\mathbf Q}_{k} = {\rm graph} \, L_{k}$ where $L_{k} \, : \, {\mathbb R}^{n} \to {\mathbb R}^{m}$ is a linear function of the form 
$L_{k}(x, y) = L_{k}(x, 0)$ for $x \in {\mathbb R}^{2}$, $y \in {\mathbb R}^{n-2}$, and   
there is a linear function $L \, : \, {\mathbb R}^{n} \to {\mathbb R}^{m}$ of the form $L(x, y) = L(x, 0)$ for $x \in {\mathbb R}^{2},$ 
$y \in {\mathbb R}^{n-2}$ such that, after passing to a subsequence of $(k)$, 
$\lim_{k \to \infty} \, \widehat{E}_{k}^{-1} L_{k}  = L.$ By direct calculation, it is then not difficult to see that 
$$\lim_{k \to \infty} \, \widehat{E}_{k}^{-2}E^{2}(\eta_{Z_{k}, 1/2 \, \#} \, T_{k}, {\mathbf Q}_{k}, {\mathbf B}_{1}(0))  = 
\left(\frac{1}{2}\right)^{-n-2}\int_{B_{1/2}(0, \z)} {\mathcal G}^{2}(w(x), q\llbracket L(x) \rrbracket) \, dx,$$  so by \eqref{decomposition-7}, we must have 
that $w \equiv q\llbracket L \rrbracket$ in $B_{1/2}(0,\z)$. But this is impossible since as we have shown above, 
${\rm graph} \, w = {\rm spt} \, {\mathbf C}_{w} \cap ({\mathbb R}^{m} \times B_{1}(0)))$ with ${\mathbf C}_{w} \in {\mathcal C}_{q, p}$ for some $p \geq 2$.  
Thus the claim is established that there is $\overline{\b}_{\star} = \overline{\b}_{\star}(n, m, q, \b, \varepsilon)$ so that if $T$, ${\mathbf C}$ satisfy the hypotheses of the theorem with $\overline{\b}_{\star}$ in place of $\b_{\star}$ then \eqref{decomposition-0-0} \eqref{decomposition-0} and \eqref{decomposition-1-1} hold.

To complete the proof of the theorem choose scales $\th_{1}, \ldots, \th_{q-1} \in (0, 1/2),$  depending only on $n$, $m$ and $q,$ such that 
$\th_{j} \geq 8\th_{j+1}$ for $j=1, \ldots, q-2$ and $\nu_{j}\th_{j}^{2\mu} < 1/2$ for $j=1, \ldots, q~-~1$, where $\mu = \mu(n, m, q) \in (0, 1)$, $\nu_{1} = \nu_{1}(n, m, q) \in (0, \infty)$  and 
$\nu_{j} = \nu_{j}(n, m, q, \th_{1}, \ldots, \th_{j-1})$, $2 \leq j \leq q-1,$ are the constants as in Lemma~\ref{excess-improvement-final}. Let 
$\eta = \eta(n, m, q, M_{0}, \th_{1}, \ldots, \th_{q-1})$, $\b = \b(n, m, q, M_{0}, \th_{1}, \ldots, \th_{q-1})$ and $\d =\d(n, m, q, M_{0}, \th_{1}, \ldots, \th_{q-1})$ be as in Lemma~\ref{excess-improvement-final} taken with $M=M_{0}$, where $M_{0} = M_{0}(n, m, q) \in [1, \infty)$ is to be determined. 
Set $\widetilde{\b} = \Lambda^{-1}\min\{\b, \b_{1}\}$ and $\mu_{\star} = \min\{\mu, \mu_{1}\},$ where $\b_{1}$ and 
$\mu_{1} = \mu_{1}$ are constants to be determined depending only on $n$, $m$ and $q$ (specified in the last paragraph of the present proof), and $\Lambda \in [1, \infty)$ is a constant to be chosen depending only on $n$, $m$ and $q$.  Suppose that the hypotheses of the theorem hold with a choice of 
$\b_{\star}$ (to be fixed depending only on $n$, $m$, $q$) such that $\b_{\star}  \in (0, \overline{\b}_{\star}(n, m, q, \widetilde{\b}, \varepsilon)),$ where 
$\overline{\b}_{\star}( \cdot)$ is as established at the beginning of the proof and $\varepsilon \in (0, 1/2)$ is to be chosen; in particular, we shall require that 

\begin{equation}\label{decomposition-7-1} 
\varepsilon <  \min\left\{\delta^{2}, \Lambda^{-1}\eta\right\}
\end{equation}

Then by \eqref{decomposition-0-0}, \eqref{decomposition-0} and \eqref{decomposition-1-1} we have for every 
$Z \in {\rm sing}^{\star}_{q} \, T$ that 

\begin{equation}\label{decomposition-8}
{\rm dist} \, (Z, {\rm spine} \, {\mathbf C}) < \varepsilon,
\end{equation}
 
\begin{equation}\label{decomposition-9}
Q(\eta_{Z, 1/2 \, \#} \, T, {\mathbf C}, {\mathbf B}_{1}(0)) < \widetilde{\b} \inf_{{\mathbf P} \in {\mathcal C}_{q, 1}} \, E(\eta_{Z, 1/2 \, \#} \, T, {\mathbf P}, {\mathbf B}_{1}(0))  \;\; {\rm and}
\end{equation}

\begin{equation}\label{decomposition-9-1}
\omega_{n}^{-1}\|\eta_{Z, 1/2 \, \#} \, T\|({\mathbf B}_{1}(0)) < q + 1/2.
\end{equation}

Consider the two cases: 

\begin{itemize}
\item[(a)] $\inf_{{\mathbf P} \in {\mathcal C}_{q, 1}} \, E^{2}(T, {\mathbf P}, {\mathbf B}_{1}(0)) < \Lambda^{-2}\eta^{2} \;\; \mbox{or}$
\item[(b)] $\inf_{{\mathbf P} \in {\mathcal C}_{q, 1}} \, E^{2}(T, {\mathbf P}, {\mathbf B}_{1}(0))  \geq \Lambda^{-2}\eta^{2}.$
\end{itemize}

Suppose (a) holds. Then,  since for any ${\mathbf P} \in {\mathcal C}_{q, 1}$ and any $Z \in {\rm sing}^{\star}_{q} \, T$ we have that 
\begin{eqnarray*}
&&E^{2}(\eta_{Z, 1/2 \, \#} \, T, {\mathbf P}, {\mathbf B}_{1}(0)) = 2^{n+2} \int_{{\mathbf B}_{1/2}(Z)} {\rm dist}^{2} \, (X, Z + {\rm spt} \, {\mathbf P}) \, d\|T\|\\ 
&&\hspace{1in}\leq 2^{n+3} \int_{{\mathbf B}_{1}(0)} {\rm dist}^{2} \, (X, {\rm spt} \, {\mathbf P}) \, d\|T\| + \widetilde{C}^{2}{\rm dist}^{2}(Z, \{0\} \times {\mathbf R}^{n-2})
\end{eqnarray*}
where $\widetilde{C} = \widetilde{C}(n, q)$, it follows from \eqref{decomposition-7-1}, \eqref{decomposition-8} and (a) that 

\begin{equation}\label{decomposition-10}
\inf_{{\mathbf P} \in {\mathcal C}_{q, 1}} \, E(\eta_{Z, 1/2 \, \#} \, T, {\mathbf P}, {\mathbf B}_{1}(0)) < (2^{(n+3)/2} + \widetilde{C}) \Lambda^{-1}\eta.
\end{equation}

In view of \eqref{decomposition-9}, \eqref{decomposition-9-1} and \eqref{decomposition-10}, for every $Z \in {\rm sing}^{\star}_{q} \, T$ we can apply Lemma~\ref{excess-improvement-final} with 
$T_{Z} \equiv \Delta_{Z \, \#} \, \eta_{Z, 1/2 \, \#} \, T$ in place of $T$ (and with $M = M_{0} = M_{0}(n, q) \in [1, \infty)$ specified below) where $\Delta_{Z}$ 
is an appropriate rotation of ${\mathbb R}^{n+m}$ that fixes $\{0\} \times {\mathbf R}^{n-2}$ and takes a plane ${\mathbf P}_{Z} \in {\mathcal C}_{q, 1}$ attaining the infimum in \eqref{decomposition-10} to ${\mathbf P}_{0}$. If for some $Z \in {\rm sing}^{\star}_{q} \, T \cap {\mathbf B}_{1/4}(0)$ conclusion~(A) of Lemma~\ref{excess-improvement-final} holds with $T_{Z}$ in place of $T$, then ${\rm sing}^{\star}_{q} \, T \subset \{Y \in {\mathbb R}^{n+m} \, : \, {\rm dist} \, (Y, \{0\} 
\times {\mathbb R}^{n-2}) < \delta^{2}\}\cap {\mathbf B}_{1/2}(0) \setminus {\mathbf B}_{\delta/2}(0,y)$  for some 
$(0,y) \in \{0\} \times {\mathbb R}^{n-2} \cap {\mathbf B}_{1/2}(0)$, from which it is easy to see that the conclusion of the present theorem holds with $\Sigma = \emptyset$ and $\Gamma = {\rm sing}^{\star}_{q} \, T$. So we can assume that for each $Z \in {\rm sing}^{\star}_{q} \, T \cap {\mathbf B}_{1/4}(0)$, conclusion~(B) of Lemma~\ref{excess-improvement-final} holds with $T_{Z}$ in place of $T$. 
Arguing as in \cite[Theorem~3.1]{MW} (which in turn is based on \cite[Lemma~14.1]{Wic14}), applying Lemma~\ref{excess-improvement-final} iteratively after choosing $\Lambda= \Lambda(n, m, q)$ sufficiently large, we can now see that there is a 
fixed $\b_{\star} = \b_{\star}(n, m, q) \in (0, 1/2)$ so that if the hypotheses of the present theorem hold with this choice of $\b_{\star}$, then for each $Z \in {\rm sing}^{\star}_{q} \, T \cap {\mathbf B}_{1/4}(0)$ there exists $K_{Z} \in \{1,2,3,\ldots\} \cup \{\infty\}$ such that, letting $\Gamma^{(Z)}_{0}$ be the identity map on $\mathbb{R}^{n+m},$  $\mathbf{C}^{(Z)}_{0} = \mathbf{C}$ and $\sigma_{0}^{(Z)} = 1$: 
\begin{enumerate}[itemsep=2mm,topsep=0mm]
\item[(i)]  for each positive integer $k \leq K_Z$ there exist a rotation $\Gamma^{(Z)}_k$ of $\mathbb{R}^{n+m},$ a cone $\mathbf{C}_k^{(Z)} \in \bigcup_{p=2}^q \mathcal{C}_{q,p},$ $\theta^{(Z)}_k \in \{\theta_1,\theta_2,\ldots,\theta_{q-1}\},$ and a radius $\sigma^{(Z)}_k$ given by $\sigma^{(Z)}_k = \theta^{(Z)}_k \sigma^{(Z)}_{k-1}$ such that 
\begin{gather}
    \label{decomposition-11a} |\Gamma^{(Z)}_{k-1}(e_i) - \Gamma^{(Z)}_{k}(e_i)| \leq C Q((\Gamma^{(Z)}_{k-1})^{-1}_{\#} T_Z,\mathbf{C}^{(Z)}_{k-1},\mathbf{B}_{\sigma^{(Z)}_{k-1}}(0)) \text{ for } i = 1,\ldots,m, \\
    \label{decomposition-11b} |\Gamma^{(Z)}_{k-1}(e_{m+j}) - \Gamma^{(Z)}_{k}(e_{m+j})| \leq \frac{C Q((\Gamma^{(Z)}_{k-1})^{-1}_{\#} T_Z,\mathbf{C}^{(Z)}_{k-1},\mathbf{B}_{\sigma^{(Z)}_{k-1}}(0))}{E((\Gamma^{(Z)}_{k-1})^{-1}_{\#} T_Z,\mathbf{P}_0,\mathbf{B}_{\sigma^{(Z)}_{k-1}}(0))} \text{ for } j = 1,\ldots,n, \\
    \label{decomposition-11c} {\rm dist}_{\mathcal H}(\op{spt}\mathbf{C}^{(Z)}_{k-1} \cap \mathbf{B}_1(0), \op{spt}\mathbf{C}^{(Z)}_{k} \cap \mathbf{B}_1(0)) \leq C Q((\Gamma^{(Z)}_{k-1})^{-1}_{\#} T_Z,\mathbf{C}^{(Z)}_{k-1},\mathbf{B}_{\sigma^{(Z)}_{k-1}}(0)) , \\
    \label{decomposition-11d} Q((\Gamma^{(Z)}_{k})^{-1}_{\#} T_Z,\mathbf{C}^{(Z)}_{k},\mathbf{B}_{\sigma^{(Z)}_{k}}(0)) \leq (\theta^{(Z)}_k)^{\mu/2} \,Q((\Gamma^{(Z)}_{k-1})^{-1}_{\#} T_Z,\mathbf{C}^{(Z)}_{k-1},\mathbf{B}_{\sigma^{(Z)}_{k-1}}(0)) , 
\end{gather}
and for any $n$-dimensional plane $P \subset \mathbb{R}^{n+m}$ containing $\{0\} \times \mathbb{R}^{n-2}$ 
\begin{align}\label{decomposition-11e} 
   &\left( (\sigma^{(Z)}_k)^{-n-2} \int_{\mathbf{B}_{\sigma^{(Z)}_k}(0)} {\rm dist}^2(X,P) \,d\|(\Gamma^{(Z)}_{k})^{-1}_{\#} T_Z\|(X) \right)^{1/2} \\ \geq\,& C_1 {\rm dist}_{\mathcal H}(\op{spt}\mathbf{C}^{(Z)}_{k-1} \cap \mathbf{B}_1(0), P \cap \mathbf{B}_1(0)) - C_2 Q((\Gamma^{(Z)}_{k-1})^{-1}_{\#} T_Z,\mathbf{C}^{(Z)}_{k-1},\mathbf{B}_{\sigma^{(Z)}_{k-1}}(0)) , \nonumber 
\end{align}
where $C, C_1, C_2 \in (0,\infty)$ are constants depending only on $n,m,q$; 

\item[(ii)]  either $K_Z = \infty$ or $K_Z < \infty$ and for some point $(0,z) \in \{0\} \times \mathbb{R}^{n-2} \cap \mathbf{B}_{\sigma^{(Z)}_{K_Z}/2}(0)$
\begin{equation}\label{decomposition-11f}
    \mathbf{B}_{\delta \theta_{q-1} \sigma^{(Z)}_{K_Z}/2}(0,z) \cap 
    \{ X : \Theta(T_Z,X) \geq q \} = \emptyset . 
\end{equation}
\end{enumerate}
When $K_Z < \infty$, we shall take $K_Z$ to be the smallest positive integer for which \eqref{decomposition-11f} holds true for some $(0,z) \in \{0\} \times \mathbb{R}^{n-2} \cap \mathbf{B}_{\sigma^{(Z)}_{K_Z}/2}(0)$.

The argument to reach this conclusion proceeds by induction on $k$ exactly as in the proofs of \cite[Theorem~3.1]{MW}, \cite[Lemma~14.1]{Wic14}, assuming by induction the validity of \eqref{decomposition-11a}--\eqref{decomposition-11e} for indices $1, \ldots, k-1$ in place of $k$ (cf.~\cite[(14.2)--(14.6)]{Wic14}), and applying in the inductive step Lemma~\ref{excess-improvement-final} in place of \cite[Lemma~13.3]{Wic14}; in particular, the constant $M_{0} = M_{0}(n,  q) \in [1, \infty)$ is chosen and fixed (cf.\ \cite[pp.\ 950-951]{Wic14}) so that, subject to the above induction hypotheses, we have 
$$E^{2}\left(\eta_{0, \sigma_{k-1}^{(Z)} \, \#} \, T_{Z}, {\mathbf P}_{0}, {\mathbf B}_{1}(0)\right) \leq M_{0} \inf_{{\mathbf P} \in {\mathcal C}_{q, 1}} \, 
E^{2}\left(\eta_{0, \sigma_{k-1}^{(Z)} \, \#} \, T_{Z}, {\mathbf P}, {\mathbf B}_{1}(0)\right).$$ 

Now define,  for $j=1, 2, \ldots$, 
$E_{j} = \{Z \in {\rm sing}^{\star}_{q} \, T \cap {\mathbf B}_{1/4}(0) \, : \,  K_{Z} = j \}$
and define 
$\Sigma = E_{\infty} = \{Z  \in {\rm sing}^{\star}_{q} \, T \cap {\mathbf B}_{1/4}(0) \, : \, K_Z = \infty\} .$ 
Let $Z \in \Sigma$.  In view of \eqref{decomposition-11a}--\eqref{decomposition-11d}, there exists a rotation $\Gamma^{(Z)}$ of $\mathbb{R}^{n+m}$ and a cone $\mathbf{C}^{(Z)} \in \bigcup_{p=2}^q \mathcal{C}_{q,p}$ such that \begin{gather}
    \label{decomposition-12a} |\Gamma^{(Z)}_{k}(e_i) - \Gamma^{(Z)}(e_i)| \leq C (\sigma^{(Z)}_k)^{\mu/2} \,Q(T_Z,\mathbf{C},\mathbf{B}_1(0)) \text{ for } i = 1,\ldots,m, \\
    \label{decomposition-12b} |\Gamma^{(Z)}_{k}(e_{m+j}) - \Gamma^{(Z)}(e_{m+j})| \leq \frac{C (\sigma^{(Z)}_k)^{\mu/2} \,Q(T_Z,\mathbf{C},\mathbf{B}_1(0))}{E(T_Z,\mathbf{P}_0,\mathbf{B}_1(0))} \text{ for } j = 1,\ldots,n, \\
    \label{decomposition-12c} {\rm dist}_{\mathcal H}(\op{spt}\mathbf{C}^{(Z)}_{k} \cap \mathbf{B}_1, \op{spt}\mathbf{C}^{(Z)} \cap \mathbf{B}_1) \leq C (\sigma^{(Z)}_k)^{\mu/2} \,Q(T_Z,\mathbf{C},\mathbf{B}_1(0)) , \\
    \label{decomposition-12d} Q((\Gamma^{(Z)}_{k})^{-1}_{\#} T_Z,\mathbf{C}^{(Z)},\mathbf{B}_{\sigma^{(Z)}_{k}}(0)) \leq C (\theta^{(Z)}_k)^{\mu/2} \,Q(T_Z,\mathbf{C},\mathbf{B}_{\sigma^{(Z)}_{k-1}}(0)) 
\end{gather}
and \eqref{thm1 concl2} holds true.  Let $Z_1,Z_2 \in \Sigma$ with $\sigma = |Z_1 - Z_2| < 1/8$.  Choose $k$ such that $\sigma^{(Z_2)}_{k+1} < 2\sigma \leq \sigma^{(Z_2)}_{k}.$  Set $\widetilde{T} = \big( \eta_{0,\sigma^{(Z_2)}_k} \circ (\Gamma^{(Z_2)}_k)^{-1} \big)_{\#} T_{Z_2}$, $\widetilde{Z} = \big( \eta_{0,\sigma^{(Z_2)}_k/2} \circ (\Gamma^{(Z_2)}_k)^{-1} \circ \Delta_{Z_2} \big)(Z_1 - Z_2)$ and $\widetilde{V}_{\widetilde Z} = \big( \eta_{0,\sigma^{(Z_2)}_k} \circ (\Gamma^{(Z_2)}_k)^{-1} \big)_{\#} T_{Z_1}.$  By \eqref{decomposition-12a}, \eqref{decomposition-12b} and \eqref{decomposition-9} we have that $\|\Gamma^{(Z_2)}_k - \Gamma^{(Z_2)}\| \leq C \widetilde{\beta}$ and thus by provided $\widetilde{\beta}$ is sufficiently small, by \eqref{decomposition-11f} 
\begin{equation*}
    \mathbf{B}_{\delta \sigma}(0,z) \cap 
    \{ X : \Theta(\widetilde{T}_{\widetilde Z},X) \geq q \} \neq \emptyset 
\end{equation*}
for all $\sigma \in (0,1]$ and $(0,z) \in \{0\} \times \mathbb{R}^{n-2} \cap \mathbf{B}_{\sigma/2}(0).$  Hence by argument of the proof of \cite[Theorem~3.1]{MW}, which involves iteratively applying Lemma~\ref{excess-improvement-final} as above with $\widetilde{T}_{\widetilde Z}$ and $\mathbf{C}^{(Z)}$ in place of $T_Z$ and $\mathbf{C}$, we then obtain that $\Sigma \subset L$ for some embedded $C^{1, \mu}$ submanifold $L$ of ${\mathbf B}_{1/2}$ with ${\mathcal H}^{n-2}(L) \leq 2\omega_{n-2}\left(\frac{1}{2}\right)^{n-2}$ (in fact $L = {\rm graph} \, \varphi$ for some function $\varphi \, : \, B^{n-2}_{1/2}(0) \equiv \{0\} \times {\mathbb R}^{n-2} \cap {\mathbf B}_{1/2}(0) \to {\mathbb R}^{m+2}$ with $|\varphi|_{1, \mu; B_{1/2}^{n-2}(0)}  \leq C\eta$, $C = C(n, m, q) \in (0, \infty)$).  By~\cite[Lemma~2.6]{KrumWica}, \eqref{decomposition-12d} and \eqref{thm1 concl2}, for each $Z \in \Sigma$ 
\begin{equation*}
    \lim_{\rho\rightarrow 0^+} \op{dist}_{\mathcal H}(\op{spt} T_Z \cap \mathbf{B}_{\rho}(0),\op{spt} \mathbf{C}^{(Z)} \cap \mathbf{B}_{\rho}(0)) = 0 .
\end{equation*}
Thus, after changing the multiplicities of the planes of $\mathbf{C}^{(Z)}$ if necessary, $\mathbf{C}^{(Z)}$ is the unique tangent cone to $T$ at $Z$ and by Lemma~\ref{one plane lemma} there exists $\rho_Z \in (0,1/4]$ such that $Q(T_Z,\mathbf{C}^{(Z)},\mathbf{B}_{\sigma}(0)) \leq C(n,m,q) \,E(T_Z,\mathbf{C}^{(Z)},\mathbf{B}_{\sigma}(0))$ for all $\sigma \in (0,2\rho_Z]$.  Thus for each $Z \in \Sigma$ and $\sigma \in (0,\rho_Z]$, we can again iteratively apply Lemma~\ref{excess-improvement-final} with $\widetilde{T} = \big( \eta_{0,\sigma} \circ (\Gamma^{(Z)})^{-1} \big)_{\#} T_{Z}$ and $\mathbf{C}^{(Z)}$ in place of $T_Z$ and $\mathbf{C}$, we obtain the estimate \eqref{thm1 concl3}.  

To see the rest of the conclusions of the theorem (still in case (a)), i.e.\ to see that ${\rm sing}^{\star}_{q} \, T  \setminus \Sigma \subset \cup_{j} {\mathbf B}_{\rho_{j}}(Y_{j})$ with $\sum_{j} \rho_{j}^{n-2} < 1 - \gamma^{\star}$ for some $\gamma^{\star} = \gamma^{\star}(n, m, q) \in (0, 1)$,  note that ${\rm sing}^{\star}_{q} \, T  \setminus \Sigma \subset \left(\cup_{j=1}^{\infty} E_{j} \right) \cup ({\mathbf B}_{1/2}(0) \setminus {\mathbf B}_{1/4}(0))\cap {\rm sing}^{\star}_{q} \, T$, and that for each positive integer $j$ and each $Z \in E_{j}$, by (ii) above, 
there is a point $Y \in Z + (\{0\} \times {\mathbb R}^{n-2}) \cap {\mathbf B}_{\sigma_{j}^{(Z)}/2}(0)$ such that 
\begin{equation}\label{decomposition-14}
{\mathbf B}_{
\delta\theta_{q-1} \sigma^{(Z)}_{j}/4}(Y) \cap {\rm sing}^{\star}_{q} \, T = \emptyset;  
\end{equation}
moreover, by (i), we can apply \eqref{decomposition-8} with any point $\in {\rm sing}^{\star}_{q} \, \eta_{0, \sigma_{k}^{(Z)} \, \#} \, T_{Z}$ in place of $Z$ and ${\mathbf C}_{k}$ in place of ${\mathbf C}$ for any $k \in \{1, \ldots, j-1\}$ to see that 
\begin{equation}\label{decomposition-15}
{\rm sing}^{\star}_{q} \, T \cap {\mathbf B}_{\rho}(Z) \subset \{X \, : \, {\rm dist} \, (X, Z + \{0\} \times {\mathbb R}^{n-1}) < C\epsilon \rho\}
\end{equation} 
for each $Z \in E_{j}$ and each $\rho$ with $\sigma^{(Z)}_{j} \leq \rho < 1,$ where $C = C(n, m ,q) \in (0, \infty)$.  
We can now reach the desired conclusion by choosing $\epsilon \in (0, 1/2)$ sufficiently small depending only on $n$, $m$, $q$ and arguing exactly as in the last part of the proof of \cite[Theorem~1]{Sim93}, with the help of the covering theorem 
\cite[Theorem~2.7]{Sim93} and with \eqref{decomposition-14}, \eqref{decomposition-15} in place of \cite[5.2(13)]{Sim93}, \cite[5.2(12)]{Sim93} respectively.

Finally, to establish the theorem in case (b),  choose scales $\th_{1}^{(1)}, \ldots, \th_{q-1}^{(1)} \in (0, 1/2),$  depending only on $n$, $m$ and $q,$ such that 
$\th_{j}^{(1)} \geq 8\th_{j+1}^{(1)}$ for $j=1, \ldots, q-2$ and $\nu_{j}^{(1)}(\th_{j}^{(1)})^{2\mu_{1}} < 1/2$ for $j=1, \ldots, q-1$, where 
$\mu_{1} = \mu_{1}(n, m, q, \eta) \in (0, 1)$, $\nu_{1}^{(1)} = 
\nu_{1}^{(1)}(n, m, q, \eta) \in (0, \infty)$ and 
$\nu_{j}^{(1)} = \nu_{j}(n, m, q, \eta, \th_{1}, \ldots, \th_{j-1})$, $2 \leq j \leq q-1,$ are the constants as in Lemma~\ref{excess-improvement-final-1}. Let $\b_{1} = \b_{1}(n, m, q, \eta, \th_{1}^{(1)}, \ldots, \th_{q-1}^{(1)})$, $\d_{1} = \d_{1}(n, m, q, \eta, \th_{1}^{(1)}, \ldots, \th_{q-1}^{(1)})$ be as in Lemma~\ref{excess-improvement-final-1}. Here $\eta = \eta(n, m, q) \in (0, 1)$ is chosen and fixed as in case (a) discussed above. The argument for case (b) then proceeds similarly to case (a) (and is in fact closer to that of \cite[Theorem~1]{Sim93}), with Lemma~\ref{excess-improvement-final-1} playing the role of Lemma~\ref{excess-improvement-final}. The proof of the theorem is thus complete. 
\end{proof}

Combining Theorem~\ref{decomposition} with \cite[Theorem~1.1]{KrumWica}, we obtain uniqueness of tangent cones to $n$-dimensional area minimizing currents 
at ${\mathcal H}^{n-2}$ a.e.\ point, and for any integer $q \geq 2$, rectifiability and local finiteness of measure of the set of density $q$ singular points at which the current does not rapidly decay to a (unique) tangent plane.  

\begin{theorem}\label{structure thm} 
Let $q$ be an integer $\geq 2$. There exists numbers $\alpha = \alpha(n, m, q) \in (0, 1)$ and $C = C(n, m, q) \in (0, \infty)$ such that the following holds true.  Let $T$ be an $n$-dimensional locally area minimizing rectifiable current in an open subset $U \subset {\mathbb R}^{n+m}$.  Let ${\rm sing}_{q} \, T = {\rm sing}\,T \cap \{ X : \Theta(T,X) = q \}$.  Then there exists an open set $V_q$ in $\mathbb{R}^{n+m}$ with ${\rm sing}_{q} \, T \subset V_q$ such that 
\begin{equation*}
    {\rm sing}\, T \cap V_q \cap \{ X : \Theta(T,X) \geq q \} = \mathcal{B}_q \cup \mathcal{S}_q 
\end{equation*}
where $\mathcal{B}_q$ and $\mathcal{S}_q$ are locally compact sets such that:
\begin{enumerate}[itemsep=2mm,topsep=0mm]
    \item[{\rm (i)}]  $\mathcal{B}_q \cap \mathcal{S}_q = \emptyset$; 
    
    \item[{\rm (ii)}]  $\mathcal{B}_q$ is relatively closed in $V_q$;

    \item[{\rm (iii)}]  $\mathcal{S}_q$ is locally $(n-2)$-rectifiable (and in particular $\mathcal{S}_q$ has locally finite $\mathcal{H}^{n-2}$-measure in the sense that for every $Z \in \mathcal{S}_q$ there exists $\rho > 0$ such that $\mathcal{H}^{n-2}(\mathcal{S}_q \cap \mathbf{B}_{\rho}(Z)) < \infty$);

    \item[{\rm (iv)}]  for every $Z \in \mathcal{B}_q$, there is an $n$-dimensional plane $P_Z$ such that the unique tangent cone to $T$ at $Z$ is equal to $q \llbracket P_Z \rrbracket$ (with the orientation induced by $T$) and there is a number $\sigma_Z > 0$ such that 
    \begin{eqnarray*}
        &&\rho^{-n-2} \int_{B_{\rho}(Z)} {\rm dist}^{2}(X,  Z + P_{Z})\, d\|T\|(X)\\ 
        &&\hspace{1in}\leq C \left(\frac{\rho}{\sigma}\right)^{2\alpha} \sigma^{-n-2} \int_{B_{\sigma}(Z)} {\rm dist}^{2}(X,  Z + P_{Z})\, d\|T\|(X) 
    \end{eqnarray*}
    for all $\rho, \sigma$ with $0 < \rho \leq \sigma \leq \rho_{Z}$;

    \item[{\rm (v)}]  at $\mathcal{H}^{n-2}$-a.e.~$Z \in \mathcal{S}_q$, $T$ has a unique tangent cone $\mathbf{C}_Z$ of the form $\mathbf{C}_Z = \sum_{j=1}^{p_Z} q_j^{(Z)} \llbracket P_j^{(Z)} \rrbracket$ where $p_Z \geq 2$ and $q^{(Z)}_{1}, \ldots, q^{(Z)}_{p} \geq 1$ are integers and $P^{(Z)}_{1}, \ldots, P^{(Z)}_{p}$ are distinct $n$-dimensional oriented planes such that there is an $(n-2)$-dimensional subspace $L_Z$ with $P^{(Z)}_{i} \cap P^{(Z)}_{j} = L_Z$ for every $i \neq j$ and there is a number $\rho_Z > 0$ such that 
    \begin{eqnarray*}
        &&\rho^{-n-2} \int_{B_{\rho}(Z)} {\rm dist}^{2}(X,  Z + {\rm spt} \, {\mathbf C}_{Z})\, d\|T\|(X)\\ 
        &&\hspace{1in}\leq C \left(\frac{\rho}{\sigma}\right)^{2\alpha} \sigma^{-n-2} \int_{B_{\sigma}(Z)} {\rm dist}^{2}(X,  Z + {\rm spt} \, {\mathbf C}_{Z})\, d\|T\|(X) 
    \end{eqnarray*}
    for all $\rho, \sigma$ with $0 < \rho \leq \sigma \leq \rho_{Z}.$
\end{enumerate}
\end{theorem} 

\begin{proof} Let $\b_{\star} = \b_{\star}(n, m ,q) \in (0, 1)$ and $\gamma_{\star} = \gamma_{\star}(n, m, q) \in (0, 1)$ be the constants as in Theorem~\ref{decomposition}. Let $R = R(n, m, q, \b_{\star}, \b_{\star}) \in [2, \infty)$, $\d = \d(n, m, q, \b_{\star}, \b_{\star}) \in (0, 1)$, 
$\eta = \eta(n, m, q, \b_{\star}, \b_{\star}) \in (0, 1)$ and $\alpha = \alpha(n, m, q, \b_{\star}, \b_{\star}) \in (0, 1)$ be as in \cite[Theorem~1.1]{KrumWica} taken with 
$\epsilon = \b = \b_{\star}$.  
For each $Z \in {\rm sing}_{q}\,T$ there is $\sigma_{Z} >0$ such that the hypotheses of \cite[Theorem~1.1]{KrumWica} (taken with $\epsilon = \b = \b_{\star}$) are satisfied with $T_{Z} \equiv \eta_{Z, \sigma_{Z} \, \#} \, T$ in place of $T$.  
Let $\mathcal{B}_Z$ and $\mathcal{S}_Z$ is the sets $\mathcal{B}$ and $\mathcal{S}$ given by \cite[Theorem~1.1]{KrumWica} with $T_Z$ in place of $T$.
Set $V_q = \bigcup_{Z \in {\rm sing}_Q\,T} \mathbf{B}_{\sigma_Z}(Z)$, $\mathcal{S}_q = \bigcup_{Z \in {\rm sing}_Q\,T} \eta_{Z,\sigma_Z}^{-1}\mathcal{S}_Z,$ 
and 
\begin{equation*}
    \mathcal{B}_q = {\rm sing}_q\,T \cap V_q \cap \{ X : \Theta(T,X) \geq q \} \setminus \mathcal{S}_q .
\end{equation*}
Clearly $\mathcal{B}_q \cap \mathcal{S}_q = \emptyset$ (as in conclusion~(i)).  Since each $\mathcal{B}_Z$ is relatively closed in $\mathbf{B}_1(0)$ and $\{ X : \Theta(T,X) \geq q \}$ is closed by the upper-semi continuity of density, $\mathcal{B}_q$ is relatively closed in $V_q$ (as in conclusion~(ii)) and thus both $\mathcal{B}_q$ and $\mathcal{S}_q$ are locally compact.

We conclude from part ($II$) of \cite[Theorem~1.1]{KrumWica} and \cite[Corollary~5.3]{KrumWica} that conclusion~(iv) of the present theorem holds.  For each $Z \in {\mathcal S},$ using $\mathcal{B}_Z$ being relatively closed in $\mathbf{B}_1(0)$ and part ($I$) of \cite[Theorem~1.1]{KrumWica} 
we can choose $\rho_{0}>0$ such that $\{X \, : \, \Th \, (T_{Z}, X) \geq q\} \cap {\mathbf B}_{\rho_{0}}(0) = {\mathcal S}_{Z} \cap {\mathbf B}_{\rho_{0}}(0),$  
${\mathcal S}_{Z} \cap \overline{{\mathbf B}_{\rho_{0}}(0)}$ is compact, and for every $Y \in {\mathcal S}_{Z} \cap {\mathbf B}_{\rho_{0}}(0)$ and every $\rho \in (0, \rho_{0}]$, either there is a cone 
${\mathbf C}_{Y, \rho} \in \cup_{p=2}^{q} {\mathcal C}_{q, p}$ such that 
\begin{equation}\label{structure-1}
Q(\eta_{Y, \rho \, \#} \, T_{Z}, {\mathbf C}_{Y, \rho}, {\mathbf B}_{1}(0)) < \b_{\star} \inf_{{\mathbf P} \in {\mathcal C}_{q, 1}} \, E(\eta_{Y, \rho \, \#} \, T_{Z}, {\mathbf P}, {\mathbf B}_{1}(0)), \;\; \mbox{or}
\end{equation}
\begin{equation}\label{structure-2}
\{X \, : \, \Th \, (\eta_{Y, \rho \, \#} \, T_{Z}, X) \geq q\} \cap {\mathbf B}_{1}(0) \subset \{X \, : \, {\rm dist}(X, L) \leq \beta_{\star}\} \cap {\mathbf B}_{1}(0)
\end{equation}
for some $(n-3)$-dimensional subspace $L$ of ${\mathbb R}^{n+m}.$ We can now reach conclusions~(iii) and (v) of the present theorem, with $\rho_{Z} = \rho_{0}\sigma_{Z}/2$, by arguing exactly as in the proof of \cite[Theorem~$2^{\prime}$]{Sim93}, using \eqref{structure-1}, \eqref{structure-2}  in places where that argument depends on \cite[Theorem~2.4]{Sim93}, and Theorem~\ref{decomposition} in places where it uses \cite[Theorem~1]{Sim93}. Specifically, if \eqref{structure-1} holds true with $Y=0$ and $\rho=1$, then we can apply Theorem~\ref{decomposition} with $\eta_{0,\rho_0\#} T_Z$ in place of $T$ to obtain  
\begin{equation*}
	\{X \in \mathbf{B}_{\rho_0/2}(0) : \Theta(T_Z,X) \geq q \} \subset \mathcal{L} \cup \bigcup_{j=1}^{\infty} \mathbf{B}_{\rho_j}(Y_j) \tag{$\star$} 
\end{equation*}
where $\mathcal{L}$ is a properly embedded $(n-2)$-dimensional $C^{1,\nu_{\star}}$-submanifold of $\mathbf{B}_{\rho_0/2}(0)$ for some $\mu_{\star} = \mu_{\star}(n,m,q) \in (0,1)$ and $\{\mathbf{B}_{\rho_j}(Y_j)\}$ is a countable collection of balls with $\rho_j < \rho_0/2$ and $\sum_{j=1}^{\infty} \rho_j^{n-2} \leq (1 - \gamma_{\star}) \,\rho_0^{n-2}$.  If instead \eqref{structure-2} holds true with $Y=0$ and $\rho=1$, then by a standard covering argument $(\star)$ holds true with $\mathcal{L} = \emptyset$ (see \cite[Theorem~2$^{\prime}$]{Sim93}).  We know that by \cite[Theorem~1.1]{KrumWica}, for each $j$ either there is a cone $\mathbf{C}_{Y_j,\rho_j} \in \bigcup_{p=2}^q \mathcal{C}_{q,p}$ such that \eqref{structure-1} holds true or \eqref{structure-2} holds true.  (Notice that the more elementary argument of \cite[Lemma~2.4]{Sim93}, which is based on the monotonicity formula for area and a compactness argument and which for its conclusions crucially relies on the assumption that 
the stationary varifolds considered belong to a multiplicity 1 class, gives us the weaker statement that for $\varepsilon_{\star} = \varepsilon_{\star}(n,m,q) \in (0,1)$ suitably small, there is an $n$-dimensional area-minimizing cone $\widetilde{\mathbf C}_j$, 
\emph{possibly supported on a plane}, such that $Q(\eta_{Y_j,\rho_j\#} T_Z, \widetilde{\mathbf C}_j, \mathbf{B}_1(0)) < \varepsilon_{\star}$ and $\{ X : \Theta(\eta_{Y_j,\rho_j\#} T_Z,X) \geq q\} \cap \mathbf{B}_1(0) \subset \{ X : \op{dist}(X, \op{spine} \widetilde{\mathbf C}_j) < \varepsilon_{\star} \}$; this is insufficient for continuing to apply Theorem~\ref{decomposition} or the necessary covering argument.)  For each $j$ such that \eqref{structure-1} holds true, we can again apply Theorem~\ref{decomposition}  with $\eta_{Y_j,\rho_j\#} T_Z$ in place of $T$, whereas for each $j$ such that \eqref{structure-2} holds true, we can apply a standard covering argument.  Iteratively applying this procedure as in \cite[Theorem~2$^{\prime}$]{Sim93} gives us conclusions~(iii) and (v). \end{proof}

Finally, we deduce the following result as an immediate consequence of Theorem~\ref{structure thm}.

\begin{corollary}\label{unique-tangent-and-rectifiability}
If $T$ is an $n$-dimensional locally area minimizing rectifiable current in an open set $U \subset {\mathbb R}^{n+m}$, then for 
${\mathcal H}^{n-2}$ a.e.\ point $Z \in {\rm spt} \, T$, the current $T$ has a unique tangent cone 
${\mathbf C}_{Z}$ of the form ${\mathbf C}_{Z} = \sum_{j=1}^{p} q_{j}\llbracket P_{j} \rrbracket$ where $p$, $q_{1}, \ldots, q_{p}$ are integers $\geq 1$ and 
$P_{1}, \ldots, P_{p}$ are distinct $n$-dimensional planes such that if $p \geq 2$ (i.e.\ if ${\mathbf C}_{Z}$ is not supported on a single plane) then there is an $(n-2)$-dimensional subspace $L$ with 
$P_{i} \cap P_{j} = L$ for every $i \neq j$.  Furthermore, we have that $${\rm sing} \, T = {\mathcal B} \cup {\mathcal S}$$ 
where: 
\begin{enumerate}[itemsep=2mm,topsep=0mm]
\item [{\rm(i)}] ${\mathcal B} \cap {\mathcal S} = \emptyset$;
\item[{\rm(ii)}] for every compact set $K \subset U$ 
\begin{equation*}
    {\mathcal S} \cap K = \bigcup_{j=1}^N {\mathcal L}_j
\end{equation*}
for some positive integer $N$ and for some pairwise disjoint, locally compact, locally $(n-2)$-rectifiable sets ${\mathcal L}_1,\ldots,{\mathcal L}_N$ (and in particular each set ${\mathcal L}_j$ has locally finite $\mathcal{H}^{n-2}$-measure); and 

\item[{\rm(iii)}] every point $Z \in {\mathcal B}$ is a branch point of $T$ where $T$ has a unique tangent cone supported on an $n$-dimensional plane $P_{Z}$;  moreover, for every compact set $K \subset U$ and every $Z \in {\mathcal B} \cap K$, the current $T$ ``decays rapidly'' to $P_{Z}$ in the sense that there are 
numbers $\alpha_{K} = \alpha(K, T) \in (0, 1)$ and $C_{K} = C(K, T) \in (0, \infty)$ such that the estimate in Theorem~\ref{structure thm}(a), with 
$\alpha = \alpha_{K}$ and $C = C_{K},$ holds for 
some $\sigma_{Z}>0$ (depending on $Z$) and all $\rho, \sigma$ with $0 < \rho \leq \sigma \leq \sigma_{Z}$. 
\end{enumerate} 
\end{corollary}

\begin{proof}
Recall from~\cite[Theorem~2.26]{Almgren} that the set $E$ of all singular points of $T$ where there is no tangent cone of the form $\sum_{j=1}^{p} m_{j}\llbracket P_{j} \rrbracket$ for some integer $p \geq 1$, positive integers $m_{1}, m_{2}, \ldots, m_{p}$ and distinct $n$-dimensional planes $P_{1}, \ldots, P_{p}$ has Hausdorff dimension $\leq n-3$.  
By Theorem~\ref{structure thm}, for each integer $q \geq 2$ there exists an open set $V_{q}$ such that ${\rm sing}_{q} T \subset V_{q}$ and 
\begin{equation*}
    {\rm sing}\,T \cap V_{q} \cap \{ X : \Theta(T,X) \geq q \} = \mathcal{B}_{q} \cup \mathcal{S}_{q} 
\end{equation*}
where $\mathcal{B}_{q}$ and $\mathcal{S}_{q}$ are pairwise disjoint, locally compact sets such that $T$ decays rapidly to a unique tangent plane at every point $Z \in \mathcal{B}_{q}$ and $\mathcal{S}_{q}$ is locally $(n-2)$-rectifiable.  
We claim that the conclusion of the theorem holds true with $\mathcal{B} = \bigcup_{q=2}^{\infty} \mathcal{B}_{q}$ and $\mathcal{S} = E \cup \bigcup_{q=2}^{\infty} \mathcal{S}_{q}$.
Let $K \subset U$ be a compact set.  For each compact subset $K \subset U$, the set $\{\Th \, (T, Z) \, : \, Z \in K\}$ is bounded, and hence there is a finite set  $\{q_{1}, q_{2} , \ldots, q_{N}\}$ of integers $\geq 2$, where 
$N = N(T, K),$  such that $\{\Th \, (T, Z) \, : \, Z \in {\rm sing} \, T \cap K \setminus E\} = \{q_{1}, q_{2}, \ldots, q_{N}\}$.  Set 
$\alpha_{K} = \min \{\alpha(n, m, q_{j}) \, : \, j=1, \ldots, N\},$ $C_{K} = \max\{C(n, m, q_{j}) \, : \, j=1, \ldots, N\}$ where $\alpha(n, m, q),$ $C(n, m, q)$ are as in Theorem~\ref{structure thm}(a).  
Set $q_0 = 1$, $q_{N+1} = \infty$, and $V_{q_0} = \mathcal{S}_{q_0} = \emptyset$.  For each $j = 0,1,2,\ldots,N$, let 
\begin{equation*}
    \Gamma_j = \{ X \in K : q_j \leq \Theta(T,X) < q_{j+1} \} \cap \mathcal{S}_{q_j}  
\end{equation*}
so that $\Gamma_j$ has locally finite $\mathcal{H}^{n-2}$-measure in $V_{q_j}$ and let 
\begin{equation*}
    \widetilde{\Gamma}_j = \{ X \in K : q_j \leq \Theta(T,X) < q_{j+1} \} \setminus V_{q_j} 
\end{equation*}
so that $\widetilde{\Gamma}_j \subseteq E$ and thus $\dim_{\mathcal H} \widetilde{\Gamma}_j \leq n-3$.  By the upper semi-continuity of density, $\Gamma_j$ is the intersection of $\mathcal{S}_{q_j}$, an open set, and a closed set and hence $\Gamma_j$ is locally compact.  Similarly, $\widetilde{\Gamma}_j$ is the intersection of an open set and closed set and hence $\widetilde{\Gamma}_j$ is also locally compact.  Noting that $\mathcal{S} \cap K = \bigcup_{j=0}^N (\Gamma_j \cup \widetilde{\Gamma}_j)$, the assertions of the present theorem can be verified with the help of Theorem~\ref{structure thm}.
\end{proof}

\bigskip
\hskip-.2in\vbox{\hsize3in\obeylines\parskip -1pt 
  \small 
Brian Krummel
School of Mathematics \& Statistics 
University of Melbourne
Parkville,VIC  3010, Australia

\vspace{4pt}
{\tt brian.krummel@unimelb.edu.au}} 
\vbox{\hsize3in
\obeylines 
\parskip-1pt 
\small 
Neshan Wickramasekera
DPMMS 
University of Cambridge 
Cambridge CB3 0WB, United Kingdom
\vspace{4pt}
{\tt N.Wickramasekera@dpmms.cam.ac.uk}
}

\end{document}